\documentclass[11pt,a4paper]{report}

\usepackage{amsmath,amssymb,amsthm,graphics,amscd}
\usepackage{longtable}
\usepackage{cite,xypic}
\usepackage{color}
\usepackage{graphicx}
\usepackage{epsf}
\usepackage{fancyhdr}
\usepackage[top=2.1cm, bottom=2.1cm, left=2.1cm, right=2.1cm]{geometry}

\long\def\symbolfootnote[#1]#2{\begingroup\def\thefootnote{\fnsymbol{footnote}}
\footnote[#1]{#2}\endgroup}
  \renewenvironment{thebibliography}[1]{
    \begin{oldthebibliography}{#1}
      \setlength{\parskip}{0ex}
      \setlength{\itemsep}{0ex}
  }
  {
    \end{oldthebibliography}
  }
  
\setlongtables
\begin{document}

\newcounter{rownum}
\setcounter{rownum}{0}
\newcommand{\ab}{\addtocounter{rownum}{1}\arabic{rownum}}
\newcommand{\im}{\mathrm{im}}
\newcommand{\x}{$\times$}

\newtheorem{lemma}{Lemma}[section]
\newtheorem{theorem}[lemma]{Theorem}
\newtheorem{corollary}[lemma]{Corollary}
\newtheorem{conjecture}[lemma]{Conjecture}
\newtheorem{prop}[lemma]{Proposition}
\theoremstyle{remark}
\newtheorem{remark}[lemma]{Remark}
\theoremstyle{definition}
\newtheorem{defn}[lemma]{Definition}

\theoremstyle{theorem}
\newtheorem*{TG1}{Theorem G1}
\newtheorem*{CG2}{Corollary G2}
\newtheorem*{CG3}{Corollary G3}
\newtheorem*{T1}{Theorem 1}
\newtheorem*{T2}{Theorem 2}
\newtheorem*{C2}{Corollary 2}
\newtheorem*{C3}{Corollary 3}
\newtheorem*{T4}{Theorem 4}
\newtheorem*{C5}{Corollary 5}
\newtheorem*{C6}{Corollary 6}
\newtheorem*{C7}{Corollary 7}
\newtheorem*{C8}{Corollary 8}
\newtheorem*{claim}{Claim}

\renewcommand{\labelenumi}{(\roman{enumi})}
\newcommand{\Hom}{\mathrm{Hom}}
\newcommand{\Ext}{\mathrm{Ext}}
\newcommand{\Aut}{\mathrm{Aut}}
\newcommand{\soc}{\mathrm{Soc}}

\newenvironment{changemargin}[1]{%
  \begin{list}{}{%
    \setlength{\topsep}{0pt}%
    \setlength{\topmargin}{#1}%
    \setlength{\listparindent}{\parindent}%
    \setlength{\itemindent}{\parindent}%
    \setlength{\parsep}{\parskip}%
  }%
  \item[]}{\end{list}}

\parindent=0pt
\addtolength{\parskip}{0.5\baselineskip}

\title{$G$-complete reducibility and the exceptional algebraic groups\vspace{40pt}}
\author{\vspace{100pt}David I. Stewart}
\date{A THESIS\\\vspace{30pt}Presented to the Department of Mathematics\\ of Imperial College London,\\in partial fulfillment of the requirements\\ for the degree of 
Doctor of Philosophy}

\maketitle
\begin{center}{\bf\Large Declaration of orginality}\end{center}
I hereby declare that the thesis here presented is my own work. All work appearing in it is my own except that which is referenced.\newpage

\begin{center}{\bf\Large Abstract}\end{center} Let $G=G(K)$ be a simple algebraic group defined over an algebraically closed field $K$ of characteristic $p>0$. A subgroup $X$ of $G$ is said to be $G$-completely reducible if, whenever it is contained in a parabolic subgroup of $G$, it is contained in a Levi subgroup of that parabolic. A subgroup $X$ of $G$ is said to be $G$-irreducible if $X$ is in no parabolic subgroup of $G$; and $G$-reducible if it is in some parabolic of $G$. In this thesis, we consider the case that $G$ is of exceptional type. When $G$ is of type $G_2$ we find all conjugacy classes of closed,  connected, reductive subgroups of $G$. When $G$ is of type $F_4$ we find all conjugacy classes of closed, connected, reductive $G$-reducible subgroups $X$ of $G$. Thus we also find all non-$G$-completely reducible closed, connected, reductive subgroups of $G$. When $X$ is closed, connected and simple of rank at least two, we find all conjugacy classes of $G$-irreducible subgroups $X$ of $G$. Together with the work of Amende in \cite{Ame05} classifying irreducible subgroups of type $A_1$ this gives a complete classification of the simple subgroups of $G$. 

Amongst the classification of subgroups of $G=F_4(K)$ we find infinite collections of subgroups $X$ of $G$ which are maximal amongst all reductive subgroups of $G$ but not maximal subgroups of $G$; thus they are not contained in any maximal reductive subgroup of $G$. The connected, semisimple subgroups contained in no maximal reductive subgroup of $G$ are of type $A_1$ when $p=3$ and of semisimple type $A_1^2$ or $A_1$ when $p=2$. Some of those which occur when $p=2$ act indecomposably on the $26$-dimensional irreducible representation of $G$. 

These results extend the complete reducibility results of \cite[Thm 1]{LS96}. 

We use this classification to prove a number of corollaries relating the subgroup structure of $G=F_4$ to the representation theory of subgroups of $G$ on the $26$-dimensional Weyl module for $G$. For instance, we prove that two isomorphic subgroups of $G$ with composition factors which correspond under the isomorphism are contained in associated parabolics.

We also use this classification to find all subgroups of $G=F_4$ which are generated by short root elements of $G$, by utilising and extending the results of \cite{LS94}.
\newpage
\begin{center}{\bf\Large Acknowledgements}\end{center}
\begin{quote}Rationality went down the drain donkey's years ago and hasn't been seen since.
\begin{flushright}-- Harold Pinter, {\it Moonlight}\end{flushright}\end{quote}

I feel remarkably lucky to have been supervised by Martin Liebeck, who not only lavished time on me and this document but also introduced me to some fascinating mathematics. He must also be thanked for an incalculable amount of time writing references for me and for aiding and abetting my vandalism of a Schubert piano trio. I am very grateful to him.

I also need to thank my examiners: Steve Donkin and Sasha Ivanov, firstly for finding this document passable; and secondly, for suggesting improvements and helping me removing mistakes.

I would also like to thank the following people for helpful discussions and/or moral support: Owen Jones, James Newton, Kevins McGerty and Buzzard, Ed Segal, Donna Testerman, Gary Seitz and the cast and crew of PTGCs 2007 thru 2009.

Finally, I wish to thank my parents and GMW, Pharaoh of Stevenage. They know who they are.

This document was produced under the financial support of the EPSRC.
\tableofcontents

\chapter{Introduction and results}\label{intro}
This thesis concerns the subgroup structure of semisimple algebraic groups over algebraically closed fields. Let $G$ be a simple algebraic group over an algebraically closed field $K$ of characteristic $p$. When $p$ is $0$ the maximal closed subgroups of positive dimension in $G$ were classified in a paper of Dynkin in 1952, \cite{Dyn52}. Seitz \cite{Sei87} initiated the study in positive characteristic, and following this, Liebeck and Seitz \cite{LS04} completed the determination of the maximal closed subgroups of positive dimension in the case where $G$ is of exceptional type. The list is surprisingly small: parabolic subgroups, reductive subgroups of maximal rank, and just a small (finite) number of conjugacy classes of further reductive maximal subgroups.
 
For many reasons one wants to extend this work to study all reductive subgroups of simple algebraic groups (not just maximal ones). For example, a result of Richardson \cite{Ric77} characterises those varieties $G/H$ which are affine as being precisely those for which $H^\circ$ is a closed reductive subgroup of $G$. One might think that finding the connected reductive subgroups $H^\circ$ is simply a matter of working down from the maximal subgroups, but one quickly realises that there is a major obstacle to this -- namely, dealing with reductive subgroups in parabolics. If $H$ is a reductive subgroup of a parabolic $P = QL$ ($Q$ the unipotent radical, $L$ a Levi subgroup), one can continue working down through the maximal subgroups of $L$ if $H$ is contained in a conjugate of $L$, but, as we describe later, it is much less easy to find those $H$ which are not in any conjugate of $L$. This leads to the following definition of a $G$-completely reducible subgroup, originally introduced by Serre \cite{Ser98}.
 
Let $H$ be a subgroup of a reductive algebraic subgroup $G$. Then $H$ is said to be {\it $G$-completely reducible ($G$-cr)} if, whenever it is contained in a parabolic subgroup $P$ of $G$, it is contained in a Levi subgroup $L$ of $P$. It is said to be {\it $G$-irreducible ($G$-ir)}, when it is contained in no parabolic of $G$. It is said to be {\it $G$-reducible ($G$-red)} if it is not $G$-ir. It is said to be {\it $G$-indecomposable ($G$-ind)} when it is not contained in any Levi subgroup of $G$.

These definitions generalise the notions of a representation of a group $\phi: X\to GL(V)$ being completely reducible, irreducible or indecomposable respectively since in the case where $G=GL(V)$ the image of $X$, $\phi X$ is a subgroup of $GL(V)$ which is $G$-cr, $G$-ir or $G$-ind, resp. 

An analogue of Maschke's Theorem due to Dynkin tells us that if the characteristic of the field $K$, $p$ is $0$ then all subgroups of $G(K)$ are $G$-cr. Following the usual paradigm, the obvious question is: {\it can we ensure that a subgroup $H$ of $G$ is $G$-cr if we know that $p$ is in some sense big?}

Indeed, for the simple algebraic groups, classical and exceptional, respectively, McNinch \cite{McN98} and Liebeck-Seitz \cite{LS96} have very specific answers. We give the coarsest versions of each of the main results in these papers. For $G$ a classical group, let $V$ be the natural representation of $G$ and $l$ the rank of a reductive subgroup $H$. Then by McNinch, $V$ is semisimple as an $H$ module provided $p\geq \dim V/l$. A further result of Serre's says that if additionally, $p>n_G(V)$, for some number $n_G(V)$ depending only on $G$ we know that $H$ is $G$-cr. Meanwhile, if $G$ is an exceptional group, then by Liebeck-Seitz, $H$ is $G$-cr provided $p> N(H,G)$, where $N(A,B)$ is a prime number always less than or equal to $7$ depending only on the type of $A$ and the type of $B$. For example, $N(A_1,E_6)=5$.

In this thesis, we extend results on the subgroup structure of the exceptional algebraic groups $G_2$ and $F_4$ to characteristics excluded by Liebeck-Seitz. In particular, we find all non-$G$-cr semisimple subgroups of $G$ when $G=G_2$ or $F_4$. (To distinguish between the main results for $G_2$ and the remaining results of the paper, we refer to Theorems G1 etc.)

In case $G = G_2(K)$ our results are particularly easy to state. The only candidate for such a non-$G$-cr subgroup is of type $A_1(K)$ (since it has to lie in a parabolic). On the face of it there could be infintely many classes of non $G$-cr subgroups, but we prove there are at most 2;

\begin{TG1}Let $X\cong A_1(K)$ be a subgroup of a parabolic subgroup in $G=G_2(K)$.

If $p>2$ then $X$ is $G$-cr (two classes). 

If $p=2$ then either $X$ is $G$-cr (two classes) or $X$ is conjugate to precisely one of $Z_1$ and $Z_2$.\end{TG1}

Here $Z_1$ is embedded in $A_1\tilde A_1$ via the map $x\mapsto (x,x)$ and $Z_2$ is embedded in $SL_3$ via the 3-dimensional symmetric square representation of the natural module for $SL_2$. We get easily the following

\begin{CG2}\label{CG2} All connected reductive subgroups of $G=G_2(K)$ are $G$-cr unless $p=2$, in which case there are precisely two classes of non $G$-cr subgroups.\end{CG2}

Thus we observe that the number $N(A_1,G_2)$ in \cite{LS96} can be reduced from 3 to 2. Using this result and the classification of maximal subgroups of $G_2$ from \cite{LS04} we are able to give a complete list of the connected reductive subgroups of $G_2$. Firstly recall some notation: 

Let $B$ be a Borel subgroup of a reductive algebraic group $G$, containing a maximal torus $T$ of $G$ with character group $X(T)$. We can identify $X(T)$ with $\mathbb Z^r$ for $r$ the rank of $G$. The dominant weights $X^+(T)$ are given by $(a_1,\dots a_r)\in \mathbb Z_r$ with $a_i\geq 0$ for each $i$. Recall that for each dominant weight $\lambda\in X(T)^+\subseteq X(T)$ for $G$, the space $H^0(\lambda):=H^0(G/B,\lambda)=\mathrm{Ind}_B^G(\lambda)$ is a $G$-module with highest weight $\lambda$ and with socle $\mathrm{Soc}_G H^0(\lambda)=L(\lambda)$, the irreducible $G$-module of highest weight $\lambda$. The Weyl module of highest weight $\lambda$ is $W(\lambda)\cong H^0(-w_0\lambda)^*$ where $w_0$ is the longest element in the Weyl group.  For brevity we write $L((a_1,a_2,\dots,a_r))=a_1a_2\dots a_r$; if $G$ is of type $A_1$ we write $L(n)=n$.  When $0\leq a_i<p$ for all $i$, we say that $\lambda$ is a restricted weight and we write $\lambda \in X_1$. Recall that any module $V$ has a Frobenius twist $V^{[n]}$ induced by raising entries of matrices in $GL(V)$ to the $p^n$th power. Steinberg's tensor product theorem states that $L(\lambda)=L(\lambda_0)\otimes L(\lambda_1)^{[1]}\otimes\dots\otimes L(\lambda_n)^{[n]}$ where $\lambda_i\in X_1(T)$ and $\lambda=\lambda_0+p\lambda_1+\dots+p^n\lambda_n$ is the $p$-adic expansion of $\lambda\in\mathbb Z_+^r$. Recall that there is a unique indecomposable tilting module $T(\lambda)$ of highest weight $\lambda$. Where we write $V_1|V_2|\dots |V_n$ we list the composition factors $V_i$ of an $X$-module. For a direct sum of $X$-modules, we write $V_1+V_2$. Where a module is uniserial, we will write $V_1/\dots/V_n$ to indicate the socle (and radical) series: here the head is $V_1$ and the socle $V_n$.

\begin{CG3} Let $X$ denote a closed, connected semisimple subgroup of $G=G_2(K)$. Then up to conjugacy, $(X,p,V_7\downarrow X)$ is precisely one entry in the following table where $V_7\downarrow X$ denotes the restriction of the seven-dimensional Weyl module $W_G(\lambda_1)$ to $X$.
\begin{center}
\begin{tabular}{|c|c|c|}
\hline
$X$ & $p$ & $V_7\downarrow X$\\
\hline
$A_2$ & any & $10\oplus 01\oplus 0$\\
$\tilde A_2$ & $p=3$ & $11$\\
$A_1 \tilde A_1$ & any & $1\otimes \tilde 1\oplus 0 \otimes \tilde W(2)$\\
$\bar L_0$ & any & $1\oplus 1\oplus 0^3$ \\
$\tilde L_0$ & any &  $1\oplus 1\oplus W(2)$ \\
$Z_1$ & $p=2$ & $T(2)\oplus W(2)$  \\
$Z_2$ & $p=2$ & $W(2)\oplus W(2)^* \oplus 0$ \\
$A_1\hookrightarrow A_1 \tilde A_1; (1,\tilde 1)\downarrow X=1^{[r]}\otimes 1^{[s]}$ $r\neq s$ & any & $(1^{[r]}\otimes1^{[s]})\oplus W(2)^{[s]}$ \\
$A_1\hookrightarrow A_2$, irred & $p>2$ & $2\oplus 2\oplus 0$\\
$A_1$, max & $p\geq 7$ & $6$\\
\hline
\end{tabular}
\end{center}
\end{CG3}

In the table, $\tilde L_0$ is the long $A_1$ Levi subgroup and $\tilde L_0$ is the short. The subgroup denoted $\tilde A_2$ exists only when $p=3$ and is generated by the short root subgroups of $G$. (The above table appears not to contain the irreducible $A_1\leq \tilde A_2$. It is shown later that this  subgroup is conjugate to the subgroup $A_1\hookrightarrow A_1\tilde A_1$ where $r=1$, $s=0$.)

A subgroup is called a {\it subsystem subgroup} if it is generated by root subgroups corresponding to a closed subset of the roots. For $G_2$, the subsystems up to conjugacy are the two Levi subgroups $A_1$, $A_2$, ($\tilde A_2$ if $p=3$) and $A_1\tilde A_1$. Thus we have that in any characteristic, {\it all} connected reductive $G_2$-reducible subgroups of $G_2$ are contained in subsystem subgroups. 

The situation when $G=F_4$ is not as neat. In particular, one can have reductive non-$G$-cr subgroups of $G$ not contained in any reductive subgroup which is maximal in $G$. Thus there are (maximal reductive) subgroups which are maximal amongst all reductive subgroups of $G$ but are not themselves (reductive) maximal subgroups of $G$. Even so, we are able to give a complete description of all reductive $G$-reducible subgroups; where such subgroups are not properly contained in reductive subgroups, we are able to give a complete description of these in terms of generators.

The following provides a summary of our findings related to $G$-complete reducibility of reductive subgroups $X$ of $G$. For the remainder of this section, $G=F_4(K)$. We start with the case where $X$ is simple. Recall that $G$ has subsystem subgroups of type $B_4$, $C_4$ ($p=2$), $D_4$, $\tilde D_4$ $(p=2)$, $A_2\tilde A_2$, $B_2B_2$ ($p=2$), where $C_4$ (resp $\tilde D_4$) is the image under the graph automorphism of $B_4$ (resp. $D_4$) and $\tilde D_4$ and $\tilde A_2$ are generated by short root groups.

\begin{T1}{\rm (A).} Let $X$ be a closed, connected, simple subgroup of $G$ of rank at least $2$. Suppose $X$ is non-$G$-cr. Then $p=2$ or $3$ and $X$ is in a subsystem subgroup of $G$. Moreover $X$ is conjugate to precisely one of the nine subgroups in the table below,

\begin{center}\begin{tabular}{c|c|l}\hline
$X$ of type & $p$ & Description\\
\hline
$A_2$ & $p=3$ & $X\leq B_4\leq GL_9$ via $V_9\downarrow X=(10, 01)$\\
$A_2$ & $p=3$ & $X\leq A_2\tilde A_2$ by $(V_3,V_3)\downarrow X=(10, 10)$
\\\hline
$B_2$ & $p=2$ & $X\leq D_3$\\
$B_2$ & $p=2$ & $X\leq \tilde D_3$\\
$B_2$ & $p=2$ & $X\leq B_2B_2$ by $(V_5,V_5)\downarrow X=(10, 10)$\\
$G_2$ & $p=2$ & $X\leq D_4$\\
$G_2$ & $p=2$ & $X\leq \tilde D_4$\\
$B_3$ & $p=2$ & $X\leq D_4$\\
$B_3$ & $p=2$ & $X\leq \tilde D_4$. \end{tabular}\end{center}

(Note that these embeddings are uniquely determined up to conjugacy in $G$. There are three classes of $B_3$ in $D_4$ which are conjugate under triality, but the Weyl group of $G$ induces such a triality automorphism on $D_4$.)

{\rm (B).} Let $X$ be a closed, connected subgroup of type $A_1$. Suppose $X$ is non-$G$-cr. Then $p=2$ or $3$. 

If $p=3$, then $X$ is conjugate to just one of the following subgroups:
\begin{enumerate}\item $X\hookrightarrow B_3$ by $V_7\downarrow X=T(4)+ 0$;
\item $X\hookrightarrow C_3$ by $V_6\downarrow X=T(3)$;
\item $X\hookrightarrow B_4$ by $V_9\downarrow X=T(4)^{[r]}+ 2^{[s]}$, $rs=0$;
\item $X\hookrightarrow A_2\tilde A_2$ by $(V_3, \tilde V_3)\downarrow X=(2,\tilde 2)$;
\item $X\hookrightarrow A_1C_3$ by $(1^{[r]},T(3)^{[s]})$, $rs=0$; or
\item $X$ in a one-dimensional variety $\Omega\cong \mathbb A^1(K)\setminus \{0\}$ parameterising subgroups of type $A_1$ in no maximal reductive subgroup of $G$. These subgroups are given by generators described in \S\ref{120}. \end{enumerate}

If $p=2$ then $X$ is conjugate to just one of the following subgroups:
\begin{enumerate}\item $X\hookrightarrow A_2$ or $\tilde A_2$ by $V_3\downarrow X=W(2)$;
\item $X\hookrightarrow \tilde A_1^2\leq B_2$ or $A_1^2\leq B_2$ by $x\mapsto (x,x)$;
\item $X\hookrightarrow A_3$ or $\tilde A_3$ by $V_4\downarrow X=T(2)$;
\item $X\hookrightarrow A_1^4$ or $\tilde A_1^4$ by $x\mapsto (x,x,x,x)$;
\item $X\hookrightarrow A_2\tilde A_2$ by $(V_3,V_3)\downarrow X=(2/0,0/2)$ or $(2/0,2/0)$.
\item $X$, not one of the above, in a $2$-dimensional variety parameterising subgroups of type $A_1$ in a $B_4$ or $C_4$ subsystem subgroup of $G$, with embeddings described in \S\ref{a1sp2inb3c3b2} or \S\ref{seca1ina1a1char2};
\item $X$ in a $2$-dimensional variety parameterising subgroups of type $A_1$ not in any reductive maximal subgroup of $G$ given by generators in \S\ref{seca1ina1a1char2}.
\end{enumerate}
\end{T1}

The following theorem deals with non-simple, semisimple subgroups.

\begin{T2}Let $X$ be a non-simple semisimple closed, connected subgroup of $G$. Then either $X$ is $G$-cr or precisely one of the following holds:
\begin{enumerate}\item $p=3$ and $X$ is of type $A_1A_1$; moreover, $X$ is conjugate to exactly one of the following two subgroups:
\begin{enumerate}\item $X\hookrightarrow D_3\tilde A_1\leq B_4$ by $(V_6,V_2)\downarrow X=(T(4),1)$
\item $X\hookrightarrow A_1C_3$ by $(V_2,V_6)\downarrow X=(1,T(3))$\end{enumerate}
\item $p=2$ and $X$ is of type $A_1A_2$ in an $A_2\tilde A_2$ subsystem. Moreover, $X$ is conjugate to exactly one of the following two subgroups:
\begin{enumerate}\item $X\leq A_2\tilde A_2$ by $(V_3,\tilde V_3)\downarrow X=(2/0,10)$
\item $X\leq A_2\tilde A_2$ by $(V_3,\tilde V_3)\downarrow X=(10,2/0)$
\end{enumerate}
\item $p=2$ and $X$ is of type $A_1^2$ with $X$ in a subsystem subgroup of type $A_2\tilde A_2$. Moreover $X$ is conjugate to precisely one of the following four subgroups:
\begin{enumerate}\item $X$ is in an $A_1\tilde A_2$ subsystem by $(V_2,\tilde V_3)\downarrow X=(1,2/0)$
\item $X$ is in an $A_2\tilde A_1$ subsystem by $(V_3,\tilde V_2)\downarrow X=(2/0,1)$
\item $X\leq A_2\tilde A_2$ by $(V_3,\tilde V_3)\downarrow X=(2/0,2/0)$, 
\item $X\leq A_2\tilde A_2$ by $(V_3,\tilde V_3)\downarrow X=(0/2,2/0)$\end{enumerate}
\item $p=2$ and $X$ is of type $A_1^2$, in a $B_4$ or $C_4$ subsystem subgroup, but not conjugate to either (iii)(a) or (b) above; embeddings of these are given in \S\ref{a1a1s}.
\item $p=2$ and $X$ is of type $A_1^3$ or $A_1B_2$ with $X$ in a subsystem subgroup of type $B_4$ or $C_4$, with embeddings described in \S\ref{a1a2ora1b2}, \S\ref{a1a1a1sec}.
\item $p=2$ and $X$ is of type $A_1^2$ with $X$ in no reductive maximal subgroup of $G$. Moreover $X$ is explicitly known, parameterised by some $2$-dimensional variety, as described explicitly in \ref{a1a1s}\end{enumerate}
\end{T2}

Most of the time, then we find that there are only finitely many conjugacy classes of subgroups of a given type, and these are all properly contained in proper reductive overgroups. The following corollary summarises what happens.
\begin{C3}(A) Let $X$ be any connected reductive subgroup of $G$ with $X/Z(X)$ not of type $A_1$ or $A_1A_1$. Then $X$ is contained in a reductive maximal subgroup of $G$.

(B) When $p=2$ or $3$, there exist infinitely many classes of reductive subgroups not contained in any  reductive maximal subgroup of $G$. Any such reductive subgroup $X$ is of type $A_1$ when $p=3$, whereas if $p=2$, $X$ is of type $A_1A_1$ or $A_1$.

(C) When $p=2$ there is a one-dimensional collection of such $A_1$s acting indecomposably on the $26$-dimensional irreducible module for $G$.

(D) There are subgroups of $G$ which are maximal amongst connected reductive subgroups of $G$ yet are not maximal subgroups of $G$. These occur when $p=3$ and $X$ is of type $A_1$ and when $p=2$ and $X$ is of type $A_1^2$.\end{C3}

The above theorems extend the $G$-complete reducibility results of \cite[Thm 1]{LS96} to cover characteristics excluded by the hypotheses there.

The following theorem deals with simple subgroups in no parabolic of $G$.

\begin{T4} \label{girr}Let $X$ denote a closed, connected simple subgroup of $G$ of rank at least $2$. If $X$ is $G$-irreducible, then one of the following holds \begin{enumerate}\item $X$ is maximal, ($B_4$, $C_4$ $(p=2)$, $G_2$ $(p=7)$);
 \item $X$ is a subsystem subgroup of maximal rank, ($B_4$, $C_4$ $(p=2)$, $D_4$, $\tilde D_4$ $(p=2)$);
 \item $X$ is of type $A_2$ embedded into a subgroup of type $A_2\tilde A_2$ as $(10^{[r]},01^{[s]})$, or $(10^{[r]}.10^{[s]})$ where if $p=3$ then we exclude $(10,01)$; or
 \item $p=2$ and $X$ is of type $B_2$ embedded into a subsystem subgroup of type $B_2 B_2$ as $(10,01)$, or if $r\neq s$,  as $(10^{[r]},10^{[s]})$, $(10^{[r]},01^{[s]})$ or $(01^{[r]},01^{[s]})$. \end{enumerate}\end{T4}

Adding to the above the work of Amende \cite{Ame05}, who classifies all $G$-ir subgroups of type $A_1$, this gives a complete classification of semisimple subgroups contained in parabolics and simple $G$-ir subgroups:

\begin{C5}Let $X$ be \begin{itemize}\item any simple, closed subgroup of $G$ of rank at least $2$; or
\item any connected semisimple subgroup of $G$ contained in a parabolic subgroup of $G$.\end{itemize} Then $X$ is conjugate to just one item in the Table in \ref{tor} below.\end{C5}

Two parabolic subgroups are said to be {\it associated} if they have conjugate Levi subgroups. A check on the classification of $G$-reducible subgroups of $G$ shows the following to be true:

\begin{C6} Let $X\cong Y$ be connected reductive subgroups of $G$ contained in parabolic subgroups of $G$ and let $P_X$ and $P_Y$ be minimal parabolic subgroups containing $X$ and $Y$ respectively. Then if the composition factors $V_{26}\downarrow X$ and $V_{26}\downarrow Y$ correspond under the isomorphism $X\cong Y$, then $P_X$ and $P_Y$ are associated parabolics of $G$.\end{C6}

{\it Remark. } The above Corollary extends the obvious statement in the case $G=SL_{n+1}(K)$: The dimensions of the composition factors of the natural $n+1$-dimensional module for $G$ determine the Levi type of a minimal parabolic $P$ of $G$ containing $X$. If the composition factors are isomorphic then these parabolics must therefore be associated.

One imagines there could be many possibilities for employment of our classification. One such is the following corollary, an extension of the results of \cite{LS94}. Recall that a short (long) root element is a member $x_\alpha(t)$ of any conjugate of the root group corresponding to the short (long) root $\alpha$ of $G$. In \cite{LS94} all simple, closed, connected  subgroups $X$ of all simple algebraic groups which contain long root elements are found. We find those containing short root elements when $G=F_4(K)$. As $X$ is simple, it is generated by any one of its conjugacy classes. Thus if it contains a short root element of $G$, it is generated by them.

\begin{C7} Let $X$ be a connected, closed, simple subgroup of $G$. Then $X$ is generated by short root elements if and only if exactly one of the following holds;
\begin{enumerate}\item $X$ is a subsystem subgroup of type $B_4$, $C_4$ $(p=2)$, $D_4$ ($p>2$), $\tilde D_4$ $(p=2)$, $A_3$ $(p>2)$, $\tilde A_3$ $(p=2)$, $B_3$, $C_3$, $\tilde A_2$, $B_2$ or $\tilde A_1$;
\item $p=2$ and $X$ is non-$G$-cr, conjugate to the one of subgroups $A_1\hookrightarrow \tilde A_2$ by $V_3\downarrow X=W(2)$, $B_2\leq \tilde D_3$, $B_3\leq \tilde D_4$ or $G_2\leq \tilde D_4$.
\item $p>2$ and $X$ is a subgroup of type $A_1$ in a subsystem subgroup of type $A_1^2$ corresponding to long roots projecting non-trivially to each factor.
\end{enumerate}
\end{C7}

More eccentric corollaries are also possible:

\begin{C8}The following statements are equivalent:
\begin{enumerate}\item All closed subgroups of $G$ of type $A_1$ are contained in closed subgroups of type $A_1^2$;
\item The characteristic $p$ of $K$ is $2$ or $5$.\end{enumerate}\end{C8}

\chapter{Overview}\label{overview}
Let $G$ be a simple algebraic group defined over an algebraically closed field $K$ of characteristic $p>0$. Let $X$ be a connected reductive subgroup of $G$. Here we give an overview of the work contained in this thesis.

In Chapter \ref{general} we cover the general theory we need for the proofs of the main results. We start in \S\ref{notation} with a section of notation. 

Much of our work revolves around finding $G$-reducible subgroups $X$ of $G$. That is, we are in the  following situation: 

We have a reductive subgroup $X$ inside a parabolic subgroup $P$ of $G$. Let $L$ be a Levi subgroup of $P$ so that $P=LQ$ is a semidirect product of $L$ with the unipotent radical of $P$. As $X$ is reductive, it can have no intersection with $Q$. If $\bar X$ denotes the image of $X$ under the canonical projection $\pi:P\to L$ we have $\bar X$ isomorphic to $X$. In particular $XQ=\bar XQ$ and $X\cap Q=\{1\}$. So $X$ is a complement to $Q$ in the semidirect product $\bar XQ$. By standard results, $X$ corresponds to a $1$-cocycle $\gamma:X\to Q$. Up to $Q$-conjugacy, these are in natural bijection with a set denoted by $H^1(X,Q)$ (defined in \S\ref{nonab}). Now $Q$ is non-abelian in general, and for general non-abelian $Q$ finding $H^1(X,Q)$  would be a very hard computation. However, by \ref{abslem}, $Q$ admits a certain filtration $Q=Q(1)\geq Q(2)\geq\dots$ so that successive quotients, called `levels', have the structure of $L$-modules. This allows us to calculate $H^1(X,Q)$ via the better understood abelian cohomology groups $H^1(X,Q(i)/Q(i+1))$.

In \S\ref{nonab}, we describe the ingredients of non-abelian cohomology that we need for these calculations. To begin we remind the reader of the filtration of unipotent radicals of parabolics by modules for Levi subgroups. This filtration has the useful property that it is central. This will permit us to apply results from the theory of Galois cohomology, in particular an exact sequence, \ref{les}, relating  cohomology in degree $0$, $1$ and $2$.

In general $H^1(X,Q)$ is not a vector space; just a pointed set. Our first task is to create an approximation of $H^1(X,Q)$ by way of a surjection  $\rho$ from a certain vector space $\mathbb V$ to $H^1(X,Q)$; this happens in Definition \ref{rho} and Proposition \ref{surj}. The map $\rho$ is sometimes strictly partially defined. This is due to the existence of non-split extensions between $X$ and subgroups of the radical. Such non-split extensions arise by non-trivial values of second cohomology. We make this link more explicit in Definition \ref{blockers}.

Now, we would like to use the map $\rho$ as an explicit parameterisation of the elements of $H^1(X,Q)$. Unfortunately $\rho$ is very far from being injective in general, though \ref{nofixedpoints} gives us a sufficient condition for $\rho$ to have that property. When $\rho$ is an isomorphism, this often leads to an easy observation that $H^1(X,Q)=K$. This has an immediate corollary that there is exactly one non-$G$-cr complement to $Q$ in $XQ$ up to $G$-conjugacy (see \ref{onedimh1}).

Indeed, a lot of the work revolves around establishing conjugacy between subgroups in parabolics. In \ref{subsystemlift} we see how this task can be made easier when certain collections of root groups in unipotent radicals of parabolics are contained in subsystem subgroups of $G$; this allows us sometimes to conclude that complements are contained in subsystems of $G$.

One of the last things we discuss in \S\ref{nonab} is how our pointed set $H^1(X,Q)$ behaves under restriction. For instance, by a classical result of Cline, Parshall and Scott we can show that $H^1(X,Q)\cong H^1(B,Q)$ for $B$ a Borel subgroup of $X$ with $Q$ considered as an algebraic $B$-group.

In \S\ref{conjofcomps} we discuss conjugacy amongst complements in parabolics in a much more general setting. Once one has calculated $H^1(X,Q)$, one has established $Q$-conjugacy amongst complements to $Q$ in $XQ$. But one needs to see how these $Q$-conjugacy classes fuse under $G$-conjugacy. While Propostion \ref{nomoreconj} and Theorem \ref{nomoreconj2} contain answers which may appear to make these question significantly more complicated, they can actually be employed quite easily in practice. (Often, in simple cases one finds that $G$-conjugacy is induced by $QZ(L)$-conjugacy.)

As explained above, we will break up the calculation of $H^1(X,Q)$ into calculations of abelian cohomology. In \S\ref{abelian}, we discuss results of abelian cohomology which we will employ to calculate $H^1(X,Q)$. This largely constitutes results on $1$- and $2$-cohomology of $G$-modules.

In Chapter \ref{g2} we find all reductive subgroups of $G_2$. Let $G=G_2(K)$. Most of the work involved is finding the non-$G$-cr subgroups of $G$ and so this is our first opportunity to apply the general machinery that we have set up in the earlier chapters. As $G_2$ is quite a small group, we find that the situation is rather easier than $F_4$ and so this chapter is a very good starting place for anyone interested in understanding the situation in $F_4$.

In Chapter \ref{f4} we begin the meat-and-potatoes calculations which produce our main theorems. Let $G=F_4(K)$.

We start in \S\ref{aux} with a collection of miscellaneous results which we will need repeatedly during our calculations.

Recall that our foremost task is to classify all $G$-reducible subgroups of $G=F_4(K)$. Thus we have $X\leq P\leq G$ for some parabolic subgroup $P=LQ$ of $G$. As $X$ is a complement to the image of $X$ under projection to a Levi subgroup of $P$, we first wish to find all such images up to conjugacy. Thus we may assume $P$ is minimal subject to containing $X$. Thus $\bar X(=\pi X)$ is in no parabolic subgroup of $L$. Hence for each standard Levi subgroup $L$ of $G$ we want a list of $L$-ir subgroups of $L$. This happens in \S\ref{lir}. We observe that this gives a full list of $G$-reducible, $G$-completely reducible subgroups of $G$ up to conjugacy in $G$. It remains to find the non-$G$-cr subgroups. 

In \S\ref{cur} we list the action of each $L$-ir subgroup of $L$ on the so-called levels $Q(i)/Q(i+1)$ of the unipotent radical $Q$ of the standard parabolic subgroup $P$ containing $L$, where $H^1(\bar X,Q(i)/Q(i+1))\neq 0$ for some $i$. By \cite[Thm 1]{LS96} we know that this happens only if $p=2$ or $3$, making our check simpler. Now, by the results of \S\ref{abelian} we can see by inspection the values of $H^1(X,Q(i)/Q(i+1))$ and thus the dimension of the vector space $\mathbb V$. This constitutes our first approximation to $H^1(X,Q)$ for relevant pair $(X,Q)$. If $H^1(X,Q)=0$ then $X$ must be $G$-cr as it is $Q$-conjugate to its projection to the Levi. Thus we need only consider pairs $(X,Q)$ such that there is a non-zero value of $H^1(X,Q(i)/Q(i+1)$ for some $i$.

From \S\ref{thm1} onwards, we are committed to case-by-case analysis. The simple subgroups of rank at least two are easier and so we tackle this first in \S\ref{geq2}. Then we move onto subgroups of type $A_1$ when $p=3$ in \S\ref{a1p3}; and when $p=2$ in \S\ref{a1p2}. Note that in many cases in the latter section the situation when $X$ is of type $A_1$ is virtually identical to the situation when $X$ is of type $A_1^2$ or $A_1^3$; and one passes between the two types by embeddings of the form  $A_1\hookrightarrow A_1^n$. We often find it preferable to prove results about the non-simple subgroups and deduce results for subgroups of type $A_1$. Thus we often defer to proofs from \S\ref{xnonsimp}, in which Theorem 2 is also proved, and we find all non-simple, non-$G$-cr subgroups of $G$.

In \S\ref{girsubgsec} we find all simple $G$-irreducible subgroups of $G$ of rank at least $2$ and prove Theorem 4. In some sense this relies on Theorem 1. The strategy is simply to work down through the maximal subgroups of $G$, obtaining a list of all candidates of subgroups and then eliminating some  conjugate subgroups and eliminating those which are in parabolics. From Theorem 1, we already have  descriptions of many subgroups inside parabolics and so this task is not so arduous. 

As we have explained there is a close connection between actions on the $26$-dimensional Weyl module $V_{26}$ for $G$ and the subgroup structure of $G$. We make use of this representation in several ways. Firstly, throughout the proof we are in constant dialogue with the table in the Appendix \S\ref{tor}, which gives restriction of $V_{26}$ to semisimple subgroups of $G$. By the end of the proof, we will have shown that any $G$-reducible semisimple subgroup of $G$ or any $G$-irreducible semisimple subgroup containing no factors of type $A_1$ is conjugate to precisely one item in this table, thus proving Corollary 5.

(We remark that the thesis of Amende \cite{Ame05} contains a determination of all $G$-irreducible subgroups of type $A_1$ of all exceptional groups $G$ not of type to $E_8$; together with our work there remains outstanding a further classification of $G$-ir subgroups of type $A_1^n$ of $F_4$. To do this would appear to involve duplication of a large proportion of that thesis, so we have not attempted such a classification here.)

We also use the representation $V_{26}$ for many technical calculations. We have been provided with $48$  $26\times 26$ matrices, {\tt F4gens}, in the GAP algebra system\footnote{We thank Prof. Frank L\"ubeck for the provision of these matrices}. Each of these has a single  indeterminate {\tt T} and each corresponds to some root group $x_\alpha(T)$ of $G$. Since $G$ is generated by its root groups, one can use these matrices to construct subgroups of $G$.  Typically we will have a description of a subgroup $X$ of $G$ by the generators of $X$. These generators will be given in the form of a product of root group elements of $G$. Therefore we may realise $X$ inside GAP. We may wish to show, for example, that this subgroup $X$ is in fact isomorphic to an algebraic group of type $A_1$. Since we have generators and relations for $A_1$ in \ref{a1rels}, we can ask GAP to perform the necessary matrix multiplications to check that the generators of $X$ given satisfy the relations given there. For example, let $p=2$ and suppose that one has a generator $x$ of $X$ and one believes that $X\cong A_1$ with $x$ mapping to a root group element $x_\alpha(t)$ of $A_1$. Since $p=2$ then, one needs $x$ to be an involution. To check this, one can input $x$ as a product of matrices of the set {\tt F4gens} in GAP and then check that the relation $x^2=1$ holds. Such calculations are used primarily in establishing $1$-cohomology sets $H^1(X,Q)$.

At the end of the thesis, in  \S\ref{corproof} we tie up the proofs of the other corollaries listed in the introduction.

\chapter{General Theory}\label{general}
\section{Notation}\label{notation}
When discussing roots or weights, unless a group is of type $G_2$ we take long fundamental roots on the left of the Dynkin diagram. Thus the diagram of $F_4$ is $\circ-\circ=>\circ-\circ$, with highest root $2342$. For $G_2$ we use the diagram $\circ<\equiv\circ$ with highest root $32$.

Let $B$ be a Borel subgroup of a reductive algebraic group $G$, containing a maximal torus $T$ of $G$. Recall that for each dominant weight $\lambda\in X^+(T)$ for $G$, the space $H^0(\lambda):=H^0(G/B,\lambda)=\mathrm{Ind}_B^G(\lambda)$ is a $G$-module with highest weight $\lambda$ and with socle $\mathrm{Soc}_G H^0(\lambda)=L(\lambda)$, the irreducible $G$-module of highest weight $\lambda$. The Weyl module of highest weight $\lambda$ is $W(\lambda)\cong H^0(-w_0\lambda)^*$ where $w_0$ is the longest element in the Weyl group. We identify $X(T)$ with $\mathbb Z^r$ for $r$ the rank of $G$ and for brevity write $L(a_1,a_2,\dots,a_r)=a_1a_2\dots a_r$; if $G$ is of type $A_1$ we write $L(n)=n$.  When $0\leq a_i<p$ for all $i$, we say that $\lambda$ is a restricted weight and we write $\lambda \in X_1(T)$. Recall that any module $V$ has a Frobenius twist $V^{[n]}$ induced by raising entries of matrices in $GL(V)$ to the $p^n$th power. Steinberg's tensor product theorem states that $L(\lambda)=L(\lambda_0)\otimes L(\lambda_1)^{[1]}\otimes\dots\otimes L(\lambda_n)^{[n]}$ where $\lambda_i\in X_1(T)$ and $\lambda=\lambda_0+p\lambda_1+\dots+p^n\lambda_n$ is the $p$-adic expansion of $\lambda\in\mathbb Z_+^r$.

We have also the notion of the $d$th Frobenius untwist $V^{[-d]}$ of a module $V$, that is the unique module $W$ such that $W^{[d]}=V$ where such a module exists.

The right derived functors of $\Hom(V,*)$ are denoted by $\Ext_G^i(V,*)$ and when $V=K$, the trivial $G$-module, we have the identity $\Ext_G^i(K,*)=H^i(G,*)$ giving the Hochschild cohomology groups. We are interested in the cases where $i=0,1$ or $2$ in this thesis. Recall that $H^0(G,V)\cong V^G$ is the space of $G$-fixed points of $V$.

Recall that $F_4(K)$ has a $26$-dimensional Weyl module $V_{26}$. When $p\neq 3$, $V_{26}$ is the irreducible representation of high weight $0001$. When $p=3$, $V_{26}$ has a one dimensional radical, with a $25$-dimensional irreducible quotient of high weight $0001$. We will often talk of restrictions of $V_{26}$ to various subgroups $X$ of $F_4(K)$. Also when $p=2$ we recall that $F_4$ admits a graph morphism $\tau$ which is an automorphism of $F_4(K)$ considered as an abstract group. Typically we will use this last fact to reduce our computations; for instance, when $p=2$ any subgroup of the $C_3$-parabolic is conjugate under $\tau$ to a subgroup of the $B_3$-parabolic, so classifying the latter classifies the former.

We will often want to consider restrictions of simple $G$-modules to reductive subgroups $X$ of $G$. Where we write $V_1|V_2|\dots |V_n$ we list the composition factors $V_i$ of an $X$-module. For a direct sum of $X$-modules, we write $V_1+V_2$. Where a module is uniserial, we will write $V_1/\dots/V_n$ to indicate the socle and radical series: here the head is $V_1$ and the socle $V_n$. On rare occasions we use $V/W$ to indicate a quotient. It will be clear from the context which is being discussed. 

Recall also the notion of a tilting module as one having a filtration by modules $V(\mu)$ for various $\mu$ and also a filtration by modules $H^0(\mu)$ for various $\mu$ (equiv. dual Weyl modules). For each $\lambda\in X(T)$ we have also the indecomposable tilting module $T(\lambda)$ of high weight $\lambda$. Each $T(\lambda)$ is an indecomposable module. Two key properties of tilting modules which we use are that $H^1(X,T(\lambda))=0$ for all $\lambda\in X(T)$ (see \ref{h1tilting} below). Also, that the tensor product of any two tilting modules is again tilting: \cite[1.2]{Don93}. As we are considering very low weight representations in general, it is possible to spot that a module is a $T(\lambda)$; for instance when $p=2$, the natural Weyl module for $B_n$ has a $1$-dimensional radical, so its structure is  $V(\lambda_1)=L(\lambda_1)/K$. It is then the case that giving the Loewy series for a module $K/L(\lambda_1)/K$ uniquely characterises it as a tilting module $T(\lambda_1)$.

Lastly, if $X=Y\times Z$ is a direct product of two simple groups then we will write an irreducible module $V_X$ for $X$ in the form $V_X=V_Y\otimes V_Z$. If $V_Y\cong L(\lambda)$ and $V_Z\cong L(\mu)$ then we will abbreviate this to $V_X=(L(\lambda),L(\mu))$.

We will often need to speak about specific parabolic subgroups of $F_4$. Usually we will just give the type of $L'$ indicating which factors are short by a tilde over that factor; in most cases this determines $P$ up to conjugacy. For subgroups of type $A_1\tilde A_1$, we denote the three non-conjugate standard parabolics by $P_{13}$, $P_{24}$, and $P_{23}$ obtained by deleting the stated nodes from the Dynkin diagram.

Recall the split short exact sequence $\{1\}\to Q\to P\stackrel{\pi}{\to} L\to \{1\}$, where $\pi$ is the projection of $P$ to $L$. Where we have a subgroup $X$ of $P$ we will often consider the image $\bar X:=\pi X$ of $X$ under this projection.

There are many concrete calculations later made with root elements and Weyl group elements of $G$. Assuming some choice of Borel has been made, we have, following \cite{Car89} root group elements, $x_{\alpha}(t)$;  representatives of Weyl group elements, $n_{\alpha}(t):=x_\alpha(t)x_\alpha(-t^{-1})x_\alpha(t)\in N_G(T)$, and elements of the torus, $h_{\alpha}(t):=n_\alpha(t)n_\alpha(-1)\in T$. Sometimes we find it visually clearer to give a product of many root elements in vector notation, writing $x_{\alpha}(t)=\left[\begin{array}{c}\alpha\\ t\end{array}\right]$. We will also abbreviate $n_{\alpha}(1)$ to $n_\alpha$. Then note that $x_{\beta}(t)^{n_\alpha}=x_{w_\alpha.\beta}(t)$.
\section{Non-abelian cohomology} \label{nonab}
 \renewcommand{\labelenumi}{(\roman{enumi})}
 \subsection{A filtration of unipotent radicals}

The constructions developed in this section are valid for any reductive algebraic group $G$ and provide a practical way to enumerate complements in parabolic subgroups by means of non-abelian cohomology.
 
 We recall the notions of height, shape and level of a root from \cite{ABS90}. Take a root system $\Phi$ for $G$ with fixed base of simple roots $\Pi$. Let $J\subset \Pi$ be a subset of the simple roots and define the parabolic subgroup $P_J$ by $P_J=\langle B, x_{-\alpha}(t): \alpha\in J\rangle$. Let $\Phi_J=\mathbb{Z}J\cap \Phi$. Fix a root $\beta\in\Phi^+ - \Phi_J$. We write $\beta=\beta_J+\beta_J'$ where $\beta_J=\sum_{\alpha_i\in J} c_i\alpha_i$ and $\beta_J'=\sum_{\alpha_i\in\Pi-J} d_i\alpha_i$.  Define
\begin{align*} \text{height}(\beta)&=\sum c_i+\sum d_i\\
\text{shape}(\beta)&=\beta_J'\\
\text{level}(\beta)&=\sum d_i.\end{align*}
Now define $Q(i):=\langle x_\beta(t): t\in K,\text{level}(\beta)\geq i\rangle$ and define $V_S=\langle x_\beta(t):t\in K, \text{shape}(\beta)=S\rangle$.

\begin{lemma}\label{abslem} Let $G=G(K)$ be a simple algebraic group. For each $i\geq1$, $Q(i)/Q(i+1)$ is central in $Q/Q(i+1)$ and has the structure of a $KL$-module with decomposition $Q(i)/Q(i+1)=\prod V_S$, the product over all shapes $S$ of level $i$. Furthermore, each $V_S$ is a $KL$-module with highest weight $\beta$ where $\beta$ is the unique root of maximal height and shape $S$.\end{lemma}
\begin{proof} This is the main result, Theorem 2, of \cite{ABS90}, noting the Remark 1 at the end of the paper which gives the result even in the case $G(K)$ is special.\end{proof}

We will always refer to this specific filtration of a unipotent radical $Q$ as $Q(i)$. 
 
  \subsection{Definitions of non-abelian cohomology sets}
\begin{defn}Let $Q$ and $X$ be an algebraic groups. Then we say $Q$ is an algebraic $X$-group if $X$ acts morphically on $Q$; i.e. there is a homomorphism $X\to\Aut_k(Q)$.\end{defn} 

We will now define a set $H^1(X,Q)$, where $X$ is an algebraic group and $Q$ is an algebraic $X$-group. Later we will be taking $Q$ the unipotent radical of some parabolic subgroup of $G$ and $X$ some $L$-irreducible subgroup of $L$, a Levi subgroup of $P$.
 
Firstly we define $H^0(X,Q)=\{q\in Q:q^x=q,$ for all $x\in X\}$ to be the group of fixed points of $Q$ under the action of $X$.

Now let $\gamma$ be a morphism from $X\to Q$. Then we say $\gamma$ satisfies the $1$-cocycle condition if $\gamma(xy)=\gamma(x)^y$. The set of such morphisms is denoted $Z^1(X,Q)$. If $q$ is an element of $Q$ then define a map $\bar q:X\to Q$ by $\bar q(x)=q^{-x}q$. We denote the set of such maps by $B^1(X,Q)$. It is easy to see that $\bar q$ is a cocycle, so that $B^1(X,Q)\subseteq Z^1(X,Q)$.  We may define an equivalence relation on $Z^1(X,Q)$ by $\gamma\sim\delta$ if and only if $\delta(x)=q^{-x}\gamma(x)q$ for some $q$ in $Q$; in general $Z^1(X,Q)$ is simply a pointed set (with a zero). Then we can give meaning to the quotient $Z^1(X,Q)/B^1(X,Q)$ by $Z^1(X,Q)/\sim$ and we denote this pointed set by $H^1(X,Q)$. It is easy to check that there is a one-to-one correspondence of elements $\gamma$ of $Z^1(X,Q)$ to closed complements $Y$ to $Q$ in the semidirect product $XQ$ via $Y=\{x\gamma(x):x\in X\}$. In particular, $X$ is $Q$-conjugate to $Y$ if and only if $\gamma\in B^1(X,Q)$. For the details of this correspondence see \cite[1.5]{LS96}. For the definition of $H^2(X,A)$ where $A$ is abelian, refer to \cite[\S5.6]{Ser02}. (We shall not need the precise definition here.)

One of the key facts about the groups $H^i(X,Q)$ is the following result:

\begin{prop}[See {\cite[Prop 38]{Ser02}}]\label{les}Let $\{1\to A\to B\to C\to 1\}$ be a short exact sequence of (algebraic) $X$-groups such that the image of $A$ is contained in the centre of $B$. Then there exists an exact sequence
\[1\to H^0(X,A)\to H^0(X,B)\to H^0(X,C)\to H^1(X,A)\to H^1(X,B)\to H^1(X,C)\to H^2(X,A),\]
where maps $H^i(X,A)\to H^i(X,B)$ are induced by the inclusion $A\hookrightarrow B$ and maps $H^i(X,B)\to H^i(X,C)$ are induced by the projection $B\to C$.\end{prop}

\subsection{Approximating $H^1(X,Q)$ with the map $\rho$}
Now we describe the main procedure we use to calculate $H^1(X,Q)$.

\begin{defn}\label{rho} Let $X$ be an algebraic group defined over $K$ and let $Q$ a unipotent algebraic $X$-group defined over $K$ such that there exists a filtration $Q=Q_1\geq Q_2\geq\dots$, with each $Q_i/Q_{i+1}$ a $(K)X$-module which is  central in $Q/Q_{i+1}$. Define $\mathfrak Q$ to be the associated graded group, $\mathsf{gr}\ Q_i\cong Q_1/Q_2\oplus Q_2/Q_3\oplus\dots$ endowed with the structure of the direct sum of the $X$-modues $Q_i/Q_{i+1}$. Call $\mathbb V=\bigoplus_i H^1(X,Q_i/Q_{i+1})\cong H^1(X,\mathfrak Q)$ the {\it pre-cocycle space of $X$ with respect to the filtration $Q_i$}. We define a partial map $\rho:\mathbb V\to H^1(X,Q)$ as follows.
 
For an element $\mathbf v=[[\gamma_1],[\gamma_2],\dots]\in \mathbb V$ we will define a cocycle class $[\Gamma_\mathbf v]=\rho(\mathbf v)$ in $H^1(X,Q)$ inductively. Pick representatives $\gamma_i\in [\gamma_i]$ for each $i$. When $i=1$, we define $\Gamma_1:X\to Q/Q_2=Q_1/Q_2$ simply by $\Gamma_1=\gamma_1$ and define $\rho_1$ by $\rho_1(\mathbf v)=[\Gamma_1]$. 

For $i=2$ we let $[\Gamma_2']$ be some pre-image of $[\Gamma_1]$ under the map $H^1(X,Q/Q_3)\to H^1(X,Q/Q_2)$ induced by the quotient map, where one exists. If a preimage does not exist then we declare $\rho$ undefined at $\mathbf v$. Otherwise define $\Gamma_2:X\to Q/Q_3$ by $\Gamma_2(x)=\Gamma_2'(x)\gamma_2(x)$ and define $\rho_2:\mathbb V\to H^1(X,Q)$ by $\rho_2(\mathbf v)=[\Gamma_2]$. As $\gamma_2(x)$ is central in $Q/Q_3$ it is clear that $\Gamma_2$ is a cocycle.

For subsequent $i$, continue in the obvious way, i.e. let $[\Gamma_i']$ be a preimage of $[\Gamma_{i-1}]$ under the map $H^1(X,Q/Q_{i+1})\to H^1(X,Q/Q_i)$ where one exists. If such a preimage does not exist, then $\rho$ is undefined at $\mathbf v$. Otherwise define $\Gamma_{i}:X\to Q/Q_{i+1}$ by $\Gamma_{i}(x)=\Gamma_i'(x)\gamma_{i}(x)$, and set $\rho_{i}(\mathbf v)=[\Gamma_{i}]$. Finally, let $\Gamma_\mathbf v=\lim_{i\to \infty} \Gamma_i$; let $\rho=\lim_{i\to\infty}\rho_i$.
\end{defn}

We note that the construction of $\rho$ above is highly dependent on the choices made, so that for a different set of choices of pre-image $[\Gamma_i']$, one can have elements of $\mathbb V$ mapping under $\rho$ to different elements of $H^1(X,Q)$. However, whichever choices are made, the following proposition still holds:

\begin{prop} \label{surj}The map $\rho:\mathbb V\to H^1(X,Q)$ is surjective.\end{prop}
\begin{proof} This is best described in the language of complements. We show $\rho_i:\mathbb V\to H^1(X,Q/Q_{i+1})$ is surjective for each $i$ by induction.

Note that $\rho_1:\mathbb V\to H^1(X,Q/Q_2)$ is trivially surjective. 

Assume that $\rho_{i-1}:\mathbb V\to H^1(X,Q/Q_{i})$ is surjective. That is, up to $Q$-conjugacy, an exhaustive set of complements to $Q/Q_{i}$ in the semi-direct product $XQ/Q_{i}$ up to $Q/Q_{i}$-conjugacy is given by $X^{i-1}_\bold v=\{x\Gamma_{i-1}(x):x\in X\}$ where $[\Gamma_{i-1}]=\rho_{i-1}(\mathbf v)$ as defined above and $\mathbf v$ ranges over the elements of $\mathbb V$. Let $Y$ be a complement to $Q$ in the semidirect product $XQ$. Then $YQ_i/Q_i$ is a complement to $Q/Q_i$ in the semidirect product $XQ/Q_i$. Replacing $Y$ by a $Q$-conjugate, we have by the inductive hypothesis that $YQ_i/Q_i=X_\bold v^i$ for some $\mathbf v=[[\gamma_1],\dots,[\gamma_{i-1}],0,0,\dots]\in \mathbb V$.

Now observe that the complement $YQ_{i+1}/Q_{i+1}$ to $Q/Q_{i+1}$ in $XQ/Q_{i+1}$ corresponds to some cocycle which is a pre-image of the map $H^1(X,Q/Q_{i+1})\to H^1(X,Q/Q_i)$. Since some preimage exists, $\rho_{i}(\bold v)$ is defined, and equals $[\Gamma_i']$, say. Let $\Gamma_i'$ correspond to a complement $Z$, say to $Q/Q_{i+1}$ in $XQ/Q_{i+1}$; thus $Z=\{x\Gamma_i'(x):x\in X\}$. Since $YQ_{i+1}/Q_{i+1}=ZQ_{i+1}/Q_{i+1}$, $Y$ is a complement to $Q_i/Q_{i+1}$ in $ZQ_i/Q_{i+1}$. Thus $YQ_{i+1}/Q_{i+1}$ corresponds to some cocycle in some class in $H^1(Z,Q_i/Q_{i+1})\cong H^1(X,Q_{i}/Q_{i+1})$ with the isomorphism defined by $\gamma(x\Gamma_i'(x))\mapsto \gamma(x)$. But this cocycle class is some $[\gamma_{i}]\in H^1(X,Q_{i}/Q_{i+1})$ from the construction of $\mathbb V$. Replacing $Y$ by some $Q$-conjugate, we have that $YQ_{i+1}/Q_{i+1}=X^i_\mathbf v=\{x\Gamma_i'\gamma_i(x):x\in X\}=\{x\Gamma_i(x):x\in X\}$ and thus corresponds to some class of cocycles $\rho_{i}(\bold w)$ for some $\mathbf w=[[\gamma_1],\dots,[\gamma_i],0,\dots]\in \mathbb V$.
\end{proof}

\begin{remark} It is possible to prove the above using the standard technique of twisting by cocycles, thus avoiding the need to talk about complements, but we feel it makes the argument easier to understand.\end{remark}

We will of course apply this to the $L$-module filtration of the unipotent radical $Q$ described in \ref{abslem}.

The following definition is quite technical; it gives a name to those modules $V$ in levels of the unipotent radical which have non-trivial values of $H^2(X,V)$, have a module $W$ in an earlier level with non-trivial $H^1(X,W)$ and such that this remains the case for any permissible Frobenius untwist of the action of $X$ on $Q$.

\begin{defn}[Frobenius twists of algebraic $X$-groups]\label{untwists} Let $X$ be an algebraic group defined over $\mathbb F_p$ and $Q$ an algebraic $X$-group where the action is defined over $\mathbb F_p$; i.e. the corresponding homomorphism $X\to\Aut_k(Q)$ is defined over $\mathbb F_p$.

We may define a new action of $X$ on $Q$ via $x\star q=F_r(x)\circ q$ where $F_r$ is the $r$th Frobenius morphism $F_r:X\to X$. Denote the algebraic $X$-group obtained in this way by $Q^{[r]}$, which is isomorphic to $Q$ as an algebraic group over $\mathbb F_p$. We say that $Q^{[r]}$ is the {\it $r$th Frobenius twist of $Q$}.

Recall that when $K=\bar K$, one may equally consider an algebraic group as \begin{itemize}\item an affine variety over $K$ with a group structure whose morphisms are morphisms of varieties, or \item as a $K$-group functor $G$ from $\mathbf{Comm\ K}$-$\mathbf{Alg}$ to $\mathbf{Grp}$, such that $G$ is an affine scheme over $K$ when considered as a $K$-functor.\end{itemize} The details of the second construction can be found in \cite[I.1 \& I.2]{Jan03}. The advantage in this construction is that one can consider Frobenius kernels as normal infinitesimal algebraic groups.

Let $X$ be defined over $\mathbb F_p$ and let $X_r$ denote the kernel of the map $F_r$, the $r$th Frobenius kernel of $X$ considered as a normal infinitesimal subgroup scheme of $X$. Assume $X_r$ acts trivially on $Q$. Then the action of $X$ on $Q$ factors through the quotient $X/X_r$. By \cite[I.9.5]{Jan03} there is an isomorphism $\phi: X/X_r\to X$ such that the composition with the canonical map $X\to X/X_r$ is the Frobenius map $F_r$ on $X$. Therefore there exists an action of $X$ on $Q$ giving $Q$ the structure of an algebraic $X$-group, $P$ say, with $P^{[r]}=Q$. We say $P=Q^{[-r]}$ is the {\it $r$th Frobenius untwist of the algebraic $X$-group $Q$}.\end{defn}

\begin{defn} Let $X$ be an algebraic group and $Q$ a unipotent algebraic $X$-group with the action defined over $\mathbb F_p$ such that there exists a filtration $Q_i$ for $Q$ with each $Q_i/Q_{i+1}$ an $X$-module which is central in $Q/Q_{i+1}$. Call an indecomposable direct summand $M$ of $Q_i/Q_{i+1}$ a {\it potential blocker} if for all $r\geq 0$ such that there exists an untwist $Q^{[-r]}$ of $Q$, we have $H^2(X,M^{[-r]})\neq 0$ and there is some $j<i$ and some indecomposable direct summand $N$ of $Q_j/Q_{j+1}$ such that $H^1(X,N)\cong H^1(X,N^{[-r]})\neq 0$.\end{defn}

\begin{remark} If $X$ is reductive with no factors of type $C_n$ then $H^1(X,N)\cong H^1(X,N^{[1]})$ by \cite[7.1]{CPSV77}. In this special case then, we can can restate the definition as follows:
\begin{quote}Call an indecomposable direct summand $M$ of $Q_i/Q_{i+1}$ a {\it potential blocker} if for all $r\geq 0$ such that there exists an untwist $Q^{[-r]}$ of $Q$, we have $H^2(X,M^{[-r]})\neq 0$ and there is some $j<i$ and some indecomposable summand $N$ of $Q_j/Q_{j+1}$ such that $H^1(X,N)\neq 0$.\end{quote}\end{remark}

\begin{lemma} \label{blockers}Provided there is no potential blocker $W$ in any level $Q_i/Q_{i+1}$ of $Q$ then the maps $H^1(X,Q)\to H^1(X,Q/Q_i)$ are surjective. In other words, all complements to $Q/Q_i$ in $XQ/Q_i$ lift to complements of $Q$ in $XQ$.\end{lemma}
\begin{proof}Assume for a contradiction that $i$ is minimal such that $H^1(X,Q/Q_{i+1})\to H^1(X, Q/Q_i)$ is not surjective. Then $H^1(X,Q/Q_i)\neq 0$ and so some module $V$ in some level $Q_j/Q_{j+1}$ with $j<i$ has $H^1(\bar X,V)\neq 0$. Since $Q_i/Q_{i+1}$ is central in $Q/Q_{i+1}$, we have by \ref{les} an exact sequence 
\[H^1(X,Q/Q_{i+1})\to H^1(X, Q/Q_i)\to H^2(X,Q_i/Q_{i+1}).\]
 As the first map is not surjective, the kernel of the second map is not all of $H^1(X,Q/Q_i)$ and so $H^2(X,Q_i/Q_{i+1})\neq 0$. But this means there is some indecomposable direct summand $W\leq Q_i/Q_{i+1}$ which has $H^2(X,W)\neq 0$. 
 
Now, we have assumed that $W$ is not a potential blocker. Thus there is some untwisted action of $X$ on $Q$, say $R^{[r]}=Q$ with \begin{enumerate}\item $H^1(X,V^{[-r]})\cong H^1(X,V)$ for each indecomposable direct summand $V\leq Q_i/Q_{i+1}$ with $i<j$ and $H^1(X,V)\neq 0$, and \item $H^2(X,W^{[-r]})=0$\end{enumerate}. Now (i) implies that $H^1(X,R/R_i)\cong H^1(X,Q/Q_i)$; but (ii) is true for any indecomposable direct summand $W$ of $Q_i/Q_{i+1}$ satisfying $H^2(X,W)\neq 0$ and so $H^2(X,R_i/R_{i+1})=0$. Thus the map $H^1(X, R)\to H^1(X, R/R_i)$ is surjective by applying the above exact sequence to $R$. But then the map $H^1(X,Q/Q_{i+1})\to H^1(X,Q/Q_i)$ would be surjective too, by twisting cocycles in the latter by the Frobenius morphism $F_r$. This is a contradiction. In particular $H^2(X,W^{[-r]})\neq 0$. So some indecomposable direct summand $W$ of $Q_i/Q_{i+1}$ is a potential blocker, a contradiction.\end{proof}

Even when there are potential blockers in $Q$ we can often modify the construction of $\rho$ to bypass them by spotting that all root groups in shapes $V_S$ admitting non-zero values of $H^1(X,V_S)$ are contained in some subgroup of $G$ without the potential blockers.

\begin{lemma}\label{subsystemlift} Suppose the only shapes $S$ giving rise to $X$-modules $V_S$ in $Q$ such that $H^1(X,V_S)\neq 0$ consist of root groups lying in some (possibly subsystem) subgroup $H\geq X$. Then any complement to $Q$ in $XQ$ is $Q$-conjugate to a subgroup of $H$. Moreover, each map $H^1(X,Q)\to H^1(X,Q/Q_{i+1})$ is surjective if the map $H^1(X,Q\cap H/Q_i\cap H)\to H^1(X,Q\cap H/Q_{i+1}\cap H)$ is surjective.\end{lemma}
\begin{proof}Let $Q_H=Q\cap H$ and ${Q_H}_i=Q_i\cap H$. We will prove that each cocycle class $H^1(X,Q)$ contains a cocycle $\gamma$ such that $\im(\gamma)$ is contained in $Q_H$. To do this we refine the construction of $\rho=\lim_{i_\to \infty}\rho_i$ in \ref{rho} so that $\rho_i(\mathbf v)=[\Gamma_i]$ is such that $\im(\Gamma_i)\subset Q_H/{Q_H}_{i+1}$.

As $H^1(X,Q_i/Q_{i+1})=\bigoplus_{\mathrm{level}(S)=i}H^1(X,V_S)\cong H^1(X,{Q_H}_i/{Q_H}_{i+1})$ we may represent each element $\mathbf v\in \mathbf V$ as $\mathbf v=[[\gamma_1],[\gamma_2],\dots]$ with $\im(\gamma_i)\in {Q_H}_i/{Q_H}_{i+1}$.

Under this choice, we have automatically that $\Gamma_1:X\to Q/Q_2$ defined by $\Gamma_1=\gamma_1$ has image contained in $H$. As in \ref{rho}, set $\rho_1(\mathbf v)=\Gamma_1$.
For $i=2$ let $[\Gamma_2']$ be a preimage under the map $H^1(X,Q_H/{Q_H}_3)\to H^1(X,Q_H/{Q_H}_2)$ if one exists; it is clear that if such a pre-image exists for this map then $[\Gamma_2']$ is also a preimage of $[\Gamma_1]$ under the map $H^1(X.Q/Q_3)\to H^1(X,Q/Q_2)$. Again, as in \ref{rho}, set $\Gamma_2=\Gamma_2'\gamma_2$ and $\rho_2(\mathbf v)=[\Gamma_2]$. Then $\Gamma_2$ has image contained in $H$.

We continue in the obvious way. 

The proof of \ref{surj} goes through with our refinement, showing that each element of $H^1(X,Q)$ has a representative with image contained in $Q\cap H$. This has the desired interpretation in terms of complements.\end{proof}

\begin{corollary}\label{bypassblocks} With the notation of the lemma, if each module ${Q_H}_i/{Q_H}_{i+1}$ for $X$ contains no potential blockers, then the conclusion of Lemma \ref{blockers} holds..\end{corollary}
\begin{proof}As the filtration ${Q_H}_i$ contains no potential blockers, the proof of \ref{blockers} shows that each map $H^1(X,Q_H)\to H^1(X, Q_H/{Q_H}_i)$ is surjective. 

Thus by Lemma \ref{subsystemlift} each map $H^1(X,Q)\to H^1(X,Q/Q_i)$ is surjective.\end{proof}

\subsection{Preliminary results on $\rho$}

In small characteristics, it is rare that $\rho:\mathbb V\to H^1(X,Q)$ is an isomorphism, but the following Lemma, gives a sufficient condition for this to hold.

\begin{lemma}\label{nofixedpoints} Let $X$ be an algebraic group and $Q$ be a connected, unipotent algebraic $X$-group with a closed $X$-stable subgroup $Q_2$, such that $Q_2$ is central in $Q$ and both $Q/Q_2$ and $Q_2$ are modules for $X$. If $H^0(X,Q/Q_2)=0$ and the map $\pi:H^1(X,Q)\to H^1(X,Q/Q_2)$ is surjective then $H^1(X,Q)=H^1(X,Q/Q_2)\oplus H^1(X,Q_2)$.\end{lemma}
\begin{proof}Recall definition \ref{rho}. Since $\pi$ is surjective, the map $\rho:\mathbb V=H^1(X,Q_2)\oplus H^1(X,Q/Q_2)\to H^1(X,Q)$ is defined everywhere and is surjective by \ref{surj}. We wish to prove that it is an injection as well as a surjection.

The image under $\rho$ of the element $\mathbf v=[[\gamma_1],[\gamma_2]]\in \mathbb V$ is $\rho(\mathbf v)=[\Gamma_\mathbf v]=[\Gamma_1'\gamma_2]$ where $[\Gamma_1']$ is a pre-image of $[\gamma_1]$ under the map  $\pi:H^1(X,Q)\to H^1(X,Q/Q_2)$. Let $\mathbf w=[[\delta_1],[\delta_2]]$ and suppose that $\rho(\mathbf w)=\rho(\mathbf v)$, that is, $\Gamma_\mathbf v^q=\Gamma_\mathbf w$ for some $q\in Q$. Let $q'=qQ_2\in Q/Q_2$ so that we get in particular $\gamma_1^{q'}=\delta_1$. Thus $\gamma_1$ and $\delta_1$ are cohomologous. But, following through the definition of \ref{rho} this means that $\gamma_1=\delta_1$ since the construction of $\rho(\mathbf v)$ and $\rho(\mathbf w)$ was made with a unique choice of cocycle $\gamma_1\in[\gamma_1]$. Thus $q'$ acts trivially on $\gamma_1$, and so the coboundary $\bar{q'}$ is trivial, i.e. $\bar{q'}(x):=q'-q'^x=0$ for all $x\in X$, and so $q'=q'^x$ and $q'\in H^0(X,Q/Q_2)$. Since the latter is trivial, $q'=1$ in $Q/Q_2$ and so $q\in Q_2$. But then $\Gamma_\mathbf v^q=\Gamma_1'(\gamma_2^q)=\gamma_1\delta_2$. So $\gamma_2$ is $Q_2$-conjugate to $\delta_2$ and again because $\rho$ is defined with a unique choice of cocycle, we have $\gamma_2=\delta_2$. Thus we have shown that $\mathbf v=\mathbf w$ and so $\rho$ is injective.\end{proof}

It often follows from the above that $H^1(X,Q)\cong K$; we can take advantage of this situation.

\begin{lemma}\label{onedimh1} Let $\bar X\leq L$ be a closed subgroup of a proper parabolic subgroup of $P$ of $G$, with Levi decomposition $P=LQ$. Suppose that $\rho:\mathbb V\to H^1(\bar X,Q)$ is surjective on restriction to a $1$-dimensional subspace $V$ of $\mathbb V$. Then up to $G$-conjugacy, there is at most one non-$G$-cr complement to $Q$ in $\bar XQ$.

If, moreover $H^1(\bar X,Q)$ is non-trivial then up to $G$-conjugacy, there is exactly one non-$G$-cr complement to $Q$ in $\bar XQ$.\end{lemma}
\begin{proof} Since $P$ is proper, $Z(L)$ is at least one-dimensional and so has an action as a $K^\times$ set of $L$-isomorphisms on each level $Q(i)/Q(i+1)$ by $v\mapsto \lambda v$. This action is non-trivial on each root group in $Q$ as $C_G(Z(L))=L$. Thus it pulls back to a non-trivial $\bar X$-equivariant $K^\times$-action on $V\cong K$. Hence, up to $G$-conjugacy, one sees that $\rho$ is surjective on restriction to the two element set $V/K^\times=\{0,1\}$. So we see that there is at most one non-$G$-cr complement to $Q$.

If $H^1(\bar X,Q)\neq 0$ then there is a complement $X$ to $Q$ in $\bar XQ$ which is not $Q$-conjugate to $\bar X$. By \cite[5.9(ii)]{BMRT09}, $X$ is not $G$-conjugate to $\bar X$ and so $X$ is non-$G$-cr.\end{proof}

We want results on restriction of $1$-cohomology. For this we need the notion of twisting by cocycles. 

\subsection{Twisting}\label{twist}
Let $X$ be an algebraic group. Let $R$ be a $X$-set and let $Q$ be $X$-group with an action on $R$ which commutes with the action of $X$. i.e.
\[(r^q)^g=(r^g)^{q^g}\hspace{30pt}\text{for all }r\in R,\ q\in Q,\ g\in X\]

This happens for instance if $R$ and $Q$ are both $X$-groups with $R\leq Q$.

Now fix an arbitrary $1$-cocycle $\gamma\in Z^1(X,Q)$ and define $r * g=r^{g\gamma(g)}$.

This is a new action of $X$  on $R$ as \[r *(gh)=r^{gh\gamma(gh)}=r^{gh\gamma(g)^h\gamma(h)}=r^{g\gamma(g)h
h\gamma(h)}=r *g*h.\]

We call this the $*$-action of $X$ on $R$ with respect to $\gamma$. The set $R$ with the $*$-action is again a $X$-set, denoted $R_\gamma$ and is called a {\it twisted form} of $R$. We say that $R_\gamma$ is obtained by twisting $R$ by $\gamma$.

Observe that if $\rho:Q\to S$ is a homomorphism of $X$-groups with the image of $\gamma$ under $\rho\circ\_$ being $\beta\in Z^1(X,S)$ then we get a map $\rho_\beta:Q_\gamma\to S_\beta$.

In particular, if $R$ is a $X$-stable normal subgroup of $Q$ with $Q/R$ given the natural $X$-action then for any cocycle $\gamma\in Z^1(X,Q)$, we have a well-defined twisted form $(Q/R)_\gamma$ of $Q/R$. (Here we denote the image of $\gamma$ in $Z^1(X,Q/R)$ by $\gamma$, also.)

Now we note further that if $R\leq Z(Q)$ then the $*$-action of $X$ on $R$ coincides with the usual action of $X$ on $R$. So in this case, if $\{1\}\to R \to Q \to S\to \{1\}$ is a short exact sequence of $X$-groups, then it is clear that so is $\{1\} \to R \to Q_\gamma \stackrel{\rho}{\to} S_\gamma \to \{1\}$. (The map $\rho$ obviously commutes with the action of $X$.)

Thus for any $\gamma$, we get from Proposition \ref{les}, a new exact sequence of cohomology:
{\small \[1\to R^X \to Q_\gamma^X \to S_\gamma^X \to  H^1(X,R) \to H^1(X,Q_\gamma) \to H^1(X,S_\gamma)\to H^2(X,R_\gamma),\]}where by the proposition below, $H^1(X,Q)\cong H^1(X,Q_\gamma)$ and $H^1(X,S_\gamma)\cong H^1(X,S)$.

\begin{prop}[{\cite[Prop. 35 bis]{Ser02}}]\label{bijofh1}Let $R$ be a $X$-group and $\gamma\in Z^1(X,R)$. Then the map \[\theta_\gamma:H^1(X,R_\gamma)\to H^1(X,R);\ [\delta]\mapsto[\delta\gamma],\]

where $\delta\gamma$ denotes the map $g\mapsto\gamma(g)\delta(g)$, is a well-defined bijection, taking the trivial class in $H^1(X,R_\gamma)$ to the class of $\gamma$ in $H^1(X,R)$\end{prop}

\subsection{Commutative diagrams}
We now give some commutative diagrams that we shall need in the proof of the main result. Note that if $Y\leq X$ then for each $0\leq i\leq 2$ we get maps $H^i(X,Q)\to H^i(Y,Q)$ by restriction.

\begin{prop}\label{exactseqsrest} Let $Y$ be a subgroup of $X$ acting on a short exact sequence $\{1\}\to R\to Q\to S\to\{1\}$, such that the image of $R$ in $Q$ is central. Then: \begin{enumerate}\item Restriction to $Y$ of the exact sequence of cohomology from Proposition \ref{les} gives rise to the following commutative diagram, where the rows are exact and the vertical arrows are restrictions: 
\begin{align*}\minCDarrowwidth19pt \begin{CD} 1@>>> R^X @>>> Q^X @>>> S^X @>\delta_X>>  H^1(X,R)\\
@VVV @VVV @VVV @VVV @VVV \\
1@>>> R^Y @>>> Q^Y @>>> S^Y @>\delta_Y>>  H^1(Y,R) \end{CD}\hspace{80pt}\\
\hspace{80pt}\begin{CD} @>>> H^1(X,Q)@>>> H^1(X,S)@>>> H^2(X,R)\\
@. @VVV @VVV @VVV \\
@>>> H^1(Y,Q)@>>> H^1(Y,S)@>>> H^2(Y,R).\end{CD}
\end{align*}

\item Let $\gamma\in Z^1(X,Q)$ be a $1$-cocycle, and let $\{1\}\to R\to Q_\gamma\to S_\gamma\to \{1\}$ be the short exact sequence from \S\ref{twist}. If $\beta\in Z^1(Y,Q)$ denotes the restriction of $\gamma$ to $Y$, then we have the following commutative diagram
\[\begin{CD} \{1\}@>>> R @>>> Q_\gamma @>>> S_\gamma @>>> \{1\}\\
@VVV @VVV @VVV @VVV @VVV \\
\{1\}@>>> R @>>> Q_\beta @>>> S_\beta @>>> \{1\}.\end{CD},\] where the vertical arrows are restrictions. Moreover,
\item restriction from $X$ to $Y$ gives rise to the following commutative diagram, where the rows are exact and the vertical arrows are restrictions:
\begin{align*}\minCDarrowwidth19pt \begin{CD} 1@>>> R^X @>>> Q_\gamma^X @>>> S_\gamma^X @>\delta_X>>  H^1(X,R)\\
@VVV @VVV @VVV @VVV @VVV \\
1@>>> R^Y @>>> Q_\beta^Y @>>> S_\beta^Y @>\delta_Y>>  H^1(Y,R) \end{CD}\hspace{80pt}\\
\hspace{80pt}\begin{CD} @>>> H^1(X,Q_\gamma)@>>> H^1(X,S_\gamma)@>>> H^2(X,R)\\
@. @VVV @VVV @VVV \\
@>>> H^1(Y,Q_\beta)@>>> H^1(Y,S_\beta)@>>> H^2(Y,R).\end{CD}
\end{align*}
\end{enumerate}
\end{prop}\label{restrict}
\begin{proof} A moment's thought is required to see that (i) is true. For instance, if $qR\in Q/R\cong S$ with $(qR)^g=qR$ for all $g$ then $(\text{res}^X_Y\circ \delta_X)(qR)(b)=[q^{-\bullet}q](b)=q^{-b}q=\delta_Y(qR)=\delta_Y(\text{res}^X_Y qR)$ for all $b\in Y$.

(ii) is clear.

(iii) follows immediately from (i) and (ii).
\end{proof}

We also need to relate some of maps of (i) and (iii) together.

\begin{prop}[{\cite[p47]{Ser02}}]\label{thetarho}
Let $X$ be a group and $\rho:Q\to S$ any $X$-homomorphism of $X$-groups with $\gamma\in Z^1(X,Q)$ and $\beta$ its image in $Z^1(X,S)$ . Then we have the following commutative diagram:
\begin{center}
$\begin{CD}H^1(X,Q)@>\rho\circ\_ >> H^1(X,S)\\
@A\theta_\gamma AA @A\theta_\beta AA\\
H^1(X,Q_\gamma) @>\rho_\gamma\circ\_ >> H^1(X,S_\beta),\end{CD}$\end{center} where the vertical maps are the bijections of Proposition \ref{bijofh1}.\end{prop}
\begin{proof} Again this is clear: take $[\delta]\in H^1(X,Q_\gamma)$. Then $\theta_\gamma[\delta]=[\delta\gamma]\in H^1(X,Q)$. In turn \[\rho\circ_\_([\delta\gamma])=[(\rho\circ\delta)(\rho\circ\gamma)]=[(\rho\circ\delta)\beta]\] and the latter is clearly equal to $\theta_\beta(\rho_\gamma([\delta]))$. \end{proof}

\begin{prop}\label{thetarest}Let $Y$ be a subgroup of $X$ and $Q$ a $X$-group, with $\gamma\in Z^1(X,Q)$ and $\beta\in Z^1(Y,Q)$ its restriction to $Y$. Then we have the following commutative diagram

\[\begin{CD} H^1(X,Q_\gamma) @>\theta_\gamma>> H^1(X,Q)\\
@VVV @VVV\\
H^1(Y,Q_\beta) @>\theta_\beta>> H^1(Y,Q)\end{CD}\]
where the vertical maps are restriction.\end{prop}
\begin{proof} On elements 
\[\xymatrix@!0{ 
[\delta]\ar@{>}[rrrr] \ar@{>}[dd]& & & &[\delta\gamma]\ar@{>}[dd]\\
\\
[\delta|_Y] \ar@{>}[rrrr] & & & &[\delta|_Y\beta]=[\delta|_Y\gamma|_Y]}\]\end{proof}

Putting together the last three propositions:
\begin{prop}\label{cuboid}With the hypotheses of Proposition \ref{restrict}, we have the following commutative partial cuboid
\[\xymatrix@!0{ 
& & H^1(X,R) \ar@{>}[rrrr]\ar@{>}'[dd][dddd]
& & & & H^1(X,Q) \ar@{>}'[dd][dddd] \ar@{>}[rrrr]
& & & & H^1(X,S) \ar@{>}[dddd]
\\ \\
H^1(X,R)\ar@{>}[rrrr]\ar@{>}[dddd]
& & & & H^1(X,Q_\gamma) \ar@{>}[uurr]\ar@{>}[dddd]\ar@{>}[rrrr] 
& & & & H^1(X,S_\gamma) \ar@{>}[uurr]\ar@{>}[dddd]
\\ \\
& & H^1(Y,R) \ar@{>}'[rr][rrrr] 
& & & & H^1(Y,Q) \ar@{>}'[rr][rrrr] 
& & & & H^1(Y,S)
\\ \\
H^1(Y,R) \ar@{>}[rrrr]
& & & & H^1(Y,Q_\beta) \ar@{>}[rrrr]\ar@{>}[uurr] 
& & & & H^1(Y,S_\beta) \ar@{>}[uurr]
},\]
where rightward arrows are part of four exact sequences running through the central vertical square. 
\end{prop}
\begin{proof}The front and back faces are subdiagrams of Proposition \ref{exactseqsrest}(i),(iii); the top and bottom faces commute by Proposition \ref{thetarho}; the vertical squares commute by Proposition \ref{thetarest}.
\end{proof}

\subsection{Versions of the five lemma}

The following proposition is the five lemma . We will need to alter it so as to use it in our context where we are dealing with non-abelian cohomology and therefore only pointed sets in general. 

\begin{lemma}[The five lemma. {\cite[2.72]{Rot09}}]\label{fivelem} Consider a commutative diagram of abelian groups with 
exact rows. 
\[\begin{CD} A_1@>>> A_2@>>> A_3@>>>A_4@>>>A_5\\
@Vh_1VV @Vh_2VV @Vh_3VV @Vh_4VV @Vh_5VV\\
B_1@>>> B_2@>>> B_3@>>>B_4@>>>B_5\end{CD}\]

\begin{enumerate}\item If $h_2$ and $h_4$ are surjective and $h_5$ is injective, then $h_3$ is surjective. 
\item If $h_2$ and $h_4$ are injective and $h_1$ is surjective, then $h_3$ is injective. 
\item If $h_1$, $h_2$, $h_4$, and $h_5$ are isomorphisms, then $h_3$ is an isomorphism. \end{enumerate}
\end{lemma}

Using the idea of twisting form the previous section, we can get a very similar result in our situation.

\begin{theorem}\label{fivelemmaforh1}Let $\{1\}\to R\to Q\to S \to \{1\}$ be a short exact sequence of $X$-groups, such that the image of $R$ is central in $Q$. Let $Y\leq X$. Then in the following diagram,
\begin{center}
$\begin{CD} S^X @>\delta_X>>  H^1(X,R)@>\iota_X>> H^1(X,Q)@>\pi_X>> H^1(X,S)@>\Delta_X>> H^2(X,R)\\
@V h_1VV @V h_2VV @V h_3VV @V h_4VV @V h_5VV \\
 S^Y @>\delta_Y>>  H^1(Y,R)@>\iota_Y>> H^1(Y,Q)@>\pi_Y>> H^1(Y,S)@>\Delta_Y>> H^2(Y,R).\end{CD}$\end{center}
 the following hold:
\begin{enumerate}\item If $h_2$ and $h_4$ are surjective and $h_5$ is injective, then $h_3$ is surjective. 
\item If $h_2$ and $h_4$ are injective and the restriction maps $S_\gamma^X\to S_\beta^Y$ are surjective for any $\gamma\in Z^1(X,S)$ with $\gamma|^X_Y=\beta$, then $h_3$ is injective.
\item If the hypotheses of (i) and (ii) hold, then $h_3$ is an isomorphism. \end{enumerate}\end{theorem}
\begin{proof} Assume the hypotheses of (i) and take $\gamma\in H^1(Y,Q)$. We need to produce a pre-image $\beta$, say, such that $h_3(\beta)=\gamma$. We get started by the diagram chase used to prove the five lemma for diagrams of abelian groups.

Let $\delta:=\pi_Y(\gamma)$. As the bottom row is exact, $\Delta_Y(\delta)=1$. As $h_4$ is surjective, we have a pre-image $\epsilon$ with $h_4(\epsilon)=\delta$.

Since the last square commutes, $h_5(\Delta_X(\epsilon))=1$ and since $h_5$ is injective, $\Delta_X(\epsilon)=1$.

Now, the top row is exact, so $\epsilon\in\ker \Delta_X$ and hence $\epsilon\in\im\pi_X$; say, $\epsilon=\pi_X(\eta)$, say. Let $h_3(\eta)=\theta$.

Since the penultimate square commutes, we have $\sigma(\theta)=\delta$.

The picture is now as follows:

\begin{center}
$\begin{CD}  * @>\delta_X >>  ?,*@>\iota_X>> ?,\eta@>\pi_X>> ?,\epsilon@>\Delta_X>> 1\\
@V h_1 VV @V h_2 VV @V h_3 VV @V h_4 VV @V h_5 VV \\
 *@>\delta_Y>>  ?,*@>\iota_Y>>{\gamma},{\theta}@>\pi_Y>> \delta@>\Delta_Y>> 1,\end{CD}$\end{center}
 where we would like to establish the $?$s and are not interested in $*$s.
 
In the case of abelian groups, one would continue the proof of the five lemma by taking the difference $\gamma-\theta$; observing that this maps under $\pi_Y$ to $1$ and continuing the diagram chase. Since we cannot do this in the case of pointed sets we use twisting by $\eta$ and continue into the partial cuboid of Propostion \ref{cuboid}.

\[\xymatrix@!0{ 
& & {*}\ar@{>}[rrrr]^{\iota_X}\ar@{>}'[dd][dddd]_{h_2}
& & & & ?,\eta \ar@{>}'[dd][dddd]_{h_3} \ar@{>}[rrrr]^{\pi_X}
& & & & \epsilon \ar@{>}[dddd]_>>>>>>>{h_4}
\\ \\
?_2 \ar@{>}[rrrr]^>>>>>>{\iota_X'}\ar@{>}[dddd]_>>>>>>>{h_2}
& & & & ?_3,1 \ar@{>}[uurr]^{\theta_\eta}\ar@{>}[dddd]_>>>>>>>{h_3'}\ar@{>}[rrrr]^>>>>>>{\pi_X'}
& & & & 1 \ar@{>}[uurr]^{\theta_\epsilon}\ar@{>}[dddd]_>>>>>>>{h_4'}
\\ \\
& & {*}\ar@{>}'[rr]^{\iota_Y}[rrrr] 
& & & & \gamma,\theta \ar@{>}'[rr]^{\pi_Y}[rrrr] 
& & & & \delta
\\ \\
?_1,1\ar@{>}[rrrr]^{\iota_Y'}
& & & & \gamma',1 \ar@{>}[rrrr]^{\pi_Y'}\ar@{>}[uurr]^{\theta_\theta} 
& & & & 1 \ar@{>}[uurr]^{\theta_\delta}
},\]

Here we are using the fact that the bijection in Proposition \ref{bijofh1} takes the neutral element in $H^1(X,Q_\eta)$ to the element $\eta$ in $H^1(X,Q)$ and the fact that the partial cuboid is commutative.

Now in the front bottom row, as $\gamma'\in\ker\pi_Y'$, we have $\gamma'\in\im\ \iota_Y'$. Thus we may put  $\kappa$ in place of $?_1$. Since $h_2$ is a surjection, we may put $\lambda$ in place of $?_2$. Then we replace $?_3$ with $\mu=\iota_X'(\lambda)$ and by the fact that the front left square commutes,  $h_3'(\mu)=\gamma'$. Finally if we replace $?$ with $\nu:=\theta_\eta(\mu)$ the commutativity of the central vertical square gives us our preimage of $\gamma$.

For (ii) the picture in the partial cuboid is as follows:
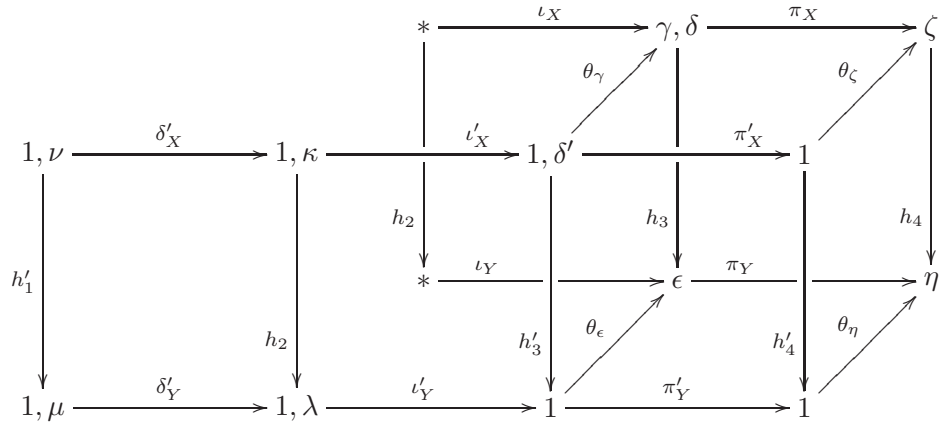
\[\xymatrix@!0{ 
& &
& & & & {*}\ar@{>}[rrrr]^{\iota_X}\ar@{>}'[dd][dddd]_{h_2}
& & & & \gamma,\delta \ar@{>}'[dd][dddd]_{h_3} \ar@{>}[rrrr]^{\pi_X}
& & & & {\zeta} \ar@{>}[dddd]_>>>>>>>{h_4}
\\ \\
1,\nu \ar@{>}[rrrr]^{\delta_X'}\ar@{>}[dddd]_{h_1'}
& & & & 1,\kappa \ar@{>}[rrrr]^>>>>>>{\iota_X'}\ar@{>}[dddd]_>>>>>>>{h_2}
& & & & 1,\delta' \ar@{>}[uurr]^{\theta_\gamma}\ar@{>}[dddd]_>>>>>>>{h_3'}\ar@{>}[rrrr]^>>>>>>{\pi_X'}
& & & & 1 \ar@{>}[uurr]^{\theta_\zeta}\ar@{>}[dddd]_>>>>>>>{h_4'}
\\ \\
& &
& & & & {*}\ar@{>}'[rr]^{\iota_Y}[rrrr] 
& & & & \epsilon \ar@{>}'[rr]^{\pi_Y}[rrrr] 
& & & & \eta
\\ \\
1,\mu \ar@{>}[rrrr]^{\delta_Y'}
& & & & 1,\lambda\ar@{>}[rrrr]^{\iota_Y'}
& & & & 1 \ar@{>}[rrrr]^{\pi_Y'}\ar@{>}[uurr]^{\theta_\epsilon} 
& & & & 1 \ar@{>}[uurr]^{\theta_\eta}
}\]
 where we have twisted by $\gamma$. Here a preimage $\kappa$  of $\delta'$ under $i'_X$ exists since $\delta'$ is in the kernel of $\pi'_X$. Simlarly $\mu$ is a preimage of $\lambda$ under $\delta'_Y$; and $\nu$ is a preimage of $\mu$ under $h_1'$.
 
 This shows that $\delta'=\iota'_X(\delta'_X(\nu))$ and hence is equal to $1$ since the composition of these two maps is trivial. Thus $\gamma=\delta$ since $\theta_\gamma$ is a bijection.
 
(iii) is now obvious.
\end{proof}

A simpler result holds in degree zero.

\begin{prop}\label{fivelemmaforh0}Let $1\to R\to Q\to S\to 1$ be a short exact sequence of $X$-groups, (with $R$ not necessarily central in $Q$) and let $Y$ be a subgroup of $X$. Then in the following diagram:
\[\minCDarrowwidth19pt \begin{CD} 1@>>> R^X @>>> Q^X @>>> S^X @>\delta_X>>  H^1(X,R)\\
@V h_1 VV @V h_2 VV @Vh_3 VV @Vh_4 VV @Vh_5 VV \\
1@>>> R^Y @>>> Q^Y @>>> S^Y @>\delta_Y>>  H^1(Y,R).\end{CD},\]
\begin{enumerate}\item If $h_2$ and $h_4$ are surjective and $h_5$ is injective, then $h_3$ is surjective.
\item If $h_2$ and $h_4$ are injective then $h_3$ is injective.
\item If $h_2$ and $h_4$ are isomorphisms  and $h_5$ is an injection then $h_3$ is an isomorphism.\end{enumerate}
\end{prop}
\begin{proof} The usual proof of the five lemma goes through in this case. Where one would take the `difference' of two elements $g$ and $h$ in an abelian group  one uses the element $gh^{-1}$. The proof then works as normal.\end{proof}

\subsection{Restriction maps on $1$-cohomology}
\begin{prop}\label{h1res} Let $X$ be an algebraic group defined over a field $K$ and let $Y$ be a closed subgroup of $X$. Let $Q$ be a connected, unipotent algebraic $X$-group such that $Q$ admits a filtration $Q=Q_1\geq Q_2\geq\dots$ with $Q_{i}/Q_{i+1}$ central in $Q/Q_{i+1}$ and successive quotients have the structure of rational $X$-modules. If the restriction maps $H^n(X,Q_i/Q_{i+1})\to H^n(Y,Q_i/Q_{i+1})$ are isomorphisms for $n=0,1,2$ and $H^0(X,(Q_i/Q_{i+1})_\gamma)\to H^0(Y,(Q_i/Q_{i+1})_\gamma)$ is an isomorphism for any $\gamma\in Z^1(X,Q)$, then $H^1(X,Q)\cong H^1(Y,Q)$.\end{prop}
\begin{proof} 
We firstly need to ensure the hypotheses of \ref{fivelemmaforh1}(ii) hold.

\begin{claim}For any $\gamma\in H^1(X,Q/Q_{i+1})$, the map $H^0(X,(Q/Q_i)_\gamma)\to H^0(Y,(Q/Q_i)_\gamma)$ is an isomorphism.\end{claim}
{\it Proof of claim:}
If $i=1$ then this is trivial so assume we have proved the proposition up to $n-1$. Now, in \ref{fivelemmaforh0} let $`R'=(Q_i/Q_{i+1})$, $`Q'=(Q/Q_{i+1})_\gamma$ and $`S'=(Q/Q_i)_\gamma$. By \ref{exactseqsrest}(ii) $1\to R\to Q\to S$ is indeed a short exact sequence of $X$-groups with the image of $R$ in the centre of $Q$. By hypothesis, $H^0(X,R)\to H^0(Y,R)$ and $H^1(X,R)\to H^1(X,R)$ are isomorphisms; by induction $H^0(X,S)\to H^0(Y,S)$ is an isomorphism; and so \ref{fivelemmaforh0}(iii) implies that $H^0(X,Q)\to H^0(X,Q)$ is an isomorphism. Thus the claim follows by induction.

Equipped with the claim, a similar induction on $i$ using \ref{fivelemmaforh1}(iii) yields the proposition.
\end{proof}

We apply the above proposition to the situation of reductive algebraic groups.

\begin{corollary}\label{h1resb} Let $X$ be a closed reductive subgroup of the reductive group $G$ contained in a parabolic $P=LQ$ of $G$ and let $B$ be a Borel subgroup of $X$. Suppose that $Q$ is a connected, unipotent algebraic $X$-group such that $Q$ admits a filtration $Q=Q_1\geq Q_2\geq\dots$ with $Q_i/Q_{i+1}$ is central in $Q/Q_{i+1}$ and successive quotients have the structure of rational $X$-modules. Then the restriction map $H^1(X,Q)\to H^1(B,Q)$ is an isomorphism. \end{corollary}
\begin{proof} The statement follows from the previous proposition since we have an isomorphism $H^n(X,V)\to H^n(B,V)$ for any rational $KX$-module by \cite[2.1]{CPSV77}.\end{proof}

\begin{corollary}Let $X$ be a closed reductive subgroup of $G$ contained in a parabolic $P=LQ$ of $G$ and let $B$ be a Borel subgroup of $X$. Then $X$ is $G$-conjugate to a subgroup of $L$ if and only if $B$ is.\end{corollary}
\begin{proof} This follows easily from the above corollary and the fact that $X$ is $Q$-conjugate to a subgroup of $L$ if and only if it is $G$-conjugate to a subgroup of $L$ from \cite[5.9(ii)]{BMRT09}.\end{proof}

The following will be useful as a construction at points in this thesis.

\begin{defn}[Canonical lift]\label{clift} Let $Q$ be the unipotent radical of a parabolic subgroup of a reductive algebraic group $G$. Choose some regular ordering `$\leq$' on the positive roots of $G$; that is, some ordering compatible with the height function on roots. We define a new ordering `$\prec$' on the roots of $Q$. Suppose $\alpha$ is a root in level $i$ and $\beta$ is a root in level $j$. Then $\alpha\prec\beta$ if $i<j$ and $\beta\prec\alpha$ if $j<i$. If $i=j$ then $\alpha\prec\beta$ if and only if $\alpha<\beta$.

Suppose $qQ(j)\in Q(i)/Q(j)$. Write $q$ as a product of elements in the root groups of $Q$ in increasing  '$\leq$'-order and then use the commutator formula to write $qQ(j)=\prod_\alpha x_\alpha(t_\alpha)Q(j)$ where the product is over roots $\alpha$ with $i\leq\mathrm{level}(\alpha)<j$ and $\alpha$ appears earlier in the product than $\beta$ if and only if $\alpha\prec\beta$.  We define an element $\hat q\in Q(i)=\prod_\alpha x_\alpha(t_\alpha)$. We refer to $\hat q$ as the {\it canonical lift} of $q$.\end{defn}

\section{Abelian cohomology results}\label{abelian}

Firstly, we have the following well known result originally due to Cline, see e.g. \cite[2.2]{AJL83}:

\begin{prop}\label{extfora1} Let $X=SL_2(K)$. Let $\lambda, \mu$ be positive integers with $p$-adic expansions $\lambda=\sum a_ip^i$, $\mu=\sum b_ip^i$. Then Ext$^1_X(\lambda,\mu)\neq 0$ if and only if there exists $k\geq v_p(\lambda+1)$ (where $v_p(\lambda+1)=\text{max }\{i:p^i$ divides $\lambda+1\}$) such that \[a_i=b_i\text{ for } i\neq\{k,k+1\}, a_k+b_k=p-2\text{ and } a_{k+1}-b_{k+1}=\pm 1.\] Further Ext$^1_X (\lambda,\mu)$ has dimension at most 1.\end{prop}
When $\lambda=0$ we recover

\begin{corollary}\label{h1a1} Let $X=SL_2(K)$. If $H^1(X,L(\mu))\neq 0$ then $H^1(X,L(\mu))=K$ and $\mu$ is some twist of the irreducible module $(p-2)\otimes 1^{[1]}$ for $X$.\end{corollary}

We come across certain $SL_2$-modules quite regularly and need to calculate their first cohomology. Recall the notation $V_1/V_2/\dots/V_n$ indicating a uniserial module with head $V_1$ and socle $V_n$. Recall also the notation $\frac{U+V}{W}$ indicating a module with socle $W$ and head $U+V$.

\begin{lemma}\label{someh1s}Let char $K=2$. In the category of $SL_2$-modules, for any $r,s>0$, we have\begin{enumerate}
\item $H^1(SL_2,1^{[r]}/0/1^{[r+1]})=0$, where $1^{[r]}/0/1^{[r+1]}={(W(4))^{[r-1]}}^*=H^0(4)^{[r-1]}$
\item $H^1(SL_2,0/1^{[r+1]}/0/1^{[r]})\cong H^1(SL_2,1^{[r]}\otimes(1^{[r]}/0))=0$,
\item $H^1(SL_2,1^{[r]}/0/1^{[r+1]}/0/1^{[r]})=0$, where $1^{[r]}/0/1^{[r+1]}/0/1^{[r]}\cong T(4)^{[r-1]}$,
\item $H^1(SL_2,1^{[r]}\otimes (1^{[s]}/0))=K$ if $r\neq s$.
\item $H^1(SL_2,\begin{array}{c}1^{[r_1]}+\dots+1^{[r_n]}\\\hline 0\end{array})\cong K^n$, where the $r_i$ are all distinct strictly positive integers. (The notation for the module indicates an extension of the trivial module where the socle is just the trivial module itself.)\end{enumerate}\end{lemma}
\begin{proof} The isomorphism in (ii) is from \ref{varrest} (ii) below.

Suppose that the first three items have strictly positive values of $H^1$. Then in each of these three cases one sees that there exists a module $0/1^{[r]}/0/1^{[r+1]}$. Thus there exists a module $1^{[r+1]}/0/1^{[r]}/0$. But this is a contradiction as the image of a high weight vector $\langle v\rangle$ under $G$ generates the image of the Weyl module $W(2^{r+1})$, whose weight spaces are all 1-dimensional, thus cannot contain two trivial factors.

For (iv), first recall that $H^1(SL_2,1^{[r]}\otimes (1^{[s]}/0))\cong \Ext^1_{SL_2}(1^{[r]},1^{[s]}/0)$. We use the long exact sequence of $SL_2$-modules of which part is \[\Hom_{SL_2}(1^{[r]},1^{[s]})=\Ext^0_{SL_2}(1^{[r]},1^{[s]})\to \Ext^1_{SL_2}(1^{[r]},0) \to \Ext^1_{SL_2}(1^{[r]},1^{[s]}/0)\to \Ext^1_{SL_2}(1^{[r]},1^{[s]})\]
Clearly the first term is $0$ as $r\neq s$. From \ref{extfora1} we see that the last term is also $0$. Thus $K\cong \Ext^1_{SL_2}(1^{[r]},0) \cong \Ext^1_{SL_2}(1^{[r]},1^{[s]}/0)$.

Item (v) is similar.
\end{proof}

Sometimes we can spot that a module is tilting, and so the following well-known result is useful.
\begin{lemma} \label{h1tilting} Let $T(\lambda)$ be a tilting module for an algebraic group $X$. Then $H^1(X,T(\lambda))=0$\end{lemma}
\begin{proof} This follows from \cite[II.4.13 (2)]{Jan03} as $T(\lambda)$ has a good filtration.\end{proof}

In order to establish potential blockers when $\bar X$ is of type $A_1$, we use the following theorem from \cite{Ste09}.

\begin{lemma}\label{h2sl2} Let $V=L(r)^{[d]}$ be any Frobenius twist (possibly trivial) of the irreducible $SL_2$-module $L(r)$ with highest weight $r$ where $r$ is one of \begin{align*}&2p\\
&2p^2-2p-2\ (p>2) \\
&2p-2+(2p-2)p^{e}\ (e>1)\end{align*}
Then $H^2(G,V)=K$. For all other irreducible $G$-modules $V$, $H^2(G,V)=0$.
 \end{lemma}
 
If $\bar X$ is non-simple, we can use the following result to calculate cohomology groups.
 
 \begin{lemma}[K\"unneth Formula, See {\cite[10.85]{Rot09}}] \label{kunneth} Let $X=Y\times Z$ be a direct product of linear algebraic groups acting on a module $V=V_Y\otimes V_Z$. Then \[H^n(X,V)\cong \bigoplus_{i+j=n} H^i(Y,V_Y)\otimes H^{n-i}(Z,V_Z).\]\end{lemma}

\begin{corollary}\label{kuncor}Let $X=Y\times Z$ be a direct product of simple algebraic groups acting on an irreducible module $V=V_Y\otimes V_Z$. Then, up to swapping $Y$ and $Z$, if $H^1(X,V)\neq 0$, we have $V_Y=K$ is the trivial module for $Y$ and $H^1(X,V)\cong H^1(Z,V_Z)$.\end{corollary}

\section{More representation theory}\label{morerep}

In finding all $G$-irreducible subgroups one can employ the classification of maximal subgroups due to Seitz and Liebeck-Seitz recursively. Once one has a subgroup $X$ contained in some classical subgroup $M$ of $G$ the following gives a criterion for $X$ to be contained in no parabolic of $M$.

\begin{lemma}[{\cite[p32-33]{LS96}}] \label{girclassical}Let $G$ be a simple algebraic group of classical type, with natural module 
$V = V_G(\lambda_1)$, and let $X$ be a $G$-irreducible subgroup of $G$. 
\begin{enumerate}\item If $G = A_n$, then $X$ acts irreducibly on $V$ 
\item If $G = B_n$, $C_n$, or $D_n$ with $p\neq 2$, then $V\downarrow X = V_1 \perp\dots\perp V_k$ with the $V_i$ all non-degenerate, irreducible, and inequivalent as $X$-modules.
\item If $G = D_n$ and $p = 2$, then $V\downarrow X = V_1\perp\dots\perp V_k$ with the $V_i$ all non-degenerate, $V_2\downarrow X,\dots,V_k\downarrow X$, irreducible and inequivalent, and if $V_1\neq 0$, $X$ acting on $V_1$ as a $B_{m-1}$-irreducible subgroup where $\dim V_1 = 2m$.\end{enumerate}\end{lemma}

In most cases, even when the characteristic of $K$ is small, one can usually distinguish isomorphic, yet non-conjugate reductive subgroups $X\leq G$ by their actions on small (e.g. natural) representations of $G$. Typically such a situation may arise when one considers reductive complements in parabolic subgroups of $G$.

\begin{lemma}\label{comp factors} Let $X$ be a reductive algebraic group, $Q$ any unipotent $X$-group and $Y$ a complement to $X$ in the semidirect product $XQ$. Suppose $V$ is a rational $XQ$ module. Then the composition factors of $V\downarrow X$ correspond to the composition factors of $V\downarrow Y$ under an isomorphism $X\to Y$.\end{lemma}
\begin{proof} Let $V_i$ be a composition factor of $V\downarrow XQ$. By Lie-Kolchin, $Q$ stabilises a vector in $V_i$, $w$ say, and $KXw=KXQw=V_i$. But now if $v=kxw\in V_i$ is any vector for some $kx\in KX$ then $Qv=kQxw=kxQw=kxw=v$ so $Q$ acts trivially on $V_i$. Thus $V_i$ is an irreducible $X$- and $Y$-module with high weights corresponding under $X\stackrel{\sim}{\to} Y$. \end{proof}

The following lemma becomes immensely useful when passing between subgroups of type $A_1^n$ and $A_1^m$ where cohomologically, the situation is identical.
\begin{lemma}[{\cite[1.4]{LS98}}]\label{samesubs}Let $X$ be a linear algebraic group over $K$ and let $Y$ 
be a closed subgroup of $X$. Suppose $V$ is a finite-dimensional rational $KX$-module 
satisfying the following conditions: \begin{enumerate}
\item every $X$-composition factor of $V$ is $Y$-irreducible; 
\item for any $X$-composition factors $M$, $N$ of $V$, the restriction map $\Ext^1_X(M, N)\to \Ext_Y^1(M, N)$ is injective; 
\item for any $X$-composition factors $M,\ N$ of $V$, if $M\downarrow Y\cong N\downarrow Y$ then $M\cong N$ (as $X$-modules).\end{enumerate}
Then $X$ and $Y$ fix exactly the same subspaces of $V$. 
\end{lemma}
\begin{proof}The statement is as \cite[Prop 1.4]{LS98}, but with $Y$ allowed to be any closed subgroup of any linear algebraic group $X$; the proof goes through the same.\end{proof}

\begin{lemma}\label{varrest} In the category of modules for $SL_2$, with $K=\bar K$ and char $K=2$, the following identities hold:
\begin{enumerate} \item If $r\neq s$ and $r,s>0$ then $1^{[r]}\otimes (1^{[s]}/0)\cong (1^{[r]}\otimes 1^{[s]})/1^{[r]}$
\item $2\otimes (2/0)\cong 0/4/0/2$
\item $(0/2)\otimes (2/0)\cong 2/0/4/0/2+0=T(4)+0$
\item If $SL_2 \leq SL_3$ such that the natural module for $SL_3$, $10\downarrow SL_2=0/2$ then $11\downarrow SL_2=2/0/4/0/2$.
\end{enumerate}
\end{lemma}
\begin{proof} 
For (i) observe that there is a submodule $1^{[r]}$ and a quotient $1^{[r]}\otimes 1^{[s]}$, so we must show that these are not direct summands. But $\Hom(1^{[r]},1^{[r]}\otimes (0/1^{[s]}))=\Hom(1^{[r]}\otimes 1^{[r]},0/1^{[s]})=\Hom(0/1^{[r+1]}/0,0/1^{[s]})=0$, thus there is no submodule $1^{[r]}$ in $1^{[r]}\otimes (0/1^{[s]})$ and so there is no quotient $1^{[r]}$ in $1^{[r]}\otimes (1^{[s]}/0)$.

 For (ii), observe that there $2\otimes (2/0)$ contains a $0/4/0$ quotient and a $2$ submodule. The options are then $2+0/4/0$, $0/(4/0+2)$, $0/4/(0+2)$ or $0/4/0/2$. Now $\Hom(0,2\otimes (2/0))=\Hom(2,2/0)=0$ so there is no trivial submodule. This rules out all but the last option.

For (iii) observe that $2/0\otimes 0/2$ must contain a submodule $0/4/0/2$. SInce it is self-dual, it must contain a submodule $2/0/4/0/2$. This is $8$-dimensional. Thus there is a further trivial composition factor, which being a submodule and a quotient, must split off as a direct summand. That the module $2/0/4/0/2$ is isomorphic to $T(4)$ is clear as the exact sequence $0\to 4/0/2 \to 2/0/4/0/2 \to 2/0 \to 0$ gives a Weyl filtration of $T(4)$ and the dual gives a dual Weyl filtration.

For (iv) recall that $10\otimes 10^*=10\otimes 01=11+00$ so $2/0\otimes 0/2=(11+00)\downarrow A_1$. The assertion is now obvious from (iii).\end{proof}

\section{Conjugacy results for complements in parabolics}\label{conjofcomps}

The following lemma will be needed for the subsequent proposition. The statement and proof are taken from part of \cite[Theorem 3.1]{BMR05}

\begin{lemma}\label{bmr} Let $P$ and $Q$ be two parabolic subgroups of a reductive algebraic group $G$ minimal with respect to containing a closed subgroup $H$. Fix Levi subgroups $L_P$ of $P$ and $L_Q$ of $Q$ both containing a maximal torus $T$. Then $L_P=L_Q$.\end{lemma}
\begin{proof} Since $(P\cap Q)R_u(Q)$ is a parabolic subgroup of $G$ contained in $Q$ and $H\leq P\cap Q$, the minimality of $Q$ implies that $Q=(P\cap Q)R_u(Q)$.Thus $P$ contains a Levi subgroup $M_Q$ of $Q$. Similarly, $Q$ contains a Levi subgroup $M_P$ of $P$. We claim that $P\cap Q$  contains a common Levi subgroup of both $P$ and $Q$. Let $L_P$ and $L_Q$ as in the hypotheses.  Now we have a decomposition \cite[2.8.7]{Car93}:
\[P\cap Q=(L_P\cap L_Q)(L_P\cap R_u(Q))(R_u(P)\cap L_Q)(R_u(P)\cap R_u(Q)).\]
Moreover, $R_u(P\cap Q)$ is the product of the last three factors. Since $M_P$ is reductive, $M_P\cap R_u(P\cap Q)$ is trivial. Thus $M_P$ is isomorphic to a subgroup of $L_P\cap L_Q$. Since $\dim (L_P\cap L_Q)\leq \dim L_P=\dim M_P$ we get $L_P\leq L_Q$. Likewise, $L_Q\leq L_P$.
\end{proof}

We use the above lemma to digest $G$-conjugacy amongst complements in parabolics to something more manageable.

\begin{prop} \label{nomoreconj} Let $X$ and $Y$ be closed subgroups of a parabolic subgroup $P$ of a reductive algebraic group $G$ and let $\pi:P\to L$ be the projection to a Levi subgroup $L$ of $P$ containing a maximal torus $T$. Let $\bar X:=\pi X$, $\bar Y:=\pi Y$ and suppose $\bar X=\bar Y$ with $\bar X$ being $L$-ir. Then the following are equivalent:
\begin{enumerate}\item[(A)] $X$ and $Y$ are $G$-conjugate
\item[(B)] $X$ and $Y$ are $QN_G(L)Q$-conjugate.
\item[(C)] $X$ and $Y$ are $QL'(N_G(T)\cap N_G(L))Q$-conjugate.
\item[(D)]$X$ and $Y$ are $Q\left[\left(L'(N_G(T)\cap N_G(L))\right)\cap N_G(\bar X)\right]Q$-conjugate.
\end{enumerate}
\end{prop}

\begin{proof} Clearly (B), (C) and (D) imply (A).

Fix a Borel subgroup $B=TU$ and assume that $X$ and $Y$ are $G$-conjugate. We may assume that $P$ is a standard parabolic $P=BN_JB$ for some $J\subset \Pi$. Suppose $X^g=Y$ for  $g\in G$. Then the Bruhat decomposition gives $g=huwu_w$ with $h\in T,\ u\in U,\ w\in N_G(T),\ u_w\in U_w$ where $U_w$ is the group generated by root groups of $U$ which are sent into $U^-$ under the action of $w$. Now $hu\in B\leq P$. Decompose this as $ql$ for unique $l\in L$, $q\in Q$. Similarly let $u_w=m_2q_2$.

We have $Y=X^g=(X^q)^{lwm_2q_2}$. Let $X'=X^{ql}\leq P$ and $Y'=Y^{{q_2}^{-1}{m_2}^{-1}}\leq P$. Thus we have $X'^{w}=Y'$ and we will show that $w\in N_G(L)$. Now $L^{w}$ is a Levi subgroup of the parabolic $P^w$. Also since $T\leq L$ and $T=T^w$, we have $T\leq L^w$. Now $X'\leq P$ so $Y'=X'^w\leq P^w$, showing that $P$ and $P^w$ are two parabolic subgroups minimal with respect to containing $Y'$, with fixed Levi subgroups $L$ and $L^w$, respectively, both containing a common maximal torus $T$. Thus by \ref{bmr}, $L=L^w$. So $w\in N_G(L)$, indeed $lwm_2\in N_G(L)$. So $g'=lwm_2\in N_G(L)$, so $g\in QN_G(L)Q$. This proves that (A) $\Rightarrow$ (B).

Additionally, as $w\in N_G(L)$, we may write $lwm_2=l'tw$ for some $l't\in L'Z(L)=L$. Now $tw\in N_G(T)\cap N_G(L)$. This proves that (A)$\Rightarrow$(C).

Moreover $\pi X'=\pi Y'=\pi{(X')^{lwm_2}}=\pi (X')^{lwm_2}$, so $g'=lwm_2\in N_G(\pi(X^q))=N_G(\bar X)$.  This proves that (A) $\Rightarrow$ (D).
\end{proof}

We can further refine (D) above provided we add in a further condition.

\begin{theorem}\label{nomoreconj2} Let $X$ and $Y$ be closed, semisimple subgroups of a parabolic subgroup $P=LQ$ such that $\bar X=\bar Y$ are $L$-irreducible and $X$ and $Y$ are complements to $Q$ in the semidirect product $\bar XQ$. Let $H$ be a maximal torus of $\bar X=\bar Y$ and $T\geq H$ a maximal torus of $L$. Suppose that $H\leq X$ and $H\leq Y$. Suppose further that $N_G(L')=L'C_G(L')$.

Then $X$ and $Y$ are $G$-conjugate, if and only if they are $MZ(L)C_Q(H)(C_G(L')\cap N_G(T))C_Q(H)$-conjugate, where $M$ is a transversal of $\bar X$ in $N_{L'}(\bar X)$ normalising $H$.
\end{theorem}
\begin{proof}
By \ref{nomoreconj}, we have $X^g=Y$ where $g\in Q\left[\left(L'(N_G(L)\cap N_G(T))\right)\cap N_G(\bar X)\right]Q$ but by hypothesis this is $Q\left[\left(L'(C_G(L')\cap N_G(T))\right)\cap N_G(\bar X)\right]Q$. So we may write $g=qlnq_2$ with $l\in L'$ and $n\in C_G(L')\cap N_G(T)$, with also $ln\in N_G(\bar X)$. Since $n\in N_G(T)$ we may write it as $n=tzw$ with $t\in L'\cap T$, $z\in Z(L)$ and $w=\prod_\alpha n_\alpha(1)\in C_G(L')$. Also as $t\in L'$ we may absorb $t$ into $l$. So $g=lzwq_3q_2$ where $q_3=q^{lzw}$. We have $lzw\in N_G(\bar X)\cap N_G(L')$ with $z$ and $w$ centralising all of $L'$; thus $l\in N_{L'}(\bar X)$. But $\bar X$ is $L'$-irreducible and so $N_{L'}(\bar X)=\bar X.M$ where we may take $M$ to be some finite set of elements of $L'$ inducing graph automorphisms of $\bar X$ and normalising the maximal torus $H$ of $\bar X$. So we may write $g=xmzwq_3q_2$ with $x\in\bar X$, $m\in M\subset L'$, $z\in Z(L)$ and $w=\prod_\alpha n_\alpha(1)\in C_G(L')\cap N_G(T)$.

Moreover, we have $xmzw\in N_G(\bar X)$ and so $H^{xmzw}$ is a maximal torus of $\bar X$. Thus by adjusting $x$ we may take $xmzw\in N_G(H)$. Now $mz\in N_G(H)$ and $w\in N_G(T)\cap N_G(\bar X)$, so $w\in N_G(\bar X\cap T)=N_G(H)$, too. So we have arranged $H=H^x=H^m=H^z=H^w$.

Now $x\in\bar X$, so $X^x$ is $Q$-conjugate to $X$, say $X^x=X^{q_4}$; thus $X^g=X^h$ where $h=q_4mzwq_3q_2$. We may write $h=mzq_5wq_2$ with $q_5\in Q$: we can move $q_4$ past $mz$ as $q_4^{mz}\in Q$, and $wq_3=q^{lz}w$, with $q^{lz}\in Q$ as well.

Write $q_5$ in canonical form as a product of elements $\prod_r x_r(t)$ with $r\in \Phi^+$. For each $r$, decide whether $w.r$ is a root in (the Lie algebra of) $Q$ or $Q^-$; as $w\in N_G(L)$ it is clear that it must be in one of these. Now use the commutator relations to move all $r$ such that $w.r\in Q^-$ to the left so that $q_5^w=q_-q_+$ with $q_-\in Q^-$ and $q_+\in Q$.

Now we have that $H^g=H^{mzwq_-q_+q_2}=H^{q_-q_+q_2}\leq Y\leq P$; since $q_+q_2\in Q$, $H^{q_-}\leq P$. But as $H$ normalises $Q^-$ we must have $q_-\in N_{Q^-}(H)=C_{Q^-}(H)$. Thus $H^{q_+q_2}=H$. Thus, similarly we must have $q_+q_2\in C_Q(H)$.

Hence we have $h=mzwq_-q_+q_2=mzq_6wq_7$ with $m\in M$, where $M$ is any finite set of elements of $L'$ inducing graph automorphisms of $\bar X$ normalising $H$, $z\in Z(L)$, $q_6, q_7\in C_Q(H)$, where $q_6={q_-}^{w^{-1}}$, $q_7=q_+q_2$, and $w\in N_G(T)\cap C_G(L')$ as required.\end{proof}

\section{Auxiliary results}

Sometimes there will be no way to fully describe a subgroup of type $X=A_1$ (or $A_1^2$) by reference to embeddings in subsystem subgroups. The following lemmas give us concrete expressions for the images of cocycles on generators of $X$. Using the following lemmas in tandem with the construction of $\rho$ allows one to give explicit descriptions of the generators of complements in parabolics of $G$ in terms of root group elements of $G$.

Firstly we give the following lemma, originally due to Steinberg which describes generators and relations for subgroups of type $A_1(K)$.

\begin{lemma}[{\cite[12.1.1 \& Rk. p198]{Car89}}] \label{a1rels}Let $X$ be the group defined by generators $x_+(t)$ and $x_-(t)$ where $t$ runs over all elements $t\in K$ subject to the following relations:
\begin{enumerate}
\item $x_{\pm}(t_1)x_\pm(t_2)=x_\pm(t_1+t_2)$,
\item $h(t)h(u)=h(tu)$ and
\item $n(t)x_+(t_1)n(t)^{-1}=x_{-}(-t^{-2}t_1)$,
\end{enumerate}
for all $t_1,t_2\in K$ and $t,u\in K^\times$ where $n(t)=x_+(t)x_{-}(-t^{-1})x_+(t)$ and $h(t)=n(t)n(-1)$.

Then $X$ is simple of type $A_1$, isomorphic to $SL_2(K)$ or $PGL_2(K)$.\end{lemma}

This description of subgroups of type $A_1(K)$ allows us to give the non-trivial $1$-cocycles of $A_1(K)$ quite explicitly.

\begin{lemma}\label{h1a1p2} Let char $K=2$. Let $X$ be of type $A_1$ and let $W:=V^{[r]}$ be a non-trivial Frobenius twist of the standard module $V$ for $X$. Take a basis $\{e_1,e_2\}$ for $W$ so that the generators $x_+(t)$ and $x_-(t)$ act as 
\begin{align*}x_+(t).e_1&=e_1\\x_+(t).e_2&=t^{2^{r}}e_1+e_2\\x_-(t).e_1&=e_1+t^{2^{r}}e_2\\x_-(t).e_2&=e_2\end{align*}
Then the following defines a set of representatives of equivalence classes of cocycles $\gamma_k$ in $H^1(X,W)\cong K$: \[\gamma_k:X\to V\ ;\ x_+(t)\to kt^{2^{r-1}}e_1\ ,\ x_-(t)\to kt^{2^{r-1}}e_2.\]
\end{lemma}
\begin{proof}
Regard $XV$ as a semidirect product and define $x_{\pm,k}(t)=x_{\pm}(t)\gamma_k(x_\pm(t))$. Then to show that each $\gamma_k$ is a cocycle It suffices to show that the group $\langle x_{+,k}(t),x_{-,k}(t):t\in K\rangle$ isomorphic to $A_1$ since then the map $\gamma_k:X\to Q$ must indeed satisfy the cocycle condition.

We check \ref{a1rels} to see that the relations given hold. For (i), we see that \begin{align*}x_{+,k}(t)x_{+,k}(u)&=x_+(t).kt^{2^{r-1}}e_1.x_+(u).ku^{2^{r-1}}\\
&=x_+(t)x_+(u).kt^{2^{r-1}}e_1.ku^{2^{r-1}}e_1\ &\text{\ \ as $x_+(u)$ commutes with $e_1$}\\
&=x_+(t+u).k(t+u)^{2^{r-1}}\ &\text{\ \ as $x_+(t).x_+(u)=x_+(t+u)$ commutes with $e_1$}\\
&=x_{+,k}(t+u)\end{align*}
as required. The statement for $x_{-,k}$ is similar.

Secondly, we must check that (ii) holds. Let $n_{+,k}(t)=x_{+,k}(t)x_{-,k}(-t^{-1})x_{+,k}(t)$ and let $h_{+,k}(t)=n_{+,k}(t)n_{+,k}(-1)$. In the following, we move elements $x_\pm$ to the left:
\begin{align*}n_{+,k}(t)&=x_{+,k}(t)x_{-,k}(t^{-1})x_{+,k}(t)\\
&=x_+(t).kt^{2^{r-1}}e_1.x_-(t^{-1}).kt^{-2^{r-1}}e_2.x_+(t).kt^{2^{r-1}}e_1,\ \text{\ \ by defn of $x_{\pm,k}$} \\
&=x_{+}(t)x_-(t^{-1})(kt^{2^{r-1}}e_1^{x_-(t^{-1})}).kt^{-2^{r-1}}e_2.x_{+}(t)kt^{2^{r-1}}e_1\\
&=x_{+}(t)x_-(t^{-1}).(kt^{2^{r-1}}e_1+kt^{2^{r-1}}t^{-2^{r}}e_2+kt^{-2^{r-1}}e_2).x_{+}(t)kt^{2^{r-1}}e_1\\
&=x_{+}(t)x_-(t^{-1}).(kt^{2^{r-1}}e_1+(kt^{-2^{r-1}}+kt^{-2^{r-1}})e_2).x_{+}(t)kt^{2^{r-1}}e_1\\
&=x_{+}(t)x_-(t^{-1})x_+(t).(kt^{2^{r-1}}e_1+kt^{2^{r-1}}e_1)\\
&=x_{+}(t)x_-(t^{-1})x_+(t)\\
&=n_+(t)\end{align*}

Hence $h_{+,k}(t)=h_+(t)$ and (ii) holds for the former, since it holds for $h_+(t)$. For (iii) first observe that $n_+$ is an involution and so:
\begin{align*}n_{+,k}(t)x_{+,k}(t_1)n_{+,k}(t)^{-1}&=x_{+,k}(t_1)^{n_+(t)}\\
&=x_{+}(t_1)^{n_+(t)}(kt_1^{2^{r-1}}e_1)^{n_+(t)}\\
&=x_{-}(t_1t^{-2})(kt_1^{2^{r-1}}e_1+kt_1^{2^{r-1}}t^{2^{-r}}e_2))^{x_+(t)}\\
&=x_{-}(t_1t^{-2})(k(t_1t^{-2})^{2^{r-1}}e_2)\\
&=x_{-,k}(t_1t^{-2})\end{align*}
as required.

Lastly we must show that if a cocycle $\gamma_k$ is cohomologous to a cocycle $\gamma_l$, then $k=l$. But this is to say that $\gamma_k$ and $\gamma_l$ differ by a coboundary, $v_x=v-v^x$, say, for some $v=v_1e_1+v_2e_2=:[v_1,v_2]\in V$. Now on $x_+(t)$, this difference is $\gamma_k(x_+(t))-\gamma_l(x_+(t))=(k-l)t^{2^{r-1}}e_1$ but $v-v^{x_+(t)}=[v_1,v_2]-[v_1+v_2.t^{2^r},v_2]=[v_2.t^{2^r},0]$ and these are not equal for all $t\in K$ unless $v_2=0$ and $k=l$.
\end{proof}

\begin{lemma}\label{a1h1p3} Let char $K=3$. Let $X=A_1$ and let $W:=V^{[r]}$ be a non-trivial Frobenius twist of the module $V=1\otimes 1^{[1]}$ for $X$. Take a basis $\{e_1,e_2,e_3,e_4\}$ for $W$ so that a torus $T$ of $X$ has weights $4.3^r,2.3^r,-2.3^r,-4.3^r$ respectively and under this representation, the generators $x_\pm(t)$ map to the following matrices
\[x_+(t)\mapsto\left[
\begin{array}{c c c c}1&t&t^3&t^4\\
&1&0&t\\
&&1&t^3\\
&&&1\end{array}\right];\text{ and }x_+(t)\mapsto\left[
\begin{array}{c c c c}1&&&\\
t&1&&\\
t^3&0&1&\\
t^4&t^3&t&1\end{array}\right].\]
Then the following defines a set of representatives of equivalence classes of cocycles $\gamma_k$ in $H^1(X,W)\cong K$: \begin{align*}\gamma_k:X\to V\ ;\begin{cases} x_+(t)\to kt^{2{p^r}}e_1-kt^{p^r}e_2,\\
 x_-(t)\to -kt^{{p^r}}e_3+kt^{2p^r}e_4.\end{cases}\end{align*}
\end{lemma}
\begin{proof}
This is similar.\end{proof}

\begin{lemma}\label{a1h1resb} Let $X$ be a simple algebraic group of type $A_1$ defined over an algebraically closed field $K$ and let $B$ denote a Borel subgroup of $X$ with $B=HU$, where $U=\{ x(t):t\in K\}$ is the unipotent radical of $B$ and $H=\{h(t):t\in K^\times\}$ is a maximal torus of $X$. Now let $Q$ be a connected, unipotent algebraic $X$-group with a filtration $Q(i)$ such that each $Q(i)/Q(i+1)$ is central in $Q/Q(i+1)$ and has the structure of an $X$-module.

Then each class in $H^1(X,Q)\cong H^1(B,Q)$ (the isomorphism by \ref{h1resb}) has a member $\gamma$ satisfying $\gamma(h(t))=1$ for all $t\in K$.

Moreover, any map $\gamma:B\to Q$ satisfying $\gamma(h(t))=1$ is a cocycle if and only if
\begin{enumerate}\item Considering $XQ$ as a semidirect product, we have $x(t)\gamma(x(t))x(u)\gamma(x(u))=x(t+u)\gamma(x(t+u))$.
\item $\gamma(x(t)h(u))=\gamma(x(t))^{h(u)}$
\item $\gamma(h(u)x(t))=\gamma(x(t))$\end{enumerate}
\end{lemma}
\begin{proof}For the first statement we have by \cite[6.2.5]{Ric82} that $H^1(H,Q)=\{1\}$.  Thus each complement to $Q$ in $HQ$ is $Q$-conjugate to $H$. Hence for every complement to $Q$ in $XQ$ or in $BQ$ there is a $Q$-conjugate containing $H$.

Now let $C$ be a complement to $Q$ in $BQ$ determined by a cocycle $\gamma$ with $\gamma(h(t))=1$ for all $t\in K^\times$. Then $H=\{h(t):t\in K^\times\}$ is a maximal torus of $X$. We have $C=HU'$ where $U'=\{x(t)\gamma(x(t)):t\in K\}$ is the unipotent radical of $C$. Since $C\cong B$ it is easy to see that (i) and (ii) follow since $x(t).x(u)=x(t+u)$ and $x(t)^{h(u)}=x(tu^{-2})$ in $B$.  This proves the 'only if' part.

Now let $\gamma$ be a map satisfying $\gamma(h(t))=1$ and the three conditions given. Thus it follows that the cocycle condition holds on $\gamma$ restricted to $U$ and to $H$. Now, any two elements of $B$ can be written uniquely as $m=x(t)h(u)$ and $n=x(v)h(w)$, so we need $m\gamma(m)n\gamma(n)=mn\gamma(mn)$.

From (ii) we get $x(t)h(u)\gamma(x(t)h(u))=x(t)\gamma(x(t))h(u)$.

Thus \begin{align*}m\gamma(m)n\gamma(n)&=x(t)\gamma(x(t))h(u)x(v)\gamma(x(v))h(w)\\
&=x(t)\gamma(x(t))x(vu^2)\gamma(x(v))^{h(u^{-1})}h(u)h(w)\\
&=x(t+vu^2)\gamma(h(u)x(t+vu^2))h(u)h(w)\\
&=x(t+vu^2)\gamma(x(t+vu^2))h(uw)\\
&=x(t+vu^2)h(uw)\gamma(x(t+vu^2)h(uw))\\
&=mn\gamma(mn).\end{align*}\end{proof}



\chapter{Reductive subgroups of $G_2$}\label{g2}
\section{Preliminaries}
For this chapter let $G=G_2(K)$ be a simple algebraic group of type $G_2$ defined over an algebraically closed field $K$ of characteristic $p>0$. Here we prove Theorem G1, Corollary G2 and Corollary G3 stated in the Introduction \S\ref{intro}.

Recall that $G_2$ has maximal rank subgroups of type $A_1\tilde A_1$ and $A_2$ (also $\tilde A_2$ generated by all short root groups of $G$ when $p=3$). When $p=2$ we define $Z_1$ to be the subgroup of type $A_1$ obtained from the embedding 
\[A_1(K)\to A_1(K)\circ A_1(K)\leq G ;\hspace{20pt} x\mapsto (x,x).\]
Also when $p=2$, we define $Z_2$ to be the subgroup of type $A_1$ obtained from the embedding
\[A_1(K) \to A_2(K)\leq G\]
where $A_1\cong PSL_2(K)$ embeds in $A_2$ by its action on the three-dimensional space Sym$^2 V$ for $V$ the standard module for $SL_2(K)$. It is shown later that these subgroups are contained in the long root parabolic of $G$, that is, $P=\langle B, x_{-r}(t):t\in K\rangle$ where $r$ is the long simple root associated with the choice of Borel subgroup $B$.

Let $\bar L$, (resp. $\tilde L$) denote the standard Levi subgroup of the standard long root (resp. short root) parabolic subgroup of $G$ containing the Borel subgroup $B$. Let $\bar L_0$ (resp $\tilde L_0$) denote the derived subgroup of $\bar L$ (resp $\tilde L$). Observe that $\bar L_0\cong \tilde L_0\cong SL_2(K)$.

For convenience we reprint here the table given in Corollary G3 and, since we will need it in what follows, prove the restrictions of the subgroups listed are indeed as stated. 

\begin{lemma}\label{g2rest}The following table lists the restrictions of the $7$-dimensional Weyl module $V_7$ for $G$ to various subgroups $X$ of $G$.
\begin{center}
\begin{tabular}{|c|c|c|}
\hline
$X$ & $p$ & $V_7\downarrow X$\\
\hline
$A_2$ & any & $10\oplus 01\oplus 0$\\
$\tilde A_2$ & $p=3$ & $11$\\
$A_1 \tilde A_1$ & any & $1\otimes \tilde 1\oplus 0 \otimes \tilde W(2)$\\
$\bar L_0$ & any & $1\oplus 1\oplus 0^3$ \\
$\tilde L_0$ & any &  $1\oplus 1\oplus W(2)$ \\
$Z_1$ & $p=2$ & $T(2)\oplus W(2)$  \\
$Z_2$ & $p=2$ & $W(2)\oplus W(2)^* \oplus 0$ \\
$A_1\hookrightarrow A_1 \tilde A_1; x\mapsto (x^{(p^r)},x^{(p^s)})$ $r\neq s$ & any & $(1^{(p^r)}\otimes1^{(p^s)})\oplus W(2)^{(p^s)}$ \\
$A_1\hookrightarrow A_2$, irred & $p>2$ & $2\oplus 2\oplus 0$\\
$A_1$, max & $p\geq 7$ & $6$\\
\hline
\end{tabular}
\end{center}
\end{lemma}
\begin{proof} The restriction $V_7\downarrow X$ for the maximal $A_1$ when $p\geq 7$ is well known and is listed in \cite[Main Theorem]{Tes88}. 

Consider $G_2$ embedded in $D_4$ as the fixed points of the triality automorphism. We consider the restriction of the natural 8-dimensional module $V_8$ for $D_4$. Recall that $V_8\downarrow G_2=0/V_7$. For $p=2$, $V_7$ becomes reducible and $V_8\downarrow G_2=0/V_6/0$.

Recall that $\bar L_0, \tilde L_0$ are the simple, connected subgroups of the long and short Levi subgroups respectively. We first consider $V_7\downarrow \bar L_0, \tilde L_0$ and $A_1\tilde A_1$. 

We can see that $A_1\tilde{A_1}\leq A_1^4\leq D_4$. It is clear that the $A_1^4$ subsystem in $D_4$ is realised as $A_1\otimes A_1\perp A_1\otimes A_1\cong SO_4\perp SO_4$. Take the long $A_1$ to be the first of the four and the short $\tilde A_1$ to be embedded diagonally in the other three. 

Now it follows that we have $V_8\downarrow \tilde L_0=0^2\otimes 1\perp 1\otimes 1=1\oplus1\perp T(2)$ for $p=2$ and $1\oplus1\perp 2\oplus 0$ for $p>2$.  This gives $V_7\downarrow \tilde L_0=1\oplus 1\perp W(2)$ for $p=2$ and $V_7\downarrow \tilde L_0=1\oplus 1\perp 2$ for $p>2$. 

We also have $V_8\downarrow \bar L_0=1\otimes 0^2 \perp 0^2\otimes 0^2=1\oplus1\perp 0^4$. Hence $V_7\downarrow L_0=1\oplus 1\perp 0^3$. It follows also that $V_7\downarrow A_1\tilde A_1=1\otimes \tilde 1 \oplus 0 \otimes \tilde W(2)$.

Next we establish $V_7\downarrow A_2$. As the $A_2$ is a subsystem subgroup of $G_2$, it is in a subsystem of the $D_4$. It is therefore contained in an $A_3$. We can see easily that $\lambda_1$ for $D_4$ restricts to $A_3$ as $\lambda_1\oplus\lambda_3=\lambda_1\oplus\lambda_1^*$ (see e.g. \cite[13.3.4]{Car89}). Since $A_2$ sits inside $A_3$ such that the natural module for $A_3$ restricts to $A_2$ as $\lambda_1\oplus 0$ we see that $V_7\downarrow A_2=\lambda_1 \oplus \lambda_1^*\oplus 0$. 

Using this we can restrict to the irreducible $A_1\leq A_2$ for $p>2$, and to $Z_2\leq A_2$, when $p=2$. In this case the natural module for $A_2$, $\lambda_1\downarrow A_1=2$ for $p>2$ and $\lambda_1\downarrow Z_2=W(2)$. Hence $V_7\downarrow A_1=2\oplus 2\oplus 0$ and $V_7\downarrow Z_2=W(2)\oplus W(2)^*\oplus 0$.

Now we compute $V_7\downarrow X$ for $X:=A_1\hookrightarrow A_1\tilde A_1$ twisted by $p^r$ on the first factor and $p^s$ on the second. Using the decomposition above, we read off $V_7\downarrow X=1^{(p^r)}\otimes 1^{(p^s)} \oplus 2^{(p^s)}$. For $s=r=0$ when $p=2$, this gives $V_7\downarrow Z_1=T(2)\oplus 2/0$. 

Lastly let $X=\tilde A_2$ ($p=3$). One checks that a base of simple roots $\{\beta_1,\beta_2\}$ for $G$ is expressed in terms of the roots of $D_4$ as $\{\frac{1}{3}(\alpha_1+\alpha_3+\alpha_4),\alpha_2\}$. On these two elements, the weight $\lambda_1$ for $D_4$ has $\lambda_1(\beta_1)=1$ and $\lambda_1(\beta_2)=1$ implying $V_8\downarrow \tilde A_2$ has composition factors $11|00$ so that $V_7\downarrow \tilde A_2=11$.
\end{proof}

\section{Complements in parabolics: proof of Theorem G1}

Let $G=G_2(K)$ with $K$ algebraically closed of characteristic $p$ and let $X\cong A_1(K)$ be a subgroup of $G$ contained in a parabolic subgroup $P=LQ$ of $G$ with image $\bar X$ under the projection from $P$ to $L$  . Then $X$ is a complement to $Q$ in $\bar XQ$, where $\bar X\cong L'$. Recall the notation $\bar L_0$ and $\tilde L_0$ denoting the cases where $L_0$ is a long root $A_1$ and short root $A_1$ respectively. 

\begin{lemma}\label{inlong} If $X$ is non-$G$-cr, then $p=2$ and $X$ is contained in the long root parabolic subgroup of $G$.\end{lemma}
\begin{proof}
We use \ref{abslem} to calculate that the levels in $Q$ and thus determine $\mathbb V$ (see \ref{rho}). Firstly if $P$ is the short root parabolic, then we find that there are two levels with high weights $0$ and $3$ respectively. For $p>3$ they are restricted and thus irreducible. For $p=3$ they are the modules $0$ and $1^{[1]}/1$; for $p=2$ they are $0$ and $1^{[1]}\otimes 1$.

For the long parabolic one calculates that there are three levels with high weights $1$, $0$, and $1$ respectively. These are restricted and irreducible for all characteristics.

Now, the action of $\bar X$ on the levels of $Q$ is some Frobenius twist of the actions given above. Hence, using \ref{h1a1} we see that $\mathbb V=0$ unless $X$ is contained in the long parabolic of $G$. Since $\rho:\mathbb V\to H^1(\bar X,Q)$ is surjective by \ref{surj}, if $\mathbb V=0$ then $H^1(\bar X,Q)=0$. Thus if $X$ is not in the long parabolic of $G$, then it must be $G$-cr as it is $Q$-conjugate to $\bar X\leq L$.\end{proof}

From this point we assume that $p=2$, $X\leq P$ the long root parabolic, a complement to $Q$ in $\bar L_0 Q$ and $X$ is not conjugate to $\bar L_0$ . As $H^1(X,1^{[r]})$ is $1$-dimensional for all $r>0$ we may assume that $r=1$, observing that we can obtain any other complement to $Q$ by applying a Frobenius map to an appropriate complement we get for $r=1$. 

We are now in a position to prove Theorem G1, which we reprint here for convenience.

\begin{TG1}Let $X\cong A_1(K)$ be a subgroup of a parabolic subgroup in $G=G_2(K)$.

If $p>2$ then $X$ is $G$-cr (two classes). 

If $p=2$ then either $X$ is $G$-cr (two classes) or $X$ is conjugate to precisely one of $Z_1$ and $Z_2$.\end{TG1}
\begin{proof}We have that $X$ is in the long parabolic subgroup of $G$. From the remark above we have that $\mathbb V\cong K^2$ and it is easy to check from \ref{h2sl2} that there are no potential blockers.

Thus \ref{blockers} says that the map $\rho:\mathbb V\to H^1(\bar X,Q)$ is everywhere defined and surjective. So an exhaustive set of cocycle representatives in $H^1(\bar X,Q)$ is indexed $\mathbf v=[v_1,v_2]\in\mathbb V$. We will now observe that $\rho$ is surjective on a smaller set than the whole of $\mathbb V$.

The level $Q(2)/Q(3)\cong K$ is given by a single root group $\{x_{\alpha}(t):t\in K\}$, say. Now let $w$ be an element in level 1, i.e.  $w\in Q(1)/Q(2)$ and let $\hat w$ be its canonical lift (cf \ref{clift}). For each  $\lambda\in K$ define a map $c_\lambda:Q(1)/Q(2)\to Q(3)$ by $c_\lambda(w)=[x_\alpha(\lambda),\hat w]$. By the commutator relations, as $\alpha$ is in level $2$ of $Q$, the image of $w$ under $c_\lambda$ must be in $Q(3)$ as claimed. We now claim that $c_\lambda$ is an isomorphism of $\bar X$-modules. Thus we need to show that $c_\lambda$ is $\bar X$-equivariant: $c_\lambda (w^x)=[x_\alpha(\lambda),\hat w^x]=[x_\alpha(\lambda)^xv',\hat w^x]$ for some (conceivably non-trivial) $v'\in Q(3)$. But $v'$ is central in $Q$ so this is $[x_\alpha(\lambda)^x,\hat w^x]=[x_\alpha(\lambda),\hat w]^x=(c_\lambda (w))^x$. Since $c_\lambda$ is also non-zero, it is an isomorphism. By Schur's lemma, $\Hom(W,Q(2))=K$ and so the isomorphisms $c_\lambda$ together with $0$ form a vector space with $(c_\lambda+c_\mu)(w)=c_\lambda (w)+c_\mu (w)$.

Each $c_\lambda$ pulls back to a map on $\mathbb V$ which we may choose to be $\lambda:[k_1,k_2]\mapsto [0, \lambda k_1]$. Since $q.[q,\bar v]=q^v$ for any $q\in Q$ we get another map on $\mathbb V$ by $[k_1,k_2]^{\widehat{\lambda v}}=[k_1,k_2+\lambda k_1]$ induced by $Q$-conjugacy, which is obviously compatible with $\rho$.

Let $X$ be a complement to $Q$ in $\bar XQ$. Then up to $Q$-conjugacy, $X$ corresponds by \ref{exhaust} to an element $[k_1,k_2]\in \mathbb V$. Now if $k_1\neq 0$, then for some choice of $\lambda$, we get $[k_1,k_2]^{x_\alpha(\lambda)}=[k_1,0]$. Thus $\rho:\mathbb V\to H^1(\bar X,Q)$ is surjective on restriction to the set $A:=\{\mathbf v=[v_1,v_2]\in\mathbb V : v_1v_2=0\}$.

Now since $L=\bar L_0Z(L)=A_1T_1$, the torus $T_1$ normalises $Q$. It is easy to check that if $t(\lambda)\in T_1$ and if $x_{\alpha}(u)Q(2)\in Q/Q(2)$ (resp.  $x_{\alpha}(u)\in Q(3)$) then $x_{\alpha}(u)^{t(\lambda)}=x_{\alpha}(u\lambda)$ (resp. $x_{\alpha}(u)^{t(\lambda)}=x_{\alpha}(u\lambda^3)$). Thus $T_1$ induces an action on $\mathbb V$ with $(v_1,v_2)^{t_\lambda}= (\lambda v_1,\lambda^3 v_2)$. Since $\rho$ is surjective on restriction to the set $A$ above, we have $(v_1,v_2)=(0,0)$ or either $v_1=0$ and $(0,v_2)$ which is $T_1$-conjugate to $(0,1)$ or $v_2=0$ with $(v_1,0)$ $T_1$-conjugate to $(1,0)$. Since $\rho(0,0)$ is represented the zero cocycle, in $H^1(\bar X,Q)$ it corresponds to the $Q$-class of complement containing $\bar X\leq L$. Thus there are at most two conjugacy classes of non $G$-cr subgroups of type $A_2$.

Now from Lemma \ref{a2sinparabs} subgroups $Z_1$ and $Z_2$ are in parabolics and not conjugate to subgroups of Levis. Since the restrictions $V_{7}\downarrow \bar X,\ Z_1$ and $Z_2$ given in \ref{g2rest} are distinct, these are all distinct up to $G$-conjugacy.\end{proof}

In conclusion we have established that a complement $X$ to $Q$ in $\bar L_0 Q$ must be conjugate to precisely one of the subgroups $Z_1, Z_2$ or $\bar L_0$. Together with \ref{inlong}, this completes the proof of Theorem G1, and Corollary G2.

\section{Classification of semisimple subgroups of $G_2$: proof of Corollary G3}

Again we remind the reader of the statement of 

\begin{CG3}\label{CG3} Let $X$ denote a closed, connected semisimple subgroup of $G=G_2(K)$. Then up to conjugacy, $(X,p,V_7\downarrow X)$ is precisely one entry in the Table in \ref{g2rest}.\end{CG3}

In the proof of Corollary G3 we need the classification of maximal subgroups of the algebraic group $G=G_2(K)$, from \cite{LS04}.

\begin{lemma} Let M be a maximal closed connected subgroup of $G$. Then $M$ is one of the following:
\begin{enumerate}\item a maximal parabolic subgroup;\item a subsystem subgroup of maximal rank;\item $A_1$ with $p\geq 7$.\end{enumerate}\end{lemma}

{\it Proof of Corollary G3:}

Firstly, a semisimple subgroup in a parabolic of $G_2$ must be of type $A_1$ and we have determined these by Theorem G1. Secondly, the subsystem subgroups of $G_2$ are well known and can be determined using the algorithm of Borel-de Siebenthal. They are $A_2$, $A_1\tilde A_1$ and $\tilde A_2$ ($p=3$) where the $\tilde A_2$ is generated by the short roots of $G_2$.

Subgroups of maximal rank are unique up to conjugacy so to verify Corollary G3 it remains to check that we have listed all subgroups of type $A_1$ in subsystem subgroups in the table in \ref{g2rest}. If $X\cong A_1$ is a subgroup of $A_2$ or $\tilde A_2$ it must be irreducible or else it is in a parabolic; we have listed these in the table in Corollary 3. If $X\leq A_1\tilde A_1$, let the projection to the first (resp. second) factor be an isogeny induced by a Frobenius morphism $x\to x^{(p^r)}$ (resp. $x\to x^{(p^s)}$). We note some identifications amongst these subgroups:

When $p\neq 2$ and $r=s$ (without loss of generality $r=s=0$), $V_7\downarrow X=2\oplus 2\oplus 0$ which is the same as $V_7\downarrow Y$ where $Y:=A_1\hookrightarrow A_2$ where $Y$ acts irreducibly on the natural module for $A_2$. Indeed these are conjugate since $G$ acts transitively on non-singular 1-spaces (see \cite[Thm B]{LSS96}). When $p=2$ we get the subgroup $Z_1$. When $r=s+1$ and $p=3$, we have $V_7\downarrow X$ is a twist of $V_7 \downarrow Y$ where $Y$ is similarly irreducible in $\tilde A_2$, and we have actually $X$ conjugate to $Y$ up to twists: the long word in the Weyl group $w_0$ induces a graph automorphism on $\tilde A_2$ and it is easy to see that we can arrange the embedding $Y\leq \tilde A_2$ such that $Y\leq C_G (w_0)$. Now $C_G(w_0)=A_1\tilde A_1$ as there is only one class of involutions in $G$ when  $p\neq 2$ by \cite[p288]{Iwa70}. The restriction $V_7\downarrow X,Y$ then gives the identification required.

Finally one can see that all other subgroups listed in the table of Corollary G3 are pairwise non-conjugate as the restrictions of $V_7$ in the table are all distinct.

This proves Corollary G3.

\chapter{Reductive subgroups of $F_4$}\label{f4}
For the remainder of this document, we take $G=F_4(K)$. By $X$ we will typically be referring to a reductive subgroup of $G$. When $X$ is contained in a parabolic subgroup of $G$, we may also consider $\bar X$ as the image of $X$ under projection to a Levi subgroup $L$ of $P$.

\section{Preliminaries}\label{aux}
We begin with some miscellany specific to $G=F_4$ which we will have to use repeatedly throughout our analysis.

The following result due to Liebeck and Seitz, \cite{LS04}[Corollary 2(ii),p. 4], classifies the closed, connected, maximal subgroups of $G$.

\begin{lemma}\label{maxs}Let $M$ be a closed, connected, maximal subgroup of $G=F_4(K)$. Then $M$ is one of the following;
\begin{itemize}\item a maximal parabolic subgroup;
\item a subsystem subgroup of maximal rank;
\item a subgroup of type $A_1$, when $p\geq 13$ (one class);
\item a subgroup of type $G_2$, when $p=7$ (one class); or
\item a subgroup of type $A_1G_2$ when $p\neq 2$ (one class).\end{itemize}\end{lemma}

Also due to Liebeck and Seitz, the following result, \cite{LS96}[Thm1] reduces the study of non-$G$-cr subgroups of $G$ to the characteristics $2$ and $3$.
\begin{lemma}\label{p23orxgcr}Let $X$ be a closed, connected, reductive subgroup of $G=F_4(K)$ and let char $K=p>3$. Then $X$ is $G$-cr.\end{lemma}

Thus when $p>3$ the study of $G$-reducible subgroups of $G$, i.e. those contained in parabolics, is reduced to the study of $L$-irreducible subgroups of $G$ for each Levi subgroup $L$ of $G$.

If in \ref{nomoreconj}, we take $P$ to be the $B_3$ parabolic subgroup of $G$, we can be more specific about conjugacy of complements to the unipotent radical $Q$ of $P$ .

\begin{lemma}\label{B3parabconj}Let $X$ and $Y$ be closed subgroups of the $B_3$ parabolic of $G$ such that $\bar X=\bar Y$ and $\bar X$ is $B_3$-ir. Suppose $X$ is not $G$-conjugate to a subgroup of $L$. Then $X$ and $Y$ are $G$-conjugate, if and only if they are $N_L(\bar X)Q$-conjugate.\end{lemma}
\begin{proof} For the $B_3$, parabolic, $L'(N_G(L)\cap N_G(T))=L\rtimes Z_2$ where $Z_2$ is generated by a  representative $\tilde w_0$,  the reflection in the hyperplane perpendicular to the highest short root $\tilde \alpha_0\in\Phi(G)$. Thus by \ref{nomoreconj}(C), if $X^g=Y$ we may write $g=q_1lwq_2$ where $w$ is either $1$ or $\tilde w_0$, and $lw\in N_G(\bar X)$ by \ref{nomoreconj}(D). If $w=1$ we are now done, so we assume otherwise for a contraction, i.e. that $g=q_1l\tilde w_0q_2$ is impossible.

Now $X$ and $Y$ are $QN_L(\bar X)Q$-conjugate whenever $X^{q_1}$ and $Y^{q_2^{-1}}$ are, so we may replace $X$ and $Y$ by $X^{q_1}$ and $Y^{q_2^{-1}}$ with no loss of generality. So $X^g=Y$ where $g=\tilde w_0l$. Since $P^{\tilde w_0}=P^-$, $X^{\tilde w_0}\leq P^-$, where $P\cap P^-=L$ but by hypothesis we have $X^{\tilde w_0}\nleq L$, so that $X^{\tilde w_0}\nleq P$, and hence $X^{\tilde w_0l}\nleq P$. This is a contradiction as $Y\leq P$.
\end{proof}

\begin{corollary}\label{b3imp} With the hypotheses of \ref{B3parabconj}, if $N_L(\bar X)=\bar XZ(L)$, then $X$ and $Y$ are $G$-conjugate if and only if they are $QZ(L)$ conjugate.\end{corollary}
\begin{proof} Let $X^g=Y$ where we may assume $g=lq$ with $l$ in $N_L(\bar X)$ and $q\in Q$ by \ref{B3parabconj}. By hypothesis we have $l=xz$ with $x\in \bar X$ and $z\in Z(L)$. But $X^x=X^{q'}$ for some $q'\in Q$. Thus, $X^{lq}=X^{zq''}$ for appropriate $q''\in Q$. \end{proof}

%
%
%


\begin{lemma}\label{girsubsconj}Let $Z$ be a simple algebraic groups of classical type with irreducible natural module $V$ and let $X$ and $Y$ be closed $Z$-irreducible subgroups of $Z$.

Then $X$ and $Y$ are $GL(V)$-conjugate if and only if they are $Z$-conjugate.\end{lemma}
\begin{proof}As $Z\leq GL(V)$, $Z$-conjugacy implies $GL(V)$-conjugacy so we must check the other direction. Thus assume that $X$ and $Y$ are $GL(V)$-conjugate, say $X^g=Y$.

In case $Z$ is of type $A_n$, the result is trivial as $GL(V)$-conjugacy implies $SL(V)$-conjugacy. Thus we may assume that $Z$ preserves a non-degenerate bilinear form $J$ on $V$ (and is of type $B$, $C$, or $G$).

As $X$ is $G$-ir, we have from \ref{girclassical} that $V\downarrow X = V_1 \perp\dots\perp V_k$ with the $V_i$ all non-degenerate, irreducible, and inequivalent as $X$-modules.

First suppose that $X$ acts irreducibly on $V$. We have $XJX^T=J$ and $YJY^T=J$. Since $X^g=Y$ this implies that $g^{-1}XgJg^TX^Tg^{-T}=J$ and so $XJ^gX^T=J^g$. Thus $X$ preserves both $J$ and $J^g$. As $X$ acts irreducibly on $V$, $X$ preserves a unique non-degenerate bilinear form on $V$ and so $J^g=\lambda J$ by Schur's lemma. Thus $g\in O(V)$ or $Sp(V)$. In the former case we may adjust by a scalar if necessary to retrieve $g\in SO(V)$. Thus $g\in Z$ as required.

Now if $X$ does not act irreducibly on $V$ we have that $V\downarrow X=V_1\perp W$ with $W:=V_2\perp\dots\perp V_k$. As $V\downarrow X\cong V\downarrow Y$, there is some essentially unique element $g\in Z$ with $Y^g$ stabilising $V_1$ (and $W$). We may replace $Y$ with $Y^g$ with the hypotheses intact.

Now we have $X$ and $Y$ in $Cl(V_1)\times Cl(W)\leq Z$ and an element $g=(g_1,g_2)\in GL(V_1)\times GL(W)$ with $X^g=Y$. Projection to the irreducible subspace $V_1$ yields a unique element $g_1'\in CL(V_1)$ (by the above observation) and $g_2'\in CL(W)$ (by induction) with $g'=(g_1',g_2')$ in $Cl(V_1)\times Cl(W)$ and $X^{g'}=Y$.
\end{proof}
  
\begin{lemma}\label{a1sinb3} Let $L'=B_3$ and let $\bar X$ be an $L'$-irreducible subgroup of $L'$ of type $A_1$.

If $p=3$, then $\bar X\leq A_1^2\tilde A_1$ where $(V_2,V_2,V_2)\downarrow \bar X=(1^{[r]},1^{[s]},1^{[t]})$ and if $r=s$ then $r\neq t$.

If $p=2$, then either $\bar X\leq A_1^2\tilde A_1$ where $(V_2,V_2,V_2)\downarrow \bar X=(1^{[r]},1^{[s]},1^{[t]})$ and $r\neq s$; or $\bar X\leq \tilde A_1^3$ where $(V_2,V_2,V_2)\downarrow \bar X=(1^{[r]},1^{[s]},1^{[t]})$ and $r,\ s$ and $t$ are pair-wise distinct.

Let $L'=C_3$, $p=2$ and let $\bar X$ be an $L'$-irreducible subgroup of $L'$ of type $A_1$. Then either $\bar X\leq A_1^3$ where $(V_2,V_2,V_2)\downarrow \bar X=(1^{[r]},1^{[s]},1^{[t]})$ and $r,\ s$ and $t$ are pairwise distinct; or $\bar X\leq \tilde A_1^2A_1$ where $(V_2,V_2,V_2)\downarrow \bar X=(1^{[r]},1^{[s]},1^{[t]})$ and $r\neq s$.\end{lemma}

 \begin{proof}
Let $p=3$. By \ref{girsubsconj}, the conjugacy class of $\bar X$ is determined by the restriction $V_7\downarrow \bar X$. By \ref{girclassical}, the possible embeddings are $V_7\downarrow \bar X=2^{[r]}+ 2^{[t]}+ 0$, with $r\neq t$; or $1^{[r]}\otimes 1^{[s]}+ 2^{[t]}$ with $r\neq s$.

Since $V_7\downarrow A_1^2\tilde A_1=(1,1,0)+(0,0,2)$, we can realise the latter case as $\bar X\leq A_1^2\tilde A_2$ by $(V_2, V_2, V_2)\downarrow \bar X=(1^{[r]},1^{[s]},1^{[t]}$. The former is the same, but with $r=s$ so $V_7\downarrow\bar X=1^{[r]}\otimes 1^{[r]}+2^{[t]}=2^{[r]}+2^{[t]}+0$.

Now let $p=2$. we prove the statement about $C_3$ instead as the result for $B_3$ will follow by using the graph automorphism of $G$. Thus we need $V_6\downarrow \bar X$ to be a sum of non-degenerate, non-isomorphic symplectic $A_1$-modules. Thus $V_6\downarrow \bar X=1^{[r]}+1^{[s]}+1^{[t]}$ with $r,\ s$ and $t$ all distinct or $V_6\downarrow \bar X=1^{[r]}\otimes 1^{[s]}+1^{[t]}$ with $r\neq s$. The former is evidently in an $A_1^3$. For the latter, the projection to the $1^{[t]}$ summand gives a projection of $\bar X$ to a long $A_1$ root subgroup of $C_3$. The centraliser of this is a $C_2$ with the projection of $\bar X$ to this $C_2$ acting on its natural module as $1^{[r]}\otimes 1^{[s]}$, which describes an $A_1$ subgroup of $\tilde A_1^2$. Thus $\bar X\leq \tilde A_1^2A_1$ as claimed.
\end{proof}

\section{$L$-irreducible subgroups of $L'$}\label{lir}
\newcounter{lcount}
\begin{lemma}[$L$-ir subgroups $\bar X$ of $L'$] The following table lists, for each Levi subgroup $L$, all proper $L$-ir subgroups $\bar X$ of $L$ up to $L$-conjugacy, together with any constraints on the characteristic $p$ of $K$ for such subgroups to occur.

(When $p=2$ the graph automorphism $\tau$ of $G$ bijects $B_3$-irreducible subgroups to $C_3$-irreducible subgroups so we only cover those in $B_3$.)
\begin{center}\begin{longtable}{|c|c|c|}
\hline\endlastfoot \hline
$L'$ & $\bar X$ & constraints\usecounter{lcount}\\
\hline
$B_3$ & $A_2$ & $p=3$\\
 & $G_2$ &\\
 & $\tilde A_1B_2$ & $p=2$\\
 & $A_1^2\tilde A_1$ & \\
 & $\tilde A_1^3$ & $p=2$\\
 & $A_1A_1\leq A_1^2\tilde A_1$ (1)& first factor in $A_1$, second in $A_1\tilde A_1$\\
 & $A_1A_1\leq A_1^2\tilde A_1$ (2) & first factor in $A_1^2$, second in $\tilde A_1$\\
 & $A_1A_1\leq \tilde A_1^3$ & \begin{tabular}{c}$p=2$, first factor in $\tilde A_1$\\ second with different twists in $\tilde A_1^2$\end{tabular}\\
 & $A_1\leq A_1^2\tilde A_1$ & different twists in first two factors\\
 & $A_1\leq \tilde A_1^3$ & $p=2$, pairwise distinct twists\\
 & $A_1$ & $p\geq 7$, irreducible\\
$C_3$ & $A_1C_2$ & $p>2$\\
 & $A_1^3$ & $p>2$\\
 & $A_1^2\leq A_1^3$ & \begin{tabular}{c}$p>2$, first factor in $A_1$,\\ second with different twists in $A_1^2$\end{tabular} \\
 & $A_1\leq A_1^3$ & $p>2$, different twists on each factor\\
 & $A_1^2$ by $V_6\downarrow X=(2,1)$ & $p>2$\\
 & $A_1\leq A_1^2$ by $V_6\downarrow X=2^{[r]}\otimes 1^{[s]}$ & $p>2$, precisely one of $r$ and $s$ non-zero\\
 & $A_1^2\leq A_1C_2$ & $p\geq 5$, second factor irreducible in $C_2$ \\
 & $A_1\leq A_1C_2$ & $p\geq 5$,  proj to second factor irreducible\\
 & $\bar X\leq B_3$ image under $\tau$ & $p=2$\\
$A_1\tilde A_2$ & $A_1A_1$ & $p>2$, second $A_1$ irreducible in $\tilde A_2$\\
 & $A_1\leq A_1A_1$ & $p>2$\\
$A_2\tilde A_1$ & $A_1A_1$ & $p>2$, first $A_1$ irreducible in $A_2$\\
 & $A_1\leq A_1A_1$ & $p>2$\\
$B_2$ & $A_1^2$ & \\
 & $\tilde A_1^2$ & $p=2$\\
 & $A_1\leq A_1^2$ & different twists\\
 & $A_1\leq \tilde A_1^2$ & $p=2$, different twists\\
 & $A_1$ & $p\geq 5$, irreducible\\
$A_2$ & $A_1$ & $p>2$, irreducible\\ 
$A_1\tilde A_1\ (\times 3)$& $A_1$& any twists\\
$A_1$&&\\\hline\end{longtable}\end{center}\end{lemma}
\begin{proof} As all Levi subgroups of $G$ are classical, with natural representation $V_{nat}$ this is a routine check of semisimple groups $X$ of rank less than or equal to that of $L'$ having a $\dim V_{nat}$-dimensional representation which is the direct sum of non-degenerate spaces, by \ref{girclassical} and \ref{girsubsconj}.\end{proof}

\section{Cohomology of unipotent radicals}\label{cur}
Recall that we are looking for a closed, connected, semisimple subgroup $X\leq P=LQ$ with $\bar X=\pi X$, where $\pi$ is the canonical projection $P\to L$. Taking $X$ in a parabolic $P$ which is minimal subject to containing $X$, we have $\bar X$ an $L$-ir subgroup. Thus $\bar X$ is one subgroup of the above list.

In classifying all $G$-reducible subgroups of $G$, we need to calculate $H^1(\bar X,Q)$. The next lemma gives the actions of each $\bar X$ in the table above on the various $V_S$ in Lemma \ref{abslem}.

Thus for each $L$-irreducible subgroup $\bar X$ contained in a parabolic $P=LQ$, in this part we calculate the restriction of each level $Q(i)/Q(i+1)$ described in \ref{abslem} and calculate the dimension of $H^1(\bar X,Q(i)/Q(i+1))$ for each level. This gives us our first approximation to $H^1(\bar X,Q)$.

\begin{lemma}[Restrictions of shapes to $\bar X$]\label{levelsforx} The following table lists the restrictions of the levels  $Q(i)/Q({i+1})$ to the subgroups $\bar X$ in the table above for $p=2$ or $3$ where $\mathbb V>0$; i.e where some level has $H^1(\bar X,Q(i)/Q(i+1))>0$. We embolden all indecomposable direct summands $M$ such that $H^1(G,M)\neq 0$, possibly under certain conditions described below. According with these conditions we have given the non-zero dimensions of each level $i$, i.e. $\dim H^1(\bar X,Q(i)/Q(i+1))$.

We also put asterisks on all potential blockers.\\
{\small
\begin{longtable}{|c|c|c|c|c|}
\hline\hline\endlastfoot
\hline
$\bar X$ & $p$ & Level $i$ & $V$ & $\dim H^1(\bar X,Q(i)/Q(i+1))$\\
\hline\hline
$\bar X$ for $L'=B_3$&&&&\\
\hline
$B_3$ & $p=2$ & $1$ & $001$&\\
  &  & $2$ & $\mathbf{100/000}$  &$1$\\
$A_2$ & $p=3$ &  $1$ & $\mathbf{11}+00$&$1$\\
 & &$2$ & $\mathbf{11}$&$1$\\
$G_2$ & $p=2$ & $1$ & $00/10/00$&\\
&&$2$ & $\mathbf{10/00}$&$1$\\
$\tilde A_1B_2$ & $p=2$ & $1$ & $(1,01)$&\\
&& $2$ &$\begin{array}{c} \mathbf{(1^{[1]},0)+(0,10)}\\\hline \mathbf{0}\end{array}$&$2$\\
$A_1^2\tilde A_1$ & $p=2$ & $1$ & $(1,0,1)+(0,1,1)$&\\
& & $2$ & $(1,1,0)+\mathbf{(0,0,1^{[1]})/0}$&$1$\\
$\tilde A_1^3$ & $p=2$ & $1$ & $(1,1,1)$&\\
&&$2$&$\mathbf{\frac{(1^{[1]},0,0)+(0,1^{[1]},0)+(0,0,1^{[1]})}{0}}$&$3$\\
$A_1^2\leq A_1^2\tilde A_1$ (1) & $p=3$ & $1$ & $(1^{[r]},1^{[t]})+(1^{[s]},1^{[t]})$&\\
&&$2$&$\mathbf{(1^{[r]}\otimes 1^{[s]},0)}+(0,2^{[t]})$&$\leq 1$\\
 & $p=2$ & $1$ & $(1^{[r]},1^{[t]})+(1^{[s]},1^{[t]})$&\\
& & $2$ & $(1^{[r]}\otimes 1^{[s]},0)+\mathbf{(0,1^{[t+1]})/(0,0)}$&$1$\\
$A_1^2\leq A_1^2\tilde A_1$ (2) & $p=3$ & $1$ & $(1^{[r]},1^{[t]})+\mathbf{(0,1^{[s]}\otimes 1^{[t]})}$&$\leq 1$\\
&  & $2$ & $(1^{[r]},0)+(0,1^{[s]})+(0,2^{[t]})^\mathbf{*}$&\\
& $p=2$ & $1$ & $(1^{[r]},1^{[t]})+(0,1^{[s]}\otimes 1^{[t]})$&\\
& & $2$ & $\mathbf{(1^{[r]},0)+(0,1^{[s]})}+\mathbf{(0,1^{[t+1]})/(0,0)}$&$3$\\
$A_1^2\leq \tilde A_1^3$ &$p=2$ & $1$ & $(1^{[r]}\otimes 1^{[s]},1^{[t]})$&\\
&& $2$ & $\mathbf{\frac{(1^{[r+1]},0)+(1^{[s+1]},0)+(0,1^{[t+1]})}{(0,0)}}$&$3$\\
$A_1\leq A_1^2\tilde A_1$ & $p=3$ & $1$ & $\mathbf{1^{[r]}\otimes 1^{[t]}+1^{[s]}\otimes 1^{[t]}}$&$\leq 2$\\
&  & $2$ & $\mathbf{1^{[r[}\otimes 1^{[s]}}+{2^{[t]}}^\mathbf{*}$&$\leq 1$\\
&$p=2$ & $1$ & $1^{[r]}\otimes 1^{[t]}+1^{[s]}\otimes 1^{[t]}$&\\
&  & $2$ & $1^{[r]}\otimes 1^{[s]}+\mathbf{1^{[t+1]}/0}$&$1$\\
$A_1\leq \tilde A_1^3$ &$p=2$ &  $1$ & $1^{[r]}\otimes 1^{[s]}\otimes 1^{[t]}$&\\
&&$2$&$\mathbf{\frac{1^{[r+1]}+1^{[s+1]}+1^{[t+1]}}{0}}$&$3$\\
\hline
\hline
$\bar X$ for $L'=C_3$ & $p=2$ &- & from image under $\tau$ of $\bar X\leq B_3$ &\\
\hline\hline
$\bar X$ for $L'=A_1\tilde A_2$&&&&\\
\hline
$A_1\tilde A_2$ & $p=2$ & $1$ & $(1,10^{[1]})/(1,01)$&\\
& & $2$ & $(0,01^{[1]})/(0,10)$&\\
& & $3$ & $\mathbf{(1,00)}$&$1$\\
\hline\hline
$\bar X$ for $L'=A_2\tilde A_1$&&&&\\
\hline
$A_2\tilde A_1$ & $p=2$ & $1$ & $(01,1)$&\\
& & $2$ & $\mathbf{(00,1)}$ &$1$\\
&& $3$ & $(10,1^{[1]})/(10,0)$ &\\
&&$4$ & $(01,1)$&\\
\hline\hline
$\bar X$ for $L'=B_2$&&&&\\
\hline
$B_2$ & $p=2$ & $1$ & $\mathbf{10/00+01}$ &$2$\\
&& $2$ & $00+\mathbf{01}$&$1$\\
&& $3$ & $\mathbf{10/00}$&$1$\\
&& $4$ & $00$&\\
$A_1^2$ &$p=2$ & $1$ &$(1,1)+(0,0)+\mathbf{(1,0)+(0,1)}$&$2$\\
&&$2$ & $(0,0)+\mathbf{(1,0)+(0,1)}$&$2$\\
&& $3$ & $(1,1)^\mathbf{*}+(0,0)$&\\
&& $4$ & $(0,0)$&\\
$\tilde A_1^2$ & $p=2$ & - & images under $\tau$&\\
$A_1\leq A_1^2$ & $p=3$ & $1$ & $\mathbf{1^{[r]}\otimes 1^{[s]}}+0+1^{[r]}+1^{[s]}$&$\leq 1$\\
& & $2$ & $0+1^{[r]}+1^{[s]}$&\\
& & $3$ & $\mathbf{1^{[r]}\otimes 1^{[s]}}+0$&$\leq 1$\\
& & $4$ & $0$&\\
 & $p=2$ & $1$ & $1^{[r]}\otimes 1^{[s]}+0+\mathbf{1^{[r]}+1^{[s]}}$&$2$\\
& & $2$ & $0+\mathbf{1^{[r]}+1^{[s]}}$&$2$\\
& & $3$ & ${1^{[r]}\otimes 1^{[s]}}^\mathbf{*}+0$&\\
& & $4$ & $0$&\\
$A_1\leq \tilde A_1^2$ & $p=2$ & - & images under $\tau$&\\
\hline
\hline
$L'=A_2,\ \tilde A_2$ & any & - & no non-zero $H^1$s&\\
\hline
\hline
$\bar X$ for $L'=A_1\tilde A_1\leq P_{24}$&&&&\\
\hline
$A_1\tilde A_1$  & $p=2$ & $1$ & $\mathbf{(0,1)}+\mathbf{(1,1^{[1]})/(1,0)}$&$2$\\
&& $2$ & $(1,1)^\mathbf{*}+(0,0)$ &\\
&& $3$ & $\mathbf{(1,0)+(0,1)}$ &$2$\\
&& $4$ & $\mathbf{(0,1^{[1]})/(0,0)}^\mathbf{*}$&$1$\\
&& $5$ & $\mathbf{(1,0)}$ &$1$\\
$A_1$ & $p=3$ & $1$ & $1^{[r]}+2^{[s]}$&\\
& & $2$ & $\mathbf{1^{[r]}\otimes 1^{[s]}}+0$&$1$\\
&& $3$ & $1^{[r]}+1^{[s]}$ &\\
& & $4$ & ${2^{[s]}}^\mathbf{*}$&\\
&& $5$ & $1^{[r]}$ &\\
& $p=2$ & $1$ & $\mathbf{1^{[s]}}+\mathbf{(1^{[r]}\otimes 1^{[s+1]})/1^{[r]}}$&$2$\\
&& $2$ & ${1^{[r]}\otimes 1^{[s]}}^\mathbf{*}+0$&\\
&& $3$ & $\mathbf{1^{[r]}+1^{[s]}}$ &$2$\\
&& $4$ & $\mathbf{1^{[s+1]}/0}^\mathbf{*}$&$1$\\
&& $5$ & $\mathbf{1^{[r]}}$&$1$\\
\hline
\hline
$\bar X$ for $L'=A_1\tilde A_1\leq P_{13}$&&&&\\
\hline
$A_1\tilde A_1$ & $p=2$ & - & images under $\tau$&\\
$A_1$ & $p=3$ & $1$ & $1^{[r]}+\mathbf{1^{[r]}\otimes 1^{[s]}}$&$\leq 1$\\
& & $2$ & $1^{[s]}+{2^{[s]}}^\mathbf{*}$&\\
& & $3$ & $1^{[r]}\otimes 2^{[s]}$&\\
& & $4$ & $1^{[s]}$&\\
& & $5$ & $1^{[r]}$&\\
& & $6$ & $0$&\\
\hline
\hline
$\bar X$ for $L'=A_1\tilde A_1\leq P_{23}$&&&&\\
\hline
$A_1\tilde A_1$ & $p=2$ & $1$ & $\mathbf{(0,1)+(1,0)}$ &\\
&& $2$ & $(1,1)^\mathbf{*}$&\\
&& $3$ & $\mathbf{(1,1^{[1]})/(1,0)}^\mathbf{*}$&$1$\\
&& $4$ & $\mathbf{(0,1^{[1]})/(0,0)}^\mathbf{*}$&$1$\\
&& $5$ & $\mathbf{(0,1)}$ &$1$\\
&& $6$ & $(0,0)$&\\
&& $7$ & $\mathbf{(1,0)}$&$1$\\
$A_1$ & $p=3$ & $1$ & $1^{[r]}+1^{[s]}$&\\
& & $2$ & $\mathbf{1^{[r]}\otimes 1^{[s]}}+2^{[s]}$&$\leq 1$\\
& & $3$ & $1^{[r]}\otimes 2^{[s]}$&\\
& & $4$ & ${2^{[s]}}^\mathbf{*}$&\\
& & $5$ & $1^{[s]}$&\\
& & $6$ & $0$&\\
& & $7$ & $1^{[r]}$&\\
& $p=2$ & $1$ & $\mathbf{1^{[r]}+1^{[s]}}$&$2$\\
&& $2$ & ${1^{[r]}\otimes 1^{[s]}}^\mathbf{*}$&\\
&& $3$ & $\mathbf{(1^{[r]}\otimes 1^{[s+1]}/1^{[r]}})^\mathbf{*}$&$1$\\
&& $4$ & $\mathbf{1^{[s+1]}/0}^\mathbf{*}$&$1$\\
&& $5$ & $\mathbf{1^{[s]}}$&$1$\\
&& $6$ & $0$&\\
&& $7$ & $\mathbf{1^{[r]}}$&$1$\\
\hline
\hline
$L'=A_1\leq P_{234}$  & $p=2$ & $1$ & $0^2+\mathbf{1}$&$1$\\
&& $2$ & $0+\mathbf{1}$&$1$\\
&& $3$ & $\mathbf{1}+\mathbf{1}$&$2$\\
&& $4$ & $\mathbf{1}+0$&$1$\\
&& $5$ & $\mathbf{1}+0$&$1$\\
&& $6$ & $0^2$&\\
&& $7$ & $0$&\\
&& $8$ & $0$&\\
&& $9$ & $\mathbf{1}$&$1$\\
\hline
\hline
$L'=A_1\leq P_{134}$ & $p=2$ & $1$ & $0+\mathbf{1+1}$&$2$\\
&& $2$ & $\mathbf{1}+0^2$&$1$\\
&& $3$ & $0^2+\mathbf{1}$&$1$\\
&& $4$ & $0+\mathbf{1}$&$1$\\
&& $5$ & $0+\mathbf{1}$&$1$\\
&& $6$ & $0$&\\
&& $7$ & $\mathbf{1}$ &$1$\\
&& $8$ & $0$&\\
\hline
\hline
$L'=\tilde A_1\leq$ $P_{123},\ P_{124}$ & $p=2$ & - & images under $\tau$&
\end{longtable}}
Note: In the above, when $p=2$ the modules $V=1$, $10$ and $100$ for $X=A_1$, $C_2$, or $C_3$ respectively have $H^1(X,V)=0$ and $H^1(X,V^{[i]})=K$ for $i>0$. When $p=3$ a module $1^{[r]}\otimes 1^{[s]}$ has positive dimensional $H^1$ if and only if $r=s\pm 1$ by \ref{h1a1}.
\end{lemma}
\begin{proof} The calculation of levels for $\bar X=L'$ is routine using \ref{abslem}. For instance, take $L'=A_2$. Then the level 1 shapes are $**10$ and $**01$. Since the only root of the form $**01$ is $0001$ we conclude that this is a shape affording a one-dimensional, and therefore trivial representation for $\bar X$. The roots in the shape $**10$ are $\{0010,0110,1110\}$ and so the shape represents the module $10$ for $\bar X$. Where we have stated that a level is an indecomposable module for $\bar X$ we have checked the structure of this by using the commutator relations. For instance, the root $1232$ is in level $2$ when $\bar X=B_3$. The commutator relations yield that the corresponding root group is a trivial submodule for $B_3$ under conjugation. It is also easy to check that it cannot be a direct summand.

We show also how to calculate the restrictions to non-Levi subgroups: For example take $G_2\leq L'=B_3$. We have $W_{B_3}(100)\downarrow G_2=W(10)=10/00$. To calculate $L_{B_3}(01) \downarrow G_2$ observe that $G_2\leq D_4$ as the fixed points of the triality automorphism. Thus if we know the natural module $1000$ for $D_4$ restricted to $G_2$ then we know the restriction of the spin modules $0010$ and $0001$ as well, since the restrictions to $G_2$ will be isomorphic. For the first of these, observe that $1000\downarrow B_3=100+000$ when $p>2$ and $000/100/000$ when $p=2$ (as $W_{B_3}(100)$ is a submodule and the module is self-dual). Then we see that $1000$, $0010$ and $0001$ restricted to $G_2$ are $10+00$ for $p>2$ and $00/10/00$ when $p=2$. Lastly the $D_4$ modules, $0010\downarrow B_3=0001\downarrow B_3=001$ for $B_3$.

The $B_3$-ir subgroup of type $A_2\leq B_3$ when $p=3$ is in fact the short $A_2$ maximal rank subgroup $\tilde A_2\leq G_2\leq B_3$ (see \cite{Sei87}). Thus we know that $L_{G_2}(10)\downarrow A_2=11$ and so the restrictions are as given.

When $p=3$, and $\bar X\leq L'=C_3$ we remark that both levels have trivial $H^1$: one sees this directly for  $\bar X=C_3$ as $Q(2)\downarrow C_3\cong 000$ and $Q(1)/Q(2)\downarrow C_3\cong 001$ are irreducible Weyl modules. One similarly calculates this directly for $\bar X=A_1C_2$. Also, any $C_3$-irreducible subgroup of type $A_1$ has odd weights on the $6$-dimensional natural module for $C_3$ and hence has odd weights on the module $\bigwedge^3 100=001+100$. Thus it has odd weights on $001$. Hence there can be no composition factor $1\otimes 1^{[1]}$; the statement follows for products of $A_1$s by \ref{kunneth}. 

To see that there are no potential blockers when $X$ has a factor of rank $2$ or more, recall that we need some level $Q(i)/Q(i+1)$ in $Q$ to have a non-trivial value of $H^2$. For the values of $H^2$ we can refer to \cite{Ste09-1} for $A_2$ and for the remainder, we can use the dimension-shifting identity $H^i(\bar X, L(\lambda))=H^{i-1}(\bar X, H^0(\lambda)/\soc(H^0(\lambda)))$ (see \cite[II.4.14]{Jan03}). For example, $H^2(B_3,100)=H^1(B_3,000)=0$, so $100/000$ is not a potential blocker since neither composition factor admits a non-trivial $H^2$.

Our indication of the values of $H^1(\bar X,V_i)$  when $\bar X=A_1$ follows from \ref{h1a1} and \ref{someh1s}. Most of the potential blockers which occur when $\bar X$ is of type $A_1$, $A_1^2$ or $A_1^3$ can be calculated from \ref{h2sl2} using the K\"unneth formula, \ref{kunneth} to find non-trivial $H^2$s which have a non-trivial $H^1$ in a lower level of $Q$. Then one needs to observe that one cannot untwist the action of $\bar X$ leaving the $H^1$s in tact. For instance, consider the case $A_1\leq A_1\tilde A_1\leq P_{24}$. Checking \ref{h2sl2} we see that the summand $(1{[r]}\otimes 1^{[s]})/1^{[r]}$ has a non-trivial value of $H^2$. However it is not a potential blocker since it is in the lowest level and so certainly has no non-trivial $H^1$ in a lower one. On the other hand, Level 4 contains a potential blocker because, by \ref{h1a1} we require $s>1$ for a non-trivial $H^1$ from the first summand in Level 1, yet by \ref{h2sl2} $H^2(A_1,1^{[2]})=K$ and hence the condition $s>1$ gives us a non-trivial $H^2$ in Level 4.

To find the other values of $H^1(\bar X,V)$ for the remaining $\bar X$ with $V=V_S$, one can first reduce to the untwisted modules as listed in the table by \cite[7.1]{CPSV77}, unless $\bar X\cong C_n$ and one must consider modules twisted up to once. Then use the identity $H^1(\bar X, L(\lambda))=H^{0}(\bar X, H^0(\lambda)/\soc(H^0(\lambda)))$ to get the value of $H^1(\bar X,M)$ for any irreducible module $M$. For example $H^1(B_3,100)=H^0(B_3,000)=K$; additionally $H^1(B_2,01^{[1]})=H^0(B_2,00/10/00)=K$. Now for any indecomposable summand $M=\frac{W}{V}$ we can use the long exact sequence of cohomology corresponding to the short exact sequence $0\to V\to M\to W\to 0$ of $\bar X$-modules, \[0\to H^0(\bar X,V)\to H^0(\bar X,M)\to H^0(\bar X,W)\to H^1(\bar X,V)\to \dots\] Alternatively one can observe that some of the modules are tilting modules (e.g. $00/10/00$ for $G_2$) and so have trivial $1$-cohomology by Lemma \ref{h1tilting}.
\end{proof}

\begin{lemma}\label{lift}Provided $\bar X$ is not in the following list of exceptions, all complements lift. That is to say that the map $H^1(\bar X,Q)\to H^1(\bar X,Q/Q(i))$ is surjective for each $i$. 

Exceptions $\bar X=$: 
\begin{itemize}\item $A_1^2\leq A_1^2\tilde A_1\leq B_3$ (2) when $p=3$, $t>0$ and $s=t\pm 1$ (see the table in \ref{lir});
\item $A_1\leq A_1^2\tilde A_1\leq B_3$ when $p=3$, $t>0$ and $s=t\pm 1$;
\item  $L'=A_1\tilde A_1$ when $p=2$;
\item $A_1\leq A_1\tilde A_1=L'$ when $p=2$ and $r\neq s\pm 1$\end{itemize} 
\end{lemma}
\begin{proof} This follows using  \ref{blockers} from an examination of  \ref{levelsforx} provided there are no potential blockers. For $\bar X=A_1\leq A_1^2\leq B_2=L'$ when $p=2$, for $\bar X=A_1^2\leq B_2=L'$ when $p=2$ and for $\bar X=A_1\leq A_1\tilde A_1=L'$ when $p=3$ we show the surjectivity in spite of the presence of a potential blocker: 

For $A_1\leq B_2$, tackle instead the case of $A_1\leq\tilde A_1^2\leq B_2\leq L'$. The result will follow on applying the graph automorphism of $G$. In this case, we use observe that the shapes $V$ affording modules with non-zero $H^1(\bar X,V)$ are $1**0$ and $1**2$; for $B_2$ they both give modules $W(10)$, the natural module and so for $\bar X$ they give the modules $(1^{[r]}+1^{[s]})/0$ which afford non-zero $H^1$. Now the shapes $1**0$ and $1**2$ are in fact part of the $B_4$ subsystem subgroup, whereas the shape affording the potential blocker is not. Thus we have the surjectivity by \ref{subsystemlift}.

The proof for $\bar X=A_1\leq A_1\tilde A_1=L'$ is similar. Since we will need to be more specific about this situation, we defer this case to \ref{a1sinchar3}.
 \end{proof}
 
 \begin{corollary}\label{exhaust} Provided $(\bar X,L)$ is not in the above list, of exceptions $\rho$ is defined on the whole pre-cocycle space $\mathbb V$. Thus an exhaustive set of representatives of cocycles $\bar X\to Q$ is indexed by $K^r$ where $r=\dim \mathbb V$.\end{corollary}
 
 \begin{remark} Later we will show that in each of the exceptions listed above one has that $H^1(\bar X,Q)\to H^1(\bar X,Q/Q(i))$ is not surjective for some $i$; in other words, not all complements to quotients of $Q$ in a semidirect product with $\bar X$ will lift.\end{remark}

\section{Simple subgroups of parabolics: Proof of Theorem 1}\label{thm1}

Before we begin in earnest the proof of theorem 1 we would like to remind the reader of part of the strategy in the introduction. Specifically, we remind the reader that we are in dialogue with the table of restrictions in \ref{tor}. By the end, this will contain a list of all semisimple subgroups of $G=F_4(K)$ inside parabolics. Most of the restrictions are calculated there although a few are deferred to elsewhere in the proof; occasions where this happens are indicated in the table.

The section below is a typical example of what happens; we show by cohomological methods that there can be up to two non-$G$-cr subgroups of type $A_2$, that these must be in the $B_3$-parabolic of $G$  and exhibit two candidates for these, $Z_1$ and $Z_2$. We use a result of \cite{CC88} to see that these are in parabolics, then \ref{tor} to see that they are distinct from each other and from the Levi subgroup they project to. This completes the classification

\subsection{Subgroups of rank $\geq 2$}\label{geq2}

Here we begin the proof of the classification of simple closed, connected subgroups of $G$ contained in parabolics of $G$. The table in \S\ref{lir} gives a list of all $G$-reducible, $G$-cr subgroups of $G$. Hence, we assume $X$ is non-$G$-cr, with $\bar X\leq L$ such that $\bar X$ is $L$-irreducible. 

\begin{lemma}\label{a2sinparabs}Let $p=3$ and define the subgroups $Z_1:=A_2\hookrightarrow A_2\tilde A_2$ as $(10,01)$ and $Z_2:=A_2\hookrightarrow B_4$ by $V_9\downarrow Z_2=T(11)$. Then $Z_1$ and $Z_2$ are both contained in parabolic subgroups of $G$.\end{lemma}
\begin{proof}From \ref{tor} we deduce that the 25-dimensional irreducible module 0001 for $G$ has restrictions $0001\downarrow Z_1=T(11)+T(11)+11$ and $0001\downarrow Z_1=T(11)+ 11/00+00/11$. Thus both subgroups stabilise totally isotropic $2$-spaces in $V$.

In \cite[Table 2]{CC88} we are given the action of $F_4(q)$ on the on $1$-spaces of the $27$-dimensional representation $V_{27}$ for $E_6(q)$ where $F_4(q)$ stabilises the subspace $\langle e\rangle\leq V_{27}$. When $q=3^r$, the irreducible 25-dimensional module $L(0001)$, is isomorphic to $e^\perp/\langle e\rangle$. The stabilisers of subspaces of $L(0001)$ are then the stabilisers of points of $\mathbb P(V_{27})/\langle e\rangle$ which are contained in $e^\perp$. From the table we see that a subgroup of $F_4(q)$ stabilising an isotropic $1$-space in $L(0001)$ must be in a parabolic or a subgroup of $D_4(q)$. 

Let $G(q)=G^{F^q}$, be the finite group of Lie type, i.e. the fixed points of the standard $q$th Frobenius map. For any $F$-stable subgroup $H$ of $G$, \cite[2.10]{BMR05} implies that it is always possible to take a large enough value of $q$ so that $H^F$ is contained in any parabolic subgroup $P$ of $G$ if and only if $H$ is.

Thus $Z_1$ and $Z_2$ are contained in parabolic subgroups or subgroups of type $D_4$. Now all subgroups of $D_4$-ir subgroups of type $A_2$ are contained in $B_3$ subgroups of $D_4$ and these are then conjugate to the Levi $L'=B_3$ in $G$. Examining the composition factors of $V_{26}\downarrow Z_1,Z_2$ in \ref{tor} we conclude that if $Z_1$ or $Z_2$ were contained in $B_3$, the natural module for $B_3$ would restrict to $Z_1$ or $Z_2$ as $11$. Comparing the restrictions of $V_{26}$ to $Z_1$, $Z_2$ and the subgroup $A_2\leq B_3$ by $V_7\downarrow A_2=11$, we see this is impossible.\end{proof} 

\begin{prop}[Classification of $A_2$s]\label{a2s} Suppose $X$ is a subgroup of type $A_2$ which is not $G$-cr. Then $p=3$, $X$ is in a $B_3$-parabolic and there are two conjugacy classes of such.\end{prop}
\begin{proof}
If there exist complements to $Q$ in $\bar XQ$ which are not $Q$-conjugate to $\bar X$ then we must have $H^1(\bar X,Q(i)/Q(i+1))\neq 0$ for some $i$. We check \ref{levelsforx} to see that this forces $p=3$ and that $X$ is in a $B_3$-parabolic. Moreover, from \ref{levelsforx} we get $Q/Q(2)\cong 11+00$ and $Q(2)\cong 11$. Since $\dim \mathbb V=2$, we get from Corollary \ref{exhaust} an exhaustive set of cocycles $H^1(\bar X,Q)$ indexed by $[k_1,k_2]\in K^2$.

We will show that there are at most two classes of non-$G$-cr subgroups of type $A_2$, then we exhibit two distinct non-$G$-cr classes of subgroups.

Let $w\in Q/Q(2)$ be any element of $Q/Q(2)$ and let $\hat w$ be its canonical lift (cf \ref{clift}). Now let $\langle v \rangle=Q/Q(2)^{\bar X}$ denote the $1$-dimensional space of fixed points in $Q/Q(2)$ under the action of $\bar X$ and let $W$ be the submodule $11\leq Q/Q(2)$. For $\lambda\in K^\times$, define a map $c_\lambda :W\to Q(2)$ by $c_\lambda (wQ(2))=[\widehat{\lambda v},\hat w]$. We claim that $c_\lambda$ is an isomorphism of $\bar X$-modules: By the commutator relations, $[\widehat{\lambda v},\hat w]$ is confined to $Q(2)$. Moreover it is $\bar X$-equivariant: $c_\lambda (w^x)=[\widehat{\lambda v},\hat w^x]=[\widehat{\lambda v}^xv',\hat w^x]$ for some $v'\in Q(2)$. But $v'$ is central in $Q$ so this is $[{\widehat{\lambda v}}^x,\hat w^x]=[\widehat{\lambda v},\hat w]^x=(c_\lambda (w))^x$. Since $c_\lambda$ is also non-zero, it is an isomorphism. By Schur's lemma, $\Hom(W,Q(2))=K$ and so the isomorphisms $c_\lambda$ together with $0$ form a vector space with $(c_\lambda+c_\mu)(w)=c_\lambda (w)+c_\mu (w)$. 

Since $Q/Q(2)=W+K$ and $H^1(\bar X,Q/Q(2)=H^1(\bar X,W)+H^1(\bar X,K)$ with the second summand being $0$, we have that any cocycle $\gamma:\bar X\to Q/Q(2)$ has image in $W$. Now, by composing any cocycle $\gamma:\bar X\to W$ with $c_\lambda$ we get a new cocycle $c_\lambda\circ \gamma: \bar X\to Q(2)$ giving our $c_\lambda$ as an explicit description of the isomorphisms $H^1(\bar X,Q/Q(2))\to H^1(\bar X,Q(2))$. Since we have a map $H^1(\bar X,Q)\to H^1(\bar X,Q/Q(2))$, we also get maps $c_\lambda\circ:H^1(\bar X,Q)\to H^1(\bar X,Q(2))$.

Now, each $c_\lambda\circ$ pulls back under $\rho$ to an endomorphism $C_\lambda$ on $\mathbb V$, with $C_\lambda:[k_1,k_2]\mapsto [0, \lambda k_1]$: it is easy to see that this forms a commutative diagram
\[\begin{CD}\mathbb V @>C_\lambda>>\mathbb V\\
@V\rho VV @V\rho VV\\
H^1(\bar X,Q) @>c_\lambda\circ >> H^1(\bar X,Q(2)) \end{CD}\] where the right hand vertical map makes sense as the first co-ordinate of $\mathbb V=H^1(\bar X,Q/Q(2))+H^1(\bar X,Q(2))$ is $0$. Since $q.[q,\hat v]=q^v$ for any $q\in Q$ we get another map $d_\lambda:H^1(\bar X,Q)\to H^1(\bar X,Q)$ via $d_\lambda(\gamma)(x)=\gamma(x).c_\lambda(\gamma(x))=\gamma^{\widehat(\lambda v)}(x)$. So we can see that each $d_\lambda$ pulls back to an endomorphsim $D_\lambda$ of $\mathbb V$, with $D_\lambda:[k_1,k_2]\mapsto [k_1,k_2+\lambda k_1]$, induced by $Q$-conjugacy.

Let $X$ be a complement to $Q$ in $\bar XQ$. Then up to $Q$-conjugacy, $X$ corresponds by \ref{exhaust} to an element $[k_1,k_2]\in \mathbb V$. Now if $k_1\neq 0$, then for some choice of $\lambda$, we get $D_\lambda([k_1,k_2])=[k_1,0]$. Thus $\rho:\mathbb V\to H^1(\bar X,Q)$ is surjective on restriction to the set $A:=\{\mathbf v=[v_1,v_2]\in\mathbb V : v_1v_2=0\}$.

Now since $L=L'Z(L)=B_3T_1$, the torus $T_1$ normalises $Q$. It is easy to check that if $t(\lambda)\in T_1$ and if $x_{\alpha}(u)Q(2)\in Q/Q(2)$ (resp.  $x_{\alpha}(u)\in Q(2)$) then $x_{\alpha}(u)^{t(\lambda)}=x_{\alpha}(u\lambda)$ (resp. $x_{\alpha}(u)^{t(\lambda)}=x_{\alpha}(u\lambda^3)$). Thus $T_1$ induces an action on $\mathbb V$ with $(v_1,v_2)^{t_\lambda}= (\lambda v_1,\lambda^3 v_2)$. The $T_1$ orbits on the set $A$ are therefore represented by $(v_1,v_2)=(0,0),\ (0,1)$ or $(1,0)$. Since $\rho(0,0)$ is the zero cocycle, in $H^1(X,Q)$ it corresponds to a complement $X_0=\bar X\leq L$. Thus there are at most two conjugacy classes of non $G$-cr subgroups of type $A_2$.

Now from Lemma \ref{a2sinparabs} subgroups $Z_1$ and $Z_2$ are in parabolics and not conjugate to subgroups of Levis. Since the restrictions $V_{26}\downarrow \bar X,\ Z_1$ and $Z_2$ are distinct, these are all distinct up to $G$-conjugacy.\end{proof}

\begin{lemma}[Classification of subgroups of type $B_3$ and $G_2$]\label{classb3sg2s} Let $X$ be a simple subgroup of $G$ of type $B_3$ or $G_2$ and suppose $X$ is non $G$-cr. Then $p=2$ and $X$ is conjugate to just one of the following subgroups:
\begin{enumerate}
\item $B_3\leq D_4$
\item $B_3\leq \tilde D_4$
\item $G_2\leq D_4$
\item $G_2\leq \tilde D_4$
\end{enumerate}
\end{lemma}
\begin{proof}
That $p=2$ follows from \ref{levelsforx}.

Clearly $L'=B_3$ or $C_3$. Without loss, we consider $\bar X\leq L'$ for $L'=B_3$ and obtain the remainder of the results by taking the images of $X$ under the graph automorphism of $G$.

We see from \ref{levelsforx} that $\mathbb V\cong K$ and thus $H^1(\bar X,Q)$ is at most one dimensional. 

\begin{claim} $H^1(B_3,Q)=H^1(G_2,Q)=K$.\end{claim}
{\it Proof of claim}
In fact, examining \ref{levelsforx} we see that $H^1(B_3,Q/Q(2))=H^0(B_3,Q)=0$ and so we have  from \ref{nofixedpoints} that $H^1(B_3,Q)=H^1(B_3,Q/Q(2))\oplus H^1(B_3,Q(2))=K$.

We wish to show that $H^1(G_2,Q)=K$ also. (We cannot use \ref{nofixedpoints} for this as $H^0(G_2,Q/Q(2))\neq 0$.) Consider the long exact sequence of (non-abelian) cohomology from \ref{les}:
\[0\to Q(2)^{\bar X}\to Q^{\bar X}\to Q/Q(2)^{\bar X}\to H^1(X,Q(2))\to H^1(X,Q)\to H^1(X,Q/Q(2)).\]
Inserting terms that we know gives us
\[0\to K \to Q^{\bar X} \stackrel{1}{\to} K \stackrel{2}{\to} K \stackrel{3}{\to} H^1(X,Q)\to 0.\]
As the last map is zero, it follows that map 3 is surjective. We wish to show that it is also injective. This will follow if map 2 is zero. This in turn will follow if the map 1 is surjective.

Consider now the $A_2$ subsystem subgroup of $G_2$ generated by long root groups. As $L_{G_2}(10)\downarrow A_2=10+01$ when the characteristic of $K$ is $2$ (see e.g \cite{Ste09-2}), we get $Q/Q(2)\downarrow A_2=10+01+00^2$ and $Q(2)\downarrow A_2=10+01+00$. 

Now this $A_2$ subgroup of $G_2$ is conjugate in $B_3$ to the $A_2$ Levi subgroup of $B_3$. Thus we may assume that it is generated by root groups with roots of the form $**00$. Now we can pick out  the root groups in $Q/Q(2)$ centralised by the action of $A_2$; they are those corresponding to roots $0001$ and $1231$: one cannot add any root of the form $**00$ to these to get a new root, so it follows by the commutator formula that this $A_2$ centralises these root groups. Similarly, in $Q(2)$, the stable subspace is generated by the root group corresponding to $1232$.

Thus the fixed submodule of $Q/Q(2)$ for $\bar X$ is contained in the subspace generated by root groups corresponding to $0001$ and $1231$. These root groups generate a three dimensional subgroup $R$ of $Q$ corresponding to roots $0001,\ 1231$ and $1232$. Now $\bar X$ centralises the latter root group, and centralises a one dimensional subspace of the quotient of $R$ by this subgroup, so stabilises a $2$-dimensional subgroup $E$ of $R$. Since a maximal torus of $\bar X$ centralises the whole of $R$, it must be contained in the kernel of this action; but as $\bar X$ has no positive dimensional normal subgroups, $\bar X$ must in fact centralise $E$.

We have shown that the $1$-dimensional submodule of $Q/Q(2)$ in fact has a lift to a $1$-dimensional submodule of $Q$. This is to say that the map 1 above is surjective. This proves the claim.

Thus by \ref{onedimh1} there is one class of non-$G$-cr subgroups of type $G_2$ in the $B_3$-parabolic and one non-$G$-cr subgroup of type $B_3$. Let $Y$ be the subgroup $B_3\leq D_4$  and $Z$ the subgroup $G_2\leq D_4$. We wish to show that $Y$ and $Z$ are non-$G$-cr.

Now, again from \ref{levelsforx}, we conclude that the shape $S$ affording a module $V_S$ with a  non-trivial value of $H^1(\bar X,V_S)$ is the shape $***2$ whether $\bar X=B_3$ or $G_2$. Thus $H^1(\bar X, Q)\cong H^1(\bar X,Q(2))$ via the injection $Q(2)\to Q$. Consider the subsystem $B_4$ corresponding to roots $\{***0,***2\}$. Then we have $Q\cap B_4=Q(2)$, so that all the root groups in $Q(2)$ together with the Levi $B_3$ subgroup of $G$ are contained in this $B_4$.  Thus by \ref{subsystemlift} we may assume our non-$G$-cr subgroups are in $B_4$. As $Y$ and $Z$ are in $B_4$, we are done if we can show that they are complements to $Q(2)$ in $\bar XQ(2)$ and not $B_4$-conjugate to their respective $\bar X$s.

Now, the restriction of the $9$-dimensional Weyl module $V_9\cong 1000/0000$ for $B_4$ to the $B_3$ Levi  is $V_7+000^2$, where $V_7$ is the $7$-dimensional Weyl module for $B_3$. Thus if $\bar X=G_2$, $V_9\downarrow \bar X=10/00 +00^2$. But $V_9\downarrow D_4=1000+0000$ and so , $V_9\downarrow Y=000/100/000+000$; thus $Y$ is in a $B_3$-parabolic. Also we see easily that $V_9\downarrow Z=00/10/00+00$; thus $Z$ is in a $B_3$-parabolic. Since neither is $GL(V_9)$ conjugate to the respective $\bar X$ of type $G_2$ or $B_3$, neither is $B_4$-conjugate to $\bar X$. Thus $Y$ and $Z$ must be representatives of the two non-$B_4$-cr subgroups and thus representatives of the two non-$G$-cr subgroups as claimed.

To finish off, we need to know that the images of these subgroups under the graph automorphism of $G$ must be distinct. The image of the $Y=B_3\leq D_4$ under the graph automorphism must be an isogenic image of $B_3$, thus is either $B_3$ or $C_3$, but sitting in a $\tilde D_4$ means it must be a subgroup of type $B_3$; comparing the restrictions of both non-$G$-cr $B_3\leq D_4$ and $B_3\leq \tilde D_4$ in \ref{tor} shows that these are not conjugate. The image of $Z$ under $\tau$ is similarly not conjugate to $Z$ from \ref{tor} as they do not even have the same composition factors on the restriction of $V_{26}$.
\end{proof}

In fact, one can show by similar means that when $p=2$ that a subgroup $B_n$ contained in a $D_{n+1}$ subgroup of $B_{n+1}$ is always non-$B_{n+1}$-cr. With that in mind, we refer to subsystem subgroups of $G$ of type $D_3$ rather than $A_3$ in the following lemma.

\begin{lemma}[Classification of $B_2$s]\label{classb2s} Let $X$ be a simple subgroup of $G$ of type $B_2$ and suppose $X$ is non $G$-cr. Then $p=2$ and $X$ is conjugate to just one of the following subgroups
\begin{enumerate}\item $X\leq D_3$
\item $X\leq \tilde D_3$
\item $X\leq B_2B_2$ by $(10,10)$.\end{enumerate}
 \end{lemma}
\begin{proof}Suppose $X$ is of type $B_2$. We apply a similar argument to the $A_2$ case, \ref{a2s}. From \ref{levelsforx}, $L'=B_2$, $p=2$ and there are two shapes affording non-trivial $H^1$ if the action of $\bar X$ on $Q$ is untwisted; viz. $1**0$ and $1**2$, both affording modules $10/00$. There are two more if the action of $\bar X$ is twisted; viz. $0**1$ and $1**1$ affording modules $01$ where $H^1(\bar X,01^{[1]})=K$. To get all complements we may as well assume that the action is twisted. By \ref{exhaust}, complements are parameterised by $[k_1,k_2,k_3,k_4]\in K^4$, with each coordinate corresponding to shapes $1**0,0**1,1**1$ and $1**2$, affording modules $10^{[1]},01^{[1]},01^{[1]}$ and $10^{[1]}$ respectively. Recall the notation of \ref{a2s}, where $\widehat{x}\in Q$ denotes the canonical lift of $x$, an element projecting to $x\in Q/Q(i)$ under the quotient map. 

Let $\langle v\rangle=x_{0122}(t)$ be the trivial submodule (i.e. fixed points) of $\bar X$ on $Q(2)/Q(3)$. Let $W$ denote the submodule $10^{[1]}/00\leq Q(1)/Q(2)$ and $W'$ the module $10^{[1]}/00$ in $Q(3)/Q(4)$. Then for each $\lambda\in K^\times$, we get an isomorphism $c_\lambda:W\to W'$ by $c_\lambda(vQ(2))=[\widehat{\lambda v},\hat w]Q(4)$: the image is confined to $Q(3)/Q(4)$ by the commutator relations. It is also $\bar X$ equivariant by the same argument in the proof of \ref{a2s}. Now $\Hom_{\bar X}(W,W')=K$ and so the isomorphisms $c_\lambda$ form a vector space with $(c_\lambda+c_\mu)(w)=c_\lambda(w)+c_\mu(w)$.

Still following the proof of \ref{a2s}, we can pull back the $c_\lambda$ to maps $C_\lambda$ on $\mathbb V$ which we check are $C_\lambda:[k_1,k_2,k_3,k_4]\to [0,k_2,k_3,\lambda k_1]$, with a map induced by conjugacy with the element $\widehat{\lambda v}$ of $Q$ given explicitly as $[k_1,k_2,k_3,k_4]^{\widehat{\lambda v}}=[k_1,k_2,k_3,k_4+\lambda k_1]$. Thus if $k_1\neq 0$, by some choice of $\lambda$ we have $[k_1,k_2,k_3,k_4]^{\widehat{\lambda v}}=[k_1,k_2,k_3,0]$. Thus $\rho:\mathbb V\to H^1(\bar X,Q)$ is surjective on $\mathbf v\in \mathbb V$ such that $k_1k_4=0$. 

Similarly, we let $\langle w\rangle=x_{1110}(t)$ be the subspace of fixed points of $X$ on $Q/Q(2)$. Then by conjugating with $w$. we pull back to an action on $\mathbb V$ which is $[k_1,k_2,k_3,k_4]^{\widehat{ \lambda w}}=[k_1,k_2,k_3+\lambda k_2,k_4]$. So again we can deduce that $\rho:\mathbb V\to H^1(\bar X,Q)$ is surjective on $\mathbf v\in \mathbb V$ such that $k_2k_3=0$.

Then we see again that conjugation by the two dimensional torus $T_2=Z(L)$ reduces the possibilities to up to eight (up to $G$-conjugacy): $[0,0,0,0]$, $[1,0,0,0]$, $[0,1,0,0]$, $[0,0,1,0]$, $[0,0,0,1]$,  $[1,1,0,0]$, $[1,0,1,0]$ and $[0,1,0,1]$. Amongst these, however we can exhibit some further $G$-conjugacy. If $\mathbf v=[0,0,0,1]\in \mathbb V$ then $\rho v$ is a cocycle $\gamma:\bar X\to Q(3)$. Now the root groups $x_r(t)$ generating $Q(3)$ have $r\in \{1122,1222,1232,1242,1342,2342\}$ and an elementary calculation shows that these are mapped under the Weyl group representative $n_{0122}$ to $\{1000,1100,1110,1120,1220,1122\}$. This element also centralises the Levi subgroup $B_2$ so that in particular, conjugation by $n_{0122}$ induces a map $[0,0,0,1]\mapsto [1,0,0,0]$ in $\mathbb V$. Similarly $n_{1110}$ induces a map $[0,1,0,0]\mapsto [0,0,1,0]$. Putting this together we see that we have up to only four possibilities: $[0,0,0,0]$, $[1,0,0,0]$, $[0,1,0,0]$, $[1,1,0,0]$. Thus there are up to four conjugacy classed of  subgroups $B_2$. One of these is the Levi $B_2$. We now wish to show that there are three non-$G$-cr subgroup $X$ of type $B_2$.

First we show that there is a non-$G$-cr subgroup of type $C_2$ in a $C_3$ Levi of $G$. We recall the basic facts of exceptional isogenies from \cite[3.9]{BT73}. There is a map $B_3\to GL_6$ induced by factoring out the radical of the $7$-dimensional natural module for $B_3$. The image of this map is in fact a $C_3$, where the action of this image on the resulting $6$-dimensional space has the structure of a module $100^{[1]}$. The $D_3$ subgroup of $B_3$ acts on $V_7$ as $100+000$; thus the radical for the $V_7$ for $B_3$ is a trivial direct summand for $D_3$. 

Take $X=B_2\leq D_3$. Then it acts as $00/10/00+00$. Thus the image $Y$ say of $X$ induced by factoring out Rad $V_7$ is a $C_2$ subgroup of $C_3$ with action $00/10^{[1]}/00$ on the $6$-dimensional space. Now $V_{26}$ restricted to a Levi $C_2$ of this $C_3$ is $V_{26}\downarrow C_2=(10^{[1]})^4+00^2+00/01^{[1]}/00$. However, this contains no submodule $00/10^{[1]}/00$ so $Y$  is not conjugate to the $C_2$ Levi of $C_3$. Thus, using the graph automorphism, we see that $X$ is not conjugate to the $B_2$ Levi of $G$.

Now $X$ together with its image under the graph automorphism of $G$, $X^\tau\leq \tilde D_3$ give two of the three subgroups we require, which are non-conjugate by observing their different actions in \ref{tor}.

For the remaining non-$G$-cr subgroup of type $B_2$, let $X:=B_2\leq B_2B_2$ by $(10,10)$. Then $X$ is in a $B_4$ with $V_9\downarrow X=(10+10)/00$. Since there are two equivalent submodules modulo the radical, $X$ is in a $D_3$ parabolic subgroup of $B_4$. Thus in $F_4$, $X$ must be in a $B_2$-parabolic subgroup of $G$. Again, checking restrictions in \ref{tor}, we get that $X$ is not  conjugate to any subgroup already discovered.
\end{proof}

\subsection{Subgroups of type $A_1$, $p=3$}\label{a1p3}

We broadly follow the formula of the previous section. First we give a statement of the results; that is, a list of the conjugacy classes of non-$G$-cr subgroups of type $A_1$ in characteristic $p=3$. Where $\mathbb V$ is non-zero, we establish the image of partial map $\rho$; we find in fact that it is either fully-defined or zero. We then observe as much conjugacy between the images of $\rho$ as we can, and give descriptions of each subgroup as contained in subsystem subgroups in terms of the restrictions of modules for the subsystems {\it or} in one set of cases observe that there is no more conjugacy to be observed.

\begin{theorem}\label{a1sinchar3}\label{a1schar3} Let $p=3$ and let $X$ be a non $G$-cr subgroup of type $A_1$. Then $X$ is conjugate to precisely one of the following, where $L'$ is a minimal Levi subgroup containing $\bar X$: \begin{enumerate}
\item $X\hookrightarrow B_3$ by $V_7\downarrow X=T(4)+ 0$ ($L'=B_2$),
\item $X\hookrightarrow C_3$ by $V_6\downarrow X=T(3)$ ($L'=A_1\tilde A_1$),
\item $X\hookrightarrow B_4$ by $V_9\downarrow X=T(4)^{[r]}+ 2^{[s]}$, $rs=0$ ($L'=B_3$).
\item $X\hookrightarrow A_2\tilde A_2$ by $(V_3, \tilde V_3)\downarrow X_2=(2,\tilde 2)$ ($L'=B_3$)
\item $X\hookrightarrow A_1C_3$ by $(1^{[r]},T(3)^{[s]})$, $rs=0$ ($L'=B_3$)
\item $X$ in a one-dimensional variety $\Omega\cong \mathbb A^1(K)\setminus \{0\}$ of subgroups of type $A_1$, with no proper reductive overgroup in $G$. ($L'=B_3$)
\end{enumerate}
\end{theorem}

{\it Proof of Theorem \ref{a1sinchar3}}

Observe from \ref{levelsforx} that $L'=A_1\tilde A_1,\ B_2$ or $B_3$. Each of these possibilities for $L'$ does in fact give rise to non-$G$-cr subgroups.

\underline{$L'=B_2$} (one subgroup arising)

The embedding $\bar X\leq L$ is determined by the restriction of the natural module $V_5$ for $B_2$ to $\bar X$ by $V_5\downarrow \bar X=1^{[r]}\otimes 1^{[s]}+ 0$. As $\mathbb V\neq 0$, $r=s\pm 1$. Without loss of generality, we may assume that $r=1$, $s=0$.

We may employ the argument we used in \ref{a2s} to reduce to up to two possibilities. Specifically, we have that since there are no potential blockers, $\rho$ is fully-defined, giving a parameterisation of $H^1(\bar X,Q)$ by $[k_1,k_2]\in K^2$. As in \ref{a2s} have an action induced on $\mathbb V$ by conjugation with $x_{0122}(t)$, giving up to two possibilities up to $QZ(L)$-conjugacy; viz. $\bold k=[0,1]$, and $\bold k=[1,0]$. Further, conjugation by $n_{0122}$ pulls back to an action on $\mathbb V$ given by $[k_1,k_2]^{n_{0122}}=[k_2,k_1]$. Thus we find that there is up to one conjugacy class of non-$G$-cr subgroups of type $A_1$ contained in the $B_2$ parabolic. 

Take a subgroup inside a $B_3$ by acting on the seven-dimensional natural module $V_7$ for $B_3$ with restriction $T(4)+0=0/1\otimes 1^{[1]}/0+0$: this is in a $B_2$-parabolic of $B_3$, so in a $B_2$-parabolic of $G$; it is in nothing smaller, since the composition factors of the restriction in \ref{tor} are incompatible with being in an $A_1$- or $\tilde A_1$-parabolic. Moreover, from \ref{tor} it is not conjugate to $\bar X$ as it has a different action on $V_{26}$.

\underline{$A_1 \tilde A_1$ parabolics} (one subgroup arising)

Let the embedding $\bar X\to L$ be $(1^{[r]},1^{[s]})$. From \ref{levelsforx}, the only way to get a module of type $1\otimes 1^{[1]}$ is, without loss of generality, if we let $r=0$,$s=1$. 

\begin{claim} For each parabolic $P_{13}$, $P_{23}$, $P_{24}$, the root groups $x_\alpha(t)$  corresponding to shapes $S$ such that $H^1(\bar X,V_S)\neq 0$ together with $\bar X$ are contained in subgroups of type $C_3$.\end{claim}
{\it Proof of claim:} This is obvious for $P_{13}$ as the shape $S$ affording a non-trivial $H^1(\bar X,V_S)$ is $0*1*$, with roots $\{0010,0110,0011,0111\}$ is visibly contained in the $C_3$ Levi. Conjugating by the Weyl group representative $n_{1000}n_{0100}$ we see that the the Levi subgroup $L'$ of $P_{13}$ is mapped to that of $P_{23}$ and the shape $0*1*$ is mapped to $*11*$ as required. Thus $C_3^{n_{1000}n_{0100}}\cap P_{23}$ contains $L'$ and the shape $*11*$. Further, conjugating by $n_{0010}n_{0001}$ moves $L'\leq P_{23}$ to $L'\leq P_{24}$ and $*11*$ to $*1*1$ as required. This proves the claim.

For each parabolic, we have $\mathbb V=K$ with an action of $Z(L)$ on $Q$ by conjugacy pulling back to a $K^\times$-action on $\mathbb V$. Thus $\rho$ is surjective on the set $\{0,1\}\in \mathbb V$; in other words, there is up to one non-$G$-cr complement to $Q$. 

Consider now just $P_{13}$. Since there are no potential blockers in $Q\cap C_3$, (the $C_3$ Levi,) we see from \ref{subsystemlift} that $H^1(\bar X,Q)=K$ and from \ref{onedimh1} that there is exactly one non-$G$-cr subgroup, $Y$ say, and some conjugate of $Y$ is contained in the $C_3$ Levi. By the claim again, $Y$ is conjugate to a non-$G$-cr subgroup in $P_{24}$ and $P_{23}$, thus all of them.

Now we describe $Y$: the representation $2\otimes 1=1/1^{[1]}/1=T(3)$ is self-dual and so describes a subgroup $X$ of type $A_1$ in $D_3$ or $C_3$. Since $2\otimes 1$ has odd high weights, it is in a $C_3$. In fact, it is easy to see that it is in an $A_1\tilde A_1$-parabolic of $C_3$, and thus of $G$, but is not conjugate to $\bar X$ by considering the restriction in \ref{tor}.

\underline{$L'=B_3$} \label{a1schar3sec}

Since $p=3$, we get from \ref{a1sinb3} that the possible embeddings $\bar X\leq L'$ are given by restriction of the natural module $V_7$ for $B_3$ to $\bar X$ as $V_7\downarrow \bar X=1^{[r]}\otimes 1^{[s]}+2^{[t]}$, where if $r=s$, $t\neq r$. From \ref{levelsforx} we have $Q(2)\cong V_7$, and $Q(1)/Q(2)\cong 1^{[r]}\otimes 1^{[t]}+1^{[s]}\otimes 1^{[t]}$. There are various cases to deal with according to the values of $r$, $s$ and $t$. We need to consider all cases where there is a module with non-zero $H^1$ in $Q(2)$ or $Q(1)/Q(2)$, viz., using \ref{h1a1}, when $s=r\pm 1$, $t=r\pm 1$ or $t=s\pm 1$. Potential blockers can be read from \ref{h2sl2}: in this case, they are either $2^{[t]}$ with $t>0$ or any twist of $1\otimes 1^{[2]}$. Note that when $r=s>0$ then $1^{[r]}\otimes 1^{[s[}=2^{[r]}+0$ can yield a potential blocker. 

As there is symmetry in $r$ and $s$ induced by conjugacy (by some $n_\alpha$), we assume $r\leq s$. Thus the cases that need considering are, 

$\begin{array}{c c} \hspace{100pt}r=s,\ t=r+1  \hspace{100pt}& 
r=s,\ t=r-1\\
s=r+1,\ t=r-1
& s=r+1,\ t=r\\
 s=r+1,\ t=r+1 &
 s=r+1,\ t=r+2\\
 s=r+1,\ t>r+2 &
 s=r+1,\ t<r-1\\
s=r+2,\ t=r+1 &
s=r+2,\ t=r-1\\
s=r+2,\ t=s+1 &
s\neq r+1,\ r+2,\ t=r-1\\
s\neq r+1,\ r+2,\ t=r+1 &
s\neq r+1,\ r+2,\ t=s-1\\
s\neq r+1,\ r+2,\ t=s+1 &
\end{array}$

\begin{claim} Let $\bar X$, $\bar Y$ be subgroups of type $A_1$ embedded $B_3$-irreducibly in a subgroup $B_3$. Then if $V_{26}\downarrow \bar X\cong V_{26}\downarrow \bar Y$ as $A_1$-modules we have that $\bar X$ is $B_3$-conjugate to $\bar Y$.\end{claim}
{\it Proof of claim: }Let $V_7=L_{B_3}(100)\downarrow \bar X=1^{[r_1]}\otimes 1^{[s_1]}+2^{[t_1]}$ and $V_7\downarrow \bar X=1^{[r_2]}\otimes 1^{[s_2]}+2^{[t_2]}$ so that $L_{B_3}(001)\downarrow \bar X=1^{[r_1]}\otimes 1^{[t_1]}+1^{[s_1]}\otimes 1^{[t_1]}$, and similarly for $\bar Y$. We note we are done if we can show $(r_1,s_1,t_1)=(r_2,s_2,t_2)$ or $(s_2,r_2,t_2)$.

Then as $V_{26}\downarrow B_3=100+001^2+000^3$, we have \begin{align*}V_{26}\downarrow \bar X=1^{[r_1]}\otimes 1^{[s_1]}&+2^{[t_1]}+(1^{[r_1]}\otimes 1^{[t_1]})^2+(1^{[s_1]}\otimes 1^{[t_1]})^2+0^3\\&=1^{[r_2]}\otimes 1^{[s_2]}+2^{[t_2]}+(1^{[r_2]}\otimes 1^{[t_2]})^2+(1^{[s_2]}\otimes 1^{[t_2]})^2+0^3=V_{26}\downarrow\bar Y.\end{align*}
From this equation we see that if $r_1\neq s_1$, then the only module which is a twist of $2$ in the LHS is $2^{[t_1]}$; thus it follows that $t_1=t_2$ so (possibly swapping $r_2$ and $s_2$), $r_1=r_2$ and $s_1=s_2$. 

If $r_1=s_1$ then as $\bar X$ is $B_3$-ir, $t_1\neq r_1$, and we have \[2^{[r_1]}+2^{[t_1]}+(1^{[r_1]}\otimes 1^{[t_1]})^4+0=1^{[r_2]}\otimes 1^{[s_2]}+2^{[t_2]}+(1^{[r_2]}\otimes 1^{[t_2]})^2+(1^{[s_2]}\otimes 1^{[t_2]})^2.\] To have an equal number of trivial summands on the RHS we require $r_2=s_2$. Then we deduce either $(r_1,t_1)=(r_2,t_2)$ or $(r_1,t_1)=(t_2,r_1)$. In the first case we are done. If $r_1=s_1=t_2$ and $r_2=s_2=t_1$ then $V_7\downarrow \bar X=V_7\downarrow \bar Y$ in any case. This proves the claim.

Up to twists, we may choose to set the lowest of $r,s,t$ to be equal to $0$. Thus the cases are $(r,s,t)=$
\begin{align*}&(0,0,1) 
&(1,1,0)\\
&(1,2,0) 
&(0,1,0) \\
 &(0,1,1) 
&(0,1,2)\\ 
&(0,1,t),\ t>2 
&(r,r+1,0),\ r>1\\ 
&(0,2,1) 
&(1,3,0) \\
&(1,s,0),\  s>3
&(0,s,1),\ s>2\\
&(0,s,s-1),\ s>2
&(0,s-1,s),\ s>3\\
\end{align*}
Here we will see the first time that $\rho$ is not always defined on the whole of $\mathbb V$, specifically the case where $(r,s,t)=(0,0,1)$.

\subsubsection{(r,s,t)=(0,0,1)}\label{001}

Here we have $Q(1)/Q(2)\cong 1\otimes 1^{[1]}+1\otimes 1^{[1]}$ and $Q(2)=2+2^{[1]}+0$. Thus $\mathbb V=H^1(\bar X,Q(1)/Q(2))\cong K^2$. We show that $\rho:\mathbb V\to H^1(\bar X,Q)$ is defined precisely on the one-dimensional variety $\{(k,\pm k):k\in K\}\leq \mathbb V$.

We have from \ref{a1sinb3} that $\bar X\leq A_1^2\tilde A_1$. The latter is a subsystem of $B_3$ corresponding to roots \[\{\mp 1220,\pm 1000,\pm 0010\}\] where we have picked $\mp 1220$ to be consistent with the Borel-Siebenthal algorithm. Thus up to conjugacy, we take $\bar X$ to be generated by the elements 
\begin{align*}
x_{\alpha_+}(t)&=
  \left[\begin{array}{c}-1220\\t\end{array}\right]
 \left[\begin{array}{c}1000\\t\end{array}\right]
  \left[\begin{array}{c}0010\\t^3\end{array}\right]\\
  x_{\alpha_-}(t)&=
   \left[\begin{array}{c}1220\\t\end{array}\right]
 \left[\begin{array}{c}-1000\\t\end{array}\right]
  \left[\begin{array}{c}-0010\\t^3\end{array}\right],
\end{align*}
where we recall the notation $ \left[\begin{array}{c}\alpha\\t\end{array}\right]:=x_\alpha(t)$.

Let $Q(1)/Q(2)\downarrow \bar X=V_1+V_2$, where $V_1$ (resp. $V_2$) is generated over $K\bar X$ by a high weight vector $x_{0121}(t)$ (resp. $x_{0011}(t)$). We claim that 
 \[ V_1=\left\{\left[\begin{array}{c}0111\\t_1\end{array}\right]
 \left[\begin{array}{c}1111 \\t_2\end{array}\right]
  \left[\begin{array}{c}0121 \\t_3\end{array}\right]
    \left[\begin{array}{c}1121\\t_4\end{array}\right]Q(2):t_i\in K\right\}\] and \[ V_2=\left\{
    \left[\begin{array}{c}1221\\t_1\end{array}\right]
 \left[\begin{array}{c}0001 \\t_2\end{array}\right]
  \left[\begin{array}{c}1231 \\t_3\end{array}\right]
    \left[\begin{array}{c}0011\\t_4\end{array}\right]Q(2):t_i\in K\right\}\]
    
To see this, observe that one cannot add any root from the set $\{0111,1111,0121,1121\}$ to the roots $\mp 1220$. It follows from the commutator relations that the first factor $A_1\leq A_1^2\tilde A_1$ centralises these root groups. Similarly adding $\mp 1220$ to the set $\{1221,0001,1231,0011\}$ gives roots in that set or non-roots. It follows that the first factor preserves the decomposition $V_1+V_2$. Similarly one sees that the other factors preserve the decomposition $V_1+V_2$. Hence $\bar X\leq A_1^2\tilde A_1$ preserves this decomposition as well. In fact one checks from the commutator relations that one gets an isomorphism to $V_1$ and $V_2$ from the explicit representation of $\bar X$ described in \ref{a1h1p3}. So in the notation of that lemma, we get $V_1$ by $e_1\mapsto x_{1121}(1)$, $e_2\mapsto x_{0121}(1)$, $e_3\mapsto x_{1111}(1)$ and $e_4\mapsto x_{0111}(1)$; one has a similar isomorphism for $V_2$. This will allow us to write down the cocycles explicitly.

Let $B$ be a Borel subgroup of $\bar X$ with maximal torus $H$. We now calculate $H^1(B,Q)$ by calculating the map $H^1(B,Q)\to H^1(B,Q(1)/Q(2))$. This is best done in terms of complements. We have by \ref{a1h1p3} that a collection of representative of the complements corresponding to cocycles of $H^1(B,Q(1)/Q(2))$ is $B_{k,l}=\{x_{+,k,l}(t):t\in K\}H$ where 

\[x_{+,k,l}(t)=
  \left[\begin{array}{c}-1220\\t\end{array}\right]
 \left[\begin{array}{c}1000\\t\end{array}\right]
  \left[\begin{array}{c}0010\\t^3\end{array}\right]
  \left[\begin{array}{c}1121\\kt^2\end{array}\right]
 \left[\begin{array}{c}0121\\-kt\end{array}\right]
   \left[\begin{array}{c}0011\\lt^2\end{array}\right]
 \left[\begin{array}{c}1231\\-lt\end{array}\right]Q(2)\]

so that $x_{+,k,l}(t)=x_{\alpha_+}(t)\gamma_{k,l}(x_{\alpha_+}(t))$ with \[\gamma_{k,l}(x_{\alpha_+}(t))=  \left[\begin{array}{c}1121\\kt^2\end{array}\right]
 \left[\begin{array}{c}0121\\-kt\end{array}\right]
   \left[\begin{array}{c}0011\\lt^2\end{array}\right]
 \left[\begin{array}{c}1231\\-lt\end{array}\right]Q(2).\]
 
\begin{claim}The cocycle $\gamma_{k,l}$ lifts to a cocycle $\gamma:\bar X\to Q$ if and only if $k=\pm l$.\end{claim}
{\it Proof of claim: }
We must find functions $f_i(t)$ for $i\in\{0122, 1122,1222,1232,1242,1342,2342\}$ (the roots in $Q(2)$) such that 
\begin{align*}\gamma(x_{\alpha_+}(t))=  \left[\begin{array}{c}1121\\kt^2\end{array}\right]
 \left[\begin{array}{c}0121\\-kt\end{array}\right]
   \left[\begin{array}{c}0011\\lt^2\end{array}\right]
    \left[\begin{array}{c}1231\\-lt\end{array}\right]
 \left[\begin{array}{c}0122\\f_1(t)\end{array}\right]
 \left[\begin{array}{c}1122\\f_2(t)\end{array}\right]\\
 \left[\begin{array}{c}1222\\f_3(t)\end{array}\right]
 \left[\begin{array}{c}1232\\f_4(t)\end{array}\right]
 \left[\begin{array}{c}1242\\f_5(t)\end{array}\right]
 \left[\begin{array}{c}1342\\f_6(t)\end{array}\right]
 \left[\begin{array}{c}2342\\f_7(t)\end{array}\right] \end{align*}
 
 and $\gamma$ satisfies the cocycle condition of \ref{a1h1resb}. In particular, letting $x_+(t)=x_{\alpha_+}(t)\gamma(x_{\alpha_+}(t))$ we require that $x_+(t).x_+(u)=x_+(t+u)$, i.e. $x_+(t).x_+(u).x_+(t+u)^{-1}=I$. 
  
In principle we could use the commutator formula to investigate the restrictions this places on the $f_i(t)$. It is more convenient however to use the 26-dimensional representation of $G$. Recall from the Overview \S\ref{overview} that we have, in GAP, a set of matrices {\tt F4gens}; that is, for each root $\alpha$ of $G$, a $26\times 26$ matrix $X_\alpha(T)$, with the set of all $X_\alpha(t)$ satisfying (in particular) the commutator relations of $G$. Then one can input $x_+(t)$ as a $26\times 26$ matrix $M$ which is a product of some of the $X_\alpha(T)$. For our purposes $M$ will contain $7$ indeterminates $T_i$ corresponding to the $f_i(t)$ and one further, $T_0$ corresponding to $t$.

As a useful first piece of data about the $f_i(t)$, note that we require that $x_+(t)^3=$id for all $t\in K$ as char $K=3$. Asking GAP to cube the matrix, $M$ and insisting the result be the identity matrix we immediately discover that $f_3(t)=0$, $f_6(t)=0$  and $k^2=l^2$. This justifies the `only if' part of the claim.

So let us consider the case where $k=\pm l$. By simply setting values of the functions $f_i(t)$ to $0$, we find that we can lift the cocycle $\gamma_{k,\pm k}$ to a cocycle $\gamma:B\to Q$ if we set \[\gamma(x_{\alpha_+}(t))= \left[\begin{array}{c}1121\\k t^2\end{array}\right]
 \left[\begin{array}{c}0121\\-k t\end{array}\right]
 \left[\begin{array}{c}0011\\\pm k t^2\end{array}\right]
 \left[\begin{array}{c}1231\\\pm -k t\end{array}\right]\]
and \[\gamma(h(t))=1\] and extending to the whole of $B$ by the cocycle condition.

We check that the stipulations of \ref{a1h1resb} hold. One similarly uses the description of $x_+(t)$ as a $26\times 26$ matrix. This proves the claim.

For each $(k,k)\in \mathbb V$, let the complement corresponding to $\rho(k,k)$ be $X_{k^+}$ and the complement corresponding to $\rho(k,-k)$ be $X_{k^-}$. 

\begin{claim} Each $X_{k^+}$ and $X_{k^-}$ is conjugate to $X_{1^+}$. \end{claim}
{\it Proof of claim: }
We observe that we are done if we can show this up to $P$-conjugacy in the semidirect product $\bar XQ/Q(2)$: if $X_{k^\pm}Q(2)$ is $P$-conjugate to $X_{l^\pm}Q(2)$ then so is any lift since $H^1(\bar X,Q(2))=0$. Now $Z(L)$ acts on $Q/Q(2)$ as scalars, and $Q/Q(2)$ is a module for $\bar X$ so $Z(L)$ acts by scalars on $H^1(\bar X,Q/Q(2))$. Thus it follows that $X_{k^+}$ is conjugate to $X_{1^+}$ (respectively $X_{k^-}$ is conjugate to $X_{1^-}$). 

Lastly, a direct check shows that $X_{1^+}$ is conjugate to $X_{1^-}$ by the element $h_{\alpha_1}(-1)h_{\alpha_3}(1)$. This proves the claim.

Finally, we find a more natural description of this subgroup:

\begin{claim} The subgroup $X_{1^+}$ is conjugate to the subgroup $Y\cong A_1\leq A_1C_3$ by $(V_2, V_6)\downarrow Y=(1,1\otimes 2)$.\end{claim}
{\it Proof of claim: }Since $2\otimes 1=1/1^{[1]}/1=T(3)$, $Y$ is certainly in a parabolic of $A_1C_3$ since the projection of $Y$ to $C_3$ is in a parabolic of $C_3$ by \ref{girclassical}. Thus $Y$ is in a parabolic of $G$. It is not minimally contained in an $A_2A_1$ or $C_3$ parabolic of $G$ since there are no non-$G$-cr subgroups in such parabolic subgroups by examination of \ref{levelsforx}. It is also not in any smaller rank parabolic since by \ref{tor} such parabolics stabilise subspaces of $V_{26}$ which have larger dimension than that stabilised by $Y$. Thus it is in a $B_3$-parabolic and so is a complement to some $B_3$-ir subgroup of type $A_1$. Now one checks the composition factors of this subgroup in \ref{tor} to see that it must be conjugate to a complement in $\bar XQ$. But it is not conjugate to $\bar X$ since it has a different action on $V_{26}$, again from \ref{tor}. This proves the claim.

Thus we have shown that there is just one non-$G$-cr subgroup up to conjugacy in this case and that it is the subgroup described in the claim. We see that this is one of the subgroups listed in the statement of \ref{a1sinchar3}. This concludes the analysis of the case $(r,s,t)=(0,0,1)$.

\subsubsection{(r,s,t)=(1,1,0)}

Observe that the natural module $V_7$ for $B_3$ restricts to $\bar X$ as $2+2^{[1]}+0$. Thus $\bar X$ is conjugate to the subgroup $\bar X$ in case $(r,s,t)=(0,0,1)$ by \ref{girsubsconj}.

\subsubsection{(r,s,t)=(0,2,1)}\label{021}

Here we have $H^1(\bar X,Q(2))=0$, $H^1(\bar X, Q/Q(2))=K^2$. We will use the same method as \ref{001}

This time $\bar X$ is generated by 

\begin{align*}
x_{\alpha_+}(t)&=
  \left[\begin{array}{c}-1220\\t\end{array}\right]
 \left[\begin{array}{c}1000\\t^6\end{array}\right]
  \left[\begin{array}{c}0010\\t^3\end{array}\right]\\
  x_{\alpha_-}(t)&=
   \left[\begin{array}{c}1220\\t\end{array}\right]
 \left[\begin{array}{c}-1000\\t^6\end{array}\right]
  \left[\begin{array}{c}-0010\\t^3\end{array}\right]
\end{align*}

and

\begin{align*}\gamma(x_{\alpha_+}(t))=  \left[\begin{array}{c}0121\\kt^6\end{array}\right]
 \left[\begin{array}{c}1111\\-kt^3\end{array}\right]
   \left[\begin{array}{c}0011\\lt^2\end{array}\right]
    \left[\begin{array}{c}1231\\-lt\end{array}\right]
 \left[\begin{array}{c}0122\\f_1(t)\end{array}\right]
 \left[\begin{array}{c}1122\\f_2(t)\end{array}\right]\\
 \left[\begin{array}{c}1222\\f_3(t)\end{array}\right]
 \left[\begin{array}{c}1232\\f_4(t)\end{array}\right]
 \left[\begin{array}{c}1242\\f_5(t)\end{array}\right]
 \left[\begin{array}{c}1342\\f_6(t)\end{array}\right]
 \left[\begin{array}{c}2342\\f_7(t)\end{array}\right] \end{align*}.
 
 That $(x_{\alpha_+}(t)\gamma(x_{\alpha_+}(t)))^3=1$ implies that $f_6(t)=-klt^{-5}$ and $f_3(t)=l^2t^{-3}$.
 
 Then the identity $(x_{\alpha_+}(t)\gamma(x_{\alpha_+}(t)))(x_{\alpha_+}(u)\gamma(x_{\alpha_+}(u)))=(x_{\alpha_+}(t+u)\gamma(x_{\alpha_+}(t+u)))$ yields that the equation $l^2t^6-l^2t^3u^3-l^2u^3=0$ holds for all $t$ and $u$. Therefore we must have $l=0$.
 
 Lastly we find that if $l=0$, then we can find a lift for this complement for any $k$. One such is given by setting $f_i$ identically zero. Then we check the cocycle condition from \ref{a1h1resb} holds.
 
 Again the action of the torus $T_1\leq Z(L)$ reduces this (one-dimensional collection) to up to one  non-$G$-cr subgroup up to $G$-conjugacy,
 
We claim that there exists such a non-$G$-cr subgroup in an $A_1C_3$ subsystem. This can be seen in two ways, either taking an explicit description of the generators $x_+(t)$ and $x_-(t)$ and showing by direct calculation that the roots lie in an $A_1C_3$, or alternatively by arguing as before that it is the image of the representation $(1,(2\otimes 1)^{[1]})=(1,T(3)^{[1]})$ of $A_1$ in $A_1C_3$. 

Again, this subgroup is listed in the statement of \ref{a1sinchar3} and so the analysis for this case is concluded.

\subsubsection{(r,s,t)=(0,s,s-1), $s>2$, (1,s,0) $s>3$,\ (1,3,0)}\label{1s0}

These cases can be seen to be similar to \ref{021}. We have $\mathbb V\cong H^1(\bar X,Q/Q(2))\cong K$, in all cases and we calculate using {\tt F4gens} that the map $H^1(\bar X,Q)\to H^1(\bar X,Q/Q(2))$ is surjective. Then by \ref{onedimh1} we get one non-$G$-cr complement for each case. (The case $s=3$ should be dealt with separately as $H^2(\bar X,Q)\neq 0$. In spite of this, we find that $\rho$ is everywhere defined and there is again one non-$G$-cr complement.)

In the $(0,s,s-1)$ case, this subgroup is described as the representation $(1,T(3)^{[s-1]})$ of $A_1$ in $A_1C_3$. 

In the $(1,s,0)$ case, with $s>3$ or $s=3$. this subgroup is described as the representation $(1^{[s]},T(3))$ of $A_1$ in $A_1C_3$.

\subsubsection{(r,s,t)=(0,1,1), (0,s,1) $s>2$, (0,s-1,s) $s>3$}\label{0s1}

For these cases one uses again similar arguments and calculations as in \ref{021}. This time one finds that $\rho$ is nowhere defined, i.e. the map $H^1(\bar X,Q)\to H^1(\bar X,Q/Q(2))$ is zero, so that no cocycles $H^1(X,Q/Q(2))$ will lift. Thus there is just one conjugacy class of complements to $Q$ in $\bar XQ$ and so and all complements are $G$-cr.

\subsubsection{(r,s,t)=(0,1,0)}

We have $Q(1)/Q(2)\cong 1\otimes 1^{[1]}+2+0$ and $Q(2)\cong 1\otimes 1^{[1]}+2$.

\begin{claim} There is an isomorphism $H^1(A_2,11)\to H^1(A_1,1\otimes 1^{[1]}\oplus 2)$ given by restriction to the irreducible $A_1$ subgroup in $A_2$. \end{claim}
{\it Proof of claim: }The standard module $10$ for the $A_2$ restricts to the $A_1$ as a high weight $2$ module. Now $10\otimes 01=00/11/00=T(11)$ contains the non-trivial quotient $E_1:=00/11$ giving an extension of $11$ by the trivial module, hence a SES \[0\to 11\to E_1\to K\to 0.\] Now $(10\otimes 01)\downarrow A_1=2\otimes 2=0/(1\otimes 1^{[1]})/0+2=T(4)+2$. and so we have a SES \[0\to 1\otimes 1^{[1]}+2\to E_2\to K\to 0.\] Hence restriction to the $A_1\leq A_2$ sends a non trivial extension of the module $11$ for $A_2$ by a trivial module to a non-trivial extension of the module $1\otimes 1^{[1]}+2$ for $A_1$ by the trivial. Since the map in the claim is then a non-zero map of one-dimensional vector spaces, it must be an isomorphism. This proves the claim.

Using the claim above we see that the hypotheses of \ref{h1res} go through, viz. \begin{align*}H^2(A_2,Q(2))=0&\to 0=H^2(\bar X,1\otimes 1^{[1]}),\\
H^1(A_2,Q(2))&\to H^1(\bar X,Q(2)),\\
H^1(A_2, Q(1)/Q(2))&\to H^1(\bar X,Q(1)/Q(2)),\\
H^0(A_2,Q(1)/Q(2))&\to H^0(\bar X,Q(1)/Q(2))\end{align*} are all isomorphisms by restriction and thus so is $H^1(A_2,Q)\to H^1(\bar X,Q)$. Now it follows that any non-$G$-cr complement to $\bar X$ is in some non-$G$-cr subgroup $A_2$; thus, by \ref{a2sinparabs} and \ref{a2s}, in a subsystem $A_2\tilde A_2$, irreducible in each factor or in a $B_4$ via $V_9\downarrow X=T(4)+2$.

From \ref{tor} we find that the actions of these $A_1$s on the $V_{26}$ are distinct. Thus we have found two more of the subgroups listed in the statement of \ref{a1schar3}.

\subsubsection{(r,s,t)=(r,r+1,0) (r$>$1), (0,1,t) (t$>$2)}

Here $H^1(\bar X,Q(1)/Q(2))=0$, $H^1(\bar X, Q(2))=K$ and $H^0(\bar X,Q(1)/Q(2))=0$. Thus by \ref{nofixedpoints} $H^1(\bar X,Q)=K$ and so by \ref{onedimh1} there is one non-$G$-cr complement $X$ to $Q$. Since $L'Q(2)\leq B_4$, we have that $X$ is in $B_4$ by \ref{subsystemlift}.

Now we claim that for case $(r,r+1,0)$, (resp. case $(0,1,t)$), $X$  is determined by the action on the natural module for $B_4$, $V_9\downarrow X=0/1^{[r]}\otimes 1^{[r+1]}/0+2=T(4)^{[r]}+2$, (resp. $T(4)+2^{[t]}$). Certainly each is in a $B_3$ parabolic of $B_4$ with $\bar X$ as the image under the projection to the Levi $B_3$. Hence, for each $r$ and $t$, each subgroup $X$ is a complement to $Q$ in $\bar XQ$ and is not $G$-conjugate to $\bar X$ since it has a different action on $V_{26}$ by \ref{tor}.

We check that these subgroups are included in the statement of \ref{a1schar3}. This concludes this case.

\subsubsection{(r,s,t)=(0,1,2)}

Here $H^1(\bar X,Q/Q(2)\cong K\cong H^1(\bar X,Q(2))$, and so $\mathbb V=K^2$. Using the methods of \ref{021} one finds that $\rho$ is undefined on $[k_1,k_2]$ with $k_1\neq 0$. In other words, the non-trivial complements to $Q(1)/Q(2)$ in $\bar XQ(1)/Q(2)$ do not lift.

Now $H^1(\bar X,Q(2))=K=H^1(\bar X,Q)$ and so we get one class of non-$G$-cr complements in a $B_4$ by \ref{subsystemlift} and \ref{onedimh1}. It is easy to see that this is given by $V_9\downarrow X=T(4)+2^{[2]}$.

We check that this subgroup is included in the statement of \ref{a1schar3}.

\subsubsection{(r,s,t)=(1,2,0)}\label{120}

Here $Q(1)/Q(2)\cong 1\otimes 1^{[1]}+1\otimes 1^{[2]}$ and $Q(2)\cong 1^{[1]}\otimes1^{[2]}+2$.

We have $H^1(\bar X,Q(1)/Q(2))=H^1(\bar X,Q(2))=K$, and $H^0(\bar X,Q/Q(2))= H^2(\bar X,Q(2))=0$. Thus by \ref{nofixedpoints}, $H^1(\bar X,Q)=K^2$.

Let $X$ and $Y$ be two complements determined by different elements of $H^1(\bar X,Q)$. From \ref{B3parabconj} we have that $X$ and $Y$ are $G$-conjugate if and only if they are $N_L(\bar X)Q$-conjugate. Now $N_L(\bar X)=\bar X.\langle\sigma\rangle.Z(L)$ where $\sigma$ is an involution acting on $O_4\times O_3\leq B_3$ by diag $(-1,1)$ on the two orthogonal factors. 

Now $Z(L)$  acts on $Q(1)/Q(2)$ and $Q(2)$ as $K^\times$, so we have that up to $G$-conjugacy, the set of complements to $Q$ is isomorphic to the variety $\mathbb P^1(K)$.

Now we observe that if $k_1,k_2\neq 0$ then the complement corresponding to $[k_1,k_2]$ is in no maximal reductive subgroup of $G$. A calculation with the MeatAxe shows that there are only two possibilities for the restriction $V_{26}\downarrow X$:

\[\begin{array}{c}1\otimes 1^{[2]}\\\hline 1\otimes 1^{[2]}+1^{[1]}\otimes 1^{[2]}\\\hline 1\otimes 1^{[2]}\end{array}+\begin{array}{c}0\\\hline 1\otimes 1^{[1]}\\\hline 0\end{array}+2+0\]

or

\[\begin{array}{c}1\otimes 1^{[2]}\\\hline 1\otimes 1^{[2]}\\\hline 1\otimes 1^{[2]}\end{array}+\begin{array}{c}0\\\hline 1\otimes 1^{[1]}+1^{[1]}\otimes 1^{[2]}\\\hline 0\end{array}+2+0\]

Now it is easy to see from \ref{tor} that $X$ is not preserving the same dimensions of direct summands as any of the maximal reductive subgroups of $G$. For instance, the $B_4$ preserves a $1+9+16$ direct sum decomposition, but $X$ does not, thus $X\nleq B_4$. The other cases are similar.

We have found the subgroups claimed in the statement of \ref{a1schar3} and proved that they are in not properly contained in any propoer reductive subgroup of $G$. Thus we have concluded this case 

Since all cases of $(r,s,t)$ have been analysed, this completes the proof of \ref{a1schar3}.

\subsection{Subgroups of type $A_1$, $p=2$}\label{a1p2}

We now begin the classification of subgroups of type $A_1$ in characteristic $2$. Firstly we need not consider any parabolics with the Levi having a factor of type $A_2$ since $A_1$ does not embed $A_2$-irreducibly in this characteristic. Also, since there is a graph automorphism of $F_4$ when $p=2$ most of our analysis will be done in the cases $L'=B_3$, $B_2$, $A_1$ and $A_1\tilde A_1$; we will deduce results for $C_3$ using the graph automorphism.

We warn the reader that we will jump between calculations in $B_n$ and $C_n$ at will, using the exceptional isogenies $\phi:B_n'\to C_n$ and $\psi:C_n\to B_n$ where $\phi\circ\psi$, (resp. $\psi\circ\phi$) induces a standard Frobenius morphism on $C_n$ (resp. $B_n$). The reason for this is that it is often possible to distinguish subgroups up to conjugacy in $GL(V_7)$ where $V_7$ is the natural (Weyl) module for $B_3$, yet it is easier to establish $C_3$-irreducibility by considering restrictions of its natural (irreducible) module $V_6$.

\subsubsection{$L'=B_3, C_3$ or $B_2$}\label{a1sp2inb3c3b2}

We show in this section that the classification of $A_1$ subgroups minimally contained in parabolics of type $B_2$, $B_3$ and $C_3$ is essentially the same as the classification of $A_1^2$, $A_1^3$ and $A_1^3$ subgroups minimally contained in those parabolic subgroups, respectively.

We give here the means of passing from a classification of the latter to the former. An actual classification of the latter is done later in the thesis.

\begin{lemma} \label{a1sina1a1a1s}Let $X$ be a non-$G$-cr subgroup of type $A_1$ contained subject to minimality in a $B_3$-parabolic (resp. $C_3$-parabolic) of $G$. Then $X$ is contained in exactly one non-$G$-cr subgroup $Y$ of type $A_1^3$. Moreover, \begin{enumerate}\item If $Y$ is a complement to $Q$ in $\tilde A_1^3Q$ (resp. $A_1^3Q$) then $X\hookrightarrow Y$ by $(V_2,V_2,V_2)\downarrow X=(1^{[r]},1^{[s]},1^{[t]})$ for some distinct $r,\ s,\ t$.
\item If $Y$ is a complement to $Q$ in $A_1^2\tilde A_1Q$ (resp. $\tilde A_1^2A_1Q$) then $X\hookrightarrow Y$ by $(V_2,V_2,V_2)\downarrow X=(1^{[r]},1^{[s]},1^{[t]})$ for some distinct $r$ and $s$.\end{enumerate}.
Hence the restriction of $V_{26}\downarrow X$ can be computed by restriction of some entry in \ref{tor}.\end{lemma}
\begin{proof}The statement for $C_3$ follows from that for $B_3$ by applying the graph automorphism. Thus we may assume that $X$ is contained in the $B_3$-parabolic of $G$.

From \ref{a1sinb3} we know that $\bar X\leq \tilde A_1^3$ or $A_1^2\tilde A_1$ satisfying the respective constraints on $r,\ s$ and $t$ in the lemma.  

Using the K\"unneth formula, \ref{kunneth}, we claim that under these conditions, the hypotheses of \ref{h1res} hold, and so $H^1(\bar X,Q)\cong H^1(\tilde A_1^3,Q)$ or $H^1(\bar X,Q)\cong H^1(A_1^2\tilde A_1,Q)$. For example take $\bar X\leq A_1^2\tilde A_1$ with twists $(r,s,t)$. Then $Q(2)\downarrow A_1^2\tilde A_1=(1,1,0)+(0,0,1^{[1]}/0)$ and $Q(2)\downarrow \bar X=1^{[r]}\otimes 1^{[s]}+1^{[t+1]}/0$. We see using \ref{kunneth} that $H^1(\bar X,1^{[r]}\otimes 1^{[s]})\cong H^1(A_1^2\tilde A_1,(1,1,0))\ (\cong  0)$ and that $H^1(\bar X,1^{[t+1]}/0)\cong H^1(A_1^2\tilde A_1,(0,0,1^{[1]}/0))\ (\cong K)$. So $H^1(A_1^2\tilde A_1,Q(2))\cong H^1(\bar X,Q(2))$. Similarly one has $H^i(A_1^2\tilde A_1,Q(2))\cong H^i(\bar X,Q(2))$ and $H^i(A_1^2\tilde A_1,Q/Q(2))\cong H^i(\bar X,Q/Q(2))$ for $i=0,1,2$ as required. The same statement holds when $\bar X\leq \tilde A_1^3$. Thus we have the claim.

From the isomorphisms $H^1(\bar X,Q)\cong H^1(\tilde A_1^3,Q)$ or $H^1(\bar X,Q)\cong H^1(A_1^2\tilde A_1,Q)$, it follows that a complement $X$ is contained in exactly one complement $Y\cong A_1^3$ to $Q$ in $A_1^3Q$ up to $Q$-conjugacy. But by \ref{b3imp} we see that the same $G$-conjugacy is induced amongst $Q$-complements in $A_1^3Q$ and $\bar XQ$; viz. by $QZ(L)$. Since $Z(L)$ pulls back to the same action on $\mathbb V$ in each case, we have that $X$ is contained in exactly one complement $Y\cong A_1^3$.

The fact that one can compute the restrictions of $V_{26}\downarrow X$ from $V_{26}\downarrow Y$ follows since the hypotheses of \ref{samesubs} hold for the subgroups as given. Since we have included the restrictions $V_{26}\downarrow Y$  in \ref{tor}, the last statement of the lemma is justified.
\end{proof}

\begin{remark}We establish a complete classification of non-$G$-cr subgroups of type $A_1^3$ in \S\ref{a1a1a1sec} later.\end{remark}
 
\begin{lemma}\label{a1sinb2}  $X$ be a non-$G$-cr subgroup of type $A_1$ contained subject to minimality in a $B_2$-parabolic subgroup of $G$. Then $X$ is contained in exactly one non-$G$-cr subgroup $Y$ of type $A_1^2$. Moreover, if $Y$ is a complement to $Q$ in $A_1^2Q$ or $\tilde A_1^2Q$ then $X\hookrightarrow Y$ by $(V_2,V_2)\downarrow X=(1^{[r]},1^{[s]})$ with $r\neq s$.\end{lemma}
\begin{proof} We know from \ref{lir} that $\bar X\leq \tilde A_1^2$ or $A_1^2$, by $(V_2,V_2)\downarrow X=(1^{[r]},1^{[s]})$ with $r\neq s$. Assume the former first.

By the K\"unneth formula, \ref{kunneth}, we check that under these conditions, the hypotheses of \ref{h1res} hold, and so $H^1(X,Q)\cong H^1(A_1^2,Q)$ and so $X$ is contained in exactly one complement $Y\cong A_1^2$ to $Q$ in $A_1^2Q$ up to $Q$-conjugacy. Moreover if $\mathbb V_{\bar X}$ is the pre-cocycle space of $\bar X$ and $\mathbb V_{A_1^2}$ is the pre-cocycle space of $A_1^2$ then $\mathbb V_{\bar X}\cong \mathbb V_{A_1^2}$. 

Now, If $X'=X^g$ is a complement to $Q$ which is $G$-conjugate to $X$ containing a maximal torus of $X$, then we claim from \ref{nomoreconj2} that we may take $g\in C_G(L')$. To prove the claim observe that in the notation of \ref{nomoreconj2} we have $N_L(\bar X)=\bar X$ giving $M=1$, while $C_G(L')\cong B_2 Z(L)=\langle x_{0122}(t),x_{1110}(t),1,n_{1110},n_{0122},Z(L)\rangle$, with $C_Q(H)=\langle x_{1232}(t),x_{2342}(t):t\in K\rangle\leq C_Q(H)$. 

Then it is easy to calculate that conjugacy by these elements pulls back to an action on $\mathbb V_{\bar X}$ which is isomorphic to that on $\mathbb V_{A_1^2}$. (We calculate these actions explicitly in \ref{a1a1inb2} below.) Thus we have that $X$ is contained in exactly one complement $Y\cong A_1^2$.

The statement for $\tilde A_1^2Q$ follows from the statement for $A_1^2Q$ by applying the graph automorphism of $G$.\end{proof}
\begin{remark}We establish a complete classification of non-$G$-cr subgroups of type $A_1^2$ in $B_2$ parabolics in \S\ref{a1a1s} later.\end{remark}

\subsubsection{Subgroups of type $A_1$, $L'=A_1$}\label{seca1ina1char2}

\begin{prop}\label{a1ina1char2} Let $X$ be a non-$G$-cr subgroup of type $A_1$ contained in a long $A_1$ parabolic of $G$. Then $X$ is in a $B_4$ subsystem subgroup of $G$ and $X$ is conjugate to precisely one of the following $4$ subgroups where we give the action of $X$ on the natural module $V_9$ for $B_4$:
\begin{center}\begin{tabular}{c|c|c}\hline Item & Description & $V_9\downarrow X_\bold k$\\
\hline
1 &$X\leq A_2$ with $V_3\downarrow X=Sym^2(V)$ &$0/2+ 2/0+ 0^3$\\
2 &$X\leq \tilde A_1^2\leq B_2$ as $x\mapsto (x,x)$ & $2+2/0+0^4$\\
3 & $X\leq A_3$ as $V_4\downarrow X=T(2)$ & $0/2+2/0+0^3$\\
4 & $X\leq A_1^4$ as $x\mapsto (x,x,x,x)$ & $T(2)+ T(2)+0$.

\end{tabular}\end{center}
In the above we have written $2$ in place of the 2-dimensional irreducible modules $L(2)\cong 1^{[1]}$.

Let $X$ be a subgroup of type $A_1$ contained in a short $\tilde A_1$ parabolic of $G$. Then $X$ is in a $C_4$ subsystem subgroup of $G$ and $X$ is conjugate to precisely one of the following $4$ subgroups where we give the action of $X$ on the natural module $V_8$ for $C_4$:
\begin{center}\begin{tabular}{c|c|c}\hline Item & Description & $V_8\downarrow X_\bold k$\\
\hline
1' &$X\leq \tilde A_2$ with $V_3\downarrow X=Sym^2(V)$ &$0/2+ 2/0+ 0^2$\\
2' &$X\leq A_1^2\leq C_2$ as $x\mapsto (x,x)$ & $2+2+0^4$\\
3' & $X\leq A_1^2\leq B_2\leq \tilde D_3$ as $x\mapsto (x,x)$ & $2/0+0/2+0^2$\\
4' & $X\leq \tilde A_1^4$ as $x\mapsto (x,x,x,x)$ & $T(2)+ T(2)$.
\end{tabular}\end{center}
\end{prop}

{\it Proof of Proposition \ref{a1ina1char2}:}

It clearly suffices to prove the statement about long $A_1$ parabolics and the other will follow by applying the graph automorphism of $G$.

\begin{lemma}\label{a1ina1gens}Let $\bold k= [k_1,\dots,k_7]$ vary over all elements of $K^7$. We claim the groups $X_{\bold k}$ generated by the following elements constitute a complete collection of complements to the unipotent radical $Q$ in $\bar XQ\leq P$ where $P$ is the $A_1$-parabolic of $G$ corresponding to node 1:

{\footnotesize \begin{align*}
x_+(t):&=\left[\begin{array}{c}1000\\t^2\end{array}\right]
\left[\begin{array}{c}1100\\k_1t\end{array}\right]
\left[\begin{array}{c}1110\\k_2t\end{array}\right]
\left[\begin{array}{c}1111\\k_3t\end{array}\right]
\left[\begin{array}{c}1120\\k_4t\end{array}\right]
\left[\begin{array}{c}1121\\k_5t\end{array}\right]
\left[\begin{array}{c}1122\\k_6t\end{array}\right]
\left[\begin{array}{c}2342\\k_7t\end{array}\right]\\
x_-(t):&=\left[\begin{array}{c}-1000\\t^2\end{array}\right]
\left[\begin{array}{c}0100\\k_1t\end{array}\right]
\left[\begin{array}{c}0110\\k_2t\end{array}\right]
\left[\begin{array}{c}0111\\k_3t\end{array}\right]
\left[\begin{array}{c}0120\\k_4t\end{array}\right]
\left[\begin{array}{c}0121\\k_5t\end{array}\right]
\left[\begin{array}{c}0122\\k_6t\end{array}\right]
\left[\begin{array}{c}1342\\\chi t+k_7t\end{array}\right]
\end{align*}}
for all $t\in K$ where $\chi=k_1k_4k_6+k_1k_5^2+k_2^2k_6+k_3^2k_4$.
\end{lemma}
\begin{proof}
Let $\bar X$ be generated by $x_{\pm 1000}(t^{2^r})$. From \ref{levelsforx}, the pre-cocycle space $\mathbb V$ is $0$ if $r=0$ and $K^7$ if $r>0$. Without loss of generality, we may assume that $r=1$, since for any $r\geq 1$ $\mathbb V\cong K^7$ and so one can obtain all complements for the case $r>1$ from the case $r=1$ by applying Frobenius to those obtained from the latter.  Referring to \ref{levelsforx} we see that there are no potential blockers and so $\bold k= [k_1,\dots,k_7]$ parameterises the elements of $H^1(\bar X,Q)$ via $\rho$.

Using the $26$-dimensional representation of $G$, one can input  $x\pm(t)$ in terms of {\tt F4gens} into GAP and check the relations \ref{a1rels} to see that for each $\bold k= [k_1,\dots,k_7]$ the group $X_\bold k$ generated by the given elements is isomorphic to $A_1$. It is also evidently contained in $\bar X Q$ so it is a complement to $Q$ in $\bar X Q$. Indeed by direct calculation one finds that the torus $H=\{h_+ (c):c\in K^\times\}$ of $X_\bold k$ is the same as the torus $\{h_{1000}(c):c\in K^\times\}$ of $\bar X$.

We claim that the set of cocycles $\gamma_\bold k$ corresponding to the $X_\bold k$ has been constructed compatibly with the map $\rho$. It is enough to see this on restriction to the Borel subgroup $\langle H, x_+(t)\rangle$ of $\bar X$ by \ref{h1resb}. But now we see that a valid construction of the map $\rho_B:\mathbb V=K^7\to H^1(B,Q)$ can be obtained by letting each $\Gamma_i'=\hat{\Gamma_{i-1}}$; that is, for each preimage under the map $H^1(B,Q/Q(i)\to H^1(B,Q/Q(i+1))$ taking the naive `canonical lift' (\ref{clift}) of each cocycle.

Let me give an example: we know that $H^1(B,Q/Q(2)=K=H^1(B,Q(2)/Q(3))$. We give a map $\rho_1:K\to H^1(B,Q/Q(2))$ in the obvious way, that is $k\mapsto [\gamma_k]$ where $\gamma_k(x_{1000}(t^2))=x_{1100}(kt)$ and $\gamma_k(H)=1$. Now we take preimages of each element of $H^1(B,Q/Q(2))\to H^1(B,Q/Q(3))$ simply by setting $\Gamma_k(x_{1000}(t^2))=\widehat{\gamma_k(x_{1000}(t^2))}Q(3)/Q(3)=x_{1100}(kt)$. Then we get representatives of all cocycle classes in $H^1(B,Q/Q(3))$ by $\Gamma_{k,l}=\Gamma_k\mu_l$ where the $\mu_l$ are cocycles running over all classes in $H^1(B,Q(2)/Q(3))$. Specifically, $\Gamma_{k,l}(x_{1000}(t^2))=x_{1100}(kt)x_{1110}(lt)$. This is clearly a valid construction of the map $\rho$.\end{proof}

In the following we will identify $X_\bold k$ and $\gamma_\bold k$ with $\bold k\in H^1(\bar X,Q)$. Where $X_\bold k^g=X_{\bold k'}$ for some $\bold k'$ we write $\bold k^g=\bold k' $
\begin{lemma}\label{someconj}
In the following table we indicate the action of various elements of $G$ by conjugacy on the set of cocycles defined above:
\begin{center}\begin{tabular}{c|l}$g$ & $\bold k^g$\\\hline
$n_{0010}(1)$ & $[k_4,k_2,k_5,k_1,k_3,k_6, k_7]$\\
$n_{0001}(1)$ & $[k_1,k_3,k_2,k_6,k_5,k_4, k_7]$\\
$x_{0010}(c)$ & $[k_1,k_2+ck_1,k_3,k_4+c^2k_1,k_5+ck_3,k_6,k_7]$\\
$x_{0001}(c)$ & $[k_1,k_2,k_3+ck_2,k_4,k_5+ck_4,k_6+c^2k_4,k_7]$\\
$x_{0011}(c)$ & $[k_1,k_2,k_3+ck_1,k_4,k_5+ck_2,k_6+c^2k_1,k_7]$\\
$x_{1242}(c)$ & $[k_1,k_2,k_3,k_4,k_5,k_6,k_7+ck_1]$\\
$h_{0010}(c)$ & $[k_1c^{-2},k_2,k_3c^{-1},k_4c^2,k_5c,k_6,k_7]$\\
$h_{0001}(c)$ & $[k_1,k_2c^{-1},k_3c,k_4c^{-2},k_5,k_6c^2,k_7]$\\
$h_{1242}(c)$ & $[k_1c^{-1},k_2,k_3,k_4c,k_5c,k_6c,k_7c]$\\
\end{tabular}\end{center}
Additionally, \[[0,k_2,0,0,0,0,k_7]^{n_{1242}(1)}=[k_7,k_2,0,0,0,0,0].\]
\end{lemma}
\begin{proof} Notice that each element $g\in C_G(\bar X)$. Using the generators given in \ref{a1ina1gens} this is a routine application of the commutator formula and standard Coxeter group calculations. It can also be verified using the matrix representation of $G$.\end{proof}

\begin{lemma}\label{downtoafew}
Up to conjugacy by the above elements, \begin{align*} \bold k=&[1,0,0,0,0,0,0],\\&[0,1,0,0,0,0,0],\\&[1,0,0,1,0,0,0],\text{ or }\\&[1,0,0,1,0,1,0]\end{align*}\end{lemma}
\begin{proof}
Fix $\bold k=[k_1,\dots,k_7]\neq [0,\dots,0]$. We successively replace $\bold k$ with $\bold k^g$ for $g$ an element described in \ref{someconj} to reduce $\bold k$ to one of the above possibilities.

First assume one of $k_2,\ k_3$ or $k_5$ is non-zero; we see that we can replace $\bold k$ by a conjugate so that $k_3=k_5=0$: If $k_5\neq 0$ then replacing $\bold k$ by its conjugate $\bold k^{n_{0010}}$ we may assume $k_3\neq 0$. If $k_3\neq 0$ then replacing $\bold k$ by $\bold k^{n_{0001}}$ we may assume $k_2\neq 0$. Now by replacing $\bold k$ by $\bold k^{x_{0001}(k_3/k_2)}$ we may assume that $k_3=0$. Similarly by replacing $\bold k$ with $\bold k^{x_{0011}(k_5/k_2)}$ we may assume that $k_5=0$. The result is that $\bold k=[*,*,0,*,0,*,*]$.

Now also assume that one of $k_1$, $k_4$ or $k_6$ is non-zero; we claim that we can replace $\bold k$ by a conjugate so that $k_2=k_3=k_5=0$. If $k_1\neq 0$, then replace $\bold k$ by its conjugate by $x_{0010}(c)$ to set $k_2=0$. If $k_4\neq 0$, then we can replace $\bold k$ by its conjugate by $n_{0001}$ to have $k_1\neq 0$ and $k_2\neq 0$. Similarly conjugation by $x_{0010}(c)$ gives $k_2=0$. Now if $\bold k$ is of the form $[0,*,0,0,0,*,*]$ we replace it by its conjugate by $n_{0010}n_{0011}$ to have $\bold k=[*,0,0,0,*,0,*]$. Now replace $\bold k$ by its conjugate by $x_{0010}(c)$ for appropriate $c$ so that we have $\bold k=[*,*,0,*,*,0,*]$ with the first, second, fourth and fifth entry non-zero. Then replace by its conjugate by $x_{0001}(c)$ to have $\bold k=[*,*,0,*,0,0,*]$ and by its conjugate by $x_{0010}(c)$ again to leave $\bold k=[*,0,0,*,0,0,*]$ as required. 

Additionally we have that if one of $k_1,\ k_4$ and $k_6$ is non-zero, then we may assume it is $k_1$. If two are non-zero then we may assume it is $k_1$ and $k_4$. Lastly, if $k_1\neq 0$ replacing $\bold k$ by its conjugate by $x_{1242}(k_1/k_7)$ we can assume $k_7=0$. Replacing $\bold k$ by its conjugate $\bold k^{x_{0010}(k_2/k_1)}$ we may assume that $k_2=0$. 

From the previous two paragraphs we deduce: if one of $k_1,k_4$ or $k_6$ is non-zero, then we may assume that $k_2=k_3=k_5=k_7=0$.

If $k_1=k_4=k_6=0$ then the argument of the first paragraph shows that we may additionally assume $k_3=k_5=0$. Now if $k_7\neq 0$ then replacing $\bold k$ by $\bold k^{n_{1242}}$ gives $k_1\neq 0$, a contradiction.

Thus we have reduced $\bold k$ to one from the following list:
\begin{align*} \bold k=&[k_1,0,0,0,0,0,0],\\&[0,k_2,0,0,0,0,0],\\&[k_1,0,0,k_4,0,0,0],\\&[k_1,0,0,k_4,0,k_6,0]\end{align*}

Now it is easy to see that conjugating by the various elements $h_\alpha(c)$ with suitable $c\in K^\times$ we can replace $\bold k$ by one of those in the statement of the lemma.
\end{proof}

We must now perform a similar calculation for the $A_1$-parabolic of $G$ corresponding to node 2.

\begin{lemma}Any complement to $\bar X=A_1$ in the $A_1$-parabolic of $G$ corresponding to node 2 is $G$-conjugate to a subgroup of the parabolic of $G$ corresponding to node 1.\end{lemma}
\begin{proof} A similar argument to \ref{a1ina1gens} shows that as $\bold k$ ranges over all values of $K^7$ that the groups $X_\bold k$ generated by the following form a representative collection of complements:

{\footnotesize \begin{align*}
x_+(t):&=\left[\begin{array}{c}0100\\t^2\end{array}\right]
\left[\begin{array}{c}1100\\k_1t\end{array}\right]
\left[\begin{array}{c}0110\\k_2t\end{array}\right]
\left[\begin{array}{c}0111\\k_3t\end{array}\right]
\left[\begin{array}{c}1220\\k_4t\end{array}\right]
\left[\begin{array}{c}1221\\k_5t\end{array}\right]
 \left[\begin{array}{c}1222\\k_6t\end{array}\right]
\left[\begin{array}{c}1342\\k_7t\end{array}\right]\\
x_-(t):&=\left[\begin{array}{c}-0100\\t^2\end{array}\right]
\left[\begin{array}{c}1000\\k_1t\end{array}\right]
\left[\begin{array}{c}0010\\k_2t\end{array}\right]
\left[\begin{array}{c}0011\\k_3t\end{array}\right]
\left[\begin{array}{c}1120\\k_4t+k_1k_2^2t\end{array}\right]
\left[\begin{array}{c}1121\\k_5t+ k_1k_2k_3t \end{array}\right]\\
&\hspace{170pt}
\left[\begin{array}{c}1122+k_1k_3^2t\\k_6t\end{array}\right]
\left[\begin{array}{c}1242\\k_7t+(k_2^2k_6+k_3^2k_4)t\end{array}\right]
\end{align*}}
for all $t\in K$.

Observe that if $k_1=0$ then conjugating any $X_\bold k$ by the element $n_{1000}n_{0100}$ gives a subgroup of $P_1$. Assume that $\bold k$ is not conjugate to one with $k_1=0$.

Observe the following conjugacy amongst the $\bold k$:

\begin{center}\begin{tabular}{c|l}$g$ & $\bold k^g$\\\hline
$x_{0001}(c)$ & $[k_1,k_2,k_3+ck_2,k_4,k_5+ck_4,k_6+c^2k_4,k_7]$\\
$n_{0001}(c)$ & $[k_1,k_3,k_2,k_6,k_5,k_4,k_7]$\\
$x_{0120}(1)$ & $[k_1,k_2,k_3,k_4+ck_1,k_5,k_6, k_7+ck_6]$
\end{tabular}\end{center}
Additionally, \[[k_1,0,k_3,k_4,k_5,k_6,k_7]^{n_{0010}n_{0100}n_{0010}}=[k_4,0,k_3,k_1,k_5,k_7,k_6]\]

If $k_2\neq 0$ then replacing $\bold k$ by its conjugate by $x_{0001}(k_3/k_2)$ we may assume that $k_3=0$. Now, replacing $\bold k$ by its conjugate by $n_{0001}$ we may assume $k_2=0$.

Replace $\bold k$ by its conjugate by $x_{0120}(k_4/k_1)$ to assume that $k_4=0$ (retaining also  $k_2=0$). Finally replacing $\bold k$ by its conjugate by $n_{0010}n_{0100}n_{0010}$ we deduce that $k_1=0$, a contradiction.
\end{proof}

To finish the proof of Proposition \ref{a1ina1char2} we must just observe that each of the values of $\bold k$ in Lemma \ref{downtoafew} is conjugate to just one of the items 1-4 in the statement of the Lemma.

To do this, note that Proposition \ref{a1ina1char2} describes 4 subgroups which are in long $A_1$ parabolics of $B_4$.

Now examining the restrictions $V_{26}\downarrow X$ for each of these in \ref{tor}, we see that there can be no conjugacy between any of the items 1-4 in the statement of the proposition, except possibly we may have item 2 conjugate to item 3. We will show that item 2 is not in an $A_3$, and then we will be done since item 3 is by definition. 

From the possibilities in \ref{downtoafew}, one can see that there must be at least two non-$G$-cr subgroups in $A_3$s by looking at the roots involved in the generators. We claim there are precisely two:

Any subgroup $X$ contained in an $A_3$ is determined up to conjugacy in the $A_3$ by the restriction of the natural module $V_4=L_{A_3}(100)\downarrow X$. The composition factors are $2|0^2$ and so we have either $V_4\downarrow X=2+0^2$, $V_4\downarrow X=0/2/0=T(2)$ or $V_4\downarrow X=2/0+0$ or $V_4\downarrow X=0/2+0$. The first of these gives a Levi $A_1$ and so is $G$-cr. The last two of these are in $A_2$s and are conjugate by the element of the Weyl group which induces a graph automorphism on this $A_2$. Thus there are at most two; so there must be exactly two. Since $V_4\downarrow X=2/0=0$ gives item 1 in the proposition, and this is not conjugate to item 3, which is defined by $V_4\downarrow X=T(2)$, we are done.

This concludes the proof of Proposition \ref{a1ina1char2}.

\begin{remark} One can in fact see that item i in Proposition \ref{a1ina1char2} corresponds to item i in Lemma \ref{downtoafew} by further examination of the roots involved.\end{remark}

\begin{remark} It is worth noticing some slightly unusual facts about conjugacy. Most of the time two complements in a $B_4$ subsystem will be conjugate in $G$ if and only if they are conjugate in $B_4$. However, the subgroup $X=A_1\leq A_1^2\tilde A_1$ defined by $x\mapsto (x,x,x)$ has restriction $V_9\downarrow X=T(2)+2/0+0^2$. In $B_4$ then, it is not $B_4$-conjugate to the subgroup of item 3 of \ref{a1ina1char2} but it is in fact $G$-conjugate.

Also we notice that item 1 in the lemma is $GL(V_9)$-conjugate to item 3 but it is not $G$-conjugate, so in particular it is not $B_4$-conjugate.
\end{remark}

\subsubsection{Subgroups of type $A_1$, $L'=A_1\tilde A_1$}\label{seca1ina1a1char2}
 
The situation in the three different $A_1\tilde A_1$ parabolics is complicated. Many from one parabolic are conjugate to those in another, while there are several infinite families in one parabolic which are not $G$-conjugate. 

Our method here is similar to \S \ref{seca1ina1char2}. That is, to get a list of candidates up to $G$-conjugacy in each parabolic, take the union of all of them, pare this down by using \ref{nomoreconj2} to achieve a list of distinct conjugacy classes and then, where possible, find reductive overgroups.

Ultimately this work results in a classification of conjugacy classes by elements of $\mathbb V=K^6$, but not in any easily stated form.

\begin{prop}\label{a1ina1a1char2} Let $X$ be a non-$G$-cr subgroup of type $A_1$ contained subject to minimality in an $A_1\tilde A_1$-parabolic subgroup of $G$. Let $\bar X$ be the projection of $X$ to $A_1\tilde A_1$ such that the natural module $(V_2, V_2)\downarrow \bar X=(1^{[r]}, 1^{[s]})$.

Then $X$ is known up to conjugacy in $G$. Moreover, $X$ is conjugate to precisely one item in the table in \ref{classofa1a1s} below.
\end{prop}

{\it Proof of Proposition \ref{a1ina1a1char2}:}

We will concentrate initially on parabolic subgroups $P_{24}$ and $P_{23}$. Later we will use the graph automorphism to deduce the results for $P_{13}$.

By examining the list \ref{levelsforx} in combination with \ref{someh1s} we see that $\mathbb V=K^6$ if $P=P_{24}$ or $P_{23}$, provided $r\neq s+1$, and that $\mathbb V=K^5$ in each case if $r=s+1$. The next lemma shows that $\rho$ is not defined on any $[k_1,k_2\dots,k_6]\in\mathbb V$ where $k_1,k_2\neq 0$.

\begin{lemma}\label{a1a1blocking} Let $Q$ be unipotent radical of the parabolic $P_{24}$ or $P_{23}$, and let $r\neq s+1$. The map $H^1(\bar X,Q/Q(3))\stackrel{\sigma}{\to} H^1(\bar X, Q/Q(2))=K^2$ has image \[\{(x,y)\in H^1(\bar X, Q/Q(2)):xy=0\}.\]\end{lemma}
\begin{proof} We follow the same method as used previously in \S\ref{001}. 

For $P_{24}$, observe $Q/Q(2)\downarrow \bar X= 1^{[s]}+(1^{[r]}\otimes 1^{[s]})/1^{[r]}$. For $P_{23}$, $Q/Q(2)\downarrow \bar X=1^{[r]}+1^{[s]}$.

For $P_{24}$, $1^{[s]}$ is generated by $\langle \left[\begin{array}{c}0001\\k\end{array}\right],
\left[\begin{array}{c}0011\\k\end{array}\right] \rangle Q(2)$ as $k$ runs over $K$. The $1^{[r]}$ submodule is generated by $\langle  \left[\begin{array}{c}0110\\k\end{array}\right], \left[\begin{array}{c}1110\\k\end{array}\right]\rangle Q(2)$ as $k$ runs over $K$.

For $P_{23}$, $1^{[r]}$ is $\langle  \left[\begin{array}{c}0100\\t\end{array}\right], \left[\begin{array}{c}1100\\t\end{array}\right]\rangle Q(2)$ as $t$ runs over $K$.  $1^{[s]}$ is $\langle \left[\begin{array}{c}0010\\t\end{array}\right],
\left[\begin{array}{c}0011\\t\end{array}\right] \rangle Q(2)$ as $t$ runs over $K$. 

By \ref{h1a1p2} we have a complete collection of representatives of complements to $Q/Q(2)$ in $\bar XQ/Q(2)$ given by $X_\bold kQ(2)$ for $\bold k=[k_1,k_2,0,0,0,0]\in \mathbb V$ and $X_\bold k$ is generated by

 \begin{align*}x_+(t):&=\left[\begin{array}{c}1000\\t^{2^r}\end{array}\right]
\left[\begin{array}{c}0010\\t^{2^s}\end{array}\right]
\left[\begin{array}{c}1110\\k_1 t^{2^{r-1}}\end{array}\right]
\left[\begin{array}{c}0011\\k_2 t^{2^{s-1}}\end{array}\right]Q(2)\\
x_-(t):&=\left[\begin{array}{c}-1000\\t^{2^r}\end{array}\right]
\left[\begin{array}{c}-0010\\t^{2^s}\end{array}\right]
\left[\begin{array}{c}0110\\k_1 t^{2^{r-1}}\end{array}\right]
\left[\begin{array}{c}0001\\k_2 t^{2^{s-1}}\end{array}\right]Q(2)\end{align*} 

when $P=P_{24}$ and

 \begin{align*}x_+(t):&=\left[\begin{array}{c}1000\\t^{2^r}\end{array}\right]
\left[\begin{array}{c}0001\\t^{2^s}\end{array}\right]
\left[\begin{array}{c}1100\\k_1 t^{2^{r-1}}\end{array}\right]
\left[\begin{array}{c}0011\\k_2 t^{2^{s-1}}\end{array}\right]Q(2)\\
x_-(t):&=\left[\begin{array}{c}-1000\\t^{2^r}\end{array}\right]
\left[\begin{array}{c}-0001\\t^{2^s}\end{array}\right]
\left[\begin{array}{c}0100\\k_1t^{2^{r-1}}\end{array}\right]
\left[\begin{array}{c}0010\\k_2t^{2^{s-1}}\end{array}\right]Q(2)\end{align*} 

when $P=P_{23}$.

To show that the map $\sigma$ has the image as stated, we instead show that if $k_1,k_2\neq 0$ then the complements to $Q/Q(2)$ in $\bar XQ/Q(2)$ given above will not lift to complements to $Q/Q(3)$ in $\bar XQ/Q(3)$.

In order to show this, for $P_{24}$, we define \[x_+'(t)=x_+(t).\left[\begin{array}{c}0111\\f_1(t)\end{array}\right]
\left[\begin{array}{c}1111\\f_2(t)\end{array}\right]
\left[\begin{array}{c}0121\\f_3(t)\end{array}\right]
\left[\begin{array}{c}1121\\f_4(t)\end{array}\right]
\left[\begin{array}{c}1220\\f_5(t)\end{array}\right]
Q(3)\]

and for $P_{23}$,  \[x_+'(t)=x_+(t).\left[\begin{array}{c}0110\\f_1(t)\end{array}\right]
\left[\begin{array}{c}0111\\f_2(t)\end{array}\right]
\left[\begin{array}{c}1110\\f_3(t)\end{array}\right]
\left[\begin{array}{c}1111\\f_4(t)\end{array}\right]Q(3)\]

Using the $26$-dimensional matrix representation of $G$ we see in either case that the relation $x_+(t)x_+(u)=x_+'(t+u)$ cannot hold for any choice of the $f_i$ unless $k_1$ or $k_2=0$ a contradiction. If $k_1k_2=0$ then we can take all $f_i(t)=0$ to give a lift.
\end{proof}

For $P_{24}$ (resp. $P_{23}$), the shape $*1*0$ (resp. $*12*$) gives rise to the second coordinate $k_2$ (resp. third coordinate $k_3$) of $\bold k\in\mathbb V$ when $r\neq s+1$. When $r=s+1$, the module $V_S$ for this shape actually affords a trivial $H^1$. Thus we may relabel elements of $\mathbb V$ in the case $r=s+1$  as $[0,k_2,k_3,k_4,k_5,k_6]$, (resp. $[k_1,k_2,0,k_4,k_5,k_6]$). This will give us a consistent notation for the lemma which follows.

\begin{lemma}\label{noneinp23}When $P=P_{23}$, all complements to $Q$ in $\bar XQ$ are $G$-conjugate to subgroups of $P_{24}$ or $P_{13}$.\end{lemma}
\begin{proof} We know that any complement is conjugate to some $X_\bold k$ corresponding to a cocycle $\rho(\bold k)$ where $\bold k=[k_1,k_2,k_3,k_4,k_5,k_6]$. By \ref{a1a1blocking} we have $k_1k_2=0$. If $k_2=0$, then by inspection of the roots involved in the generators given in the proof of \ref{a1a1blocking} we see that $X_\bold k\leq L'\langle x_{1100}(t), x_{0100}(t):t\in K\rangle Q(2)$. It is easy to check that all roots involved in this expression are sent into $P_{24}$ under conjugation by $n_{0010}n_{0001}$. The graph automorphism of $G$ induces a transposition on the roots in $Q/Q(2)$ given explicitly as $(1100,0011)(0100,0010)$ and so if $k_1=0$ we see that $X_\bold k$ is conjugate to a subgroup of $P_{13}$.\end{proof}

\begin{lemma}\label{a1a1gens} Let $P=LQ$ be an $A_1\tilde A_1$-parabolic subgroup of $G$ and let $X\cong A_1$ be a complement to $Q$ in $\bar XQ$. Then either $X$ is conjugate to one of the following collection of subgroups $X_\bold k\leq P_{24}$ given by generators 
{\footnotesize \begin{align*}
x_+(t):&=\left[\begin{array}{c}1000\\t^{2^r}\end{array}\right]
\left[\begin{array}{c}0010\\t^{2^s}\end{array}\right]
\left[\begin{array}{c}1110\\k_1t^{2^{r-1}}\end{array}\right]
\left[\begin{array}{c}0011\\k_2t^{2^{s-1}}\end{array}\right]
\left[\begin{array}{c}1122\\k_3t^{2^{r-1}}\end{array}\right]
\left[\begin{array}{c}1231\\k_4t^{2^{s-1}}\end{array}\right]
\left[\begin{array}{c}1242\\k_5t^{2^{s}}\end{array}\right]
\left[\begin{array}{c}2342\\(k_6+\chi)t^{2^{r-1}}\end{array}\right]\\
x_-(t):&=\left[\begin{array}{c}-1000\\t^{2^r}\end{array}\right]
\left[\begin{array}{c}-0010\\t^{2^s}\end{array}\right]
\left[\begin{array}{c}0110\\k_1t^{2^{r-1}}\end{array}\right]
\left[\begin{array}{c}0001\\k_2t^{2^{s-1}}\end{array}\right]
\left[\begin{array}{c}0122\\k_3t^{2^{r-1}}\end{array}\right]
\left[\begin{array}{c}1221\\k_4t^{2^{s-1}}\end{array}\right]
\left[\begin{array}{c}1222\\k_5t^{2^s}\end{array}\right]
\left[\begin{array}{c}1342\\k_6t^{2^{r-1}}\end{array}\right]
\end{align*}} where $\chi=k_1^2k_3$, and $\bold k=[k_1,\dots,k_6]$ runs over all elements of $K^6$ subject to $k_1k_2=0$ and  $k_1=0$ if $r=s+1$;

or $X$ is conjugate to one of the subgroups $X_\bold m\leq P_{13}$ given by generators
{\footnotesize \begin{align*}
x_+(t):&=\left[\begin{array}{c}0001\\t^{2^r}\end{array}\right]
\left[\begin{array}{c}0100\\t^{2^{s+1}}\end{array}\right]
\left[\begin{array}{c}0122\\m_1t^{2^{r}}\end{array}\right]
\left[\begin{array}{c}1100\\m_2t^{2^s}\end{array}\right]
\left[\begin{array}{c}1111\\m_3t^{2^{r-1}}\end{array}\right]
\left[\begin{array}{c}1342\\m_4t^{2^s}\end{array}\right]
\left[\begin{array}{c}1221\\m_5t^{2^s}\end{array}\right]
\left[\begin{array}{c}1232\\m_6t^{2^{r-1}}\end{array}\right]\\
x_-(t):&=\left[\begin{array}{c}-0001\\t^{2^r}\end{array}\right]
\left[\begin{array}{c}-0100\\t^{2^{s+1}}\end{array}\right]
\left[\begin{array}{c}0120\\m_1t^{2^{r}}\end{array}\right]
\left[\begin{array}{c}1000\\m_2t^{2^{s}}\end{array}\right]
\left[\begin{array}{c}1110\\m_3t^{2^{r-1}}\end{array}\right]
\left[\begin{array}{c}1242\\m_4t^{2^{s}}\end{array}\right]
\left[\begin{array}{c}1121\\m_5t^{2^s}\end{array}\right]
\left[\begin{array}{c}1231\\m_6t^{2^{r-1}}\end{array}\right]
\end{align*}} as $\bold m=[m_1,\dots,m_6]$ runs over all elements of $K^6$ subject to $m_1m_2=0$, and if $r=s+1$ $m_1=0$.
\end{lemma}
\begin{proof}
Consider $X_\bold k\leq P_{24}$ initially. 

Using the $26\times 26$ matrices {\tt F4gens} one can input the given generators and check that they satisfy the relations necessary to be isomorphic to generate a group of type $A_1$. Further, using \ref{h1a1p2}, one sees that each $X_\bold k$ has been constructed as prescribed in \ref{rho}. Thus by \ref{surj} we have that the $X_\bold k$ give an exhaustive set of representatives of conjugacy classes of complements to $Q$ in $\bar XQ$.

Now, applying the graph automorphism to the $X_\bold k$ to convert root group elements in the generators we get the collection $X_\bold m\leq P_{13}$ given above which must also be an exhaustive set of representatives of conjugacy classes of complements to $Q$ in $\bar XQ\leq P_{13}$. 

By \ref{noneinp23}, this establishes a complete set of representatives of conjugacy classes of  complements.\end{proof}

\begin{lemma}\label{a1a1conjugacies} Suppose $r\neq s$. Then subgroups $X_\bold k, X_{\bold{k'}}\leq P_{24}$ are $G$-conjugate if and only if they are $C_G(L')$-conjugate, where $C_G(L')\cong A_1\tilde A_1$ corresponds to the  subsystem given by roots $\{\pm 1220,\pm 1232\}$ of $G$.

Subgroups $X_\bold m, X_{\bold{m'}}\leq P_{13}$ are $G$-conjugate if and only if they are $C_G(L')$-conjugate, where $C_G(L')\cong A_1\tilde A_1$ corresponds to the  subsystem given by roots $\{\pm 0121,\pm 2342\}$ of $G$.

Moreover, the following table generates all conjugacy in $G$ amongst the $X_\bold k$ and $X_\bold m$, where only some $\bold k$ and $\bold m$ are valid; (other values of $\bold k$ and $\bold m$ would be sent outside the parabolic by the given conjugating element).

\begin{tabular}{|c|c|c|}
\hline
$x\in C_G(L')$ & Valid $\bold k$ & $\bold k'=\bold k^x$\\
\hline
$x_{1220}(t)x_{1232}(u)$ & $[k_1,0,k_3,k_4,k_5,k_6]$ & $[k_1,0,k_3,k_4,k_5,k_6+tk_3]$\\
& $[0,k_2,k_3,k_4,k_5,k_6]$ & $[0,k_2,k_3,k_4+tk_2,k_5+tk_2^2,k_6+tk_3]$\\
$h_{1220}(\lambda)h_{1232}(\mu)$ & $[k_1,0,k_3,k_4.k_5,k_6]$ & $[\lambda k_1,0,\lambda^{-1}\mu^2 k_3,\lambda\mu k_4,\mu^2 k_5,\lambda\mu^2 k_6]$\\
& $[0,k_2,k_3,k_4,k_5,k_6]$ & $[0,\lambda^{-1}\mu k_2,\lambda^{-1}\mu^2 k_3,\lambda\mu k_4,\mu^2 k_5,\lambda\mu^2 k_6]$\\
$n_{1220}$ & $[0,k_2,k_3,k_4,k_5,k_6]$ & $[0,k_4,k_6,k_2,k_5,k_3]$\\
$n_{1232}$ & $[k_1,0,0,0,0,0]$ & $[k_1,0,0,0,0,0,0]$\\
$n_{1220}n_{1232}$ & $[k_1,0,0,0,0,0]$ & $[k_1,0,0,0,0,0,0]$\\\hline\end{tabular}

\begin{tabular}{|c|c|c|}
\hline
$x\in C_G(L')$ & Valid $\bold m$ & $\bold m'=\bold m^x$\\
\hline
$x_{0121}(t)x_{2342}(u)$ & $[m_1,0,m_3,m_4,m_5,m_6]$ & $[m_1,0,m_3,m_4,m_5,m_6+tm_3]$\\
& $[0,m_2,m_3,m_4,m_5,m_6]$ & $[0,m_2,m_3,m_4+t^2m_2,m_5+tm_2,m_6+tm_3]$\\
$h_{0121}(\lambda)h_{2342}(\mu)$ & $[m_1,0,m_3,m_4.m_5,m_6]$ & $[\lambda^2 m_1,0,\lambda^{-1}\mu m_3,\lambda^2\mu m_4,\mu m_5,\lambda\mu m_6]$\\
& $[0,m_2,m_3,m_4,m_5,m_6]$ & $[0,\lambda^{-2}\mu m_2,\lambda^{-1}\mu m_3,\lambda^2\mu m_4,\mu m_5,\lambda\mu m_6]$\\
$n_{0121}$ & $[0,m_2,m_3,m_4,m_5,m_6]$ & $[0,m_4,m_6,m_2,m_5,m_3]$\\
$n_{2342}$ & $[m_1,0,0,0,0,0]$ & $[m_1,0,0,0,0,0,0]$\\
$n_{0121}n_{2342}$ & $[m_1,0,0,0,0,0]$ & $[m_1,0,0,0,0,0,0]$\\\hline\end{tabular}

If $r=s$ then $G$-conjugacy amongst the $X_\bold k$ and $X_\bold m$ is generated by the items above, and the following further set of elements:

\begin{tabular}{|c|c|c|}
\hline
$x\in C_G(L')$ & Valid $\bold k$ & $\bold k'=\bold k^x$\\
\hline
$x_{1111}(t)x_{0121}(t)$ & $[k_1,0,k_3,k_4,k_5,k_6]$ & $[k_1,0,k_3,k_4+tk_1,k_5,k_6]$\\
& $[0,k_2,k_3,k_4,k_5,k_6]$ & $[0,k_2,k_3,k_4,k_5,k_6]$\\\hline
$x\in C_G(L)$ & Valid $\bold m$ & $\bold m'=\bold m^x$\\\hline
$x_{1122}(t)x_{1220}(t)$& $[m_1,0,m_3,m_4,m_5,m_6]$ & $[m_1,0,m_3,m_4+tm_1,m_5,m_6]$\\
& $[0,m_2,m_3,m_4,m_5,m_6]$ & $[0,m_2,m_3,m_4,m_5,m_6]$\\\hline
\end{tabular}.
\end{lemma}

\begin{proof}Firstly, let $r\neq s$. We appeal to \ref{nomoreconj2}. We calculate that $X_\bold k$ all share a maximal torus $H$ and all $X_\bold m$ share a maximal torus $K$. Then we observe that $M=\{1\}$, and that $C_G(L')\cap N_G(T), Z(L)$ and $C_Q(H)(=\{x_{1220}(t)x_{1232}(u):t,u\in K\})$ are all subgroups of $C_G(L')$.

The table follows by a routine application of commutator relations and Weyl group elements to the generators $x_\pm(t)$ of $X_\bold k$. The given elements generate $C_G(L')$ and thus all $G$-conjugacy amongst the $X_\bold k$ and $X_\bold m$.\end{proof}

\begin{lemma}\label{lsinp24} With the notation as above, the subgroups $X_\bold m$ which are conjugate to a subgroup of $P_{13}$ are precisely those with $m_3^2m_4=m_2m_6^2$.\end{lemma}
\begin{proof} Assume $X_\bold m$ is conjugate to $X_\bold k$ for some $\bold m$ and $\bold k$ by $g\in G$ and let $H$ and $K$ be maximal tori of $X_\bold m$ and $X_\bold k$ respectively, contained in the Levi subgroups of $P_{13}$ and $P_{24}$ respectively. Then we may assume that $H^g=K$. Thus $g\in N_G(H)w$ where $w=n_{1000}n_{0100}n_{0010}n_{0001}$, since $H^w=K$.

Now $N_G(H)=\langle n\rangle.C_G(H)$ where $\langle n\rangle H=N_{X_\bold k}(H)$ and $C_G(H)$ is the Levi subgroup of $G$ corresponding to roots $\pm 0121$ and $\pm 2342$.

Notice from the roots given in \ref{a1a1gens} that $X_\bold m^n\leq P_{24}$ if and only if $m_2=m_3=0$. Since $X_\bold m^g\leq P_{24}$ we must have $X_\bold m^r$ is some $X_\bold{m'}$ with $m'=[m_1',0,0,m_4',m_5',m_6']$ where $r\in C_G(H)$. We now use \ref{a1a1conjugacies} to establish that one must have $m_3^2m_4=m_2m_6^2$. 

Firstly, let $m_3^2m_4=m_2m_6^2$. We wish to see that $X_\bold m$ is conjugate to to $X_\bold{m'}$ of the above form. Now if $m_2=m_3=0$ then this is obvious. If $m_4=m_6=0$ then just take $X_\bold{m'}=X_\bold m^{n_{0121}} $. If $m_2=m_4=0$ or $m_3=m_6=0$ then conjugate by $x_{0121}(t)$ for appropriate $t$ to get $m_4=m_6=0$ and then by $n_{0121}$ again. Lastly we need to deal with the case that $m_3^2m_4=m_2m_6^2\neq 0$. So let $t^2=m_4/m_2=(m_6/m_3)^2$. We see from \ref{a1a1conjugacies} that we can replace $X_\bold m$ by its conjugate $X_{\bold m}^{x_{0121}(t)}$ to have $m_4=m_6=0$. Then by replacing $X_\bold m$ by its conjugate $X_\bold m^{n_{0121}}$ we may assume $m_2=m_3=0$ as required.

Now, If $m_3m_5^2\neq m_2m_6^2$ then working through conjugacy by elements in \ref{a1a1conjugacies} one quickly convinces oneself that that there is no value of $r\in C_G(L')$ giving us $X_\bold m^r$ with the second and third entries $0$. (For example, one could check that conjugation by the element $h_{0121}(\lambda)h_{2342}(\mu)$ has the effect of multiplying both $m_3m_5^2$ and $m_2m_6^2$ by $\mu^3$.)
\end{proof}

\begin{lemma}\label{a1a1exhaust} Let $X$ be a closed non-$G$-cr subgroup of $G$ of type $A_1$ contained subject to minimality in an $A_1\tilde A_1$ parabolic of $G$. Then $X$ is conjugate to  $X_\bold k$ or $X_\bold m$ where, $\bold k$ and $\bold m$ is an entry in the list below. If an entry is `$k_i$' or '$m_i$` rather than '0', it is intended to be a non-zero value.

For each value of $\bold k$ we have listed a subsystem overgroup $H$ of $X_\bold k$; or that it is no maximal reductive subgroups; or stated by `$\leq_G P_{24}$' that $X_\bold m$ is $G$-conjugate to a subgroup of $P_{24}$; or stated that $X_\bold k$ or $X_\bold m$ is conjugate to some other $X_\bold k'$ or $X_\bold m'$, respectively.
\begin{center}{\small
\begin{tabular}{l|l|c||l|l|c}
\hline
\# & $\bold k$ & $\leq M$ & \# & $\bold k$ & $\leq M$
\\\hline
1 & $[k_1,0,k_3,k_4,k_5,0]$ & none & 12 & $[0,k_2,k_3,k_4,0.k_6]$ & none \\
2 & $[k_1,0,0,k_4,k_5,k_6]$ & none & 13 & $[0,k_2,0,k_4,0,k_6]$ & none\\
3 & $[k_1,0,k_3,k_4,0,0]$ & $C_4$ & 14 & $[0,k_2,k_3,k_4,0,0]$ & none\\
4 & $[k_1,0,k_3,0,k_5,0]$ & $B_4$ & 15 & $[0,k_2,k_3,0,0,0]$ & none \\
5 & $[k_1,0,0,k_4,0,k_6]$ & none & 16 & $[0,k_2,0,k_4,0,0]$ & $C_4$ \\
6 & $[k_1,0,0,0,k_5,k_6]$ & $B_4$ & 17 & $[0,k_2,0,0,0,k_6]$ & $A_2A_2$\\
7 & $[k_1,0,k_3,0,0,0]$ & $\tilde A_1A_3$ & 18 & $[0,k_2,0,0,0,0]$ & $A_1\tilde A_2$\\
8 & $[k_1,0,0,k_4,0,0]$ & $C_4$ & 19 & $[0,0,k_3,k_4,k_5,0]$ & $\sim$ 13\\
9 & $[k_1,0,0,0,k_5,0]$ & $B_4$ & 20 & $[0,0,k_3,k_4,0,0]$ & $\sim$ 17\\
10 & $[k_1,0,0,0,0,k_6]$ & $A_3\tilde A_1$ & 21 & $[0,0,k_3,0,k_5,0]$ & $B_4$\\
11 & $[k_1,0,0,0,0,0]$ & $B_3$ & 22 & $[0,0,k_3,0,0,0]$ & $A_2\tilde A_1$\\
\hline
23 & $[0,0,0,k_4,k_5,k_6]$ & $\sim$ 14 \\
24 & $[0,0,0,k_4,k_5,0]$ & $\sim$ 16 \\
25 & $[0,0,0,k_4,0,k_6]$ & $\sim$ 15 \\
26 & $[0,0,0,k_4,0,0]$ & $\sim$ 18 \\
27 & $[0,0,0,0,k_5,k_6]$ & $\sim$ 21 \\
28 & $[0,0,0,0,k_5,0]$ & $\tilde A_3A_1$\\
29 & $[0,0,0,0,0,k_6]$ & $\sim$ 22
\end{tabular}
\begin{tabular}{l|l|c||l|l|c}
\hline
\# & $\bold m$ & $\leq M$ & \# & $\bold m$ & $\leq M$
\\\hline
1 & $[m_1,0,m_3,m_4,m_5,0]$ & none & 12 & $[0,m_2,m_3,m_4,0,m_6]$ & \parbox{120pt}{\begin{center}$\leq_G P_{24}$ if $m_2m_6^2=m_3m_4^3$,\\ none otherwise\end{center}}\\
2 & $[m_1,0,0,m_4,m_5,m_6]$ & $\leq_G P_{24}$ & 13 & $[0,m_2,0,m_4,0,m_6]$ & none\\
3 & $[m_1,0,m_3,m_4,0,0]$ & $B_4$ & 14 & $[0,m_2,m_3,0,m_5,0]$ & $\leq_G P_{24}$\\
4 & $[m_1,0,m_3,0,m_5,0]$ & $C_4$ & 15 & $[0,m_2,m_3,0,0,0]$ & $\leq_G P_{24}$ \\
5 & $[m_1,0,0,m_4,0,m_6]$ & $\leq_G P_{24}$ & 16 & $[0,m_2,0,0,m_5,0]$ & $\leq_G P_{24}$ \\
6 & $[m_1,0,0,0,m_5,m_6]$ & $\leq_G P_{24}$ & 17 & $[0,m_2,0,0,0,m_6]$ & $A_2\tilde A_2$\\
7 & $[m_1,0,m_3,0,0,0]$ & $A_1\tilde A_3$ & 18 & $[0,m_2,0,0,0,0]$ & $\leq_G P_{24}$\\
8 & $[m_1,0,0,m_4,0,0]$ & $\leq_G P_{24}$ & 19 & $[0,0,m_3,m_4,m_5,0]$ & $\sim$ 13\\
9 & $[m_1,0,0,0,m_5,0]$ & $\leq_G P_{24}$ & 20 & $[0,0,m_3,m_4,0,0]$ & $\sim$ 17\\
10 & $[m_1,0,0,0,0,m_6]$ & $\leq_G P_{24}$ & 21 & $[0,0,m_3,0,m_5,0]$ & $\leq_G P_{24}$\\
11 & $[m_1,0,0,0,0,0]$ & $\leq_G P_{24}$ & 22 & $[0,0,m_3,0,0,0]$ & $\leq_G P_{24}$\\
\hline
23 & $[0,0,0,m_4,m_5,m_6]$ & $\leq_G P_{24}$ \\
24 & $[0,0,0,m_4,m_5,0]$ & $\leq_G P_{24}$ \\
25 & $[0,0,0,m_4,0,m_6]$ & $\leq_G P_{24}$ \\
26 & $[0,0,0,m_4,0,0]$ & $\leq_G P_{24}$ \\
27 & $[0,0,0,0,m_5,m_6]$ & $\leq_G P_{24}$ \\
28 & $[0,0,0,0,m_5,0]$ & $\leq_G P_{24}$\\
29 & $[0,0,0,0,0,m_6]$ & $\leq_G P_{24}$
\end{tabular}}
\end{center}
Where item $n$ is conjugate to item $m$ we have written $\sim m$.
\end{lemma}
\begin{proof}
We will discuss the $\bold k$, first.

We leave to the reader to check that we have given an exhaustive list of the possible values of $\bold k\in K^6$ with no repetitions subject to the conditions \begin{enumerate}\item at least one of $k_1$ and $k_2$ is zero. (From \ref{a1a1blocking})
\item If $k_3$ is non-zero then $k_6=0$. (From \ref{a1a1conjugacies})
\item If $k_2$ is non-zero then $k_5=0$. (From \ref{a1a1conjugacies})
\end{enumerate}

In the third column, we have given any further conjugacy arising from Weyl group representatives from \ref{a1a1conjugacies} leaving the remainder as a unique collection of representatives of conjugacy classes amongst the $X_\bold k$ modulo $T$-conjugacy.

The statements about overgroups follow mostly by inspection. For instance we see that when $\bold k=[1,0,1,0,1,0]$ (case 4), $X_\bold k$ is generated by root elements contained in $\pm ***0,\pm ***2$ (those which generate $B_4$).

To justify the statement that a subgroup $X_\bold k$ lies in no proper reductive overgroup, we will now examine subgroups of type $A_1$ contained in $A_1\tilde A_1$-parabolic subgroups of $G$ arising from the maximal reductive subgroups $B_4$, $C_4$ and $A_2\tilde A_2$. 

\underline{$B_4$}
Suppose $X$ is one of the items claimed to be in no proper reductive subgroup and assume that $X$ is in a subgroup of type $B_4$.
\begin{claim} $X$ is in an $A_1\tilde A_1$-parabolic of $B_4$.\end{claim}
{\it Proof of claim:}

$V_{26}\downarrow X$ has the same composition factors as $V_{26}\downarrow \bar X={1^{[r]}}^2+{1^{[s]}}^4+{1^{[r]}\otimes 1^{[s]}}^2+0/1^{[s+1]}/0+0^2$. Since $V_{26}\downarrow B_4=T(1000)+0001$, we must find a way of splitting these composition factors into two self-dual pieces of dimension 10 and 16 with an appropriate asymmetry of $s$ and $t$. Checking all possible ways of doing this, we find that the only possibility is to have $T(1000)\downarrow X={1^{[r]}}^2|1^{[s+1]}|0^4$. 

Now we have $X$ in some parabolic of the $B_4$ by \ref{girclassical}. $X$ must be in a semisimple rank two or more parabolic of $B_4$ since it is in a s.s rank two parabolic of $G$. Now $T(1000)\downarrow A_2=10+01+0^4$ and it is impossible to make this consistent with the composition factors of $T(1000)\downarrow X$. Similarly we cannot have $X$ in a $B_2$-parabolic of $B_4$. Also, since $T(1000)\downarrow A_2\tilde A_1=(10,0)+(01,0)+(00,T(2))$, if$X$ were in an $A_2\tilde A_1$-parabolic of $X$ we would have the projection of $X$ in $A_2$ with composition factors $1^{[r]}|0$ on the natural module for $A_2$, putting it in a parabolic of $A_2$. Similarly we check $X$ cannot be in a $B_3$- or $A_3$-parabolic. Thus $X$ is in an $A_1\tilde A_1$-parabolic as claimed.

Now, since all subgroups of type $B_4$ are conjugate in $G$, $X^g$ is in the standard $B_4$ generated by roots $\mp 2342,\pm 1000,\pm 0100,\pm 0010$ with $X$ in some $A_1\tilde A_1$-parabolic of this $B_4$. By conjugating by an element $b$ of $B_4$ we have $X^{gb}$ in the standard $A_1\tilde A_1$ parabolic $P'$ of $B_4$ such that $P\cap B_4=P'$. corresponding to roots $1000,0010$ of $G$ or to roots $-2342,0010$ of $G$. By direct calculation with roots under the action of the Weyl group we find that the former parabolic of $B_4$ is conjugate in $B_4$ to the subgroup $P_{24}\cap B_4$ and the latter to $P_{13}\cap B_4$.

Further by conjugating $X$ by an element $q$ of $P$ we can arrange that $H\leq X^{gbq}$. Now we have $H\leq Y:=X^{gbq}\leq B_4\cap P$ a complement to $Q$ in $P$ and thus we know that it is conjugate to $X$ by an element of the form given in \ref{nomoreconj2}, i.e. by an element $x$ of $C_G(L')$. But $X^x\nleq B_4\cap P$ for any such $x$.

\underline{$C_4$}

The result follows from the above by using the graph automorphism of $G$.

\underline{$A_2\tilde A_2$}

Suppose that $X$ is in a subgroup of type $A_2\tilde A_2$. We note that the natural module $V=V_{A_2}(\lambda_1)\otimes V_{A_2}(\lambda_2)$ for $A_2\tilde A_2$ must restrict to $X$ as $(1^{[r]}|0,1^{[s]}|0)$ for some $r$ and $s$.

Since the long word in the Weyl group of $G$ induces a graph automorphism of both $A_2$ factors, the following (finite) collection of embeddings determine and exhaust all possibilities for $X$ up to $G$-conjugacy:
\begin{enumerate}\item $V\downarrow X=(1^{[r]}/0,1^{[s]}/0)$\\
\item $V\downarrow X=(1^{[r]}/0,0/1^{[s]})$\\
\item $V\downarrow X=(1^{[r]}/0,1^{[s]}+0)$\\
\item $V\downarrow X=(1^{[r]}+0,1^{[s]}/0)$\\
\item $V\downarrow X=(1^{[r]}+0,1^{[s]}+0)$\end{enumerate}

The last is clearly $\bar X$. We can see by inspection of the roots that (iv) corresponds to case $\bold k=[0,1,0,0,0,0]$, case (iii) to $\bold k=[0,0,0,0,0,1]$, cases (ii) and (i) to $\bold k=[0,1,0,0,0,1]$ and $\bold m=[0,1,0,0,0,1]$.

Since these are not conjugate to any $\bold k$ which we have claimed to be in no proper reductive subgroup, we are done.

For $\bold m$ we have used \ref{lsinp24} to indicate when some $X_\bold m$ is conjugate to a subgroup of $P_{24}$. In the remaining five cases of $\bold m$, no $X_\bold m$ is conjugate to any subgroup of $P_{24}$ and none are pairwise conjugate.
\end{proof}

\begin{remark}\label{indecomposablea1}
We can take a finite group $Y\cong L_2(4)=A_5\leq X_\bold k$ by restricting the parameter $t$ to $t\in GF(2^2)$. We may use the root group matrices to input this subgroup into GAP and make a calculation with the MeatAxe to examine the structure of $V_{26}\downarrow Y$. For instance, when $\bold k=[1,0,1,1,1,0]$ and $r\equiv s+1$ mod $(2)$ we see that the module $V_{26}$ restricted to $Y$ is indecomposable, with the following structure:
\begin{center}\begin{tabular}{c c c c c c c c c c}
&$1^{[s]}$\\
\parbox{10pt}{\line(1,0){100}}\\
$1^{[r}\otimes 1^{[s]}$&&&0&&&&&$0$\\
\hspace{-40pt}\parbox{5pt}{\line(1,0){40}}&$+$&\hspace{-20pt}\parbox{10pt}{\line(1,0){80}}&&&&&\hspace{-10pt}\parbox{10pt}{\line(1,0){80}}\\
$1^{[r]}$&&$1^{[s]}$&$+$&$1^{[s+1]}$&$+$&&$1^{[s]}$&&$1^{[r]}$\\
\parbox{10pt}{\line(1,0){80}}&&&&\hspace{-20pt}\parbox{10pt}{\line(1,0){90}}&&&&$+$&\hspace{-40pt}\parbox{5pt}{\line(1,0){40}}\\
&$0$&&&&$0$&&&&$1^{[r]}\otimes 1^{[s]}$\\
&&&&&\parbox{10pt}{\line(1,0){100}}&\\
&&&&&&&$1^{[s]}$
\end{tabular}\end{center}

It follows that this is the structure of $V_{26}$ for all $\bold k$ of the form $[k_1,0,k_3,k_4,k_5,0]$, provided $r\neq s, s+1$: We check that the conditions of \ref{samesubs} holds for the subgroup $Y$ of $X_\bold k$ when $\bold k=[1,0,1,1,1,0]$, for all $r\neq s,s+1$. It then follows that $X_\bold k$ has the same indecomposable action on $V_{26}$ as that given above. Now take any $\bold k$ of the above form, we check that $X_\bold k$ can be regarded as being in an $A_1^3$-parabolic $P$ of $E_7$ with Levi $L$ and that $X_\bold k$ is $Z(L)$-conjugate to $X_{[1,0,1,1,1,0]}$. The $56$-dimensional module restricts to $F_4$ as $2$ copies of $V_{26}$ and $4$ copies of the trivial module. Further restriction to $X_\bold k$ gives $2$ $26$-dimensional indecomposable modules and $4$ trivials; but $V_{26}\downarrow X_\bold k$ is a direct summand of this, so the result follows.
\end{remark}

\begin{theorem}\label{classofa1a1s} Let $X$ be a closed non-$G$-cr subgroup of $G$ of type $A_1$ contained subject to minimality in an $A_1\tilde A_1$ parabolic of $G$. Let $\bar X$ be the projection of $X$ to $A_1\tilde A_1$ such that the natural module $V_2\otimes V_2\downarrow \bar X=1^{[r]}\otimes 1^{[s]}$ and suppose that $r\neq s,s+1$. Then modulo $Z(L)$-conjugacy, $X$ is conjugate to precisely one subgroup $X_\bold k$ or $X_\bold m$ where, $\bold k$ and $\bold m$ is given in the list below. Where an entry is given as `$k_i$' or '$m_i$` rather than '$0$', it is intended to be a non-zero value. The number of $Z(L)$-conjugacy classes for each case of $\bold k$ is listed. 

\begin{tabular}{l|l|c||l|l|c}
\hline
$\bold k$ & $\leq M$ & $\#\{CCls\}$ & 			$\bold k$ & $\leq M$ & $\#\{CCls\}$\\\hline
$[k_1,0,k_3,k_4,k_5,0]$ & none & $2$-dim		& $[0,k_2,k_3,k_4,0.k_6]$ & none & $2$-dim\\
$[k_1,0,0,k_4,k_5,k_6]$ & none & $2$-dim 		& $[0,k_2,0,k_4,0,k_6]$ & none & $1$-dim\\
$[k_1,0,k_3,k_4,0,0]$ & $C_4$ & $1$-dim 		& $[0,k_2,k_3,k_4,0,0]$ & none & $1$-dim\\
$[k_1,0,k_3,0,k_5,0]$ & $B_4$ & $1$-dim		& $[0,1,1,0,0,0]$ & none & $1$\\
$[k_1,0,0,k_4,0,k_6]$ & none & $1$-dim			 & $[0,1,0,1,0,0]$ & $C_4$ & $1$ \\
$[k_1,0,0,0,k_5,k_6]$ & $B_4$ & $1$-dim		 & $[0,1,0,0,0,1]$ & $A_2A_2$ & $1$\\
$[1,0,1,0,0,0]$ & $\tilde A_1A_3$ & $1$			 & $[0,1,0,0,0,0]$ & $A_1\tilde A_2$ & $1$\\
$[1,0,0,1,0,0]$ & $C_4$ & $1$			 		& $[0,0,1,0,1,0]$ & $B_4$ & $1$\\
$[1,0,0,0,1,0]$ & $B_4$ & $1$ 					& $[0,0,1,0,0,0]$ & $A_2\tilde A_1$ & $1$\\
 $[1,0,0,0,0,1]$ & $A_3\tilde A_1$ & $1$			 & $[0,0,0,0,1,0]$ & $A_3\tilde A_1$ & $1$\\
 $[1,0,0,0,0,0]$ & $B_3$ & $1$					 &&& \\\hline
 $\bold m$ & $\leq M$ & $\#\{CCls\}$ & 			$\bold m$ & $\leq M$ & $\#\{CCls\}$\\\hline
 $[m_1,0,m_3,m_4,m_5,0]$ & none & $2$-dim			& \parbox{100pt}{$[0,m_2,m_3,m_4,0,m_6]$ \\ (if $m_2m_6^2\neq m_3m_4^2$)} & none & $2$-dim\\
 $[m_1,0,m_3,m_4,0,0]$ & $B_4$ & $1$-dim			&$[0,m_2,0,m_4,0,m_6]$ & none & $1$-dim\\
  $[m_1,0,m_3,0,m_5,0]$ & $C_4$ & $1$-dim 		&  $[0,1,0,0,0,1]$ & $A_2\tilde A_2$ & $1$\\
  &&&  $[1,0,1,0,0,0]$ & $A_1\tilde A_3$ & $1$\\\hline
\end{tabular}

If $r=s+1$ then modulo $Z(L)$-conjugacy, $X$ is conjugate to precisely one subgroup $X_\bold k$ or $X_\bold m$, where $\bold k$ (respectively, $\bold m$) is some item in the list above having $k_1=0$ (respectively, $m_1=0$).

If $r=s$ then modulo $Z(L)$-conjugacy, $X$ is conjugate to precisely one subgroup $X_\bold k$ or $X_\bold m$, where $\bold k$ and $\bold m$ are given in the list below:

\begin{tabular}{l|l|c||l|l|c}
\hline
$\bold k$ & $\leq M$ & $\#\{CCls\}$ & 			$\bold k$ & $\leq M$ & $\#\{CCls\}$\\\hline
$[k_1,0,k_3,0,k_5,0]$ & $B_4$ & $1$-dim		& $[0,k_2,k_3,k_4,0.k_6]$ & none & $2$-dim\\
$[k_1,0,0,0,k_5,k_6]$ & $B_4$ & $1$-dim 		& $[0,k_2,0,k_4,0,k_6]$ & none & $1$-dim\\
$[1,0,1,0,0,0]$ & $\tilde A_1A_3$ & $1$ 		& $[0,k_2,k_3,k_4,0,0]$ & none & $1$-dim\\
 $[1,0,0,0,0,1]$ & $A_3\tilde A_1$ & $1$				& $[0,1,1,0,0,0]$ & none & $1$\\
$[1,0,0,0,1,0]$ & $B_4$ & $1$ 			 & $[0,1,0,1,0,0]$ & $C_4$ & $1$ \\
 $[1,0,0,0,0,0]$ & $B_3$ & $1$				 & $[0,1,0,0,0,1]$ & $A_2A_2$ & $1$\\
			&&	& $[0,1,0,0,0,0]$ & $A_1\tilde A_2$ & $1$\\
			 		&&& $[0,0,1,0,1,0]$ & $B_4$ & $1$\\
		&&	& $[0,0,1,0,0,0]$ & $A_2\tilde A_1$ & $1$\\
	&&	 & $[0,0,0,0,1,0]$ & $A_3\tilde A_1$ & $1$\\\hline
 $\bold m$ & $\leq M$ & $\#\{CCls\}$ & 			$\bold m$ & $\leq M$ & $\#\{CCls\}$\\\hline
  $[m_1,0,m_3,0,m_5,0]$ & $C_4$ & $1$-dim			& \parbox{100pt}{$[0,m_2,m_3,m_4,0,m_6]$\\ (if $m_2m_6^2\neq m_3m_4^2$)} & none & $2$-dim\\
 $[1,0,1,0,0,0]$ & $A_1\tilde A_3$ & $1$		& $[0,m_2,0,m_4,0,m_6]$ & none & $1$-dim\\
 &&&$[0,1,0,0,0,1]$ & $A_2\tilde A_2$ & $1$\\\hline
 
\end{tabular}

Since the action of $Z(L)$ is known explicitly by \ref{a1a1conjugacies} the above amounts to a classification of subgroups of type $A_1$ in $A_1\tilde A_1$-parabolic subgroups of $G$.
\end{theorem}
\begin{proof} This is essentially a refinement of the previous lemma. We know from \ref{lsinp24} that no $X_\bold m$ is conjugate to any $X_\bold k$. We also know that the complements $X_\bold k$ (resp. $X_\bold m$) are $G$-conjugate if and only if they are $C_G(L')$-conjugate. Now $Z(L)$ only permutes elements $\bold k$ in the same entry in the table amongst themselves and by checking \ref{a1a1conjugacies} conjugacy by any other element either centralises a $\bold k$ or $\bold m$ or changes its entry number. Since we have given all instances where an element in one entry is conjugate to another, the above list gives an exhaustive list with no coincidences up to $Z(L)$-conjugacy. \end{proof}

We remark that, with exception of the parts of the proof deferred until \S\ref{xnonsimp} we have completed the proof of Theorem 1. Lemma \ref{p23orxgcr} gives the first assertion of (A).  Proposition \ref{a2s} gives the first two items in the table, and Lemmas \ref{classb3sg2s} and \ref{classb2s} give the remainder. For (B) the first assertion is \ref{p23orxgcr} again while the remaining statement is a brief summary of the results of sections \S\ref{a1p3} and \S\ref{a1p2}.

\section{Non-simple yet semisimple subgroups of parabolics: Proof of Theorem 2}\label{xnonsimp}

Now we tackle the non-simple semisimple closed subgroups of $G$. We first remind ourselves of the possible isomorphism types that need considering from \ref{levelsforx}.

\begin{lemma}\label{xnonsimpleoptions}Let $X=\prod X_i$ be non-simple closed semisimple subgroup of $G$. Then $X$ is $G$-cr unless $p=3$ and $X$ is of type $A_1A_1$ in a $B_3$-parabolic of $G$ or $p=2$ and $X$ is of type $A_1A_2$, $A_1B_2$, $A_1^2$ or $A_1^3$.\end{lemma}
\begin{proof} We examine \ref{levelsforx} to observe that $\mathbb V=0$ for each choice of $\bar X\leq L$ for a  parabolic with Levi $L$ apart from the exceptions listed.
\end{proof}

\subsection{Subgroups of type $A_1^2$, $p=3$}
\begin{lemma}\label{xnonsimplepis3} Let $p=3$ and $X\cong A_1A_1$ be non-$G$-cr. Then $X$ is conjugate to one of the following subgroups:
\begin{enumerate}\item $X\hookrightarrow D_3\tilde A_1\leq B_4$ by $(V_6,V_2)\downarrow X=(T(4),1)$
\item $X\hookrightarrow A_1C_3$ by $(V_2,V_6)\downarrow X=(1,T(3))$
\end{enumerate}
\end{lemma}
\begin{proof} By \ref{xnonsimpleoptions}, $\bar X$ is an irreducible subgroup of a $B_3$-Levi.  Then by \ref{a1sinb3} we have $\bar X\leq A_1A_1\tilde A_1$. Without loss of generality, we have either $\bar X\hookrightarrow A_1A_1\tilde A_1$ by $(V_2,V_2,\tilde V_2)\downarrow \bar X=((1^{[r]},0),(1^{[s]},0),(0,1^{[t]}))$ or $((1^{[r]},0),(0,1^{[s]}),(0,1^{[t]}))$. Thus $Q(2)\cong V_7\downarrow \bar X =(1^{[r]}\otimes 1^{[s]},0) +(0,2^{[t]})$ with $r\neq s$ or $(1^{[r]},1^{[s]})+(0,2^{[t]})$ and $Q/Q(2)\downarrow \bar X=(1^{[r]},1^{[t]})+(1^{[s]},1^{[t]})$ or $(1^{[r]},1^{[t]})+(0,1^{[s]}\otimes 1^{[t]})$, respectively.

In the first case, if $\mathbb V\neq 0$ then without loss of generality, $(r,s,t)=(1,0,0)$. As $H^0(\bar X,Q/Q(2))=H^1(\bar X,Q/Q(2))=0$, we have by \ref{nofixedpoints} that $H^1(\bar X,Q)\cong H^1(\bar X,Q(2))=K$. So we have a $1$-dimensional class of non-$G$-cr complements to $Q$ in $\bar XQ$ and, as usual, by \ref{onedimh1}, these are all $Z(L)$-conjugate to a single non-$G$-cr complement $X$. 

In the second case, if $\mathbb V\neq 0$ then without loss of generality, $(r,s,t)=(0,1,0)$ or $(0,0,1)$ and $\mathbb V\cong H^1(\bar X,Q/Q(2))$. If $(r,s,t)=(0,1,0)$ there are no potential blockers and we conclude that $H^1(\bar X,Q)\cong K$. So we have a $1$-dimensional class of complements to $Q$ in $\bar XQ$ and, as usual, these are all $Z(L)$-conjugate to single complements $Y$.

If $(r,s,t)=(0,0,1)$ we claim that $\rho=0$ and so there are no further complements. Assume otherwise. Then the map $H^1(\bar X,Q)\to H^1(\bar X,Q/Q(2))$ is non-zero. Now we can take a subgroup $Z\cong A_1\mapsto \bar X$ by $(V_2,V_2)=(1^{[r]},1)$ giving us case \ref{0s1} from \S\ref{a1schar3sec}. Now,  observe that the following diagram commutes, where vertical maps are restriction:

\begin{center}
$\begin{CD} H^1(\bar X,Q)@>>> H^1(\bar X,Q/Q(2))\\
@VVV @VVV\\
H^1(Z,Q) @>>> H^1(Z,Q/Q(2))\end{CD}$\end{center}

But the bottom map is zero by the calculation done in \ref{0s1}, and the right hand side map is an isomorphism by restriction, so it follows that $\rho=0$.

Thus we have up to two non-$G$-cr complements to $Q$ in $\bar XQ$: one from case $(r,s,t)=(1,0,0)$ and one from case $(r,s,t)=(0,1,0)$.

We claim that these are the two subgroups listed in the statement of this lemma.

It is clear from \ref{girclassical} that each is in a parabolic subgroup of $G$. Comparing the restrictions  of the corresponding entries in \ref{tor} with those of all $B_3$-ir subgroups of $G$ we find that neither is conjugate to a subgroup of the $B_3$ Levi. Thus $X$ and $Y$ are non-$G$-cr. Since they have distinct restrictions, they are not pairwise conjugate.
\end{proof}

\subsection{Subgroups of type $A_1A_2$ and $A_1B_2$, $p=2$}\label{a1a2ora1b2} 
\begin{lemma}\label{a1a2s}Let $p=2$ and let $X\cong A_1A_2$ be non-$G$-cr. Then $X\leq A_2\tilde A_2$ and $X$ is conjugate to one of the following subgroups:
\begin{enumerate}\item $X\hookrightarrow A_2\tilde A_2$ by $(V_3,V_3)\downarrow X=(2/0,10)$
\item $X\hookrightarrow A_2\tilde A_2$ by $(V_3,V_3)\downarrow X=(10,2/0)$.
\end{enumerate}
\end{lemma}
\begin{proof}Clearly $\bar X\cong A_1\tilde A_2$ or $A_2\tilde A_1$. In each case, we observe from \ref{levelsforx} that $\mathbb V\cong K$. We get one non-$G$-cr subgroup from \ref{onedimh1} in each case. The subgroups in the lemma are in parabolics of $A_2\tilde A_2$ and so are in parabolics of $G$. By comparing the restrictions in \ref{tor} we see that they are not conjugate to $A_2A_1$ Levi subgroups, and not conjugate to each other.\end{proof}

\begin{lemma}\label{a1b2s} Let $p=2$ and let $X\cong A_1B_2$ be non-$G$-cr. Then $X$ is contained in a subgroup of type $B_4$ or $C_4$ and there is a one-dimensional collection of conjugacy classes of each. \end{lemma}
\begin{proof}Either $X$ is contained in the $B_3$ parabolic or the $C_3$ parabolic. We will assume the former and then get results for the latter by applying the graph automorphism.

Thus $\bar X=\tilde A_1B_2\leq B_3$. From \ref{levelsforx} we get that $\mathbb V\cong K^2$. As $H^0(\bar X,Q/Q(2))=H^1(\bar X,Q/Q(2))=0$, \ref{nofixedpoints} implies that $\mathbb V\cong H^1(\bar X.Q)$. Since $\bar XQ(2)\leq B_4$, by \ref{subsystemlift} we see that $X\leq B_4$ as claimed. We may parameterise the complements $X_\bold k$ to $Q$ in $\bar XQ$ as $\bold k=[k_1,k_2]$. By \ref{b3imp}, $X_\bold k$ and $X_\bold k'$ are conjugate if any only if they they are $QZ(L)$-conjugate. Thus $X_\bold k$ is conjugate to $X_\bold{k'}$ if and only if $\bold k'=\lambda\bold k$ for some $\lambda\in K^\times$.

This shows that the collection of $X_\bold k$ up to $G$-conjugacy  form a $1$-dimensional subset of $K^2$.
 \end{proof}

\subsection{Subgroups of type $A_1^3$, $p=2$}\label{a1a1a1sec}
\begin{lemma}\label{a1a1a1} Let $p=2$ and let $X\cong A_1^3$ be non-$G$-cr. Then $X$ is contained in a subsystem subgroup of type $B_4$ or $C_4$ and $X$ is conjugate to one of the subgroups in the following list  \begin{enumerate}\item $X\hookrightarrow A_1^4$ or $\tilde A_1^4$ by $(V_2,V_2,V_2,V_2)\downarrow X=((1,0,0),(0,1,0),(0,0,1),(0,0,1))$,
\item $X\hookrightarrow D_4$ or $\tilde D_4$ by $V_8\downarrow X=\begin{array}{c}0\\\hline (1^{[1]},0,0)+(0,1^{[1]},0)+(0,0,1^{[1]})\\\hline 0\end{array}$
\item $X\hookrightarrow A_1^2\tilde A_1^2$ by $(V_2,V_2,V_2,V_2)\downarrow X=((1,0,0),(1,0,0),(0,1,0),(0,0,1))$ or by $(V_2,V_2,V_2,V_2)\downarrow X=((1,0,0),(0,1,0),(0,0,1),(0,0,1))$
\item $X\hookrightarrow \tilde A_1B_3$ by $(V_2,V_7)\downarrow X=\left((1,0,0),\begin{array}{c}0\\\hline (0,1^{[1]},0)+(0,0,1^{[1]})\\\hline 0\end{array}\right)+0$ with the projection to the second two factors contained in no proper reductive subgroup of $B_3$; a $1$-dimensional collection of such. (Or analogue for $X\hookrightarrow A_1C_3\leq C_4$.)
\item $X$ is in one of a $2$-dimensional variety  of conjugacy classes not in any proper reductive subgroup of $B_4$ or $C_4$.\end{enumerate}\end{lemma}
\begin{proof} Again, $X$ is contained in a $B_3$-parabolic or a $C_3$-parabolic but these are conjugate under the graph automorphism of $G$. So we assume the former.

Then $\bar X=A_1^2\tilde A_1$ or $\bar X=\tilde A_1^3$.

In the first case, we see from \ref{levelsforx} that $\mathbb V\cong Q(2)\cong K$. As $H^0(\bar X,Q/Q(2))=H^1(\bar X,Q/Q(2))=0$ \ref{nofixedpoints} implies that $H^1(\bar X,Q)\cong H^1(\bar X,Q(2))$. Thus there is one non-trivial complement $X$ up to $QZ(L)$-conjugacy and by \ref{subsystemlift}, $X\leq B_4$. Checking \ref{tor} we see that the subgroup given as item (i) above has composition factors consistent only with being a complement to an $A_1^2\tilde A_1$ subgroup of a $B_3$-parabolic of $G$ but is not conjugate to $\bar X$. Thus $X$ is conjugate to item (i).

If $\bar X=\tilde A_1^3$ then from \ref{levelsforx} we have $\mathbb V\cong H^1(\bar X,Q)\cong K^3$. Let the corresponding complements be $X_\bold k$ for $\bold k\in K^3$. From \ref{b3imp} we deduce that $X_\bold k$ is conjugate to $X_\bold{k'}$ if and only if $\bold k=\lambda \bold k'$ for some $\lambda\in K^\times$. Thus there is a $2$-dimensional collection of conjugacy classes $X_\bold k$. 

If $\bold k=(1,1,1)$ then it is easy to see that $X_\bold k$ is contained in a non-$G$-cr $B_3$, which is in $D_4$. There is only one of these up to conjugacy inside $D_4$ as we must have \[V_8\downarrow X_\bold k=V_8\downarrow X=\begin{array}{c}0\\\hline (1^{[1]},0,0)+(0,1^{[1]},0)+(0,0,1^{[1]})\\\hline 0\end{array}\] and there is only one $D_4$ up to conjugacy inside $B_4$; this gives case (ii). 

If $\bold k=(0,0,1)$, $(0,1,0)$ or $(1,0,0)$ then we may assume that it is last of these as the $\tilde A_1^3$ factors are permuted by the Weyl group. Then clearly projection to the second two factors of $X_\bold k$ gives a subgroup of type $\tilde A_1^2$ contained in the Levi subgroup $B_2\leq B_3$. The centraliser of this $\tilde A_1^2$ in $B_4$ is a $B_2$ and so the remaining factor must lie in this $B_2$. This factor is an $A_1$ subgroup of $B_2$ lying in a $\tilde A_1$-parabolic of $B_2$. The usual method shows that there is only such and thus that it is the subgroup stated in case (iii).

If just one of the coordinates of $\bold k$ is zero, then without loss of generality it is the first. Thus the first factor is an $\tilde A_1$ Levi subgroup of $B_4$. The centraliser of this is in $B_4$ is a $B_3$ and thus these remaining subgroups $X_\bold k$ (not inside $A_1^2\tilde A_1^2$ subgroups) have projections to this $B_3$ in no proper reductive subgroup of $B_3$; this describes item (iv).

We have found all non-$G$-cr subgroups in $D_4$, and $B_3\tilde A_1$ with the correct composition factors on restriction of $V_9$. Thus the remaining $X_\bold k$ are not in any proper reductive subgroup of $B_4$. These are in the remaining case (v).

The results for $C_4$ follow by application of the graph automorphism of $G$.
\end{proof}

\subsection{Subgroups of type $A_1^2$, $p=2$}\label{a1a1s}

\begin{lemma}\label{a1a1sina1a1a1s} Let $p=2$ and let $X\cong A_1^2$ be non-$G$-cr contained subject to minimality in a $B_3$-parabolic or $C_3$-parabolic of $G$. Then $X$ is conjugate to a subgroup of exactly one of the $A_1^3$ subgroups of $G$ listed in \ref{a1a1a1} above.\end{lemma}
\begin{proof} As usual we deal with $B_3$ and deduce the results for $C_3$ using the graph automorphism of $G$.

Let $\bar X\cong A_1A_1$ be $B_3$-ir, then. It is clear that one factor must centralise a $2$-torus in $B_3$ and the other a $1$-torus. So one of the factors is a Levi $A_1$ of $B_3$. There are two cases then from \ref{a1sinb3}: (a) $\bar X\leq A_1^2\tilde A_1$ or (b) $\bar X\leq \tilde A_1^3$.

In case (a), we have, without loss of generality either: 
\begin{enumerate}\item $V_7\cong Q(2)\downarrow \bar X=(1\otimes 1^{[s]},0)+(0,1^{[1]}/0)$, $s\neq 1$, and\\
$Q/Q(2)\downarrow \bar X=(1,1)+(1^{[s]},1)$; or
\item $V_7\downarrow\bar X=(1,1^{[s]})+(0,1^{[t+1]}/0)$, and\\
 $Q/Q(2)\downarrow \bar X=(1,1^{[t]})+(0,1^{[s]}\otimes 1^{[t]})$ respectively.\end{enumerate} Unless $s=t$ in the second case, $H^0(\bar X,Q/Q(2))=H^1(\bar X,Q/Q(2))=0$ and so from \ref{nofixedpoints} we have $H^1(\bar X,Q)\cong K$. Thus there is up to one conjugacy class of non-$G$-cr complements to $\bar X$ and comparing restrictions in \ref{tor}, it must be $A_1A_1\hookrightarrow A_1^4$ by $(V_2,V_2,V_2,V_2)\downarrow \bar X=((1^{[1]},0),(1^{[s]},0),(0,1),(0,1))$ or $(V_2,V_2,V_2,V_2)\downarrow X=((1^{[t]},0),(1^{[t]},0),(1,0),(0,1^{[s]})$ respectively. It is then easy to check that each of these is conjugate to just one of the subgroups in item (i) in \ref{a1a1a1} above. It cannot be in any other since the composition factors of $X$ are different to the others.

We claim that if $s=t$ then $H^1(\bar X,Q)$ is still one-dimensional in the second case. To prove this we use the method in \ref{classb3sg2s} that we applied to see that $H^1(G_2,Q)=K$. That is, we have $\bar X$ centralising a $1$-dimensional subgroup and $1$-dimensional quotient of the $3$-dimensional subgroup of $Q$ which is $\langle x_{1111}(t),x_{0121}(t),x_{1232}(t):t\in K\rangle$ and $H^1(\bar X,Q)=K^2$. Moreover, the sequence $\{1\}\to H^0(\bar X,Q(2))\to H^1(\bar X,Q)\to H^1(\bar X,Q/Q(2))\to \{1\}$ is exact. Thus in the exact sequence in \ref{les}, the map $H^1(\bar X,Q/Q(2))\to H^1(\bar X,Q(2))$ is zero. Since $H^1(\bar X,Q(2))\to H^1(\bar X,Q)=0$ is an injection. Thus $H^1(\bar X,Q)\cong K$. From \ref{onedimh1} we have a single non-$G$-cr complement arising and it is easily seen that this is in item (i) in \ref{a1a1a1} above. It cannot be in any other since the composition factors of $X$ are different to the others.

Now if $\bar X\leq \tilde A_1^3$ then without loss of generality, $\bar X\hookrightarrow\tilde A_1^3$ by $(V_2,V_2,V_2)\downarrow \bar X=((1,0),(0,1^{[s]}),(0,1^{[t]}))$ and $s\neq t$ by \ref{girclassical} or $\bar X$ would not be $B_3$-ir. Now we observe that $H^0(\bar X,Q/Q(2))=0=H^1(\bar X,Q/Q(2))$ and so $H^1(\bar X,Q)\cong H^1(\bar X,Q(2))=K^3$ by \ref{nofixedpoints}. Indeed $H^1(\tilde A_1^3,Q(2))\to H^1(\bar X,Q(2))$ is an isomorphism by restriction. Thus restriction of the cocycles of the former to $\bar X$ induces a one-to-one mapping between $Q$-conjugacy classes of $\tilde A_1^3$ complements to $Q$ in $\tilde A_1^3Q$ and $A_1A_1$ complements to $Q$ in $\bar XQ$. $Z(L)$-conjugacy collapses these complements in the same way as in \ref{a1a1a1} and so there is a one-to-one correspondence between $G$-conjugacy classes of the complements in each case.\end{proof}

\begin{lemma} \label{a1a1inb2} Let $p=2$ and let $X$ be a non-$G$-cr subgroup of type $A_1^2$ contained subject to minimality in a $B_2$-parabolic of $G$. Then $X$ is in a $B_4$ or $C_4$ subsystem of $G$ and $X$ is conjugate to exactly one subgroup in the following list:
\begin{enumerate}\item $X\hookrightarrow A_1^4$ or $\tilde A_1^4$ by $(V_2,V_2,V_2,V_2)\downarrow X=((1,0),(1,0),(0,1),(0,1))$
\item $X\hookrightarrow A_1^2\tilde A_1\leq B_3$ or $\tilde A_1^2A_1\leq C_3$ by $(V_2,V_2,V_2)\downarrow X=((1,0),(1,0),(0,1))$ or  $(V_2,V_2,V_2)\downarrow X=((1,0),(1,0),(0,2))$, respectively
\item $X\hookrightarrow D_3$ or $\tilde D_3$ by $V_6\downarrow X=\begin{array}{c}0\\\hline (1^{[1]},0)+(0,1^{[1]})\\\hline 0\end{array}$
\item $X\hookrightarrow B_3$ by $V_7\downarrow X=\begin{array}{c}0\\\hline (1^{[1]},0)+(0,1^{[1]})\\\hline 0\end{array}+0$ with $X$ in no proper reductive subgroup of $B_3$; a $1$-dimensional collection of such. (Or analogue for $X\hookrightarrow C_3$.)
\end{enumerate}
\end{lemma}
\begin{proof}
By \ref{lir}, $\bar X\leq A_1^2$ or $\bar X\leq \tilde A_1^2$ and we may assume the latter, up to graph automorphism. Thus for the natural module for $B_2$, $V_5\downarrow \bar X=\frac{(2,0)+(0,2)}{(0,0)}$.

By \ref{levelsforx} we have four levels with only levels 1 and 3 admitting non-zero $H^1$s: $H^1(\bar X,Q(1)/Q(2))=H^1(\bar X,Q(3)/Q(4))=K^2$ so that $\mathbb V=K^4$.

As the shapes $S$ affording non-trivial values of $H^1(\bar X,V_S)$ are contained in the $B_4$ subsystem, \ref{subsystemlift} implies that  $\rho$ is defined everywhere, as there are no potential blockers in $Q\cap B_4$; and moreover $X$ is in a $B_4$ subsystem subgroup of $G$. Let $\{X_\bold k:\bold k\in \mathbb V\}$ be a set of complements sharing a common torus, corresponding to the images of $\mathbb V$ under $\rho$. As usual we will say $\bold k^g=\bold k'$ if $X_\bold k^g=X_{\bold k'}$.

Using \ref{nomoreconj2} we calculate that any two complements are $G$-conjugate, if and only if they are $C_G(L')$-conjugate. In particular if $X$ and $X'$ are $G$-conjugate complements to $Q$, then they are $B_4$-conjugate.

\begin{claim} $G$-conjugacy amongst complements $X_\bold k$ pulls back to an action of $C_G(L')$ on $\mathbb V$ generated by the following elements:

\begin{center}\begin{tabular}{|c|l|l|}
\hline
$x\in C_G(L)$ & Valid $\bold k$ & $\bold k'=\bold k^x$\\\hline
$x_{0122}(t)$ & $[k_1,k_2,k_3,k_4]$ & $[k_1,k_2,k_3+tk_1,k_4+tk_2]$\\
$x_{1110}(t)$ & $[k_1,k_2,k_3,k_4]$ & $[k_1,k_2,k_3,k_4]$\\
$n_{0122}$ & $[k_1,k_2,k_3,k_4]$ & $[k_3,k_4,k_1,k_2]$\\
$n_{1110}$ & $[0,0,k_3,k_4]$ & $[0,0,k_3,k_4]$\\
$h_{0122}(t)$ & $[k_1,k_2,k_3,k_4]$ & $[tk_1,tk_2,t^{-1}k_3,t^{-1}k_4]$\\
$h_{1110}(t)$ & $[k_1,k_2,k_3,k_4]$ & $[t^{-2}k_1,t^{-2}k_2,k_3,k_4]$\\\hline
\end{tabular}\end{center}\end{claim}
{\it Proof of claim:}
By \ref{nomoreconj2} $C_G(L')=\langle x_{0122}(t),x_{1110}(t),1,n_{1110},n_{0122},Z(L)\rangle$ as $t$ varies over $K$. Thus we need only show that each element acts as it does.

Following the argument in \ref{a2s} we find from the commutator formula that we get $\bar X$-equivarient isomorpisms  $t_\lambda: Q(1)/Q(2)\to Q(3)/Q(4)$ induced by commutating with the element $x_{0122}(\lambda)$, and thus the action of $x_{0122}$ is as given. Similarly the commutator formula allows us to calculate the other actions, noting that if one of $k_1$ and $k_2$ is non-zero, then $X_\bold k^{n_{1110}}$ is not in the parabolic $P$. This proves the claim.

Under conjugation by the above elements then, it is easy to check that if $\bold k=[k_1,k_2,k_3,k_4]$ is such that $k_3k_2\neq k_4k_1$. then $X_\bold k$ is conjugate to $X_{[1,0,0,1]}$. If $k_1=k_3=0$ with $k_2$ or $k_4$ non-zero, (resp. $k_2=k_4=0$ with $k_1$ or $k_3$ non-zero) then $X_\bold k$ is conjugate to $X_{[1,0,0,0]}$ (resp. $X_{[0,1,0,0]}$).

If $k_3k_2=k_4k_1\neq 0$ then $X_\bold k$ is conjugate to $X_{[k_1,k_2,0,0]}$ for some $k_1,k_2\neq 0$. It is easy to check that the shape associated to $k_1,\ k_2$ is $1**0$, contained in $B_3$. And so from \ref{subsystemlift} we have that $X_\bold k\leq B_3$. Now from the table above, $X_\bold k$ is conjugate to $X_{\bold k'}$ provided $\bold k=\lambda.\bold k'$. If $k_1=k_2$ we discover that $X_\bold k\leq B_2\leq D_3\leq B_3$ and otherwise, $X_\bold k$ is in no proper subsystem subgroup of $B_3$.

It remains to identify $X_{[1,0,0,0]}$; $X_{[0,1,0,0]}$ and $X_{[1,0,0,1]}$. Now the first of these (resp. second) corresponds to a cocycle which is trivial on restriction to the second factor (resp. first factor), and is in a $B_3$. Hence the other factor is contained in $C_{B_3}(\tilde A_1)=D_2$ and then for $X_\bold k$ to have the same composition factors as $\bar X$, we deduce that the first factor is embedded in the $D_2$ with the same twists, as described. Lastly we claim that $X_{[1,0,0,1]}$ is conjugate to the subgroup $X\hookrightarrow A_1^4$ or $\tilde A_1^4$ by $(V_2,V_2,V_2,V_2)\downarrow X=((1,0),(1,0),(0,1),(0,1))$. The latter has $V_9\downarrow X=\begin{array}{c}0\\\hline (2,0)\\\hline 0\end{array}+\begin{array}{c}0\\\hline (0,2)\\\hline 0\end{array}+0$ and so is not in a $B_3$ since it does not stabilise a non-degenerate $2$-space. It is also in a rank $2$ parabolic, as each factor is contained in a $B_2$ subgroup of $B_4$. A quick check on composition factors then shows that it must be in a $B_2$ parabolic of $B_4$ and moreover a complement to $\tilde A_1^2\leq B_2$. Thus it is the subgroup $X_{[1,0,0,1]}$.\end{proof}

\begin{lemma} \label{a1a1sina1a1pars}Let $p=2$ and let $X\cong A_1^2$ be non-$G$-cr contained subject to minimality in an $A_1\tilde A_1$-parabolic subgroup $P=LQ$ of $G$. Let $\bar Y\leq \bar X\leq L$ be a subgroup of type $A_1$ projecting to the factors of $\bar X$ with twists $r$ and $s$ respectively where $r,\ s>0$ and $r\neq s,\ s+1$ (such as $r=1$, $s=3$). Then $X$ contains exactly one complement $Y$ to $Q$ in the semidirect product $\bar YQ$. Moreover if $Y\leq M$ where $M$ is a reductive maximal subgroup of $G$ then $X\leq M$.

Thus the classification of conjugacy classes in $A_1$ subgroups of $G$ given in \ref{classofa1a1s} when $r\neq s,s+1$ is in one-to-one correspondence with the set of conjugacy classes of subgroups of type $A_1\tilde A_1$. Moreover,  is in a subsystem whenever its corresponding subgroup of type $A_1$ is.\end{lemma}
\begin{proof}
This is a matter of checking that the proofs of the results leading up to \ref{classofa1a1s} go through in the case that $\bar X$ is of type $A_1\tilde A_1$. 

We check for each $A_1\tilde A_1$ parabolic of $G$ that  $\dim \mathbb V=6$. Moreover, if $\bar Y$ is a subgroup of $\bar X$ of type $A_1$ and $\bar Y\hookrightarrow \bar X$ by $(V_2,V_2)\downarrow \bar X=(1^{[r]},1^{[s]})$ where $r\neq s, s+1$. then there is an isomorphism $H^1(\bar X,Q)\cong H^1(\bar Y,Q$ by restriction. It is easy to see that on the generators of the complements given in \ref{a1a1gens} this isomorphism works by splitting up each element $x_\pm(t)$ into $x_{\pm,1}(t)$ and $x_{\pm,2}(t)$ where $x_{\pm,1}(t)$ (resp. $x_{\pm, 2}(t)$ is a product of all root group elements containing an `$r$' (resp. an `$s$').

Then we check also that the conclusion of \ref{noneinp23} holds for $A_1^2$ subgroups and that the same conditions for a complement $Y$ in $P_{13}$ to be conjugate to a subgroup of $P_{24}$ hold. Lastly we establish all remaining conjugacy amongst all complements $X$ in $P$ by use of \ref{nomoreconj2}; this amounts to finding out which subgroups are $C_G(L')\cong A_1\tilde A_1$-conjugate which amounts to the same calculation as that done for $Y$.
\end{proof}

Note that this concludes the proof of Theorem 2, which is a significant abbreviation of our results in some ways. The case when $p=3$ is simply \ref{xnonsimplepis3}. When $p=2$ and $X\cong A_1A_2$, \ref{a1a2s} gives (ii). When $p=2$ and $X$ is contained in an $A_2\tilde A_2$ subsystem, it is in an $A_1\tilde A_1$-parabolic, since projection to each $A_2$ factor must be reducible. Thus one can confirm the embeddings given from \ref{classofa1a1s} in light of  Lemma \ref{a1a1sina1a1pars}, giving (iii). The statement in (iv) encapsulates results from \S\ref{a1a1s}, specifically \ref{a1a1sina1a1a1s} and \ref{a1a1inb2}, together with the cases of \ref{classofa1a1s} which arise in light of \ref{a1a1sina1a1pars} excluding those already mentioned in (iii). The statement in (v) encapsulates the results \ref{a1b2s}, which deals with $A_1B_2$s; and \ref{a1a1a1}, which deals with $A_1^3$s. Item (vi) is those cases of \ref{classofa1a1s} which arise in light of \ref{a1a1sina1a1pars} which are not covered by either (iii) or (iv).

\section{Subgroups not in parabolics: Proof of Theorem 4}\label{girsubgsec}
Let $X$ be a simple $G$-irreducible subgroup of $G$ of rank at least $2$.

We work recursively through the maximal subgroups of $G$ which are not parabolic subgroups.

Since $X$ is not contained in any parabolic subgroup of $G$, it must be contained in one of the maximal reductive subgroups $M$ of $G$ of rank at least $2$. By \cite{LS04} these are $A_1G_2$ $(p>2)$, $B_4$, $A_2\tilde A_2$, $A_1C_3$ ($p\neq 2$), $C_4$ $(p=2)$, $G_2$ $(p=7)$.

The subgroup $X$ has rank at least $2$ so if $X\leq M$ where $M=A_1G_2$ or $A_1C_3$, then $X$ has trivial projection to the $A_1$ factor of each. Then $X\leq C_G(A_1)$ and thus $X$ is in a Levi subgroup of $G$
 
Now if $M=G_2$ when $p=7$ then $X=M$ or $A_2$. Clearly $M$ is $G$-ir. If $X=A_2$ then it has a central $3$-element, thus lies in the centraliser of this $3$-element. The centraliser of this $3$-element in $G$ is an $A_2\tilde A_2$ subgroup of $G$. Thus $X\leq A_2\tilde A_2$; by checking restrictions it is easy to see that  it is embedded via $(10,10)$. We consider this subgroup separately below.

If $M=A_2\tilde A_2$ and $X\cong A_2$ is embedded via $(10^{[r]},10^{[s]})$ or $(10^{[r]},01^{[s]})$ with $r\neq s$ then we see from \ref{tor} that all embeddings listed in the table are visibly distinct and inconsistent with being in parabolics (each has two composition factors of dimension $9$). If $r=s$ and $X$ is embedded in $A_2\tilde A_2$ as $(10,01)$ then we showed in \ref{a2s} that $X=Z_1$ is in a parabolic subgroup of $G$ when $p=3$. If $p\neq 3$ then we showed in the previous section that $X$ is $G$-cr; thus it would be in a Levi subgroup of $G$, which is again impossible. Finally if $X$ is embedded in $A_2\tilde A_2$ as $(10,10)$ $V_{26}\downarrow X=20+02+10+01+11$ and again the dimensions of the composition factors are inconsistent with being in any  parabolic subgroup of $G$.

Now let $M=B_4$. From \cite{Sei87}, the maximal subgroups of $B_4$ with no factor of type $A_1$ are $D_4$ and $B_2B_2$ when $p=2$ (we may exclude all those with factors of type $A_1$ for the same reason as we excluded $A_1G_2$ above). We see from the restrictions in \ref{tor} that if $X\cong  B_2\hookrightarrow B_2B_2$ with different twsts, then $V_{26}\downarrow X$  has a composition factor of dimension $16$ making it impossible for it to be in any parabolic. For the embeddings $(10,10)$ and $(01,01)$, $X\hookrightarrow B_2B_2$, we have that the natural module for $B_4$, $V_9\downarrow B_2B_2=\frac{10+10}{0}$ or $\frac{01^{[1]}+01^{[1]}}{0}$ . Each stabilises a non trivial flag modulo the radical, thus is in a parabolic of $B_4$.

For the result when $M=C_4$, use the graph automorphism of $G$.

This completes the proof of Theorem 4.

\section{Proofs of Corollaries 3, 5, 6, 7 and 8}\label{corproof}

\subsection{Proof of Corollary 3}
Any $G$-irreducible subgroup is not contained in a parabolic subgroup and so is in a reductive maximal subgroup of $G$. Any $G$-cr subgroup is either $G$-ir or contained in a Levi, thus in a maximal subsystem subgroup. So we may assume that $X$ is non-$G$-cr for the purposes of proving this corollary.

(A) Since $X$ is connected, $Z(X)^0$ is a torus. If this torus is positive dimensional then $X$ is contained in the full centraliser of $Z(X)^0$ and so is in a Levi, thus a subsystem subgroup of maximal rank, thence a reductive maximal subgroup. So we may assume that $Z(X)^0=\{1\}$ and $X$ is semisimple. If $X$ is not of type $A_1$ or $A_1A_1$ then Theorem 1(A) and Theorem 2 show that $X$ is contained in a subsystem subgroup of $G$.

(B) This is Theorem 1(B)($p=3$)(vi) and 1(B)($p=2$)(vii) for $X$ of type $A_1$ and Theorem 2(vi) for $X$ of type $A_1^2$.

(C) This is proved in Remark \ref{indecomposablea1}.

(D) When $p=3$ and $X$ is not contained in any reductive maximal subgroup, it is of type $A_1$ as claimed. When $p=2$ and $X$ is not contained in any reductive maximal subgroup, it is of type $A_1$ or $A_1^2$ but it is shown in \ref{a1a1sina1a1pars} that each such subgroup of type $A_1$ is contained in a subgroup of type $A_1^2$.

\subsection{Proof of Corollary 5}
This statement follows from the proofs of Theorems 1, 2 and 4. One needs to go through the proof checking that the case-by-case analysis corresponds with the table in \ref{tor}.

\subsection{Proof of Corollary 6}
This follows from the validity of Corollary 5. Once one knows that each $G$-reducible subgroup is listed somewhere in Table \ref{tor}, it suffices to go through checking the composition factors of isomorphic subgroups $X$ and $Y$ on the $26$-dimensional module $V_{26}$ for $G$. By \ref{comp factors} one can reduce to $G$-cr subgroups. One finds that there is no coincidence amongst the composition factors, collected according to isomorphism type of subgroups minimally contained in non-associated parabolics.

For instance, let $X$ be a subgroup of type $A_1$. Then one finds that $X$ is minimally contained in a long $A_1$-parabolic subgroup of $G$ if and only if it has $6$ isomorphic factors of dimension $2$ and $14$ trivial factors on $V_{26}$. Similarly, let $p=2$ and let $X$ be a subgroup of type $A_1$ or $A_1^2$. Then $X$ is minimally contained in $A_1\tilde A_1$-parabolics if and only if it has the set of collected factors $\{(4,2),(2,4),(2,2),(2,1),(1,4)\}$, where $(n,m)$ indicates a maximal collection of $m$ isomorphic factors of dimension $n$. When $p>2$ the appropriate set of factors is $\{(4,2),(2,4),(2,2),(3,1),(1,3)\}$. The other cases are similar and easier.

\subsection{Subgroups generated by short root elements: Proof of Corollary 7}
Let $X$ be a closed, connected simple subgroup of $G$. Since $X$ is simple it is generated by any one of its conjugacy classes. So if $X$ contains any short root element, then it will be generated by short root elements.

We have a complete list of all simple subgroups of $G$ by collation of the Table in \ref{tor} and \cite[6.2]{Ame05}. Further, each gives the restriction $V_{26}\downarrow X$. Now, in \cite[Table 3]{Law95}, Lawther gives a list of the Jordan block sizes of the action of each class of unipotent elements on the module $V_{26}$. In that table, a short root element is a unipotent element of type `$\tilde A_1$'. We see from the table that an element has Jordan block sizes $J_2^{10}/J_1^{6}$ on $V_{26}$ when $p=2$ and $J_3/J_2^8/J_1^7$ on $V_{26}$ for $p>2$. Here $J_{n_1}^{i_1}/J_{n_2}^{n_2}/\dots$ indicates an action with $i_1$ Jordan blocks of size $n_1$, $i_2$ Jordan blocks of size $n_2$ and so forth.

Our strategy then, is simply to go through tables \ref{tor} and \cite[6.2]{Ame05} and compare each restriction with the Jordan block sizes given above. That is, for each simple subgroup $X$ and class of unipotent elements $u^X$ in $X$, we consider the Jordan block sizes of the action of $u$ on $V_{26}$ and check to see if any matches those above. For subgroups of type $A_1$, we have only one unipotent element to check, $x_\alpha(1)$; whereas for $A_2$, we much check $x_{\alpha}(1)$ and $x_{\alpha+\beta}(1)$.  For $G_2$ we must potentially check several classes. In most cases however, we can make short-cuts.

We will deal with the case $p=2$ first. Here we can employ results of Liebeck and Seitz, \cite[Theorems 2.1 and 3.1]{LS94} on subgroups generated by long root elements, by noting that a long root group is mapped to a short under the graph automorphism $\tau$. From the tables \cite[p333]{LS94} we see that  a simple group having intersection $|U\cap X|>2$ with some long root subgroup $U$ of $G$ is one of the following: \begin{enumerate}\item Subsystem subgroups of types $A_1,\ A_2,\ B_2,\ A_3,\ B_3,\ C_3,\ B_4,\ C_4$, and $D_4$, i.e. any subsystem subgroup not of type $\tilde A_1$ or $\tilde A_2$.
\item Subgroups conjugate to $B_2\leq D_3$, $B_3\leq D_4$, $G_2\leq D_4$.
\end{enumerate}

From \cite[3.1]{LS94}, the only subgroups $X$ having intersection $|U\cap X|=2$ with some long root subgroup $U$ of $G$ are conjugates of $X\cong A_1\leq A_2$ by $V_3\downarrow A_1=W(2)$.

Thus it immediately follows that the images of these subgroups under the graph automorphism are generated by short root elements. Thus the statement of Corollary 7 follows in the case $p=2$.

Now let $p>2$. We deal with subsystem subgroups first. It is clear that simple subgroups containing a short root group are generated by short root elements. This accounts for the subsystem groups $B_2$,\ $B_3$, $C_3$, $B_4$, $\tilde A_1$ and $\tilde A_2$ as stated in the corollary. We wish to show that $A_3$ and $D_4$ subsystems are also generated by short root elements. It suffices to prove this for $A_3$ since $A_3\leq D_4$. Now we claim that $X$ is conjugate to a subsystem $\tilde A_1$ of  $G$. Indeed from \ref{tor} $V_{26}\downarrow A_1^2\leq B_2=(1,1)+((1,0)+(0,1))^4+ 0^6$, so $V_{26}\downarrow X=2+1^{8}+0^7$. Since the table in \ref{tor} contains a single representative of each conjugacy class, looking down the restrictions, one observes that in fact $X$ is conjugate to $\tilde A_1$. Since $X$'s non-trivial unipotent elements are thus short, we are done.

We claim that a long $A_2$ subsystem is not generated by short root elements. There are two unipotent classes in $A_2$. One is in a Levi $A_1$ of $A_2$ and one is in an irreducible $A_1$ of $A_2$. The former can be ruled out since its only unipotent elements are long. For the latter, $V_{26}\downarrow X=2^6+0^8$, so a non-trivial unipotent element has Jordan block sizes $J_3^6/J_1^8$ on $V_{26}$. Thus it is not a short root element. (It is in fact a unipotent element of type $A_2$.) Hence, as neither of the unipotent classes in $A_2$ is short in $G$, $A_2$ is not generated by short root elements. This concludes the case of subsystem subgroups.

We have already observed that the subgroup $X=A_1\hookrightarrow A_1^2$ by $x\mapsto (x,x)$ is conjugate to a short Levi $\tilde A_1$. We claim also that any subgroup $Y=A_1\hookrightarrow A_1^2$ by $(V_2,V_2)\downarrow Y=(1^{[r]},1^{[s]})$ for $r\neq s$ is short. Indeed if $\alpha$ is a root of $A_1$, then the element $x_{\alpha}(1)$ is centralised by Frobenius, and so $Y$ contains the element $x_{\alpha}(1)$ of $X$. Thus it is also generated by short root elements.

We claim that there are no further subgroups $X$ generated by short root elements.

Recall $p>2$, and first suppose $X$ is not of type $A_1$. Then extracting a list of non-subsystem, simple subgroups from \ref{tor}, we must decide if the following subgroups are generated by short root elements:
\begin{enumerate}\item $X=A_2\leq B_3$, $(p=3)$;
\item $X=A_2\leq B_4$ via $V_9\downarrow A_2=T(11)$ $(p=3)$;
\item $X=A_2\leq A_2\tilde A_2$ via $(V_3,V_3)\downarrow X=(10^{[r]},01^{[s]})$ 
\item $X=A_2\leq A_2\tilde A_2$ via $(10^{[r]},(10^{[s]})$;
\item $X=G_2\leq B_3$;
\item $X=G_2$ max, $(p=7)$.\end{enumerate}

Since a short root element has only one Jordan block of size three, we can rule out many of these subgroups immediately. For instance, in item (i),  $V_{26}\downarrow X=11^3+00^5$ and so if $X$ had a short root element $u$, $11\downarrow u$ must contain a Jordan block of size 3. But then $V_{26}\downarrow u$ would have $3$. A contradiction. Similarly item (iv) is ruled out.

Since a unipotent elements $x_{\alpha}(1)$ is centralised by the action of a standard Frobenius map, all the subgroups in item (iii) have the same $G$-classes of unipotent elements. So we may as well just consider $Y=A_2\hookrightarrow A_2\tilde A_2$ by $(V_3,V_3)\downarrow X=(10,01)$. Suppose either this subgroup or the subgroup $X$, described in item (ii) is generated by short root elements. By \ref{tor}, $V_{26}\downarrow Y=T(11)^2+W(11)$ and $V_{26}\downarrow X=T(11)+W(11)+W(11)^*$. If $X$ contained a short root element $u$ then $T(11)\downarrow u$ contains a $J_3$ and the other summands contain $J_2$s and $J_1$s. But then $Y$ would contain an element $u$ acting on $V_{26}$ with two $J_3$s and the remainder $J_2$s and $J_1$s. Such a unipotent element does not exist by \cite[Table 3]{Law95}, a contradiction. Similarly if $Y$ had a short root element, then $W(11)\downarrow u$ would have a $J_3$ and $T(11)$ would consist of $J_2$s and $J_1$s; but  then $X$ would contain an element $u$ acting on $V_{26}$ with two $J_3$s and the remainder $J_2$s and $J_1$s.

The same argument shows that item (iv) contains no short root elements.

To complete the analysis when $X$ is not of type $A_1$, it remains to consider the maximal subgroup of type $G_2$ when $p=7$. By \cite[Table 1]{Law95} we know the actions of all unipotent elements on the seven-dimensional representation $V_7=L_{G_2}(10)$. The classes are denoted by $A_1,\ \tilde A_1,\ G_2(a_1)$ and $G_2$. We may take representatives of each of these four classes by taking a non-trivial unipotent element in each of the following four subgroups $X$ of type $A_1$, respectively: a long $A_1$ Levi, a short $\tilde A_1$ Levi, the irreducible $A_1\leq A_2$ and the maximal $A_1$ (note $p=7$). Now, the symmetric square $S^2(10)\cong T(20)=00/20/00$. So one can deduce, at least, the composition factors of $20\downarrow X$ by taking the symmetric square of its action on $V_7$ and discarding two trivial factors. For example, if $X=\tilde A_1$,  $V_7\downarrow X=2+1^2$. So $S^2(20)\downarrow X=4+3^2+2^3+1^2+0^2$, and so $20\downarrow X=4+3^2+2^3+1^2$. But the action of a unipotent element $u$ on $L_{A_1}(4)$ is a Jordan block of size $5$ and so $u$ is not a short root element. Similarly, taking a unipotent element $u$ of type $G_2$ in the maximal $X=A_1$ of $G_2$, we have $V_7\downarrow X=6$ so $S^2(6)=0/(5\otimes 1^{[1]})/0+T(8)$ where $T(8)=4/8/4$. So $20\downarrow X=5\otimes 1^{[1]}+T(8)$. The first summand gives rise to a Jordan block of size $6$ and so $u$ is not a short root element. The other two cases are similar and easier.

Now we consider the case where $X$ is of type $A_1$ and non-$G$-cr. Thus $p=3$ and $X$ is given by Theorem 1(B) with restrictions $V_{26}\downarrow X$ furnished by \ref{tor}. Notice that since $2\otimes 2=T(4)+2$ for $A_1$, and the module $2$ is projective for the subgroup $L_2(3)\leq A_1$, the tensor product will be projective too, and so will split into a direct sum of isomorphic modules of high weight $2$ on restriction to $L_2(3)$. It follows that a root group element $u$ of $A_1$ contained in $L_2(2)$ has $2\otimes 2\downarrow u=J_3^3$. Thus $T(4)\downarrow u=J_3^2$. Since a short root element has only one Jordan block of size $3$ on $V_{26}$ this rules out any $A_1$ acting with a direct factor of $T(4)$ (or any twist) of it on $V_{26}$. Thus items (i), (iii) and (iv) of Theorem 1(B)($p=3$) can be ruled out. 

Similarly, since $T(3)=2\otimes 1$ and $T(3)\otimes 1=T(4)+2^2$ (since a product of tilting modules is tilting), we can rule out items (ii) and (v).

To establish that the subgroups of item $(vi)$ do not contain short root elements, we have used the GAP realisation of these subgroups. One takes a unipotent element $u$ in such a subgroup acting on $V_{26}$ and checks the socle length by repeatedly asking the MeatAxe to calculate the socle of the action of $u$ and then creating a new module by factoring this out. The number of times one does this determines the socle length of $V_{26}\downarrow u$ and thus the length of the Jordan block. Since it has length bigger than 3 we are done in this case.

It remains to check the non-Levi, $G$-cr subgroups of type $A_1$ when $p>2$. The reducible ones are listed in \ref{tor} and one can use similar arguments see at a glance that none of these subgroups can contain short root elements except for the subgroups $A_1\leq A_1^2\leq B_2$ by $(V_2,V_2)\downarrow X=(1^{[r]},1^{[s]})$. Since we have established that each of these does in fact contain short root elements, the reducible ones are done.  

For the irreducible ones, we turn as promised to \cite[Table 6.2]{Ame05}. Using the restrictions from \cite[6.3]{Ame05} we can use the above arguments to immediately rule out all but possibly item 21, the maximal $X\cong A_1$ in $F_4$. If $p=13$, we have $V_{26}\downarrow X=T(16)=8/16/8$ and if $p=17$, $V_{26}\downarrow X=16+8$. The latter has a $J_{17}$ whereas the former has a submodule of high weight $8$ so has at Jordan block of size at least $J_9$.

Since we have now checked every conjugacy class of simple subgroups of $G$, we have confirmed the statement of Corollary 7.

\subsection{Proof of Corollary 8}
Let $X$ be a subgroup of type $A_1$, contained in no subgroup of type $A_1^2$. We must show that $p=3$ or $p>5$.

For $p=3$ we simply exhibit the subgroups of type $A_1$ described in \S\ref{120}.

If $p=7$ there is a maximal subgroup of type $G_2$ (cf \ref{maxs}). We claim that the maximal subgroup $X\cong A_1$ of this $G_2$ such that $V_7\downarrow X=6$ is in no subgroup of type $A_1^2$. (Note that $X$ does not appear in \ref{tor} as it is, in fact, $G$-irreducible.) Examination of the restriction $V_{26}\downarrow X=5\otimes 1^{[1]}+T(8)$, from \cite[6.3]{Ame05} in tandem with those of the maximal subgroups $M$ of $G$ given in \ref{tor} shows that it is is not contained in any other maximal subgroup $M$ of $G$. Thus its sole proper, connected, semisimple overgroup in $G$ is the maximal $G_2$. So in particular it is not in a subgroup of type $A_1^2$.

If $p\geq 11$ we claim the maximal subgroup $X$ of $B_4$ with $V_9\downarrow X=8$ is not in a subgroup of type $A_1^2$. From \cite[6.3]{Ame05} one has that $V_{26}\downarrow X=8+10+4+0$. Comparison with the restrictions of the maximal subgroups in \ref{tor} again shows that $X$ has no proper, connected, semisimple overgroups in $G$ apart from $B_4$.

(Note that for $p\geq 13$ one could also cite the maximal subgroup of type $A_1$.)

Now if $p=2$ or $p=5$ we must show that every subgroup $X$ of type $A_1$ is contained in a subgroup of type $A_1^2$.

First consider $p=2$. Then If $X$ is minimally contained in an $A_1$-parabolic, then it is one of the subgroups listed in \ref{a1ina1char2} each of which is obviously in a subgroup of type $A_1^2$. If $X$ is minimally contained in an $A_1\tilde A_1$-paraboilc, one simply cites \ref{a1a1sina1a1pars}. If $X$ is minimally contained in a $B_2$, $B_3$- or $C_3$-parabolic then it is in an $A_1^2$ or $A_1^3$ by \ref{a1sina1a1a1s} or \ref{a1sinb2} and so clearly in an $A_1^2$. These are all the cases where $X$ is in a parabolic that need checking as $A_1$ does not embed $A_2$-irreducibly when $p=2$.

For $G$-irreducible $A_1$s one refers to \cite[Table 6.2]{Ame05} to find $G$-irreducible $A_1$s. Up to conjugacy under the graph automorphism $\tau$ of $G$ we have that $X$ is in $C_4$ and has $V_8\downarrow X$ is one of \begin{enumerate} \item $1^{[z]}\otimes 1^{[q]}\perp 1^{[r]}\otimes 1^{[s]}$,
\item $1^{[z]}\otimes 1^{[q]}\perp 1^{[r+1]}\perp 1^{[s+1]}$,
\item $1^{[z+1]}\perp 1^{[q+1]}\perp 1^{[r+1]}\perp 1^{[s+1]}$.\end{enumerate}
Each puts $X$ in a $C_2C_2$ and thus in an $A_1^2$.

Next consider $p=5$. Here $X$ is $G$-cr by \ref{p23orxgcr} and so if $X$ is in a parabolic, it is in a Levi. Since every simply Levi subgroup centralises an $A_1$ except for $B_3$ we may assume that $X\leq B_3$. But by \ref{girsubsconj} the possible $B_3$-irreducible embeddings of $X$ into $B_3$ are given by $V_7\downarrow X=2^{[r]}+2^{[s]}+0$ and $1^{[r]}\otimes 1^{[s]}+2^{[t]}$. Each of these is in an $A_1^2\tilde A_1$ by the same proof as the $p=3$ case in \ref{a1sinb3}.

Now suppose $p=5$ and $X$ is $G$-irreducible, not contained in an $A_1^2$ subgroup of $G$. If it is  contained in a maximal subgroup of type $A_1C_3$, $A_1G_2$ or $A_2\tilde A_2$ then it is in an $A_1^2$, a contradiction. So $X$ is in $B_4$ with $V_9\downarrow X=2^{[r]}\otimes 2^{[s]}$, $4^{[r]}\perp 2^{[s]}\perp 0$, $3^{[r]}\otimes 1^{[s]}\perp 0$, $4^{[r]}\perp 1^{[s]}\otimes 1^{[t]}$, $2^{[r]}\perp 2^{[s]}\perp 2^{[t]}$, or $1^{[z]}\otimes 1^{[q}\perp 1^{[r}\otimes 1^{[t]}$ by \cite[Table 6.2]{Ame05}. But these are in subgroups of type $SO_3\otimes SO_3$, $O_5\times O_4$, $Sp_4\otimes Sp_2$, $O_5\times O_4$, $O_3\times O_3\times O_3$ or $O_4\times O_4$ respectively and in each case, $X$ is in a subgroup of type $A_1^2$.

This completes the proof of Corollary 8.

\chapter{Appendices}

\section{Table of subgroups of $G$} \label{tor}

\begin{prop}[The Restrictions]The restrictions of the 26-dimensional Weyl module $V_{26}=W_G(0001)$ for various subgroups of $F_4$ are given in the table below stating which are R(educible) and N(on-$G$-cr).

{\footnotesize
\begin{longtable}{|c|c|c|c|c|c|c|}
\hline\endlastfoot \hline
$X$ & $p$ & $V_{26}\downarrow X$ & $R$ & $N$\\
\hline
{\bf Maximal} &&&&\\\hline
$A_1G_2$ & $p\neq 2$ & $(2,10)+(W(4),00)$ &&\\
$A_2\tilde A_2$ & any & $(10,01)+(01,10)+(0,W(11))$&&\\
$B_4$ & $p>2$ & $1000+ 0001+ 0000$&&\\
$B_4$ & $p=2$ & $T(1000)+0001$&&\\
$C_4$ & $p=2$ & $0100$&&\\
$A_1C_3$ & $p\neq 2$ & $1\otimes 100+ 0\otimes W(010)$&&\\
$G_2$ & $p=7$ & $20$&&\\
$A_1$ max & $p=13$ & $T(16)$&&\\
$A_1$ max & $p\geq 17$ & $16+8$&&\\\hline\hline
{\bf Max rank, not maximal} &&&&\\\hline
$A_1C_3$ & $p=2$ & $1\otimes 100+ 0\otimes 010$.&&\\
$\tilde A_1B_3$ & $p=2$ & $\begin{array}{c}0\\\hline (2,000)+(0,100)\\\hline 0\end{array}+(1,001)$&&\\
$D_4$ & any & $1000+ 0001+0100+ 0000^2$&&\\
$\tilde D_4$ & $p=2$ & $0100$&&\\
$B_2B_2$ & $p=2$ & $\begin{array}{c}0\\\hline (10,00)+(00,10)\\\hline 0\end{array}+(01,01)$&&\\
$A_1^2B_2$ & $p>2$ & $\begin{array}{c}(0,0,10)+(1,1,00)+(1,0,01)\\+(0,1,01)+(0,0,00)\end{array}$&&\\
$A_1^2B_2$ & $p=2$ & $T((0,0,10))+(1,1,00)+(1,0,01)+(0,1,01)$&&\\
$\tilde A_1^2B_2$ & $p=2$ & $\begin{array}{c}0\\\hline (2,0,00)+(0,2,00)+(0,0,10)\\\hline 0\end{array}+(1,1,01)$\vspace{-10pt}&&\\
$A_1^4$ & any & \parbox{150pt}{\begin{center}$(1,1,0,0)+(0,0,1,1)+(1,0,1,0)+(1,0,0,1)+(0,1,1,0)+(0,1,0,1)+(0,0,0,0)^2$\end{center}}
\vspace{-20pt}&&\\
$A_1^2\tilde A_1^2$ & $p=2$ & \parbox{150pt}{\begin{center}$\begin{array}{c}0\\\hline (0,0,2,0)+(0,0,0,2)\\\hline 0\end{array}$\\$+(1,1,0,0)+(1,0,1,1)+(0,1,1,1)$\end{center}}&&\vspace{-20pt}\\
$\tilde A_1^4$ & $p=2$ & \parbox{150pt}{\begin{center}$\begin{array}{c}0\\\hline (2,0,0,0)+(0,2,0,0)+\dots\\\hline 0\end{array}$\\$+(1,1,1,1)$\end{center}}\vspace{-10pt}&&\\\hline\hline
{\bf Rank $3$, no $A_1$ factors} &&&&\\\hline
$A_3$  & any & $010+100^2+001^2+000^4$&\x&\\
$\tilde A_3$ & $p=2$ & $010^2+101$&\x&\\
$B_3$ Levi & $p>2$ & $100+001^2+000^3$&\x&\\
$B_3$ Levi  & $p=2$ & $T(100)+001^2+000^2$&\x&\\
$B_3\leq D_4$ & $p=2$ & $T(100)+001^2+000^2$&\x&\x\\
$C_3$ & any & $100^2+ W(010)$&\x&\\
$B_3\leq \tilde D_4$ & $p=2$ & $100^{[1]}/010^{[1]}/100^{[1]}$&\x&\x\\\hline\hline
{\bf Rank $2$, no $A_1$ factors} &&&&\\\hline
Long $\bar A_2$ Levi & any & $ 10^3+ 01^3+ 00^8$&\x&\\
Short $\tilde A_2$ Levi & any & $10^3+ 01^3+ W(11)$&\x&\\
$A_2\hookrightarrow B_3$ by $V_7\downarrow A_2=11$ & $p=3$ & $11^3 + 00^5$&\x&\\
$A_2\hookrightarrow B_4$ by $10\otimes 01$ & $p=3$ & $T(11)+ 11/00+00/11 + 00$&\x&\x\\
$A_2\hookrightarrow A_2\tilde A_2$ as $(10,01)$ & $p=3$ & $T(11)+ T(11) + 11/00$&\x&\x\\
$A_2\hookrightarrow A_2\tilde A_2$ as $(10^{[r]},10^{[s]})$ $r\neq s$ & any & $10^{[r]}\otimes 10^{[s]}+01^{[r]}\otimes 01^{[s]}+11^{[s]}$&&\\
$A_2\hookrightarrow A_2\tilde A_2$ as $(10,10)$ & $p>2$ & $20+01+02+ 10 + W(11)$&&\\
$A_2\hookrightarrow A_2\tilde A_2$ as $(10,10)$ & $p=2$ & $T(20)+T(02)+11$&&\\
$A_2\hookrightarrow A_2\tilde A_2$ as $(10^{[r]},01^{[s]})$ $r\neq s$ & any & $10^{[r]}\otimes 01^{[s]}+01^{[r]}\otimes 10^{[s]}+ W(11)^{[s]}$&&\\
$A_2\hookrightarrow A_2\tilde A_2$ as $(10,01)$ & $p\neq 3$ & $11^3+00^2$&&\\
\hline
$B_2$ Levi & $p>2$ & $10+ 01^4+ 00^5$&\x&\\
$B_2$ Levi & $p=2$ & $T(10)+ 01^4+ 00^4$&\x&\\
$B_2\leq D_3$ & $p=2$ & $T(10)+ 01^4+ 00^4$&\x&\x\\
$B_2\leq \tilde D_3$ & $p=2$ & $\left(\begin{array}{c}0\\\hline 01^{[1]}\\\hline 0\end{array}\right)^2+T(20)$&\x&\x\\
$B_2\hookrightarrow B_2B_2$ as $(01^{[1]},01^{[1]})$ & $p=2$ & $\begin{array}{c}0\\\hline 01^{[1]}\end{array}+\begin{array}{c}01^{[1]}\\\hline 0\end{array}+T(20)+00^2$&\x&\x\\
$B_2\hookrightarrow B_2B_2$ as $(10^{[r]},10^{]s]})$ $r\neq s$& $p=2$ & $\begin{array}{c}0\\\hline 10^{[r]}+10^{[s]}\\\hline 0\end{array}+01^{[r]}\otimes 01^{[s]}$&&\\
$B_2\hookrightarrow B_2B_2$ as $(10^{[r]},01^{]s]})$ & $p=2$ & $\begin{array}{c}0\\\hline 10^{[r]}+01^{[s+1]}\\\hline 0\end{array}+01^{[r]}\otimes 10^{[s]}$&&\\
$B_2\hookrightarrow B_2B_2$ as $(01^{[r]},01^{]s]})$ $r\neq s$& $p=2$ & $\begin{array}{c}0\\\hline 10^{[r]}+01^{[s+1]}\\\hline 0\end{array}+01^{[r]}\otimes 10^{[s]}$&&\\
\hline
$G_2\leq B_3$ Levi & $p>2$ & $10^3+00^5$&\x&\\
$G_2\leq B_3$ Levi & $p=2$ & $T(10)^3+00^2$&\x&\\
$G_2\leq D_4$ & $p=2$ & $T(10)^3+00^2$&\x&\x\\
$\tilde G_2\leq \tilde C_3$ & $p=2$ & $10^2+01$&\x&\\
$\tilde G_2\leq \tilde D_4$ & $p=2$ & $T(02)=10^{[1]}/01^{[1]}/10^{[1]}$&\x&\x\\
\hline\hline
{\bf Rank $3$ with $A_1$ factors} &&&&\\\hline
$A_1G_2\leq \tilde A_1B_3\leq B_4$ & $p=2$ & $\begin{array}{c}0\\\hline (2,0)+(0,10)\\\hline 0\end{array}+(1,T(10))$ && \\
$A_1G_2\leq A_1C_3\leq C_4$ & $p=2$ & $(1,10)+(0,01)$ && \\
$A_1C_2\leq C_3$ & $p>2$ & $(1,00)^2 + (0,10)^2+ (1,10)+ (0,01)+(0,00)$&\x&\\
$A_1C_2\leq C_3$ & $p=2$ & $(1,00)^2 + (0,10)^2+ (1,10)+T((0,01))$&\x&\\
$A_1B_2\leq C_4$ ($1$-dim) & $p=2$ & $\begin{array}{c}(2,00) +(0,10)\\\hline (2,10)+\begin{array}{c}0\\\hline 0,01^{[1]}\\\hline 0\end{array}\\\hline (2,00)+(0,10)\end{array}$&\x&\x\\
$\tilde A_1B_2\leq B_3$ & $p=2$ & $\begin{array}{c}0\\\hline (2,00)+(0,10)\\\hline 0\end{array}+(1,01)^2+(0,00)^2$&\x&\\
$\tilde A_1B_2\leq B_4$ ($1$-dim) & $p=2$ & $\begin{array}{c}0\\\hline (2,00)+(0,10)\\\hline 0\end{array}+(1,01)^2+(0,00)^2$&\x&\x\\
$A_1\tilde A_2$ Levi & any & $(0, 10)+ (0,01)+ (1, 10)+ (1,01)+ (0,11)$&\x&\\
$A_1\tilde A_2\leq A_2\tilde A_2$ by $(2/0,10)$ & $p=2$ & $(2/0,01)+(0/2,10)+(0,11)$&\x&\x\\
$A_2\tilde A_1$ Levi & any & $\begin{array}{c}(00, 1)^2+ (10, 0)+ (01,0)+ (10,1)\\ + (01,1)+ (00,W(2))+ (0,0)\end{array}$&\x&\\
$A_2\tilde A_1\leq A_2\tilde A_2$ by $(10,2/0)$ & $p=2$ & $(10,0/2)+(01,2/0)+(00,T(4))$\vspace{-10pt}&\x&\x\\
$A_1^2\tilde A_1\leq B_3$ Levi & $p>2$ & \parbox{130pt}{\begin{center}$(1,1,0)+(0,0,2)+[(1,0,1)+(0,1,1)]^2+(0,0,0)^3$\end{center}}\vspace{-20pt}&\x&\\
$A_1^2\tilde A_1\leq B_3$ Levi & $p=2$ & \parbox{130pt}{\begin{center}$(1,1,0)+T(0,0,2)+[(1,0,1)+(0,1,1)]^2+(0,0,0)^2$\end{center}}&\x&\\
$A_1^3\leq A_1^4; (x,y,z)\mapsto (x,y,z,z)$ & $p=2$ & \parbox{130pt}{\begin{center}$(1,1,0)+T(0,0,2)+[(1,0,1)+(0,1,1)]^2+(0,0,0)^2$\end{center}}&\x&\x\\
$\tilde A_1^3\leq B_3$ Levi & $p=2$ & $\begin{array}{c}0\\\hline(2,0,0)+(0,2,0)+(0,0,2)\\\hline 0\end{array}+0^2+(1,1,1)^2$&\x&\\
$A_1^3\leq A_1^2\tilde A_1^2; (x,y,z)\mapsto (x,x,y,z)$ & $p=2$ & $T(2,0,0)+\begin{array}{c}0\\\hline (0,2,0)+(2,0,0)\\\hline 0\end{array}+(1,1,1)$&\x&\x\\
$\tilde A_1^3\leq B_3\leq D_4$ & $p=2$ & $\begin{array}{c}0\\\hline(2,0,0)+(0,2,0)+(0,0,2)\\\hline 0\end{array}+0^2+(1,1,1)^2$&\x&\x\\
$A_1^3\leq \tilde A_1B_3$ n.r.o. of $A_1^2$ in $B_3$ & $p=2$ & $\begin{array}{c}0+0\\\hline(2,0,0)+(0,2,0)+(0,0,2)\\\hline 0+0\end{array}+(1,1,1)^2$&\x&\x\\
$A_1^3\leq A_1B_2\leq B_4$, $A_1B_2$ in n.r.o. of $B_4$ & $p=2$ & $\begin{array}{c}0+0\\\hline(2,0,0)+(0,2,0)+(0,0,2)\\\hline 0+0\end{array}+(1,1,1)^2$&\x&\x\\
$X\leq B_4$ no red. overgroup & $p=2$ & $\begin{array}{c}0+0\\\hline(2,0,0)+(0,2,0)+(0,0,2)\\\hline 0+0\end{array}+(1,1,1)^2$&\x&\x\\
$A_1\tilde A_1^2$ & $p=2$ & $((1,0,0)+(0,1,1))^2+(1,1,1)+(0,1,1)+0^2$&\x&\\
$A_1^3\leq \tilde A_1^4$;\ $(x,y,z)\mapsto (x,y,z,z)$ & $p=2$ & $\begin{array}{c}0\\\hline(2,0,0)+(0,2,0)+(0,0,2)^2\\\hline 0\end{array}+(1.1.T(2))$&\x&\x\\
$A_1^3\leq C_3$ & any & $\begin{array}{c}((1,0,0)+(0,1,0)+(0,0,1))^2\\+(1,1,0)+(0,1,1)+(1,0,1)+0^2\end{array}$&\x&\\
$A_1^3\leq A_1^2\tilde A_1^2$;\ $(x,y,z)\mapsto (x,y,z,z)$ & $p=2$ & $(1,1,0)+\begin{array}{c}0\\\hline (0,0,2)^2\\\hline 0\end{array}+(1.0.T(2))+(0,1,T(2))$&\x&\x\\
$\tilde A_1^3\leq B_3\leq \tilde D_4$ & $p=2$ & $\begin{array}{c}(2,0,0)+(0,2,0)+(0,0,2)\\\hline (2,2,0)+(0,2,2)+(2,0,2)+0^2\\\hline (2,0,0)+(0,2,0)+(0,0,2)\end{array}$&\x&\x\\
$A_1^3\leq A_1B_3$, $A_1^2$ in n.r.o of $B_3$ & $p=2$ & $\begin{array}{c}(2,0,0)+(0,2,0)+(0,0,2)\\\hline (2,2,0)+(0,2,2)+(2,0,2)+0^2\\\hline (2,0,0)+(0,2,0)+(0,0,2)\end{array}$&\x&\x\\
$A_1^3\leq A_1B_2\leq \tilde C_4$, n.r.o of $A_1B_2$ in $C_4$ & $p=2$ & $\begin{array}{c}(2,0,0)+(0,2,0)+(0,0,2)\\\hline (2,2,0)+(0,2,2)+(2,0,2)+0^2\\\hline (2,0,0)+(0,2,0)+(0,0,2)\end{array}$&\x&\x\\
$A_1^3$ in n.r.o of $C_4$  & $p=2$ & $\begin{array}{c}(2,0,0)+(0,2,0)+(0,0,2)\\\hline (2,2,0)+(0,2,2)+(2,0,2)+0^2\\\hline (2,0,0)+(0,2,0)+(0,0,2)\end{array}$&\x&\x\\
\hline\hline
{\bf $X$ of type $A_1A_1$} &&&&\\\hline
$A_1\tilde A_1$ Levi & $p>2$ & $(1,0)^2+ (0,1)^4+(0,2) + (1,1)^2+ 0^3$&\x&\\
$A_1\tilde A_1$ Levi & $p=2$ & $(1,0)^2+ (0,1)^4+(0,T(2)) + (1,1)^2+ 0^2$&\x&\\
$A_1^2\leq B_2$ & any & $(1,1)+((1,0)+(0,1))^4+ 0^6$&\x&\\
$A_1^2\leq A_1C_2\leq C_3$ & $p\geq 5$ & $(1,0)^2+(0,3)^2+0+(0,4)+(1,3)$&\x&\\
$A_1^2\leq C_3$ by $V_6\downarrow X=(2,1)$ & $p=3$ & $(2,1)^2+(W(4),0)+(2,2)$&\x&\\
$A_1^2\leq C_3$ by $V_6\downarrow X=(2,1)$ & $p\geq 5$ & $(2,1)^2+(4,0)+(2,2)$&\x&\\
\begin{tabular}{c}$A_1^2\leq D_3\tilde A_1\leq B_4$;\\$(V_6,V_2)\downarrow X=(T(4),1)$\end{tabular} & $p=3$ & $(T(4),0)+(0,2)+0+(1,W(3))+(1,W(3)^*)$&\x&\x\\
\begin{tabular}{c}$A_1^2\leq A_1C_3;$\\ $(V_2,V_6)\downarrow X=(1,T(3))$\end{tabular} & $p=3$ & $(1,T(3))+(0,W(4))+(0,W(4)^*)+0$&\x&\x\\
$\tilde A_1^2\leq B_2$ & $p=2$ & $\begin{array}{c}0\\\hline (2,0)+(0,2)\\\hline 0\end{array}+0^4+(1,1)^4$&\x&\\
$A_1^2\leq A_1^4$; $(x,y)\to (x,x,y,y)$ & $p=2$ & $T(2,0)+T(0,2)+0^2+(1,1)^4$&\x&\x\\
$A_1^2\leq A_1^2\tilde A_1$; $(x,y)\to (x,x,y)$ & $p=2$ & $T(2,0)+T(0,2)+0^2+(1,1)^4$&\x&\x\\
$A_1^2\leq B_3$; no red. overg. & $p=2$ & $T(2,0)+T(0,2)+0^2+(1,1)^4$&\x&\x\\
$A_1^2\leq \tilde A_1^4$; $(x,y)\to (x,x,y,y)$ & $p=2$ & $\begin{array}{c}0\\\hline (2,0)^2+(0,2)^2\\\hline 0\end{array}+(T(2),T(2))$&\x&\x\\
$A_1^2\leq \tilde A_1^2 A_1$; $(x,y)\to (x,x,y^2)$ & $p=2$ & $(0,2)^2+(T(2),0)^2+(T(2),2)+(2/0,0)+(0/2,0)$&\x&\x\\
$A_1^2\leq \tilde C_3$; no red. overg. & $p=2$ & $\left(\begin{array}{c}0\\\hline(2,0)+(0,2)\\\hline 0\end{array}\right)^2+T(2,2)$&\x&\x\\
$\begin{array}{c}A_1^2\leq A_1^3 \text{ minimally }\\\text{contained in a }B_3\text{-parab}\end{array}$& $p=2$ & see \ref{a1a1sina1a1a1s} and \ref{a1a1a1}.&\x&\x\\
\parbox{120pt}{\begin{center}$A_1^2$ minimally contained in\\ an $A_1\tilde A_1$-parabolic, not in $A_2\tilde A_2$\end{center}} & $p=2$ & see \ref{a1a1sina1a1pars} and \ref{a1ina1a1char2}&\x&\x\\

$A_1^2\hookrightarrow A_1\tilde A_2$; $(1,10)\downarrow X=(1,0/2)$ & $p=2$ & $ (1,0/2)+(1,2/0)+(0,T(4))$&\x&\x\\
$A_1^2\hookrightarrow A_2\tilde A_1$; $(10,1)\downarrow X=(0/2,1)$ & $p=2$ & $(0/2,1)+(2/0,1)+(0,T(2))+(0,1)^2$&\x&\x\\
$\begin{array}{c}A_1^2\hookrightarrow A_2\tilde A_2;\\ (10,10)\downarrow X=(2/0,0/2)\end{array}$& $p=2$ &  $\begin{array}{c}(0,0)\\\hline (2,0)+(0,2)\\\hline (2,2)\end{array}+ \begin{array}{c}(2,2)\\\hline (2,0)+(0,2)\\\hline (0,0)\end{array}+(0,T(4))$&\x&\x\\
$\begin{array}{c}A_1^2\hookrightarrow A_2\tilde A_2; \\ (10,10)\downarrow X=(2/0,2/0)\end{array}$ & $p=2$ & $\begin{array}{c}(0,2)\\\hline (2,2)\\\hline (2,0)\end{array}+ \begin{array}{c}(2,0)\\\hline (2,2)\\\hline (0,2)\end{array}+ (0,T(4))+(0,0)^2$&\x&\x\\

\hline\hline
{\bf $X$ of type $A_1$} &&&&\\\hline
Long $A_1$ Levi & any  & $1^6+ 0^{14}$&\x&\\
Short $A_1$ Levi & $p\geq 3$ & $1^8+ 2+ 0^7$&\x&\\
Short $A_1$ Levi & $p=2$ & $1^8+ T(2)+ 0^6$&\x&\\
$A_1\leq A_2$ Levi, irreducible & $p>2$ & $2^6+0^8$ & \x &\\
$A_1\leq \tilde A_2$ Levi, irreducible & $p=3$ & $2^7+1\otimes 1^{[1]}/0$ & \x &\\
$A_1\leq B_2$ Levi, irreducible & $p\geq 5$ & $4+0^5+3^4$ & \x& \\
$A_1\leq B_3$ Levi, irreducible & $p\geq 7$ & $6^3+0^3$&\x&\\
$\begin{array}{c}A_1\leq A_1C_2\leq C_3\text{ by }\\ (V_2,V_4)\downarrow X=(1^{[r]},3^{[s]}),\ rs=0\end{array}$ & $p\geq 5$ & $(1^{[r]})^2+(3^{[s]})^2+(1^{[r]}\otimes 3^{[s]})+ 4^2+2+0$&\x&\\
\begin{tabular}{c}$A_1\leq C_3$ by $V_6\downarrow X=2^{[r]}\otimes 1^{[s]}$\\
exactly one of $r,\ s$ is zero.\end{tabular}& $p=3$ & $(2^{[r]}\otimes 1^{[s]})^2+4^{[r]}+2^{[r]}\otimes 2^{[s]}$&\x&\\
$A_1\leq C_3$ Levi, irreducible & $p=7$ & $5^2+T(8)$&\x&\\
$A_1\leq C_3$ Levi, irreducible & $p\geq 11$ & $5^2+8+4$&\x&\\
\begin{tabular}{c}$A_1\hookrightarrow A_1^2\tilde A_1\leq B_3$ by $(V_2,V_2,V_2)\downarrow X=$\\$(1^{[r]},1^{[s]},1^{[t]})$ ($r\neq s$)\end{tabular}&$p\geq 3$&$1^{[r[}\otimes 1^{[s]}+2^{[t]}+(1^{[r]}\otimes 1^{[t]}+1^{[s]}\otimes 1^{[t]})^2+0$&\x&\\
\begin{tabular}{c}$A_1\hookrightarrow A_1^2\leq B_2$ by \\$(V_2,V_2)\downarrow X=(1^{[r]},1^{[s]})$ ($r\neq s)$\end{tabular} & any & $1^{[r]}\otimes 1^{[s]}+ 0^6+ (1^{[r]}+1^{[s]})^4$&\x&\\
$A_1\hookrightarrow A_1\tilde A_1$ by $(V_2.V_2)\downarrow X=(1^{[r]}.1^{[s]})$&$p\geq 3$ & 
${1^{[r]}}^2+{1^{[s]}}^4+2^{[s]}+(1^{[r]}\otimes 1^{[s]})^2+0^3$ &\x & \\
$A_1\hookrightarrow A_1\tilde A_1$ by $(V_2.V_2)\downarrow X=(1^{[r]}.1^{[s]})$&$p=2$ & 
${1^{[r]}}^2+{1^{[s]}}^4+T(2)^{[s]}+(1^{[r]}\otimes 1^{[s]})^2+0^2$ &\x & \\\hline
\begin{tabular}{c}$A_1\hookrightarrow D_3\tilde A_1\leq B_4$ by \\ $(V_6,V_2)\downarrow X=(T(4)^{[r]},2^{[s]})$\end{tabular}&$p=3$&$\begin{array}{c}T(4)^{[r]}+2^{[s]}+ (1^{[r+1]}/1^{[r]})\otimes 1^{[s]}\\+(1^{[r]}/1^{[r+1]})\otimes 1^{[s]}+  0\end{array}$&\x&\x\\

$A_1 \hookrightarrow A_1C_3$ by $(1^{[r]},T(3)^{[s]})$ & $p=3$ & $1^{[r]}\otimes T(3)^{[s]}+\frac{1^{[s]}\otimes 1^{[s+1]}}{0}+\frac{0}{1^{[s]}\otimes 1^{[s+1]}}+0$&\x&\x\\

$A_1\hookrightarrow A_2\tilde A_2$; $(2,2)$ diagonal & $p=3$ & $T(4)^2+ 2^3+ 1\otimes 1^{[1]}/ 0$&\x&\x\\
$A_1\hookrightarrow D_3\leq B_3$ by $V_6\downarrow X=T(4)$ & $p=3$ & $T(4)+(1/1^{[1]}+1^{[1]}/1)^2+ 0^4$&\x&\x \\
$A_1\hookrightarrow C_3$ by $2\otimes 1=T(3)$ & $p=3$ & $T(3)^2+ 2+ (1\otimes 1^{[1]})/0+0/(1\otimes 1^{[1]})+0$&\x&\x\\
$A_1$  no overgroup; $1$-dim of such & $p=3$ & $\begin{array}{c}1\otimes 1^{[2]}\\\hline 1\otimes 1^{[1]}\\\hline 1\otimes 1^{[2]}\end{array}+\begin{array}{c}0\\\hline 1\otimes 1^{[1]}+1^{[1]}\otimes 1^{[2]}\\\hline 0\end{array}+2+0$&\x&\x \\
\hline
$A_1\hookrightarrow A_2$ by $Sym^2 V_2$ & $p=2$ & $(0/2)^3+ (2/0)^3 + 0^8$&\x&\x\\
$A_1\hookrightarrow \tilde A_1^2\leq B_2$ diag & $p=2$ & $(0/2)+(2/0)+T(2)^4+0^4$&\x&\x\\
$A_1\hookrightarrow A_3$; $V_4\downarrow X=T(2)$ & $p=2$ & $(0/2)+(2/0)+T(2)^4+0^4$&\x&\x\\
$A_1\hookrightarrow A_1^4$ diag & $p=2$ & $T(2)^6+0^2$&\x&\x\\
$A_1\hookrightarrow \tilde A_2$ by $Sym^2 V_2$ & $p=2$ & $(0/2)^3+ (2/0)^3+ T(4)$&\x&\x\\
$A_1\hookrightarrow A_1A_1\leq B_2$ diag& $p=2$ & $T(4)+2^8+0^6$&\x&\x\\
$A_1\hookrightarrow \tilde A_1^4$ diag & $p=2$ & $\begin{array}{c}0\\\hline 2^4\\\hline 0\end{array}+T(4)+T(2)+T(2)$&\x&\x\\
$A_1\hookrightarrow \tilde A_3$; $V_4\downarrow X=T(2)$ & $p=2$ & $(0/2)^3+(2/0)^3+T(4)$&\x&\x\\
\hline
\begin{tabular}{c}$A_1$ minimally contained in\\ a $B_3$-parabolic\end{tabular} & $p=2$ & \begin{tabular}{c}see \ref{a1a1a1}: contained in some $A_1^3$ above.\\Use \ref{samesubs} to get restrictions.\end{tabular} &\x&\x\\
\begin{tabular}{c}$A_1$ minimally contained in\\ a $B_2$-parabolic\end{tabular}  & $p=2$ & \begin{tabular}{c} see \ref{a1sinb2}: contained in some $A_1^2$ above \\Use \ref{samesubs} to get restrictions.\end{tabular}  &\x&\x\\
\begin{tabular}{c}$A_1$ minimally contained in\\ an $A_1\tilde A_1$-parabolic\end{tabular} & $p=2$ & \begin{tabular}{c}see \ref{a1ina1a1char2}; for $A_1 \in A_2\tilde A_2$, use \ref{samesubs}\\ and $A_1^2\leq A_2\tilde A_2$ above\end{tabular} &\x&\x\\
\hline
\\\hline
\end{longtable}
} 
\begin{changemargin}{-1.0in} 
\end{changemargin}

Here square brackets enclose a submodule for which only the composition factors have been calculated, rather than the specific restrictions. 
\end{prop}
\begin{proof}
The restrictions $V_{26}\downarrow B_4,\ C_4\ (p=2),\ A_2\tilde A_2,\ D_4,\ G_2\ (p=7)$ and $\tilde D_4 \ (p=2)$ are given in \cite[Table 10.2]{LS04}. Restricting further to the subgroups described gives most of the entries of \ref{tor}. The calculation of some restrictions require observations made in the proof above; specifically, those where we have claimed that there is no proper reductive overgroup of $X$ in some (subsystem) subgroup $M$ containing $X$.

We now show how to obtain some of the non-routine restrictions. 

We first justify the entry $A_2\hookrightarrow B_4$ by $10\otimes 01$. Let $p=3$ and take $X:=A_2$ embedded in $GL_9(K)$ via the action of $X$ on the module $M:=10\otimes 01$. Since $M$ is the tensor product of two tilting modules, the result is also tilting. Because the highest weight of $M$ is $11$, $M$ must contain a $T(11)\cong 00/11/00$, which is $9$-dimensional and therefore the whole of $M$. Now $M$ is self-dual and so $X$ must preserve a bilinear form. Since the dimension of the module is $9$ and indecomposable, it must be in an $SO_9$ subgroup of $GL_9$; this justifies the inclusion of $A_2$ into $B_4$ in the way described. 

We have the restriction $1000\downarrow X$, so to finish, we need the restriction of the spin module $0001$ for $B_4$ to $X$. Take instead a regular $A_1$, $Y\leq X$. Then $V_{\text{nat}}$ for $B_4$ restricted to $Y$ is $2\otimes 2=T(4)+2$ and this is contained in a $D_3B_1$ subsystem of $B_4$. Using \cite[2.7]{LS96} we see that $0001\downarrow D_3B_1=100\otimes 1+001\otimes 1$. Now  $010=\bigwedge^2(100)$, and we notice that $\bigwedge^2(1^{[1]}/1)=T(4)$, so it follows that $0001\downarrow Y=(1^{[1]}/1)\otimes 1+(1/1^{[1]})\otimes 1$. As $\Hom(K,1/1^{[1]}\otimes 1)=\Hom(1,1/1^{[1]})=0$, the structure of $(1/1^{[1]})\otimes 1=2+(0/1\otimes 1^{[3]})$. The only way this can happen is if $0001\downarrow X=11/00+00/11$.

Similarly, consider the $\tilde G_2\leq C_3$ when $p=2$. We use the fact that the natural module $V_6$ for $C_3$ has restriction $V_6\downarrow G_2=10$. Then we recall that the $15$-dimensional module $\bigwedge^2 100=010+000$ for $C_3$. Thus we calculate that $010\downarrow G_2=01$ as required.

When $p=3$, consider $X\cong A_1^2\leq C_3$ by the embedding $V_6\downarrow X=(2, 1)$. Now, we know  $V_{26}\downarrow C_3=100^2+W(010)$. To find $W(010)\downarrow X$ we use the fact that $T(010)=\bigwedge^2(100)$. We calculate that $\bigwedge^2((2,1))=(T(4),0)+(2,2)$ using the fact that since $p\neq 2$ $\bigwedge^2(100)$ is a direct summand of $100\otimes 100$ which in turn is a tilting module for $X$ since it is a tensor square of $T(2,1)$. Now $W(010)$ is a codimension $1$ subspace of $T(010)$ so the only possibility is that $W(010)\downarrow X=(W(4),0)+(2,2)$ as required. It thus follows also that the (non-$G$-cr) subgroup $X\cong A_1\leq C_3$ by $V_6\downarrow X=2\otimes 1=T(3)$ has restriction $V_{26}\downarrow X=W(010)\downarrow X=(W(4),0)+(2,2)$.

We have had to do some of these calculations using the MeatAxe. Consider for example $X\cong A_1\leq \tilde A_3$ by $V_4\downarrow X=T(2)$. Routine calculations show that $V_{26}\downarrow X$ consists of the direct sum of $T(4)$ with some collection of $W(2)$s, $W(2)^*$s, $T(2)$s and $2$s. We can use the MeatAxe to establish how many of each. Since $T(2)=1\otimes 1$, we can easily construct two  $4\times 4$ matrices with indeterminants $\mathtt{T}$ representing the positive and negative root groups $x_{\pm\alpha}(\mathtt{T})$ of $X$ in $SL_4$. Then one can use these matrices to write the generators of $X$ in terms of the roots of $A_3$. Inputting these into GAP using {\tt F4gens} and taking values of $\mathtt T$ in $F_4$ we get a subgroup $L_2(4)$ in $G$ acting on $V_{26}$. Then one can easily establish the composition factors of socle layers and argue that the representation decomposes consistently with the one given.

We give some examples of calculations when $X$ is of type $A_1,\ A_1^2$ or $A_1^3$. For instance, when $p=2$ and $X\cong A_1^3$ is in no reductive subgroup of $B_4$, we see in \ref{a1a1a1} that $V=T_{B_4}(1000)\downarrow X=(2,0,0)|(0,2,0)|(0,0,2)|0^4$, where $T_{B_4}(1000)$ is in fact the natural module for $D_5$ restricted to the maximal $B_4$ subgroup of $D_5$. If the restriction is not as given, then either $V\downarrow X=T(2,0,0)+\begin{array}{c}0\\\hline (0,2,0)+(0,0,2)\\\hline 0\end{array}$ or $X$ stabilises a non-degenerate $2$-space in $V$. If the latter, $X$ in a $D_4$ of $D_5$, putting $X$ in a $B_3$ of $B_4$. If the former, $V$ has the same restriction to $X$ as to a subgroup of $B_2B_1\leq D_3D_2\leq D_5$ and the conjugacy class of $X$ is determined up to $D_5$-conjugacy by the action of $X$ on $V$.  This would put $X$ in a $D_3B_1$ subgroup of $B_4$.
\end{proof}

\section{gap programs}

For the strongly hypothetical edification of the reader we provide a couple of examples of the sorts of programs we have used in GAP to do the various calculations made throughout this paper.

The following short program {\tt subs(mat,ind,exp)} takes a matrix {\tt mat} defined over an extension of $\mathbb F_4$ by some indeterminates, a list of indeterminates {\tt ind} and a list of expressions {\tt exp} which should replace each item in {\tt ind} respectively.

{\tt
subs:=function(mat,ind,exp)\\
local i,j,res;\\
res:=NullMat(26,26)*Z(4);\\
for i in [1..26] do\\
for j in [1..26] do\\
res[i][j]:=Value(mat[i][j],ind,exp);\\
od;od;\\
res:= NullMat(26,26)*Z(4)+res;\\
return res;\\
end;}

The next instruction adds to a list {\tt gens} a list of matrices {\tt mt} with the parameter {\tt T} replaced by $Z(4)^k$ for all $k$ from $1$ to $3$ where $Z(4)$ is a generator of the finite field $GF(4)$.

\begin{verbatim}ad:=function(gens,mt)
local k;
for k in [0..3] do
Append(gens,[subs(mt,[T],[Z(4)^k])]);
od;
end;\end{verbatim}

In \S\ref{120} (and elsewhere but we consider the former here), we need to use GAP to construct subgroups in no reductive maximal subgroup of $G$. The next instruction sets {\tt l} equal to the unipotent element $x_{+}(T)$ of the $A_1$ subgroup of the $A_1^3$ subgroup of the Levi $B_3$  defined by twists $(1,2,0)$. Here, {\tt f[n]=f4gens[n]} is a $26\times 26$ matrix with indeterminate {\tt T} corresponding to the image of $x_{n}(T)$ in $GL(V_{26})$ where $n$ is a number from $1$ to $48$ indexing the roots of $G$. The program multiplies elements {\tt f[n]} together to give the image of $x_{+}(T)$ in $GL(V_{26})$.

\begin{verbatim} l:=subs(f[40],[T],[T^3])*subs(f[2],[T^9],[S])*subs(f[3],[T],[T]);\end{verbatim}

Then one would like to construct a cocycle $\gamma:x_+(T)\to Q$ so that $x_{+}(T)\gamma(x_+(T))$ is isomorphic to $x_{+}(T)$; i.e. a complement to $Q$ in $x_+(T)Q$. We know from \ref{a1h1p3} what this ought to look like on $Q/Q(2)$:

\begin{verbatim}gamma:=subs(f[5],[T],[k*T^2])*subs(f[1],[T],[k*T]);\end{verbatim}

One may then wish to check, for instance, that the property that {\tt l*gamma=IdentityMat(26,GF(4))} holds.

Similarly one defines a map $\delta:x_{+}(T)\to Q(2)$ and combines with $\gamma$ to get any  complement to $Q$. A similar process with $x_-(T)$ gives generators. Working out exactly how the cocycles should fit together can take some time. Finally, one can use the program {\tt ad} to get a set of generators {\tt gens} of the subgroup $X(GF(4))$.

The MeatAxe can be accessed via 

{\tt gm:=GModuleByMats(gens,GF(4));}

and one can start to use the various routines like {\tt MTX.BasisSocle()} to establish the module structure of $V_{26}\downarrow X$.

For reference, a generically behaved subgroup of \S\ref{120} is in fact given by the generators {\tt xtf} and {\tt yf} below.

\begin{verbatim}zl1:=subs(f[40],[T],[T^3])*subs(f[2],[T],[T^9])*subs(f[3],[T],[T]);;  
zl2:=subs(f[16],[T],[T^3])*subs(f[26],[T],[T^9])*subs(f[27],[T],[T]);;
x:=zl1*subs(f[5],[T],[T^2])*subs(f[1],[T],[T]);;   
y:=zl2*subs(f[17],[T],[T^2])*subs(f[19],[T],[-T]);;
xt:=x*subs(f[22],[T],[-T]);; 
xtf:=xt*subs(f[18],[T],[T^6])*subs(f[24],[T],[T^3]);;
yf:=y*subs(f[23],[T],[T^6])*subs(f[15],[T],[-T^3]);;  
\end{verbatim}
{\footnotesize
\bibliographystyle{amsalpha}
\bibliography{noteonf4b.bib}}
\end{document}